\documentclass[a4paper,12pt,twoside]{article}
\usepackage{graphicx}\usepackage{a4wide}
\usepackage{amssymb,amsmath,amsthm}
\usepackage{siunitx}
\usepackage{subfigure,color,colordvi,graphicx,xcolor}
\usepackage[only,llbracket,rrbracket,interleave]{stmaryrd}
\usepackage[superscript,sort,compress]{cite}
\usepackage{enumitem}

\newcommand{\bm}[1]{\text{\boldmath $#1$\unboldmath}}

\newcommand{\vect}[1]{\mathbf{#1}}
\newcommand{\mat}[1]{\mathbf{#1}}

\newcommand{\node}[1]{\mathrm{#1}}

\newcommand{\eltwo}{\ensuremath{\mathcal{L}_2}}

\newcommand{\Ga}[1][1]{\ensuremath{\Gamma_{\!\! #1}}}

\newcommand{\nsd}    {\texttt{n}_{\texttt{sd}}}
\newcommand{\npar}{\texttt{n}_{\texttt{pa}}}

\newcommand{\numel}{\texttt{n}_{\texttt{el}}}

\DeclareMathOperator{\adj}{adj \!}
\DeclareMathOperator{\Det}{det \!}

\DeclareMathOperator*{\argmin}{arg\,min}

\newcommand{\Div}{{\bm{\nabla}\!\!_{\bmu} \cdot\,}}
\newcommand{\Grad}{\bm{\nabla}\!\!_{\bmu}}

\newcommand{\bu}{\bm{u}}
\newcommand{\bs}{\bm{s}}
\newcommand{\bt}{\bm{g}_N}
\newcommand{\bn}{\bm{n}}
\newcommand{\bD}{\bm{D}}
\newcommand{\bE}{\bm{E}}
\newcommand{\Insd}{\mat{I}_{\nsd\!}}

\newcommand{\bL}{\bm{L}}
\newcommand{\hu}{\hat{u}}
\newcommand{\bhu}{\bm{\hu}}
\newcommand{\btau}{\bm{\tau}}

\newcommand{\sVh}{\mathcal{V}^h}
\newcommand{\shVh}{\mathcal{\widehat{V}}^h}
\newcommand{\sLh}{\mathcal{L}^h}
\newcommand{\Lh}{\bm{\mathcal{L}}^h}

\newcommand{\sVmu}{\mathcal{V}^{\bmu}}
\newcommand{\Vmu}{\bm{\mathcal{V}}^{\bmu}}
\newcommand{\Wmu}{\bm{\mathcal{W}}^{\bmu}}
\newcommand{\hVmu}{\bm{\mathcal{\widehat{V}}}^{\bmu}}

\newcommand{\Pk}{\mathcal{P}^{\nDeg}}

\newcommand{\de}{\delta\!}

\newcommand{\Vh}{\bm{\mathcal{V}}^h}
\newcommand{\Wh}{\bm{\mathcal{W}}^h}
\newcommand{\hVh}{\bm{\mathcal{\widehat{V}}}^h}

\newcommand{\bW}{\bm{W}}
\newcommand{\bv}{\bm{v}}
\newcommand{\hv}{\hat{v}}

\newcommand{\bhv}{\bm{\hat{v}}}

\newcommand{\bx}{\bm{x}^\bmu}
\newcommand{\bX}{\bm{x}}

\newcommand{\bmu}{\bm{\mu}}

\newcommand{\intExI}[2]{\big(#1,#2\big)_{\Omega_e \times \bI }}
\newcommand{\intBExI}[2]{\langle #1,#2   \rangle_{\partial \Omega_e \times \bI }}

\newcommand{\intBNoDxI}[2]{\langle #1,#2 \rangle_{( \partial \Omega_e \setminus \Ga[D] ) \times \bI }}
\newcommand{\intBNoDSxI}[2]{\langle #1,#2 \rangle_{( \partial \Omega_e \setminus (\Ga[D] \cup \Ga[S] )) \times \bI }}
\newcommand{\intBDxI}[2]{\langle #1,#2   \rangle_{( \partial \Omega_e \cap \Ga[D] ) \times \bI }}
\newcommand{\intBNxI}[2]{\langle #1,#2   \rangle_{( \partial \Omega_e \cap \Ga[N] ) \times \bI }}
\newcommand{\intBSxI}[2]{\langle #1,#2   \rangle_{( \partial \Omega_e \cap \Ga[S] ) \times \bI }}

\newcommand{\intE}[2]{\big(#1,#2\big)_{\Omega_e }}
\newcommand{\intI}[2]{\big(#1,#2\big)_{\bI}}
\newcommand{\intBE}[2]{\langle #1,#2   \rangle_{\partial \Omega_e}}

\newcommand{\intBNoD}[2]{\langle #1,#2 \rangle_{\partial \Omega_e \setminus \Ga[D]}}
\newcommand{\intBNoDS}[2]{\langle #1,#2 \rangle_{\partial \Omega_e \setminus (\Ga[D] \cup \Ga[S]) } }
\newcommand{\intBD}[2]{\langle #1,#2   \rangle_{\partial \Omega_e \cap \Ga[D]}}
\newcommand{\intBN}[2]{\langle #1,#2   \rangle_{\partial \Omega_e \cap \Ga[N]}}
\newcommand{\intBS}[2]{\langle #1,#2   \rangle_{\partial \Omega_e \cap \Ga[S] }}

\newcommand{\Lpgd}{\bL_{_{\texttt{PGD}}}}
\newcommand{\bupgd}{\bu_{_{\texttt{PGD}}}}
\newcommand{\ppgd}{p_{_{\texttt{PGD}}}}
\newcommand{\bhupgd}{\bhu_{_{\texttt{PGD}}}}
\newcommand{\rpgd}{\rho_{_{\texttt{PGD}}}}
\newcommand{\dLpgd}{\De \bL_{_{\texttt{PGD}}}}
\newcommand{\dbupgd}{\De \bu_{_{\texttt{PGD}}}}
\newcommand{\dppgd}{\De p_{_{\texttt{PGD}}}}
\newcommand{\dbhupgd}{\De \bhu_{_{\texttt{PGD}}}}
\newcommand{\drpgd}{\De \rho_{_{\texttt{PGD}}}}
\newcommand{\buref}{\bu_{_{\texttt{REF}}}}

\newcommand{\vecU}[1][1]{\ensuremath{\vect{U}_{\! #1}}}

\newcommand{\Gpgd}{\mat{G}_{_{\texttt{PGD}}}}
\newcommand{\fUtens}{\vect{f}_{\node{U}}}
\newcommand{\fPsi}{\bm{\psi}}
\newcommand{\ampUtens}{\sigma_{\! \node{U}}}
\newcommand{\fUtensT}{\vect{\tilde{f}}_{\node{U}}}
\newcommand{\fPsiT}{\bm{\widetilde{\psi}}}

\newcommand{\niter}{\texttt{n}_{\texttt{i}}}
\newcommand{\nsnap}{\texttt{n}_{\texttt{s}}}

\newcommand{\fL}{\bm{F}_{\!\! L}}
\newcommand{\fU}{\bm{f}_{\!\! u}}
\newcommand{\fP}{f_{\! p}}
\newcommand{\fHU}{\bm{f}_{\!\! \hu}}
\newcommand{\fR}{f_{\! \rho}}

\newcommand{\ampL}{\sigma_{\! L}}
\newcommand{\ampU}{\sigma_{\! u}}
\newcommand{\ampP}{\sigma_{\! p}}
\newcommand{\ampHU}{\sigma_{\! \hu}}
\newcommand{\ampR}{\sigma_{\! \rho}}

\newcommand{\DivX}{{\bm{\nabla}\cdot\,}}
\newcommand{\GradX}{\bm{\nabla}}
\newcommand{\nmap}{\texttt{n}_{\texttt{M}}}

\newcommand{\bJ}{\mat{J}}
\newcommand{\Jk}{\bJ^k}

\newcommand{\ndet}{\texttt{n}_{\texttt{d}}}
\newcommand{\nadj}{\texttt{n}_{\texttt{a}}}

\newcommand{\bA}{\mat{A}}
\newcommand{\mapping}{\bm{\mathcal{M}}_{\bmu}}
\newcommand{\Jaco}{\mat{J}_{\!\bmu}}

\newcommand{\nDir}{\texttt{n}_{\texttt{D}}}
\newcommand{\nNeu}{\texttt{n}_{\texttt{N}}}
\newcommand{\nSou}{\texttt{n}_{\texttt{S}}}

\newcommand{\bg}{\bm{g}}
\newcommand{\bgD}{\bg_{D}}
\newcommand{\bgN}{\bg_{N}}
\newcommand{\bgS}{\bg_{S}}

\newcommand{\bI}{\bm{\mathcal{I}}}
\newcommand{\I}{\mathcal{I}}

\newcommand{\bM}{\mat{M}}

\newcommand{\nDeg}{\ensuremath{k}}

\newcommand{\Rout}{R_{\text{out}}}
\newcommand{\Rref}{R_{\text{ref}}}
\newcommand{\Rint}{R_{\text{int}}}

\newcommand{\Dmax}{D_{\text{max}}}
\newcommand{\Dmin}{D_{\text{min}}}

\newcommand{\FD}{\texttt{F}_{\! \texttt{D}}}
\newcommand{\FDpgd}{\texttt{F}_{\! \texttt{D}_{\texttt{PGD}}}}
\newcommand{\FDref}{\texttt{F}_{\! \texttt{D}_{\texttt{REF}}}}

\newcommand{\De}{\varDelta}

\usepackage{xcolor}

\renewcommand{\emph}[1]{\textit{#1}}
\newenvironment{keywords}{\begin{quote}\emph{\textbf{Keywords:}}}{\end{quote}}

\theoremstyle{definition}

\newtheorem{remark}{Remark}

\RequirePackage{algorithm}[0.1] %Following SIAM cls
\usepackage{algorithmic} %Algorithm style, could alternatively use algpseudocode 

\usepackage{scalerel}
\usepackage{stackengine}
\stackMath
\def\hatgap{0pt}
\def\subdown{-2pt}
\newcommand\reallywidehat[2][]{
	\renewcommand\stackalignment{l}
	\stackon[\hatgap]{#2}{
		\stretchto{\scalerel*[\widthof{$#2$}]{\kern-.6pt\bigwedge\kern-.6pt}{\rule[-\textheight/2]{1ex}{\textheight}}}{0.5ex}_{\smash{\belowbaseline[\subdown]{\scriptstyle#1}}}}}
		
\makeatletter
\newcommand*\wt[2][0.2ex]{%
        \begingroup
        \mathchoice{\wt@helper{#1}{#2}{\displaystyle}{\textfont}}
                   {\wt@helper{#1}{#2}{\textstyle}{\textfont}}
                   {\wt@helper{#1}{#2}{\scriptstyle}{\scriptfont}}
                   {\wt@helper{#1}{#2}{\scriptscriptstyle}{\scriptscriptfont}}%
        \endgroup
        #2%
}
\newcommand*\wt@helper[4]{%
        \def\currentfont{\the#41}%
        \def\currentskewchar{\char\the\skewchar\currentfont}%
        \setbox\tw@\hbox{\currentfont#2\currentskewchar}%
        \dimen@ii\wd\tw@
        \setbox\tw@\hbox{\currentfont#2{}\currentskewchar}%
        \advance\dimen@ii-\wd\tw@
        \rlap{\raisebox{-#1}{$\m@th#3\kern\dimen@ii\widetilde{\phantom{#2}}$}}%
}
\makeatother
\graphicspath{{figsPDF/}{figsPNG/}{figs/}}

\begin{document}
%==========================================================================
\title{Separated response surfaces for flows in parametrised domains: comparison of a priori and a posteriori PGD algorithms}

\author{
\renewcommand{\thefootnote}{\arabic{footnote}}
			  M. Giacomini\footnotemark[1]\textsuperscript{ \ ,}\footnotemark[2]\textsuperscript{ \ ,}* , 
			  L. Borchini\footnotemark[1]\textsuperscript{ \ ,}\footnotemark[3] , 
			  R. Sevilla\footnotemark[3] \ and
             A. Huerta\footnotemark[1]\textsuperscript{ \ ,}\footnotemark[2]
}

\date{\today}
%________________________________________________________________________
\maketitle

\renewcommand{\thefootnote}{\arabic{footnote}}

\footnotetext[1]{Laboratori de C\`alcul Num\`eric (LaC\`aN), ETS de Ingenieros de Caminos, Canales y Puertos, Universitat Polit\`ecnica de Catalunya, Barcelona, Spain}
\footnotetext[2]{International Centre for Numerical Methods in Engineering (CIMNE), Barcelona, Spain.}
\footnotetext[3]{Zienkiewicz Centre for Computational Engineering, Faculty of Science and Engineering, Swansea University, Wales, UK
\vspace{5pt}\\
* Corresponding author: Matteo Giacomini. \textit{E-mail:} \texttt{matteo.giacomini@upc.edu}
}

%________________________________________________________________________
\begin{abstract}
Reduced order models (ROM) are commonly employed to solve parametric problems and to devise inexpensive response surfaces to evaluate quantities of interest in real-time. There are many families of ROMs in the literature and choosing among them is not always a trivial task. This work presents a comparison of the performance of a priori and a posteriori proper generalised decomposition (PGD) algorithms for an incompressible Stokes flow problem in a geometrically parametrised domain. This problem is particularly challenging as the geometric parameters affect both the solution manifold and the computational spatial domain. The difficulty is further increased because multiple geometric parameters are considered and extended ranges of values are analysed for the parameters and this leads to significant variations in the flow features. Using a set of numerical experiments involving geometrically parametrised microswimmers, the two PGD algorithms are extensively compared in terms of their accuracy and their computational cost, expressed as a function of the number of full-order solves required. 
\end{abstract}

%________________________________________________________________________

%%%%%%%%%%%%%%%%%%%%%%%%%%%%%%%%%%%%%%%%%%%%%%
\section{Introduction and literature review}
\label{sc:intro}
%%%%%%%%%%%%%%%%%%%%%%%%%%%%%%%%%%%%%%%%%%%%%%

Since their introduction by Box and Wilson in 1951~\cite{Box-BW-51}, response surfaces have been extensively used in computational engineering to devise the relationship between input variables or parameters and output quantities of interest~\cite{Breitkopf-BNRV-05,Breitkopf-ZBKZ-11}. This is especially interesting in the context of parametrised partial differential equations (PDEs), where the solution depends both on the spatial coordinates and on a set of user-defined parameters. The solutions of such parametric problems are defined as multidimensional manifolds and efficient strategies for their computation have been developed during the last decades using reduced order models (ROMs)~\cite{AH-CHRW:17,Gunzburger-PWG-18}. In this context, surrogate models of parametric response surfaces were devised by expressing the multidimensional quantities of interest in terms of the elements of the reduced basis constructed for the ROM, see e.g.~\cite{Alonso-LA-00,Maute-WEM-09,Breitkopf-XBCKSV-10,Carlberg-CF-11,Rozza-MQR-12,Passieux-GP-13,Zahr-ZF-15}. Interested readers are also referred to the collection~\cite{Breitkopf-Coelho} and references therein.

Response surfaces based on a posteriori ROMs are constructed starting from a series of snapshots obtained as full-order solutions of the problem under analysis, for a given set of values of the parameters. A critical aspect to devise competitive and accurate numerical strategies is the selection of the snapshots, that is, the sampling procedure in the parametric space. Several techniques were proposed to address this problem, starting from the classical Latin hypercube sampling~\cite{Conover-79-Latin} and centroidal Voronoi tessellation~\cite{Gunzburger-DFG-99} to greedy approaches based on a posteriori error estimates~\cite{Patera-GP-05,Patera-VP-05} and model-constrained adaptive sampling~\cite{Willcox-BWG-08}. In this context, special attention was also devoted to hyper-reduction techniques~\cite{Ryckelynck-05,Sorensen-CS-10}, required to achieve competitive performance in the evaluation of nonlinear quantities of interest~\cite{Carlberg-CBF-11,Amsallem-AZF-12,Carlberg-CFCA-13}.

Once the sampling points are selected, the computation of the snapshots is generally performed in parallel, exploiting the independence of each set of parameters to one another, to reduce the computational cost of the offline phase. Recently, an alternative strategy aiming to reduce the number of required full-order solutions was proposed via an incremental algorithm~\cite{Breitkopf-PBBVZ-20}. The idea is to compute snapshots sequentially and on-the-fly, corresponding to the values of the parameters identified by an appropriate error estimate. In a similar fashion, a priori model reduction strategies construct an approximation by means of a greedy algorithm which sequentially computes the terms of the reduced solution~\cite{Chinesta-Keunings-Leygue}. Although such procedure cannot be performed in parallel, a priori ROMs automatically determine the number of terms in the reduced basis and do not require prior knowledge of the solution, circumventing the sampling step.

Starting from the reduced solution obtained using either a priori, e.g. the proper generalised decomposition (PGD)~\cite{Chinesta-Keunings-Leygue,AH-CLBACGAAH-13}, or a posteriori, e.g. reduced basis (RB)~\cite{Maday-BMNP:04,Patera-GP-05,Patera-VP-05} or proper orthogonal decomposition (POD)~\cite{Iollo-ILD:00,Volkwein-KV:02}, approaches, parametric response surfaces can be efficiently devised. In this context, a critical aspect is represented by the interpolation strategy used to evaluate the quantities of interest depending on the solution manifold constructed using the ROMs. The difficulty of effectively interpolating the reduced solution in a multidimensional manifold was first addressed in~\cite{Amsallem-AF-08}. Since then, different strategies were proposed to reduce the dimensionality of the input space, e.g. via kernel principal component analysis~\cite{Chinesta-GACAC-18} and manifold learning~\cite{Arroyo-MA-13,Breitkopf-LRB-15,Breitkopf-MBLRV-18}, and to perform accurate interpolation using manifold walking~\cite{Breitkopf-RXBRV-13} and co-kriging \cite{Breitkopf-XZBVZ-18} techniques. Recently, manifold learning techniques and collocation methods inspired by sparse grids~\cite{Bungartz-Griebel} have been coupled with PGD-based separated representations of functions of interest~\cite{Ibanez-IBAACLC-17,Ibanez-IAAGCC-18}. The resulting methodologies, including sparse subspace learning~\cite{Borzacchiello-BAC-19} and sparse PGD~\cite{Ibanez-IAAGCHDC-18}, allow to concurrently devise low-dimensional descriptions of the parameter space and functional approximation of the solution manifold, leading to the so-called \emph{hybrid twins} paradigm~\cite{Chinesta-CCADE-20}.

Although both a priori and a posteriori ROMs have been utilised to solve parametrised PDEs and to devise parametric response surfaces, it is not possible to know a priori which reduction approach will perform better for a given problem. Indeed, to the best of the authors' knowledge, no comparison of these approaches in terms of their accuracy and their computational cost is available in the literature. The present work aims to provide a comparison of a priori and a posteriori model reduction techniques, with special emphasis on their cost in terms of number of calls to the full-order solver. It is worth noticing that the current development of a priori and a posteriori ROMs shows different levels of maturity. Indeed, a posteriori approaches feature an extended literature tackling various aspects critical for their efficiency, such as sampling strategies and error control. On the contrary, in the context of a priori ROMs, methodologies for the advanced treatment of the space of parameters~\cite{Modesto-MYZH-20} or the embedded control of accuracy~\cite{Chamoin-KCLP-19,Moitinho-RMDZ-20,Smetana-SZ-20} represent an active line of investigation. Hence, in order for the present comparison to be unbiased, similar versions of the a priori and a posteriori ROM algorithms, without targeted sampling or error control streatgies, are considered.

For the purpose of this comparison, the PGD framework, successfully applied in recent years to several problems~\cite{Aguado-AHCC-15,DM-MZH:15,diez2017generalized,Sibileau-SGAMD-18,Barroso-BGLMH-20}, is considered to construct both the a priori and the a posteriori ROM strategy. More precisely, this study focuses on PGD for geometrically parametrised PDEs. Previous works, see e.g.~\cite{AH-AHCCL:14,SZ-ZDMH:15,sevilla2020solution,RS-SBGH-20}, have shown that this class of parametric problems is particularly challenging since the parameters under analysis affect both the solution manifold and the computational spatial domain. In particular, the difficulty of such problems increases when more than one parameter is considered and when the parameters are responsible for extreme changes of the geometry. Hence, PDEs on geometrically parametrised domains offer complex benchmark cases for ROM strategies, even when a low number of parameters is considered.

The remainder of the manuscript is organised as follows. Section~\ref{sc:StokesPb} presents the incompressible Stokes flow problem in a geometrically parametrised domain and the full-order hybridisable discontinuous Galerkin (HDG) solver. The PGD framework is introduced in section~\ref{sc:PGD}, where the a priori and a posteriori algorithms are detailed and critically compared. An extensive set of numerical simulations for geometrically parametrised Stokes flows in the context of modelling of microswimmers is presented in section~\ref{sc:examples}. More precisely, a comparison of a priori and a posteriori PGD algorithms is performed in terms of accuracy and computational cost by means of parametric problems involving one or more geometric parameters, with different ranges of values, to study the sensitivity of the PGD-based methodologies to the variability of the multidimensional solution. Finally, section~\ref{sc:conclusion} summarises the presented results and the open lines of investigation and appendix~\ref{app:separatedPGD} presents some technical details on the separated form of the operators involved in the a priori PGD algorithm.

%%%%%%%%%%%%%%%%%%%%%%%%%%%%%%%%%%%%%%%%%%%%%%
\section{Full-order Stokes solver in geometrically parametrised domains}
\label{sc:StokesPb}
%%%%%%%%%%%%%%%%%%%%%%%%%%%%%%%%%%%%%%%%%%%%%%

The model problem for the present study is an incompressible Stokes flow in a domain with parametrised geometry. More precisely, the parametrised domain is denoted by $\Omega^{\bmu} \subset \mathbb{R}^{\nsd}$, $\nsd$ being the number of spatial dimensions, with boundary $\partial\Omega^{\bmu} {=} \Ga[D]^{\bmu} \cup \Ga[N]^{\bmu} \cup \Ga[S]^{\bmu}$, where the portions $\Ga[D]^{\bmu}$, $\Ga[N]^{\bmu}$ and $\Ga[S]^{\bmu}$ are disjoint by pairs. The $\npar$ parameters $\bmu\in\bI\subset\mathbb{R}^{\npar}$ control the geometric representation of the domain and are assumed to be independent, that is, $\bI := \I^1\times\I^2\times\dotsb\times\I^{\npar}$ with $\mu_j\in\I^j$ for $j=1,\ldots , \npar$.

By considering the parameters $\bmu$ as extra coordinates of a multidimensional problem in the higher-dimensional manifold $\Omega^{\bmu} \times \bI$, the solution of the Stokes problem is the parametric velocity-pressure pair, $\bu(\bx,\bmu)$ and $p(\bx,\bmu)$, such that
\begin{equation} \label{eq:stokesStrongMulti}
\left\{\begin{aligned}
-\Div(\nu \Grad \bu - p \Insd) &= \bs  															&&\text{in $\Omega^{\bmu} \times \bI$,}\\
\Div\bu &= 0  																					&&\text{in $\Omega^{\bmu} \times \bI$,}\\
\bu &= \bu_D  																					&&\text{on $\Ga[D]^{\bmu} \times \bI$,}\\
\bn^{\bmu} \cdot \bigl(\nu \Grad \bu -p\Insd \bigr) &= \bt        								&&\text{on $\Ga[N]^{\bmu} \times \bI$,}\\
\bu \cdot \bD^{\bmu} + \bn^{\bmu} \cdot \bigl(\nu \Grad \bu -p\Insd \bigr)\bE^{\bmu} &= \bm{0}  &&\text{on $\Ga[S]^{\bmu} \times \bI$,}
\end{aligned}\right.
\end{equation}
where $\nu {>} 0$ is the kinematic viscosity of the fluid, $\bs$ denotes the applied body forces and $\bn^{\bmu}$ is the outward unit normal vector to the boundary. As mentioned above, the external boundary $\partial\Omega^{\bmu}$ is partitioned in Dirichlet, $\Ga[D]^{\bmu}$, Neumann, $\Ga[N]^{\bmu}$, and slip, $\Ga[S]^{\bmu}$, boundaries which also depend on the parameters $\bmu$. Thus, the velocity $\bu_D$ and the pseudo-traction $\bt$ are imposed on $\Ga[D]^{\bmu}$ and $\Ga[N]^{\bmu}$, respectively, for each value in the parametric space $\bI$. Similarly, on $\Ga[S]^{\bmu}$, perfectly slip (i.e. symmetry) conditions are enforced by introducing the matrices $\bD^{\bmu} {:=} [\bn^{\bmu},\bm{0}_{\nsd \times (\nsd-1)}]$ and $\bE^{\bmu} {:=} [ \bm{0},\bm{t}_1^{\bmu},...,\bm{t}_{\nsd-1}^{\bmu}]$, where the tangential directions $\bm{t}_k^{\bmu}, \ k {=} 1,\ldots,\nsd-1$ form an orthonormal system of vectors $\{\bn^{\bmu},\bm{t}_1^{\bmu},\ldots,\bm{t}_{\nsd-1}^{\bmu}\}$, see~\cite{MG-GSH-20}.

Moreover, from the divergence-free equation in problem~\eqref{eq:stokesStrongMulti}, the compatibility condition 
\begin{equation} \label{eq:incompressibleConstMultidim}
%\langle \bu_D \cdot \bn^{\bmu} , 1 \rangle_{\Ga[D]^{\bmu} \times \bI} + \langle \bu \cdot \bn^{\bmu} , 1 \rangle_{ ( \partial \Omega^{\bmu} \setminus \Ga[D]^{\bmu} ) \times \bI}  = 0
\langle \bu_D \cdot \bn^{\bmu} , 1 \rangle_{\Ga[D]^{\bmu}} + \langle \bu \cdot \bn^{\bmu} , 1 \rangle_{\partial \Omega^{\bmu} \setminus \Ga[D]^{\bmu} }  = 0 \quad \text{for almost every $\bmu \in \bI$} ,
\end{equation}
is imposed, where $\langle \cdot,\cdot \rangle_S$ denotes the $\eltwo$ inner product defined on any surface $S \subset \partial \Omega^{\bmu}$. 

Finally, it is worth recalling that in case no Neumann boundary conditions are considered, that is $\partial\Omega^{\bmu} {=} \Ga[D]^{\bmu} \cup \Ga[S]^{\bmu}$, an additional constraint is required to guarantee the uniqueness of the computed pressure field. A common constraint, see~\cite{Jay-CG:09,Cockburn-CS:14,Nguyen-NPC:10}, enforces a zero mean value of the pressure on the domain, that is, for almost every $\bmu \in \bI$,
\begin{equation}
%\Big{\langle}  \frac{1}{|\partial \Omega^{\bmu}|} p,1 \Big{\rangle}_{\partial \Omega^{\bmu}  \times \bI} = 0.
\left(  \frac{1}{|\Omega^{\bmu}|} p,1 \right)_{\Omega^{\bmu}} = 0 ,
\end{equation}	
$\left(\cdot,\cdot\right)_D$ being the $\eltwo$ inner product in a generic subdomain $D \subset \Omega$.

\subsection{Multidimensional hybridisable discontinuous Galerkin solver}
\label{sc:multiHDG}
%-------------------------------------------------------

The multidimensional Stokes flow problem~\eqref{eq:stokesStrongMulti} is discretised using the full-order HDG solver described in~\cite{RS-SBGH-20}. The choice of the HDG framework~\cite{Jay-CG:09,Jay-CGL:09} allows to devise an LBB-compliant discretisation of the Stokes equations with high-order isoparametric formulations using equal order polynomial approximations for all the variables, see e.g.~\cite{Jay-CGNPS:11,Cockburn-CS:14,HDG-NEFEM,giacomini2018superconvergent,MG-GSH-20}. In addition, the PGD-ROM based on the HDG formulation provides an exact separation of the integrals appearing in the geometrically parametrised PDE~\cite{RS-SBGH-20} and does not rely on numerical separation techniques as discussed in~\cite{sevilla2020solution}.

In this section, the multidimensional HDG formulation of equation~\eqref{eq:stokesStrongMulti} is briefly recalled, whereas for a complete derivation interested readers are referred to~\cite{RS-SBGH-20}. First, the spatial, $\Omega^{\bmu}$, and parametric, $\I^j , j {=} 1,\ldots,\npar$, domains are subdivided in $\numel$ and $\numel^j$ disjoint subdomains, respectively, that is,
\begin{align*}
\Omega^{\bmu} & = \bigcup_{e=1}^{\numel} \Omega_e^{\bmu} \ , \ \text{such that} \ \Omega_i \cap \Omega_l = \emptyset \ \text{for $i \neq l$} ,\\
\I^j & = \bigcup_{e=1}^{\numel^j} \I_e^j \ , \ \text{such that} \ \I_i^j \cap \I_l^j = \emptyset \ \text{for $i \neq l$} .
\end{align*}
Moreover, the mesh skeleton $\Gamma^{\bmu}$ of the spatial domain is defined as
\begin{equation*}
\Gamma^{\bmu} := \left[ \bigcup_{e=1}^{\numel} \partial\Omega_e^{\bmu} \right] \setminus \partial\Omega^{\bmu} .
\end{equation*}

The full-order HDG solver for geometrically parametrised Stokes flows is devised according to the formulation introduced in~\cite{RS-SBGH-20}. More precisely, a reference domain $\Omega$, independent of the parameters $\bmu$, is introduced and a mapping 
\begin{equation} \label{eq:mappingDomain}
\begin{aligned}
\mapping\, :
\Omega\times\bI & \longrightarrow \Omega^{\bmu} \\
(\bX,\bmu) 		& \longmapsto \bx = \mapping(\bX,\bmu).
\end{aligned}
\end{equation}
is considered to transform it into the geometrically parametrised domain $\Omega^{\bmu}$, see~\cite{AH-AHCCL:14,SZ-ZDMH:15,sevilla2020solution,RS-SBGH-20}.

The HDG formulation of the Stokes equations is thus written on the reference domain $\Omega$ by applying the mapping~\eqref{eq:mappingDomain}. Following~\cite{RS-SBGH-20}, the functional spaces
\begin{align*}
\sVh(\Omega) &:= \{ v \in \eltwo(\Omega) : v \vert_{\Omega_e} \in \Pk(\Omega_e) \ \forall \Omega_e \ , \ e=1,\ldots ,\numel \} , \\
\shVh(S)&:= \{ \hv \in \eltwo(S) : \hv \vert_{\Ga[i]} \in \Pk(\Ga[i]) \ \forall \Ga[i] \subset S \subseteq \partial\Omega \cup \Gamma \}, \\
\sLh(\I^j) &:= \{ v \in \eltwo(\I^j) : v \vert_{\I_e^j} \in \Pk(\I_e^j) \ \forall \I_e^j \ , \ e=1,\ldots ,\numel^j \}, \\
\Lh(\bI) &:= \sLh(\I^1) \otimes \dotsb \otimes \sLh(\I^{\npar}),
\end{align*}
are introduced, where $\Pk(\Omega_e)$, $\Pk(\Ga[i])$ and $\Pk(\I_e^j)$ denote the spaces of polynomial functions of complete degree at most $k$ in $\Omega_e$, on $\Ga[i]$ and in $\I_e^j$, respectively. Moreover, for the sake of readability, the following scalar-valued, $\sVmu$, vector-valued, $\Vmu$ and $\hVmu$, and tensor-valued, $\Wmu$, discrete functional spaces are defined
\begin{align*}
\sVmu  & := \sVh(\Omega) \otimes \Lh(\bI), 
\quad
& \Vmu   & :=  \left[ \sVh(\Omega) \otimes \Lh(\bI) \right]^{\nsd}, \\
\hVmu & :=  \left[ \shVh(\Gamma \cup \Ga[N] \cup \Ga[S]) \otimes \Lh(\bI) \right]^{\nsd}, 
\quad
& \Wmu  & :=  \left[ \sVh(\Omega) \otimes \Lh(\bI) \right]^{\nsd \times \nsd} .
\end{align*}

As usual in the context of HDG formulations of the Stokes equations~\cite{Jay-CG:09,Nguyen-NPC:10,Nguyen-CNP:10,Jay-CGNPS:11,Cockburn-CS:14,HDG-NEFEM,giacomini2018superconvergent,MG-GSH-20}, the mixed variable $\bL {=} {-} \nu \Grad \bu$, i.e. a scaling of the gradient of velocity, the hybrid variable $\bhu$ representing the trace of the velocity on the element faces and the mean pressure $\rho$ on the boundary of the element are introduced.

First, a static condensation of the degrees of freedom inside each element is performed via the HDG local problems: for $e {=} 1,\ldots ,\numel$, velocity, $\bu_e$, pressure, $p_e$, and mixed variable, $\bL_e$, are written element-by-element in terms of the hybrid velocity $\bhu$ and the mean pressure $\rho_e$. The weak form of the HDG local problems on the spatial reference element $\Omega_e$ is: given $\bu_D$ on $\Ga[D]$ and $\bhu^h$ on $\Gamma \cup \Ga[N] \cup \Ga[S]$, find $(\bu_e^h, p_e^h, \bL_e^h) \in \Vmu \times \sVmu \times \Wmu$ such that
\begin{equation} \label{eq:localWeak}
	\begin{aligned}	
	A_{LL}(\bW,\bL_e^h;\bmu) + A_{Lu}(\bW,\bu_e^h;\bmu)  & = L_L(\bW;\bmu) + A_{L\hu}(\bW,\bhu^h;\bmu),  \\
	A_{uL}(\bv,\bL_e^h;\bmu) + A_{uu}(\bv,\bu_e^h;\bmu) &+ A_{up}(\bv,p_e^h;\bmu) \\[-1ex] 
	& = L_u(\bv;\bmu) + A_{u\hu}(\bv,\bhu^h;\bmu),  \\
	A_{pu}(v,\bu_e^h;\bmu) & = L_p(v;\bmu) + A_{p\hu}(v,\bhu^h;\bmu),  \\
	A_{\rho p}(1,p_e^h;\bmu) & = A_{\rho \rho}(1,\rho_e^h;\bmu),
	\end{aligned}
\end{equation}
for all $(\bv,v,\bW) \in \Vmu \times \sVmu \times \Wmu$. The multidimensional bilinear and linear forms appearing in equation~\eqref{eq:localWeak} are obtained by applying the mapping~\eqref{eq:mappingDomain} to the integrals defined on the geometrically parametrised elements $\Omega_e^{\bmu}$, leading to
\begin{subequations} \label{eq:LocalRef}
\begin{equation} \label{eq:bilinearLocalRef}
	\begin{aligned}	
	A_{LL}(\bW,\bL;\bmu)  	& := -\intExI{\bW}{\nu^{-1} \Det{(\Jaco) \bL } }, 	\\
	A_{Lu}(\bW,\bu;\bmu)  	& := \intExI{ \adj{(\Jaco)}\DivX \bW}{\bu}, 			\\
	A_{L\hu}(\bW,\bhu;\bmu) 	& := \intBNoDxI{\adj{(\Jaco)} \bn \cdot \bW}{\bhu}, 	\\
	A_{uL}(\bv,\bL;\bmu)  	& := \intExI{\bv}{\adj{(\Jaco)}\DivX \bL},			\\
	A_{uu}(\bv,\bu;\bmu)  	& := \intBExI{\bv}{\btau \bu},						\\
	A_{up}(\bv,p;\bmu)    	& := \intExI{\bv}{\adj{(\Jaco)} \GradX p},			\\
	A_{u\hu}(\bv,\bhu;\bmu) 	& := \intBNoDxI{\bv}{\btau \bhu}, 					\\
	A_{pu}(v,\bu;\bmu)    	& := \intExI{\adj{(\Jaco)} \GradX v}{\bu},			\\
	A_{p\hu}(v,\bhu;\bmu) 	& := \intBNoDxI{v}{\bhu \cdot \adj{(\Jaco)} \bn},	\\
	A_{\rho p}(w,p;\bmu)  	& := \intExI{w}{|\Omega_e|^{-1}p},			\\
	A_{\rho \rho}(w,\rho;\bmu)  	& := \intI{w}{\rho},						
	\end{aligned}	
\end{equation}
where $\left(\cdot,\cdot\right)_{D \times \bI}$ stands for the $\eltwo$ inner product in a generic subdomain $D \times \bI$ with $D \subset \Omega$ and $\langle \cdot,\cdot \rangle_{S \times \bI}$ denotes the $\eltwo$ inner product in any domain $S \times \bI$, with $S \subset \partial \Omega^{\bmu}$. Moreover, $\Jaco {=} \Jaco(\bX,\bmu)$ and $\Det{(\Jaco)}$ represent the Jacobian of the mapping and its determinant, respectively, whereas its adjoint is defined as $\adj{(\Jaco)} {=} \Det{(\Jaco)} \Jaco^{-1}$.
The corresponding linear forms are given by
\begin{equation} \label{eq:linearLocalRef}
	\begin{aligned}	
	L_L(\bW;\bmu) 	& := \intBDxI{\adj{(\Jaco)} \bn \cdot \bW}{\bu_D}, 									\\
	L_u(\bv;\bmu) 	& := \intExI{\bv}{\Det{(\Jaco)} \bs} + \intBDxI{\bv}{\btau \bu_D}, 	\\
	L_p(v;\bmu)   	& := \intBDxI{v}{\bu_D \cdot \adj{(\Jaco)} \bn} .
	\end{aligned}	
\end{equation}
\end{subequations}

It is worth noticing that the bilinear and linear forms introduced in equation~\eqref{eq:LocalRef} depend both on spatial, $\bX$, and parametric, $\bmu$, variables. On the one hand, the integrals obtained from the application of the mapping~\eqref{eq:mappingDomain} are defined on the spatial reference domain $\Omega$, which is independent of the parameters $\bmu$. On the other hand, the Jacobian of the mapping being a function of space and parameters, it follows that the terms $\Det{(\Jaco)}$ and $\adj{(\Jaco)}$ depend both on $\bX$ and $\bmu$, as further detailed in the following section.

\begin{remark}[Stabilisation in hybridisable discontinuous Galerkin methods]
The HDG stabilisation tensor $\btau$ is known to play an important role in the accuracy and convergence properties of the numerical approximation~\cite{Jay-CGL:09,Nguyen-NPC:09,Nguyen-NPC:09b,Nguyen-NPC:10}. In the context of Stokes flow problems, see~\cite{MG-GSH-20}, an isotropic stabilisation tensor is considered, namely $\btau := (\tau \nu/\ell) \Insd$ , where $\ell$ is a characteristic length of the domain and $\tau$ a positive scaling factor selected equal to $10$ in the present work.
\end{remark}

Second, the globally-coupled degrees of freedom, namely the hybrid velocity and the mean pressure, are computed by solving the HDG global problem whose weak form reads: find $\bhu^h \in \hVmu$ and $\rho^h \in \mathbb{R}^{\numel} \otimes \Lh(\bI)$ such that
\begin{equation} \label{eq:globalWeak}
	\begin{aligned}	
	\sum_{e=1}^{\numel} \left\{ A_{\hu L}(\bhv,\bL_e^h;\bmu) + A_{\hu u}(\bhv,\bu_e^h;\bmu) + A_{\hu p}(\bhv,p_e^h;\bmu)  \right. \hspace{5pt} & \\[-1ex]
	\left. + A_{\hu \hu}(\bhv,\bhu^h;\bmu) \right\} & = \sum_{e=1}^{\numel} \left\{ L_{\hu}(\bhv;\bmu) \right\},  \\
	A_{p \hu}(1,\bhu^h;\bmu)  & = -L_p(1;\bmu),
	\end{aligned}
\end{equation}
for all $\bhv \in \hVmu$, with the multidimensional bilinear and linear forms
\begin{equation} \label{eq:bilinearGlobalRef}
	\begin{aligned}
	A_{\hu L}(\bhv,\bL;\bmu) 		 := & \intBNoDSxI{\bhv}{\adj{(\Jaco)} \bn \cdot \bL} 	\\
	 &- \intBSxI{\bhv}{\adj{(\Jaco)} \bn \cdot \bL \bE} ,												\\
	A_{\hu u}(\bhv,\bu;\bmu)		 := & \intBNoDSxI{\bhv}{\btau \bu} - \intBSxI{\bhv}{(\btau \bu) \!\cdot\! \bE} ,	\\
	A_{\hu p}(\bhv,p;\bmu)			 := & \intBNoDSxI{\bhv}{p \adj{(\Jaco)} \bn} , \\
	A_{\hu \hu}(\bhv,\bhu;\bmu) := & - \intBNoDSxI{\bhv}{\btau \bhu} 											\\
	 &+ \intBSxI{\bhv}{\bhu \!\cdot\! \adj{(\Jaco)} \bD + (\btau \bhu) \!\cdot\! \bE} , \\
	  L_{\hu}(\bhv;\bmu)	 := & - \intBNxI{\bhv}{\bt} .
	\end{aligned}	
\end{equation}

Similarly to what observed for the local problem, the bilinear and linear forms of the global problem presented in equation~\eqref{eq:bilinearGlobalRef} are defined on the reference domain $\Omega$, independent of $\bmu$, whereas the adjoint of the Jacobian $\adj{(\Jaco)}$ incorporates the dependence on the space and on the parameters.

Interested readers are referred to~\cite{RS-SBGH-20} for the details of the derivation of the local and global problems of the multidimensional HDG method for a Stokes flow in a geometrically parametrised domain. For the sake of readability and except in case of ambiguity, the subscript $_e$ and the superscript $^h$ will be henceforth omitted.

%%%%%%%%%%%%%%%%%%%%%%%%%%%%%%%%%%%%%%%%%%%%%%
\section{Separated response surfaces based on the proper generalised decomposition}
\label{sc:PGD}
%%%%%%%%%%%%%%%%%%%%%%%%%%%%%%%%%%%%%%%%%%%%%%

In order to construct response surfaces for the real-time evaluation of quantities of interest, the parametric problem~\eqref{eq:stokesStrongMulti} needs to be efficiently computed for a large number of configurations. Nonetheless, the solution of the geometrically parametrised Stokes flow problem in a space of dimension $\nsd {+} \npar$ using the previously introduced multidimensional HDG solver is computationally unaffordable, even when only few parameters are considered. 

In this section, two PGD-based strategies are presented to construct response surfaces in terms of separated functions. Before detailing the two algorithms, the framework to construct a separated approximation using PGD is briefly recalled.

As classical in the context of ROMs for parametric PDEs~\cite{Patera-Rozza:07,Rozza:14}, the bilinear and linear forms introduced in section~\ref{sc:multiHDG} require an affine dependence on the parameters $\bmu$. For the problem involving geometrically parametrised domains introduced in section~\ref{sc:multiHDG}, the bilinear and linear forms are approximated by the sum of products of parameter-dependent functions and spatial operators independent of the parameters, as detailed in~\ref{app:separatedPGD}. As previously mentioned, the integrals in the weak forms~\eqref{eq:localWeak} and~\eqref{eq:globalWeak} of the HDG local and global problems are defined on the parameter-independent reference domain, whereas the determinant and the adjoint of the mapping incorporate the dependence on the parameters.

It is worth noticing that a fixed mesh is generated once for the spatial reference domain. Hence, no mesh quality issues due to the deformation of the mesh of the reference domain arise. Of course, if the transformation~\eqref{eq:mappingDomain} is responsible for extreme deformations of the reference configuration, the Jacobian of the isoparametric mapping may degrade, leading to the need of a larger number of integration points for the computation of the terms in equations~\eqref{eq:LocalRef} and~\eqref{eq:bilinearGlobalRef}. It is worth noticing that this is especially critical in the context of high-order curved meshes. Hence, it is still advisable to perform a preliminary study of mesh quality for the configurations corresponding to the limit of the geometric mapping, as reported for the examples in section~\ref{sc:examples}.

%-------------------------------------------------------
\subsection{Separated representation of the unknown variables}
\label{sc:separatedSol}
%-------------------------------------------------------

According to the PGD framework~\cite{Chinesta-Keunings-Leygue}, the unknown variables are written as rank-$m$ separable approximations as
\begin{equation}\label{eq:PGDapproxStandard}
	\begin{aligned}
	\bupgd^m (\bX,\bmu) &= \ampU^m \fU^m (\bX) \, \psi^m(\bmu)+ \bupgd^{m-1} (\bX, \bmu) , \\
	\ppgd^m  (\bX,\bmu) &= \ampP^m  \fP^m (\bX) \, \psi^m(\bmu) + \ppgd^{m-1}  (\bX, \bmu) , \\
	\Lpgd^m	 (\bX,\bmu) &= \ampL^m  \fL^m (\bX) \, \psi^m(\bmu) + \Lpgd^{m-1}  (\bX, \bmu) , \\
	\bhupgd^m(\bX,\bmu) &= \ampHU^m \fHU^m(\bX) \, \psi^m(\bmu) + \bhupgd^{m-1}(\bX, \bmu) , \\
	\rpgd^m  (\bX,\bmu) &= \ampR^m  \fR^m (\bX) \, \psi^m(\bmu)  + \rpgd^{m-1}  (\bX, \bmu) ,
	\end{aligned}
\end{equation}
where each term in the expansion, referred to as \emph{mode}, is the product of a spatial function and a function depending on the parameters and $\ampL^m$, $\ampU^m$, $\ampP^m $, $\ampHU^m$ and $\ampR^m$ denote the amplitudes of the corresponding modes. For the sake of simplicity, the parametric terms are assumed to be factorisable using one-dimensional functions, that is,
\begin{equation} \label{eq:phiM}
\psi^m(\bmu) = \prod_{j=1}^{\npar} \psi^m_j(\mu_j) .
\end{equation}
\begin{remark}[Factorisation of the parametric space]
In the context of PGD approximations, spatial and parametric modes are computed alternatively. In order to guarantee the computational efficiency of the PGD algorithm, it is critical for the number of degrees of freedom in the parametric problem to be considerably smaller than the number of unknowns in the spatial problem. Assumption~\eqref{eq:phiM} is commonly employed to enforce that the number of dimensions of each parametric subdomain is inferior to the number of spatial dimensions $\nsd$. Hence, the computational cost of the problem in the PGD parametric step is reduced by devising it as a sequence of lower dimensional problems. It is worth noticing that the number of modes computed by the PGD algorithms without separating the parametric modes in terms of one-dimensional functions will, in general, differ from the one computed starting from equation~\eqref{eq:phiM}. Since it is not possible to know a priori which approach will perform better for a given problem, assumption~\eqref{eq:phiM} will be henceforth considered to minimise the number of degrees of freedom involved in each parametric problem.
\end{remark}

\begin{remark}[PGD initial approximation]
No initial information is considered in the construction of the PGD approximation~\eqref{eq:PGDapproxStandard}, that is, $\bupgd^0 {=} \bm{0}$, $\ppgd^0 {=} 0$, $\Lpgd^0 {=} \bm{0} $, $\bhupgd^0 {=} \bm{0}$ and $\rpgd^0 {=} 0$. More precisely, contrary to traditional finite element-based PGD algorithms requiring the construction of initial modes to account for Dirichlet boundary conditions, see e.g.~\cite{Chinesta-Keunings-Leygue}, the HDG method employed as spatial solver in this work relies on the weak imposition of essential boundary conditions. Hence, no prior computation needs to be performed to initialise the approximations in equation~\eqref{eq:PGDapproxStandard}.
\end{remark}

More recently, a \emph{predictor-corrector} approach has been introduced~\cite{tsiolakis2020nonintrusive} by splitting the $m$-th modes to be computed into the predictions $\ampL^m  \fL^m \psi^m$, $\ampU^m \fU^m \psi^m$, $\ampP^m  \fP^m \psi^m$, $\ampHU^m \fHU^m \psi^m$ and $\ampR^m  \fR^m \psi^m$ and the corrections $\ampL^m  \dLpgd^m$, $\ampU^m \dbupgd^m$, $\ampP^m  \dppgd^m$, $\ampHU^m \dbhupgd^m$ and $\ampR^m  \drpgd^m$, namely,
\begin{equation}\label{eq:PGDapprox}
	\begin{aligned}
	\bupgd^m (\bX,\bmu) &= \ampU^m [\fU^m (\bX) \, \psi^m(\bmu) + \dbupgd^m(\bX, \bmu)] + \bupgd^{m-1} (\bX, \bmu) , \\
	\ppgd^m  (\bX,\bmu) &= \ampP^m  [\fP^m (\bX) \, \psi^m(\bmu) + \dppgd^m(\bX, \bmu)] + \ppgd^{m-1}  (\bX, \bmu) , \\
	\Lpgd^m	 (\bX,\bmu) &= \ampL^m  [\fL^m (\bX) \, \psi^m(\bmu)  + \dLpgd^m(\bX, \bmu)] + \Lpgd^{m-1}  (\bX, \bmu) , \\
	\bhupgd^m(\bX,\bmu) &= \ampHU^m [\fHU^m(\bX) \, \psi^m(\bmu) + \dbhupgd^m(\bX, \bmu)] + \bhupgd^{m-1}(\bX, \bmu) , \\
	\rpgd^m  (\bX,\bmu) &= \ampR^m  [\fR^m (\bX) \, \psi^m(\bmu) + \drpgd^m(\bX, \bmu)] + \rpgd^{m-1}  (\bX, \bmu) ,
	\end{aligned}
\end{equation}
with the corrections given by
\begin{equation}\label{eq:PGDcorrections}
	\begin{aligned}
	\dbupgd^m (\bX,\bmu) &:= \De\fU (\bX) \, \psi^m(\bmu)  + \fU^m (\bX) \, \De\psi(\bmu) + \De\fU (\bX) \, \De\psi(\bmu) , \\
	\dppgd^m  (\bX,\bmu) &:= \De\fP (\bX) \, \psi^m(\bmu)  + \fP^m (\bX) \, \De\psi(\bmu) + \De\fP (\bX) \, \De\psi(\bmu) , \\
	\dLpgd^m	 (\bX,\bmu) &:= \De\fL (\bX) \, \psi^m(\bmu)  + \fL^m (\bX) \, \De\psi(\bmu) + \De\fL (\bX) \, \De\psi(\bmu) , \\
	\dbhupgd^m(\bX,\bmu) &:= \De\fHU (\bX) \, \psi^m(\bmu)  + \fHU^m (\bX) \, \De\psi(\bmu) + \De\fHU (\bX) \, \De\psi(\bmu) , \\
	\drpgd^m  (\bX,\bmu) &:= \De\fR (\bX) \, \psi^m(\bmu)  + \fR^m (\bX) \, \De\psi(\bmu) + \De\fR (\bX) \, \De\psi(\bmu) .
	\end{aligned}
\end{equation}
It is worth noticing that, $\De$ being a variation, the last term in equation~\eqref{eq:PGDcorrections} represents a high-order variation and can thus be neglected in the computation.

\begin{remark}[Choice of the parametric function]\label{rmrk:paramFunc}
According to the single-parameter approach described in~\cite{diez2017generalized}, a unique parametric function $\psi^m$ is considered in~\eqref{eq:PGDapproxStandard} and~\eqref{eq:PGDapprox} for all the variables at the $m$-th mode. Alternative strategies for the definition of the parametric function in PGD approximations of incompressible flow problems are also explored in~\cite{diez2017generalized}.
\end{remark}

Following both the approach in~\eqref{eq:PGDapproxStandard} and~\eqref{eq:PGDapprox}, the number of terms in the PGD expansion is not known a priori. Indeed, assuming the modes up to $m {-} 1$ to be known, a greedy procedure is performed to compute the $m$-th mode. For this purpose, a nonlinear iterative algorithm, namely the alternating direction (AD) scheme, is devised to alternately compute the spatial and parametric modes of the PGD approximation. More precisely, first, the $m$-th parametric mode is assumed to be known and the corresponding spatial functions are computed by solving an HDG problem independent of the parameters. Second, the recently computed spatial function is fixed and the corresponding parametric mode is determined solving a linear system of equations. The procedure is thus repeated until a convergence criterion is fulfilled or a maximum number of iterations is achieved. In the following sections, the details of this greedy strategy are presented, highlighting the main differences between the a priori and a posteriori PGD algorithms.

%-------------------------------------------------------
\subsection{A priori proper generalised decomposition}
\label{sc:PGD-priori}
%-------------------------------------------------------

The a priori PGD relies on solving a separated version of the HDG local and global problems introduced in section~\ref{sc:multiHDG}. In this context, no a priori information is required and the computation of the PGD solution does not depend on any previously determined snapshot. In order to construct such an approximation, a separable representation of the geometric mapping and the user-defined data is required.

%-------------------------------------------------------
\subsubsection{Separated representation of geometric mapping and user-defined data}
%-------------------------------------------------------

Following~\cite{sevilla2020solution,RS-SBGH-20}, the mapping~\eqref{eq:mappingDomain} and its Jacobian are assumed to be separable, namely,
\begin{equation} \label{eq:displacementSep}
\begin{aligned}
\mapping(\bX,\bmu) &= \sum_{k=1}^{\nmap} \bM^k(\bX) \phi^k(\bmu) , \\
\Jaco(\bX,\bmu) &= \sum_{k=1}^{\nmap} \Jk(\bX) \phi^k(\bmu) .
\end{aligned}
\end{equation}
In addition, the determinant and the adjoint of the Jacobian can also be expressed in separated form, see~\cite{sevilla2020solution}, as
\begin{equation} \label{eq:DeterminantAdjointSep}
\begin{aligned}
\det(\Jaco)(\bX,\bmu) &= \sum_{k=1}^{\ndet} D^k(\bX) \theta^k(\bmu) , \\
\adj(\Jaco)(\bX,\bmu) &= \sum_{k=1}^{\nadj} \bA^k(\bX) \vartheta^k(\bmu) .
\end{aligned}
\end{equation}
In case an analytical separation of the mapping is not available, the separated forms in equation~\eqref{eq:displacementSep} may be numerically approximated, e.g. by means of a singular value decomposition or a high-order PGD projection. Recently, a more general approach for the construction of separated geometric mappings was proposed in~\cite{sevilla2020solution}, starting from a representation of the boundary parametrised through the control points of non-uniform rational B-splines and solving a linear elastic problem inspired by high-order curved mesh generation techniques~\cite{poya2016unified,HO-Meshing}. The scope of the present work being the comparison of a priori and a posteriori PGD algorithms, the mapping~\eqref{eq:displacementSep} is henceforth assumed to be analytically separable for the sake of simplicity.

Similarly, user-defined data like body forces, Dirichlet and Neumann boundary terms are assumed to be provided in separated form either analytically, that is,
\begin{equation} \label{eq:separatedData}
	\begin{aligned}	
	\bu_D& \!=\! \sum_{l=1}^{\nDir} \bgD^l(\bX) \lambda^l_D(\bmu),
	\\ 
	\bt	 & \!=\! \sum_{l=1}^{\nNeu} \bgN^l(\bX) \lambda^l_N(\bmu),
	\\ 
	\bs	 & \!=\! \sum_{l=1}^{\nSou} \bgS^l(\bX) \lambda^l_S(\bmu) ,
	\end{aligned}	
\end{equation}
or approximated using appropriate numerical separation techniques~\cite{Chinesta-Keunings-Leygue}.

%-------------------------------------------------------
\subsubsection{The a priori PGD algorithm}
%-------------------------------------------------------

To devise the a priori PGD algorithm for the geometrically parametrised Stokes flow, first the separated functions~\eqref{eq:PGDapprox} and the separated expressions~\eqref{eq:DeterminantAdjointSep}-\eqref{eq:separatedData} of mapping and data are introduced into the multidimensional HDG problems~\eqref{eq:localWeak} and~\eqref{eq:globalWeak}. The resulting equations are thus alternately projected on the tangent manifold associated with the spatial and parametric coordinates to perform the iterations of the AD algorithm.

The tangent manifolds for the vector-valued, $\Vmu$ and $\hVmu$, and scalar, $\sVmu$ and $\mathbb{R}^{\numel} \otimes \Lh(\bI)$, multidimensional discrete functional spaces, are obtained by selecting the test functions in~\eqref{eq:localWeak} and~\eqref{eq:globalWeak} as
\begin{equation}\label{eq:tangentU}
\begin{aligned}
\bv  &= \de \fU \psi^m + \ampU^m \fU^m  \de \psi ,			&
\bhv &= \de \fHU \psi^m + \ampHU^m \fHU^m  \de \psi ,		\\
v    &= \de \fP \psi^m + \ampP^m \fP^m  \de \psi ,			&
w    &= \de \fR \psi^m + \ampR^m \fR^m \de \psi, 			
\end{aligned}
\end{equation}
where the spatial test functions are such that $\de \fU \in \Vh \!:=  \left[ \sVh(\Omega)  \right]^{\nsd}$, $\de\fHU \in \hVh \!:=  \left[ \shVh(\Gamma \cup \Gamma_N \cup \Gamma_S)  \right]^{\nsd}$, $\de \fP \in \sVh(\Omega) $ and $\de \fR \in \mathbb{R}^{\numel}$, whereas the parametric test function is given by $\de \psi \in \Lh(\bI)$. In a similar fashion, the tangent manifold for the tensor-valued space $\Wmu$ is characterised by the test function
\begin{equation}\label{eq:tangentL}
\bW = \de \fL \psi^m + \ampL^m \fL^m  \de \psi ,
\end{equation}
for $\de \fL \in \Wh \!:=  \left[ \sVh(\Omega)  \right]^{\nsd \times \nsd} $.

Hence, in the AD scheme, the parametric function $\psi^m$ is first fixed in the spatial iteration. It follows that $\de \psi {=} 0$ and the multidimensional problems~\eqref{eq:localWeak} and~\eqref{eq:globalWeak} reduce to a parameter-independent spatial HDG problem. In the parametric step, the spatial functions are assumed to be known and, consequently, $\de \fU {=} \de \fP {=} \de \fL {=} \de \fHU {=} \de \fR {=} 0$. The parametric iteration thus results in $\npar$ linear systems of equations, each associated with a one-dimensional problem, which are solved sequentially. A detailed derivation of the spatial and parametric equations in the AD scheme of the a priori PGD algorithm is presented in~\cite{RS-SBGH-20}.

\begin{algorithm}
\caption{The a priori PGD algorithm}\label{alg:PGDpriori}
\begin{algorithmic}[1]
\REQUIRE{For the greedy enrichment loop, the value $\eta^\star$ of the tolerance. For the AD loop, the number of iterations $\niter$.}
\STATE{Set $m \gets 1$ and initialise the amplitude of the spatial mode $\ampHU^1 \gets 1$.}

%___Loop modes (greedy)
\WHILE{$\ampHU^m / \ampHU^1 > \eta^\star$}
\STATE{Set $q \gets 1$ and initialise the parametric prediction.}
\STATE{Solve the HDG global and local problems to compute the spatial prediction.}
             
%___Loop ADI
\WHILE{$q < \niter$}

\STATE{Solve the parametric linear system to compute the parametric correction.}
\STATE{Update the parametric prediction with the correction.}

\STATE{Solve the HDG global and local problems to compute the spatial correction.}
\STATE{Update the spatial prediction with the correction.}
             
\STATE{Increase the counter of the AD iterations $q \gets q+1$.}

\ENDWHILE

\STATE{Increase the mode counter $m \gets m+1$.}
\ENDWHILE
\end{algorithmic}
\end{algorithm}

The resulting a priori PGD strategy is reported in algorithm~\ref{alg:PGDpriori}. Given a guess for the prediction of the parametric mode, the loop for the PGD enrichment first determines a prediction of the spatial mode by solving the HDG global and local problems (Algorithm~\ref{alg:PGDpriori} - Step 4). Then, the AD scheme computes the parametric (Algorithm~\ref{alg:PGDpriori} - Steps 6-7) and spatial (Algorithm~\ref{alg:PGDpriori} - Steps 8-9) corrections solving a parametric linear system and the HDG global and local problems, respectively. The procedure in the AD scheme is then repeated until the maximum number of iterations $\niter$ is achieved. Finally, the greedy iterations stop when the ratio of the amplitude of the current mode to the one of the first mode is negligible (Algorithm~\ref{alg:PGDpriori} - Step 2).

\begin{remark}[Choice of the stopping criteria]
Alternative stopping criteria may be considered for both the AD scheme and the greedy algorithm. A common approach for the former relies on checking the relative amplitude of the computed correction with respect to the amplitude of the current mode~\cite{tsiolakis2020nonintrusive,RS-SBGH-20}. Concerning the latter, the ratio of the amplitude of the current mode to the cumulative amplitudes of the previously computed modes has also been utilised as stopping criterion, see~\cite{tsiolakis2020nonintrusive}. Alternative strategies to stop the enrichment procedure rely on the computation of measures obtained from the problem under analysis, e.g. the norm of the residual of the governing equations or the goal-oriented estimate for a quantity of interest~\cite{PD-GBCD-17,PD-GDBC-18,Chamoin-KCLP-19}.
\end{remark}

%-------------------------------------------------------
\subsection{A posteriori proper generalised decomposition}
\label{sc:PGD-posteriori}
%-------------------------------------------------------

Contrary to the a priori PGD introduced above, the a posteriori framework relies on constructing a reduced basis starting from a series of snapshots. Each snapshot is defined as a vector 
\begin{equation}\label{eq:snapshot}
\vecU[s]^T := \left[ \vect{\hu}^T , \ \bm{\rho}^T , \ \vect{u}^T , \ \vect{p}^T , \ \vect{L}^T \right]_s \quad , \quad s = 1, \ldots, \nsnap ,
\end{equation}
where $\vect{\hu}$, $\bm{\rho}$, $\vect{u}$, $\vect{p}$ and $\vect{L}$ denote the vectors of nodal values of the unknowns of problem~\eqref{eq:stokesStrongMulti}, computed using the full-order HDG spatial solver for a given set of the parameters. Hence, the size of each snapshot vector is equal to the number of degrees of freedom of the HDG global and local problems. The $\nsnap$ snapshots are thus gathered in a multidimensional tensor structure $\mat{G}$. For the case of a unique parameter, this is given by a tensor of order 2, that is, a matrix 
\begin{equation} \label{eq:tensorG}
\mat{G} = 
\begin{bmatrix}
\vecU[1] ,
&
\vecU[2] ,
&
\ldots ,
&
\vecU[\nsnap]
\end{bmatrix} ,
\end{equation}
where the rows correspond to the degrees of freedom of the HDG spatial discretisation and the columns are associated with the snapshots in the parametric interval. In case more than one parameter is considered, the snapshots structure is constructed as the tensor product of the matrix~\eqref{eq:tensorG} with each extra parametric dimension, leading to a multidimensional tensor of order $\npar {+} 1$, with one dimension for each parameter plus one dimension for the space.

The a posteriori PGD, also known as PGD separation or least-squares PGD~\cite{DM-MZH:15,PD-DZGH-18,PD-DZGH-20}, computes the separated approximation of $\mat{G}$ in the form of product of rank-one approximations~\eqref{eq:PGDapproxStandard} using a greedy approach, that is, given $m {-} 1$ modes,  the $m$-th term in the PGD expansion is computed as
\begin{equation}\label{eq:posterioriMin}
\left(\fUtens^m , \fPsi_1^m , \ldots , \fPsi_{\npar}^m \right) = \argmin \left\| \mat{G} - \Gpgd^{m-1} - \ampUtens^m \fUtensT^m \otimes \fPsiT_1^m \otimes \cdots \otimes \fPsiT_{\npar}^m \right\|_2 , 
\end{equation}
where each vector $\fUtensT^m , \fPsiT_1^m , \ldots , \fPsiT_{\npar}^m$ is sought in a subspace of $\mathbb{R}^d$ of appropriate dimension, namely the sizes of $\fUtensT^m$ and $\fPsiT_j^m, \ j {=} 1,\ldots,\npar$ being the number of degrees of freedom of the HDG spatial solution $\vect{U}$ and of the parametric discretisations in the directions $\I_j , \ j {=} 1,\ldots,\npar$, respectively. The greedy procedure presented above aims to compute, at each step, the best approximation $\ampUtens^m \fUtensT^m \otimes \fPsiT_1^m \otimes \cdots \otimes \fPsiT_{\npar}^m$ to describe the unresolved part $\mat{G} {-} \Gpgd^{m-1}$ of the target tensor $\mat{G}$. It is straightforward to observe that the $m$-th enrichment in the PGD loop tackles the approximation of the remaining residual from iteration $m {-} 1$, thus improving the overall separated approximation of $\mat{G}$. Numerical experiments, see e.g.~\cite{DM-MZH:15,PD-DZGH-20}, have shown that the resulting procedure is responsible for errors decreasing monotonically with the inclusion of new modes in the PGD approximation. A similar behaviour is observed for the amplitude of the modes: this is considered as a relative measure of the relevance of the newly computed mode and is thus employed to devise an appropriate stopping criterion. Of course, this strategy may be improved by equipping the a posteriori PGD algorithm with an appropriate error control step targeting the unknowns of the problem or a given quantity of interest.

From a practical point of view, the nonlinear problem~\eqref{eq:posterioriMin} is solved using an AD scheme. It is worth noticing that in this context, both spatial and parametric iterations are determined as rank-one approximations on a purely algebraic level and they do not require any information on the underlying multidimensional HDG discretisation. Hence, their computation relies on elementary tensorial operations, i.e. products and sums of separated objects~\cite{PD-DZGH-20}, and the resulting cost is proportional to the size of the vectors of spatial and parametric modes. 

The resulting a posteriori PGD strategy is reported in algorithm~\ref{alg:PGDposteriori}. First, a set of $\nsnap$ snapshots is constructed using the full-order HDG spatial solver (Algorithm~\ref{alg:PGDposteriori} - Step 1). Then, in each PGD enrichment iteration, the parametric mode is initialised with a user-defined guess and the AD loop alternately computes the spatial (Algorithm~\ref{alg:PGDposteriori} - Step 6) and parametric (Algorithm~\ref{alg:PGDposteriori} - Step 7) modes solving two rank-one problems at the algebraic level. The above routine is repeated until a convergence criterion (Algorithm~\ref{alg:PGDposteriori} - Step 8) is fulfilled or the maximum number of iterations $\niter$ is achieved. Similarly to the a priori PGD algorithm, the greedy enrichment loop stops when the ratio of the amplitude of the current mode to the one of the first mode is negligible (Algorithm~\ref{alg:PGDposteriori} - Step 3).

\begin{algorithm}
\caption{The a posteriori PGD algorithm}\label{alg:PGDposteriori}
\begin{algorithmic}[1]
\REQUIRE{For the greedy enrichment loop, the value $\eta^\star$ of the tolerance. For the AD loop, the value $\eta_{\sigma}$ of the tolerance on the amplitude variation and the maximum number of iterations $\niter$.}
\STATE{Compute $\nsnap$ snapshots solving the HDG global and local problems.}
\STATE{Set $m \gets 1$ and initialise the amplitude of the spatial mode $\ampHU^1 \gets 1$.}

%___Loop modes (greedy)
\WHILE{$\ampHU^m / \ampHU^1 > \eta^\star$}
\STATE{Set $q \gets 1$ and initialise the parametric mode.}
             
%___Loop ADI
\WHILE{$\varepsilon_{\sigma} > \eta_{\sigma}$ or $q < \niter$}

\STATE{Compute the rank-one spatial mode.}

\STATE{Compute the rank-one parametric mode.}

\STATE{Update the stopping criterion $\varepsilon_{\sigma} = (\ampHU^{m,q} - \ampHU^{m-1})/\ampHU^{m,q}$.}
\STATE{Increase the counter of the AD iterations $q \gets q+1$.}

\ENDWHILE

\STATE{Increase the mode counter $m \gets m+1$.}
\ENDWHILE
\end{algorithmic}
\end{algorithm}

%-------------------------------------------------------
\subsection{Devising separated response surfaces} 
\label{sc:surface}
%-------------------------------------------------------

Once the reduced solution is computed for all the variables using either the a priori or the a posteriori algorithms presented above, parametric response surfaces can be easily devised as a postprocess of the separated PGD solutions. More precisely, separated response surfaces are obtained as explicit functions of the parameters of interest. For the case of the drag force on an object of surface $\mathcal{B}$, the rank-$m$ separated approximation is given by
\begin{equation}\label{eq:dragParam}
\begin{aligned}
\FDpgd^m(\bmu) &= \int_{\mathcal{B}}{\left( -\ppgd^m(\bX,\bmu)\Insd - (\Lpgd^m(\bX,\bmu) + \Lpgd^m(\bX,\bmu)^T) \right) \bn \ d\Gamma} \\
&= \sum_{j=1}^m \texttt{D}^j \psi^j(\bmu) 
\end{aligned}
\end{equation}
where the $\texttt{D}^j$ corresponds to the drag coefficient of the $j$-th spatial mode and is obtained as
\begin{equation}\label{eq:dragSpatialMode}
\texttt{D}^j := \int_{\mathcal{B}}{\left( -\ampP^j \fP^j(\bX)\Insd - \ampL^j (\fL^j(\bX) + \fL^j(\bX)^T) \right) \bn \ d\Gamma} .
\end{equation}

It is worth noticing that the accuracy of the separated response surface of a quantity of interest directly depends upon the precision achieved by the PGD approximation of the variables utilised for its computation (e.g. pressure and gradient of velocity in the case of the drag). In this context, the HDG method used as full-order solver allows to achieve optimal convergence of order $k {+} 1$ for both the pressure and the mixed variable respresenting the gradient of velocity~\cite{giacomini2018superconvergent,RS-SBGH-20}. Thus, it provides additional accuracy in the approximation of the viscous part of the drag with respect to classical primal finite element formulations, in which this is obtained as a postprocess of the computed velocity field. To construct separated approximations assessing the accuracy in a given quantity of interest, interested readers are referred to~\cite{PD-GBCD-17,PD-GDBC-18,Chamoin-KCLP-19}, where PGD algorithms with goal-oriented error control were investigated.

%-------------------------------------------------------
\subsection{Critical comparison of a priori and a posteriori PGD algorithms} 
\label{sc:comparison}
%-------------------------------------------------------

Both the a priori and a posteriori approach introduced above have attractive properties and their performance differs depending upon the problem under analysis and the parameters of interest. As it is not possible to know which of the two methodologies will perform better for a given problem, this section offers a critical comparison of the two approaches, highlighting the main advantages and disadvantages of each method. It is worth noticing that geometric parameters are one of the more challenging problems to consider in the context of parametric PDEs as the changes induced by such parameters not only have an influence on the discretised equations but also on the computational spatial domain.

The main drawback of the a posteriori PGD approach is that the user is required to select of a set of snapshots, corresponding to the simulations of the full-order problem, for a given set of values of the parameters. In order to provide a comparison of the cost of a priori and a posteriori PGD in terms of full-order solves, in this work no problem-specific sampling is considered and the snapshots are computed in correspondence of the nodes of the parametric intervals used by the a priori algorithm. Although more advanced sampling techniques have been proposed, see section~\ref{sc:intro}, the selected points are expected to produce accurate representations of the solution in the parametric domain inheriting the good approximation properties of the utilised Fekete nodal distributions. 

Despite the vast literature on sampling methods, it is not possible to initially know the number of snapshots the a posteriori ROM will require to capture the multidimensional solution accurately. In contrast, the a priori PGD approach requires no previous knowledge of the solution and no snapshots need to be selected by the user. Instead, a set of modes is automatically constructed in the enrichment process and the required number of terms is automatically determined by the greedy algorithm according to a user-defined tolerance. 

An important advantage of the a posteriori approach is that the snapshots can be computed in parallel as they are completely independent of each other. In contrast, the a priori approach computes the modes sequentially within the enrichment process. The computation of each mode involves several calls to the spatial solver in order to obtain the solution of the nonlinear problem by using the AD scheme. 

The main drawback of the a priori PGD approach is that its standard implementation is generally intrusive with respect to the spatial solver, see~\cite{sevilla2020solution,RS-SBGH-20}. This means that access to the code is required to devise the PGD algorithm starting from the spatial solver. Despite some recent advances towards non-intrusive implementations of the a priori PGD~\cite{Ladeveze-CNLB-16,zou2018nonintrusive,tsiolakis2020nonintrusive}, this aspect still represents an important challenge for the application of the methodology in an industrial context, where the use of commercial software is preferred. On the contrary, the a posteriori approach does not require access to the code sources as it simply relies on a set of snapshots, which can be obtained using any computational code.

Concerning the two types of separated approximations introduced in section~\ref{sc:separatedSol}, it is worth mentioning that equation~\eqref{eq:PGDapproxStandard} and~\eqref{eq:PGDapprox} are equivalent. In the latter, the computation of each new mode is split into a prediction and a correction step. This approach is fostered for the a priori PGD as it allows to refine the stopping criterion of the AD scheme (Algorithm~\ref{alg:PGDpriori} - Step 5) by introducing an additional check to end the iteration loop when the amplitude of the correction is negligible with respect to the amplitude of the current mode, see~\cite{RS-SBGH-20}. This test has been omitted in the present work to perform a more transparent comparison with the a posteriori PGD algorithm in which no information is provided a priori to reduce the number of computed snapshots.

The points previously discussed are general for any parametric problem, but a crucial aspect specific to geometrically parametrised problems concerns the mesh generation process. More precisely,  for each parametric configuration of the geometry, a different mesh is required. It is worth noticing that these meshes need to have the same number of nodes and the same connectivity matrix. A common approach for a posteriori ROMs is thus to generate one mesh and morph it to obtain the mesh corresponding to each geometric configuration of interest. In this context, special attention needs to be paid to the morphing algorithm, especially in a high-order framework, as this deformation can significantly decrease the quality of the resulting meshes. An alternative approach, used in the context of a priori PGD in~\cite{AH-AHCCL:14,SZ-ZDMH:15,sevilla2020solution,RS-SBGH-20}, relies on defining a reference configuration and an appropriate geometric mapping. Henceforth, only one mesh is required to compute the snapshots for the a posteriori PGD and to perform the computation of the a priori PGD solution, for any geometric configuration of interest.

%%%%%%%%%%%%%%%%%%%%%%%%%%%%%%%%%%%%%%%%%%%%%%
\section{Numerical experiments} 
\label{sc:examples}
%%%%%%%%%%%%%%%%%%%%%%%%%%%%%%%%%%%%%%%%%%%%%%

This section presents a set of numerical examples to investigate the performance of a priori and a posteriori PGD approaches in the context of Stokes flows in geometrically parametrised domains. The problem of interest is the Stokes flow around the so-called \emph{push-me-pull-you} microswimmer, a geometry extensively studied in the context of microfluidics devices~\cite{avron2005pushmepullyou,Alouges-ADL-09}. The swimmer consists of two bladders of spherical shape that can change their volume and mutual distance, with the constraint that the total volume of the two bladders is kept constant. Several numerical tests will be presented, involving one and two parameters. First, parametric studies involving a unique parameter will be performed, with special attention on the influence of the range of parameters starting from the configurations of interest described in the literature~\cite{avron2005pushmepullyou,Alouges-ADL-09}. Then, the concurrent treatment of two parameters will be analysed. In this context, special attention will be devoted to analysing the possibility of extending the previously obtained results to multidimensional cases.

%-------------------------------------------------------
\subsection{Description of the geometry and parametrised mappings} 
%-------------------------------------------------------

\begin{figure}[!tb]
	\centering
	\includegraphics[width=\textwidth]{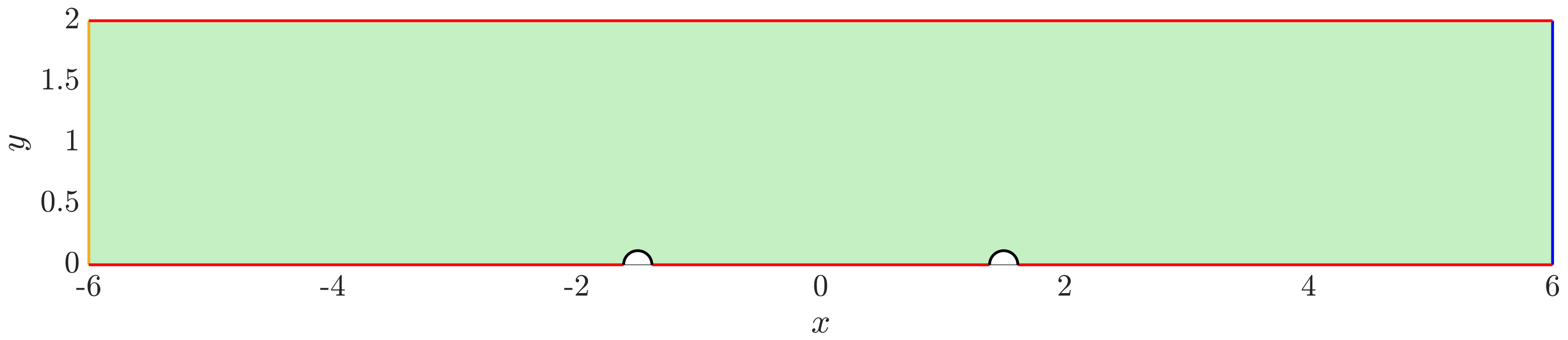}
	\caption{Computational domain for the simulation of the axisymmetric flow around the \emph{push-me-pull-you} microswimmer. The inflow boundary is highlighted in orange, the outflow boundary in blue, the slip boundary in red and the no-slip boundary in black.}
	\label{fig:swimmerGeometry}
\end{figure}

First, the axial symmetry of the problem is exploited and the computational domain is defined as $\Omega = \left([-L,L] \times [0,H]\right) \setminus \left( \mathcal{B}^+ \cup \mathcal{B}^- \right)$, where 
\begin{equation}
\mathcal{B}^{\pm} = \{ \bX \in \mathbb{R}^2 \; : \; \| \bX \pm \bX_0 \|_2 \leq \Rref \},
\end{equation} 
where $L {=} 6$, $H {=} 2$, $\bX_0 {=} (1.5,0)$ and $\Rref {=} 0.116$, as represented in figure~\ref{fig:swimmerGeometry}.

A Dirichlet boundary condition is imposed on the left portion of the boundary to simulate an inlet with unitary horizontal velocity. A homogeneous Neumann condition is applied on the right vertical boundary to simulate an outlet surface. On the surface of the two bladders, a homogeneous Dirichlet condition describes a no-slip boundary. Finally, on the remaining boundaries of the domain a perfect slip condition is enforced.

The geometry of the computational domain is described using two parameters: the parameter $\mu_1$ characterises the radius of the two spherical bladders; the parameter $\mu_2$ controls the distance between the centres of the two spheres. It is worth noticing that only one parameter is required to control the radius of the spheres because the total volume of the two bladders is kept constant. 

\begin{figure}[!tb]
	\centering
	\subfigure[$\bm{\mathcal{M}}_{\mu_1}$ \label{fig:swimmeMapRadius}]{\includegraphics[width=0.49\textwidth]{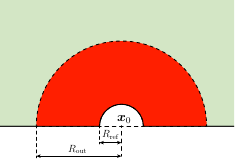}}
	\subfigure[$\bm{\mathcal{M}}_{\mu_2}$ \label{fig:swimmeMapDistance}]{\includegraphics[width=0.49\textwidth]{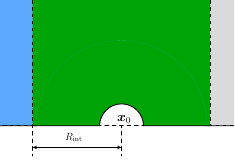}}
	\caption{Piecewise mapping (a) $\bm{\mathcal{M}}_{\mu_1}$ and (b) $\bm{\mathcal{M}}_{\mu_2}$ controlling the variation of radius and distance, respectively. Detail of the mapping in the vicinity of the sphere centred at $\bX_0$.}
	\label{fig:swimmeMap}
\end{figure}

Let $\bm{\mathcal{M}}_{\mu_1}$ be the mapping controlling the radius of the two spheres. Figure~\ref{fig:swimmeMapRadius} reports a sketch of the piecewise definition of the mapping in the vicinity of one sphere, $\Rout {=} 0.45$ being the interface between the deformable region (inside) and the fixed one (outside).
Following~\cite{sevilla2020solution,RS-SBGH-20}, the mapping $\bm{\mathcal{M}}_{\mu_1}$ is defined in the separable form~\eqref{eq:displacementSep} as
\begin{equation} \label{eq:swimmerMapRadius}
\begin{aligned}
\bM^1_1(\bX) & = \left\{ 
\begin{split}
	\frac{1}{r} \bX_0^-  & \text{\quad if $\|\bX_0^-\| \leq \Rout$	}\\
	0 & \text{\quad otherwise}
\end{split}
\right.
\quad &
\phi^1_1(\mu_1) & = \dfrac{ \Rout (R^+(\mu_1) - \Rref)}{\Rout - \Rref}, 
\\
\bM^2_1(\bX) & = \left\{ 
\begin{split}
 \bX_0^- & \text{\quad if $\|\bX_0^-\| \leq \Rout$	}\\
 0 & \text{\quad otherwise}
\end{split}
\right.
\quad &
\phi^2_1(\mu_1) & = \dfrac{\Rout- R^+(\mu_1)}{\Rout - \Rref}, 
\\
\bM^3_1(\bX) & = \left\{ 
\begin{split}
\bX_0 & \text{\quad if $\|\bX_0^-\| \leq \Rout$	}\\
0 & \text{\quad otherwise}
\end{split}
\right.
\quad &
\phi^3_1(\mu_1) & = 1, 
\\
\bM^4_1(\bX) & = \left\{ 
\begin{split}
\frac{1}{r}  \bX_0^+ & \text{\quad if $\| \bX_0^+ \| \leq \Rout$	}\\
0 & \text{\quad otherwise}
\end{split}
\right.
\quad &
\phi^4_1(\mu_1) & = \dfrac{ \Rout (R^-(\mu_1) - \Rref)}{\Rout - \Rref}, 
\\
\bM^5_1(\bX) & = \left\{ 
\begin{split}
\bX_0^+  & \text{\quad if $\| \bX_0^+ \| \leq \Rout$	}\\
0 & \text{\quad otherwise}
\end{split}
\right.
\quad &
\phi^5_1(\mu_1) & = \dfrac{\Rout- R^-(\mu_1)}{\Rout - \Rref}, 
\\
\bM^6_1(\bX) & = \left\{ 
\begin{split}
-\bX_0 & \text{\quad if $\| \bX_0^+ \| \leq \Rout$	}\\
0 & \text{\quad otherwise}
\end{split}
\right.
\quad &
\phi^6_1(\mu_1) & = 1,
\end{aligned}
\end{equation}
where $\bX_0^\pm {:=} \bX \pm \bX_0$. In this work, the interval for the first parameter is defined as $\I^1 {=} [-1,1]$. In addition, the radius of the second sphere, centred at $\bX_0$, is defined as \\ $R^+(\mu_1) {=} {-}0.0372\mu_1^2 {+} 0.0968\mu_1 {+} 0.25$, whereas the radius of the first sphere, centred at $-\bX_0$, is given by
\begin{equation}
R^-(\mu_1) = \left( \frac{1}{32} - \left[ R^+(\mu_1) \right]^3 \right)^{\tfrac{1}{3}} .
\end{equation}
Hence, for $\mu_1 {=} {-} 1$, the radii of the two spheres are $R^- {=} 0.3096$ and $R^+ {=} 0.116$, respectively, whereas their values are $R^- {=} 0.116$ and $R^+ {=} 0.3096$ for $\mu_1 {=} 1$. 

Similarly, the mapping $\bm{\mathcal{M}}_{\mu_2}$ for the distance is defined in a piecewise form as represented in figure~\ref{fig:swimmeMapDistance}. Let $\Rint {=} 0.47$ denote the location of the interface between the fixed and deformable region for the mapping affecting the distance of the spheres. The definition of the separable form of the mapping is given by 
\begin{equation} \label{eq:swimmerMapDistance}
\begin{aligned}
\bM^1_2(\bX) & = 
\begin{Bmatrix}
d(x) \\ 
0 
\end{Bmatrix} 
\qquad &
\phi^1_2(\mu_2) & = -\frac{1}{3} x_0\mu_2 , 
\\
\bM^2_2(\bX) & = \bX
\qquad &
\phi^2_2(\mu_2) & = 1,
\end{aligned}
\end{equation}
where the distance function $d(x)$ is
\begin{equation} \label{eq:fDistance}
d(x):= 
\begin{cases}
\displaystyle \frac{x+L}{x_0+\Rint-L} & \text{if }  x \in [-L,-x_0-\Rint]  \\
\displaystyle -1 & \text{if }  x \in [-x_0-\Rint,-x_0+\Rint]  \\
\displaystyle \frac{x}{x_0-\Rint} & \text{if }  x \in [-x_0+\Rint,x_0-\Rint]  \\
\displaystyle 1 & \text{if }  x \in [x_0-\Rint,x_0+\Rint]  \\
\displaystyle \frac{x-L}{x_0+\Rint-L} & \text{if }  x \in [x_0+\Rint,L]  
\end{cases} .
\end{equation}
For the parameter $\mu_2$, two different intervals are considered. The objective of this study is to evaluate the sensitivity of a priori and a posteriori PGD algorithms to variations in the amplitude of the parametric intervals, with special emphasis on the case in which large deformations of the domain are induced by the geometric parameters. On the one hand, the interval $\I^2 {=} [-2,-1]$ induces a maximum and minimum distance between the bladders equal to $\Dmax {=} 5$ and $\Dmin {=} 4$, respectively. On the other hand, the interval $\I^2 {=} [-3,2]$ is responsible for a maximum and minimum distance between the bladders equal to $\Dmax {=} 6$ and $\Dmin {=} 1$, respectively. It is worth noticing that in the latter setup, extreme deformations are introduced by the mapping, leading to complex variations of the flow features near the microswimmer.

Finally, when the two parameters are concurrently analysed, the resulting mapping is obtained as the composition of the two mappings associated with the radius and the distance. It is worth noticing that given the above piecewise definitions of $\bm{\mathcal{M}}_{\mu_1}$ and $\bm{\mathcal{M}}_{\mu_2}$, the resulting mappings are only $\mathcal{C}^0$ in the spatial domain. Hence, the meshes introduced in the following section need to be conforming with the artificial interfaces utilised to define the mappings, that is, the dashed lines in figure~\ref{fig:swimmeMap}. Of course, alternative definitions of these mappings may be devised, e.g. by imposing higher regularity across the interfaces~\cite{lovgren2009global}. Consequently, different intervals of the parameters $\mu_1$ and $\mu_2$ may arise. The choices above stem from the works~\cite{sevilla2020solution,RS-SBGH-20} and interested readers can compare the results computed for the small and large interval $\I^2$ with the ones reported in~\cite{sevilla2020solution} and~\cite{RS-SBGH-20}, respectively.

%-------------------------------------------------------
\subsection{Problem setup and comparison criteria} 
%-------------------------------------------------------

\begin{figure}[!tb]
	\centering
	\includegraphics[width=\textwidth]{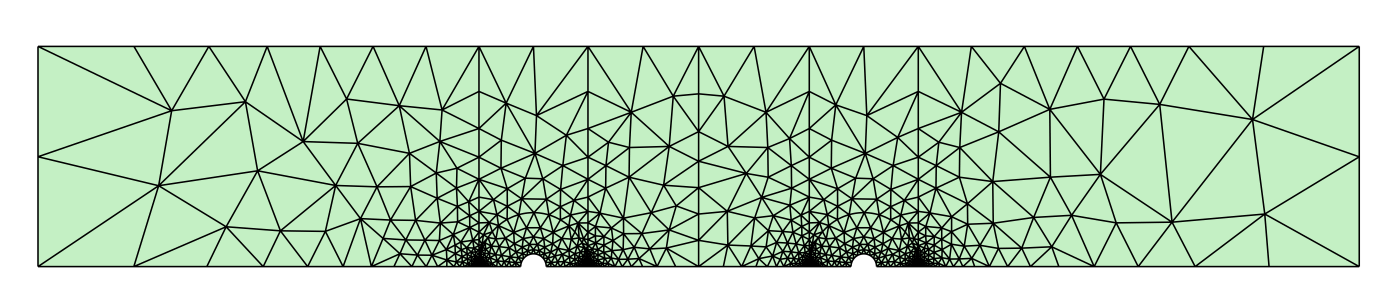}
	\caption{Computational mesh for the simulation of the axisymmetric flow around the \emph{push-me-pull-you} microswimmer.}
	\label{fig:swimmerMesh}
\end{figure}

Figure~\ref{fig:swimmerMesh} shows the computational mesh of the reference domain for the \emph{push-me-pull-you} microswimmer. The mesh has 1,426 fourth-order triangular elements.
The spatial discretisation leads to a global HDG system of 22,260 equations.

For the a priori PGD approach, 10 elements are used to discretise the parametric domain $\I^1$, whereas 20 and 100 elements are employed for the intervals $\I^2 {=} [-2,-1]$ and $\I^2 {=} [-3,2]$, respectively. The different number of elements in each parametric interval has been selected after observing that the variation in the flow induced by the first parameter is weaker than the variation induced by the second one~\cite{sevilla2020solution,RS-SBGH-20}. It is worth noticing that the set of nodes used to discretise the second parametric dimension in the first case, $\I^2 {=} [-2,-1]$, is a subset of the nodes selected for the second case, where $\I^2 {=} [-3,2]$. In all the numerical tests, a degree of approximation $k {=} 4$ is utilised for both the spatial and the parametric discretisations and non-uniform Fekete nodal distributions are employed. For the a posteriori PGD algorithm, the snapshots are computed in correspondance of the position of the nodes in the parametric space utilised for the a priori approach.

To compare the accuracy of the a priori and the a posteriori PGD algorithms, two error measures are considered. 

First, a multidimensional $\eltwo(\Omega \times \bI)$ error is defined for each variable, namely velocity, pressure and gradient of the velocity, by using a reference solution. For instance, the multidimensional $\eltwo(\Omega \times \bI)$ error measure for the velocity field is given by
\begin{equation} \label{eq:PGDerrorMulti}
E_u 
:= \left( \frac{  \int_{\bI} \int_{\Omega} ( \bupgd\!(\bX,\bmu) - \buref\!(\bX,\bmu) ) {\cdot} ( \bupgd\!(\bX,\bmu) - \buref\!(\bX,\bmu) ) d\Omega \,  d\mu }
{  \int_{\bI} \int_{\Omega}  \buref\!(\bX,\bmu) {\cdot} \buref\!(\bX,\bmu) \, d\Omega \, d\mu} \right)^{\!\! 1/2}.
\end{equation}
It is worth noticing that the evaluation of the multidimensional error~\eqref{eq:PGDerrorMulti} requires the definition of a reference solution $\buref$ at each integration point of the space of parameters. For each integration point of the space $\bI$, the reference solution is thus computed using the HDG spatial solver on a new, finer, mesh of the reference domain and with a higher order polynomial approximation. It is worth mentioning that an increased number of integration points on the spatial mesh is considered for this computation in order to guarantee that the reference solution is not affected by the introduction of the mapping in the HDG solver. It follows that the error introduced by the interpolation of the reference solution defined on the fine mesh onto the mesh used for the PGD computation is negligible and the quantity in equation~\eqref{eq:PGDerrorMulti} provides an accurate description of the error of the reduction strategy. For the example with $\I^1 {=} [-1,1]$ and $\I^2 {=} [-3,2]$, a total of 25,000 reference solutions were required to compute this error measure, since five integration points in each parametric element are utilised for $k {=} 4$.

Second, the separated response surface for the drag force $\FD$ and its error are considered to assess the accuracy of the PGD-based strategies analysed. More precisely, the $\eltwo(\bI)$ error measure for the drag force in the parametric space is defined as
\begin{subequations}
\begin{equation} \label{eq:DragErrorMulti}
E_D 
= \left( \frac{  \int_{\bI} ( \FDpgd\!(\bmu) - \FDref\!(\bmu) )^2 \,  d\mu }
{  \int_{\bI}  \FDref\!(\bmu)^2\,  d\mu} \right)^{\!\! 1/2}, 
\end{equation}
whereas the error in the quantity of interest $\FD$ as a function of the parameters is given by
\begin{equation} \label{eq:DragErrorPointwise}
\varepsilon_D (\bmu)
= \frac{  | \FDpgd\!(\bmu)  - \FDref\!(\bmu) | }{ | \FDref\!(\bmu) |} .
\end{equation}
\end{subequations}

%In addition to the errors in velocity, pressure and gradient of the velocity,  The two measures used for the accuracy of the drag force are given by
%%
%\begin{equation} \label{eq:DragErrorMulti}
%E_D 
%= \left( \frac{  \int_{\bI} ( \FDpgd - \FDref )^2 \,  d\mu }
%{  \int_{\bI}  \FDref^2\,  d\mu} \right)^{1/2}, 
%\qquad
%\varepsilon_D (\bmu)
%= \frac{  | \FDpgd- \FDref | }{ | \FDref |} .
%\end{equation}

%-------------------------------------------------------
\subsection{One geometric parameter} 
\label{sc:exOne}
%-------------------------------------------------------

In this section, two geometric mappings, affecting independently the radius of the spherical bladders and their distance, are considered. An extensive comparison of accuracy and computational cost of the a priori and the a posteriori PGD approaches is presented for these two cases and special attention is devoted to the PGD-based separated response surfaces for the drag force.

%-------------------------------------------------------
\subsubsection{Varying the radius of the spherical bladders} 
\label{sc:exMu1}
%-------------------------------------------------------

The first example involves the simulation of the Stokes flow past the \emph{push-me-pull-you} microswimmer when the domain is parametrised using $\mu_1$ and the distance between the centres of the two spheres is fixed and equal to 3. To evaluate the influence of the number of nonlinear iterations in the AD scheme of the a priori PGD, different numbers of iterations are considered, namely $\niter {=} 1,2,3,5$. In addition, to evaluate the influence of the number of snapshots used in the a posteriori PGD approach, different numbers of snapshots are employed, namely $\nsnap {=} 11,21,41$.

\begin{figure}[!tb]

	\subfigure[Mesh, $\mu_1=-1$]{\includegraphics[width=0.94\textwidth]{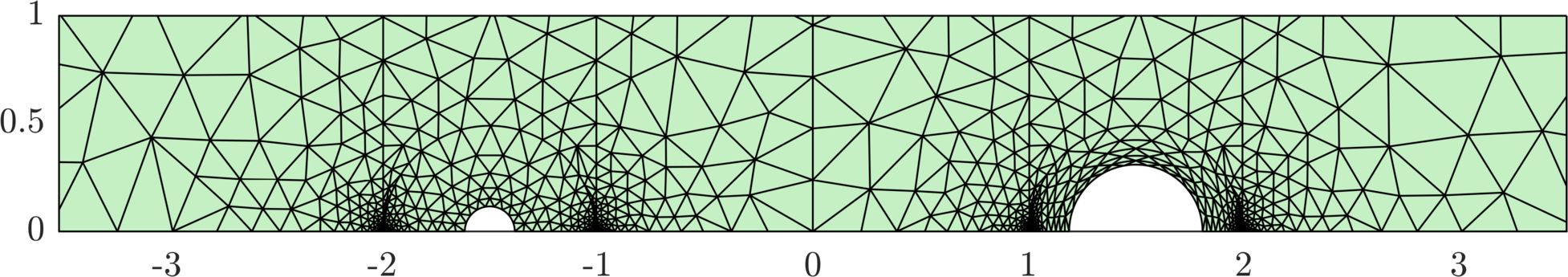}}
	
	{\centering		
	\subfigure[Quality, $\mu_1=-1$]{\includegraphics[width=\textwidth]{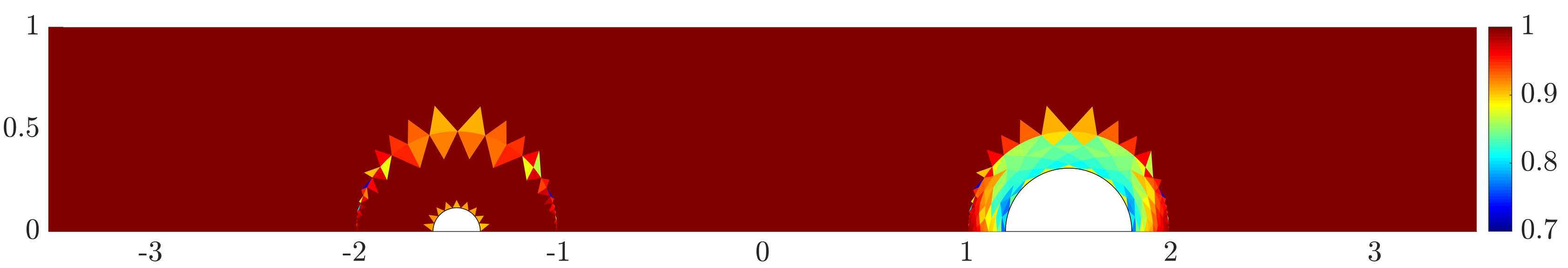}}
	}
	
	\subfigure[Mesh, $\mu_1= 1$]{\includegraphics[width=0.94\textwidth]{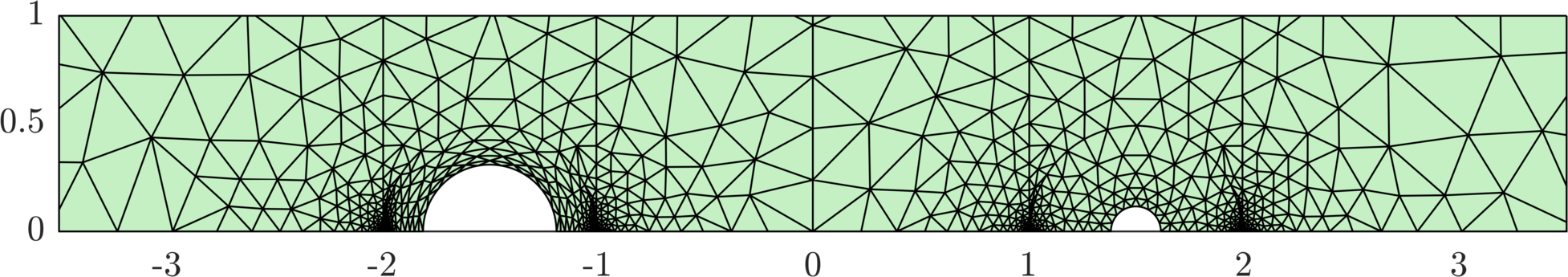}}
	
	{\centering	
	\subfigure[Quality, $\mu_1= 1$]{\includegraphics[width=\textwidth]{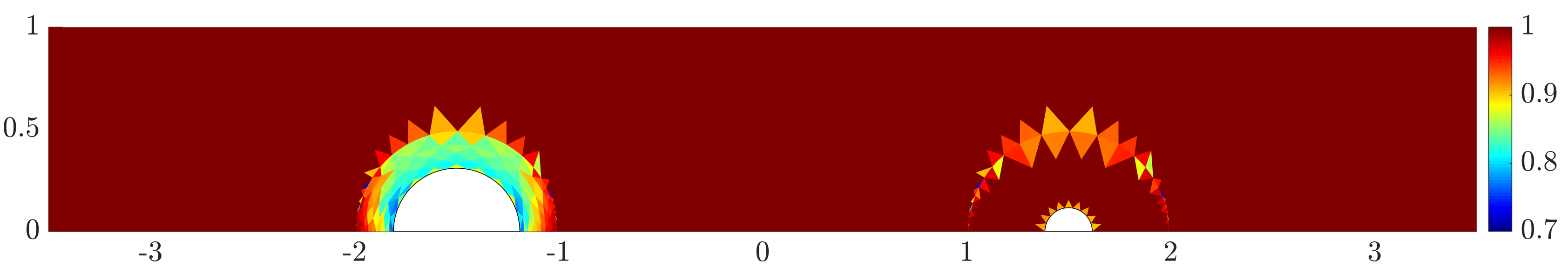}}
	}
	\caption{(a-c) Mesh configurations and (b-d) mesh quality of two deformed microswimmers for the mapping with $\mu_1$ as a geometric parameter.}
	\label{fig:swimmerMeshMappingMu1}
\end{figure}

First, the effect of the geometric mapping on the quality of the meshes is investigated. Figure~\ref{fig:swimmerMeshMappingMu1} displays the mesh configurations associated with the two extreme values of the parameter $\mu_1$ controlling the radius of the bladders. Moreover, the mesh quality, measured as the scaled Jacobian of the isoparametric mapping~\cite{poya2016unified,HO-Meshing}, is reported for the corresponding geometric configurations.
The mesh quality map shows that for the sphere with minimum radius, $0.116$, the quality is lower than one only in the elements with an edge on the boundary or on the interior interface used to define the piecewise geometric mapping described in~\cite{RS-SBGH-20}. This is due to the use of the mesh generation technique described in~\cite{poya2016unified,HO-Meshing}, where only the elements in contact with curved entities are represented with high-order polynomials. In contrast, for the sphere with maximum radius, 0.3096, all the elements in the region where the mapping is different from the identity are deformed. It is worth noticing that in all the presented cases, the majority of the elements features a mesh quality of 0.9 or higher and only few elements experience an extreme distortion reducing the value of the scaled Jacobian of the isoparametric mapping to 0.7. Hence, the mesh in figure~\ref{fig:swimmerMesh} is confirmed to be suitable for the parametric study of the influence of the radius in the microswimmer configuration.

\begin{figure}[!tb]
	\centering
	\subfigure[$\bu$] {\includegraphics[width=0.49\textwidth]{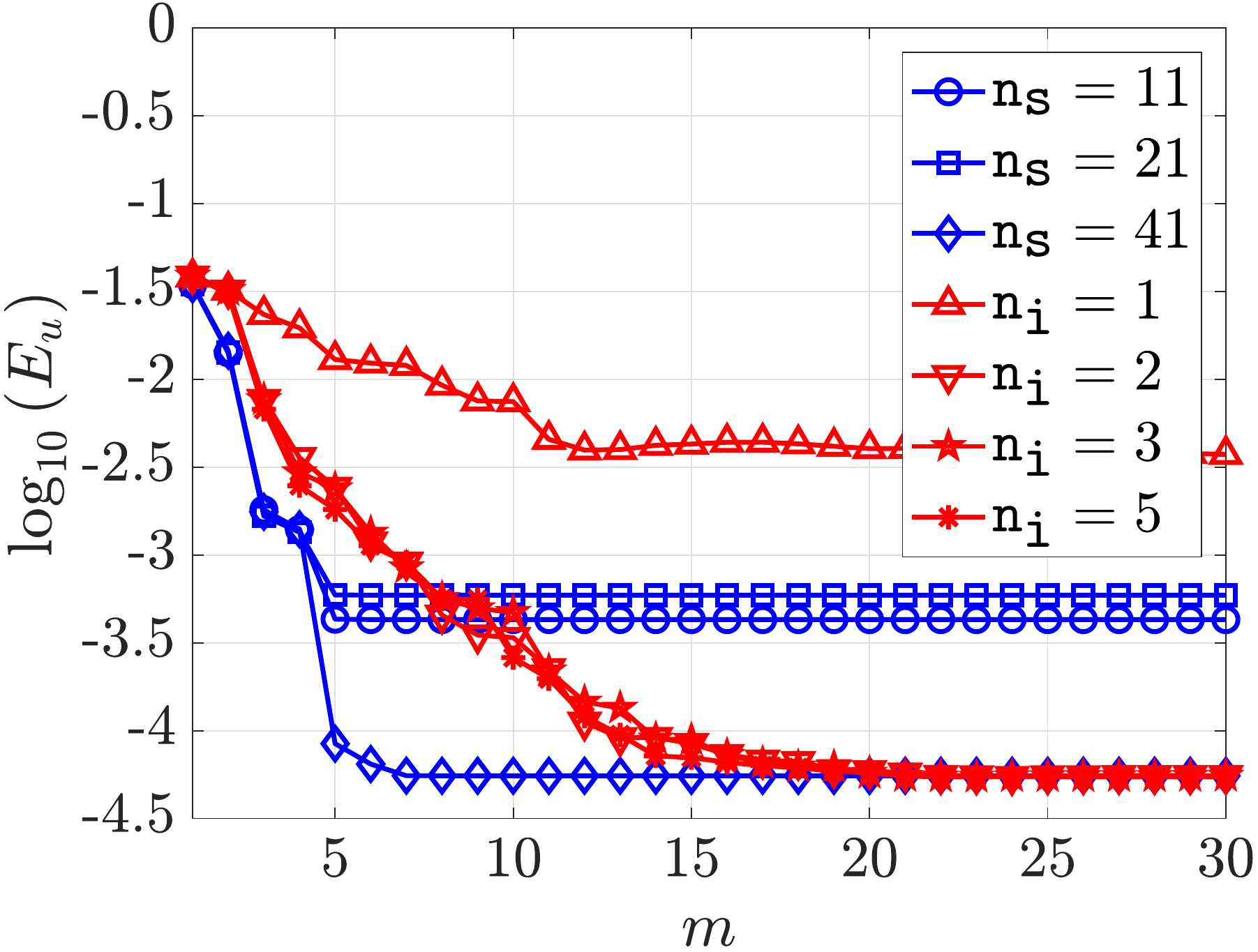}\label{fig:RadiusU}}
	\subfigure[$p$]   {\includegraphics[width=0.49\textwidth]{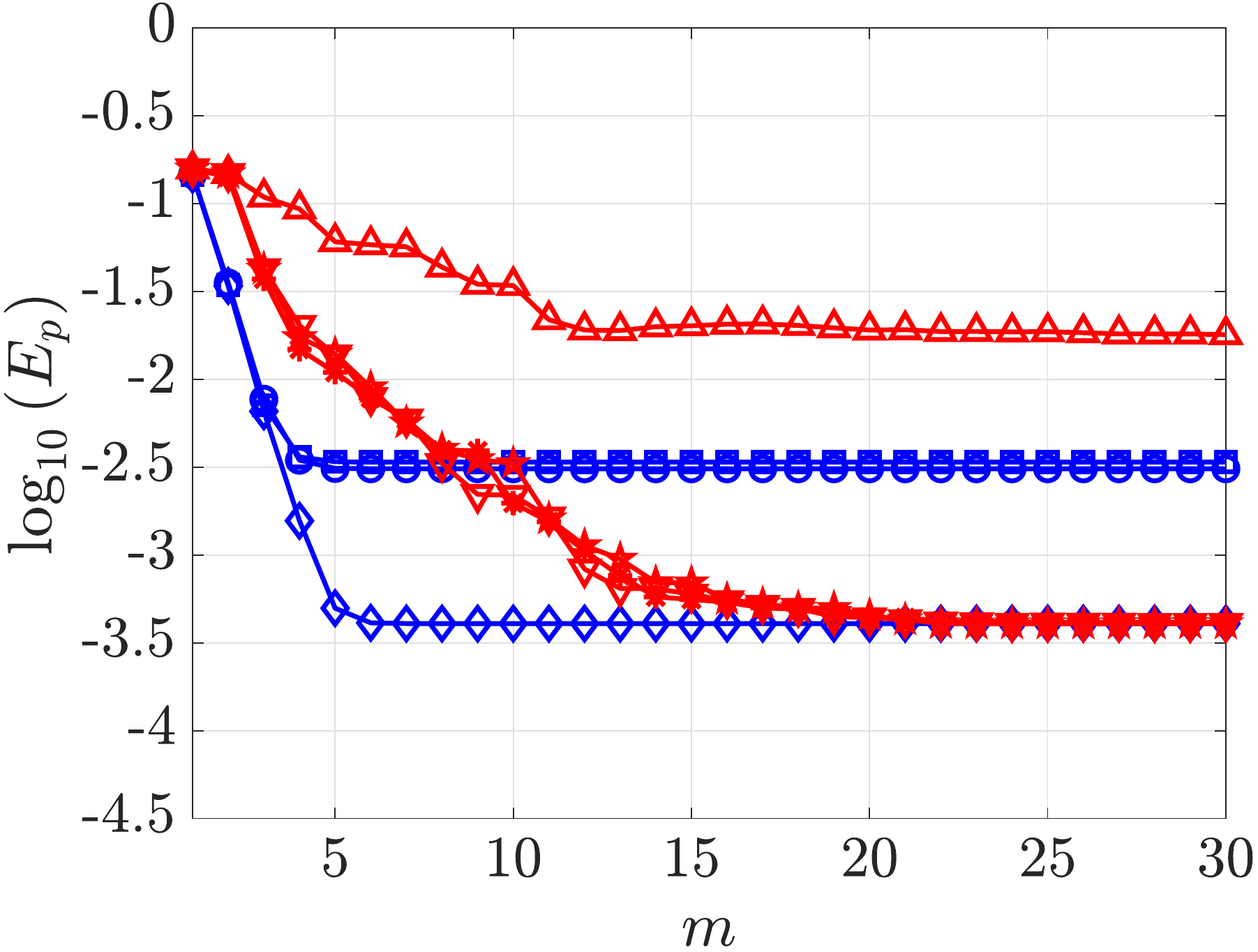}\label{fig:RadiusP}}
	\caption{Evolution of the $\eltwo(\Omega \times \I^1)$ error for (a) velocity and (b) pressure as a function of the number of PGD modes for the problem with one geometric parameter controlling the radius of the spherical bladders. The legend details the number $\nsnap$ of snapshots used by the a posteriori PGD approach (blue) and the number $\niter$ of nonlinear AD iterations used by the a priori PGD approach (red).}
	\label{fig:RadiusUP}
\end{figure}

Figure~\ref{fig:RadiusUP} shows the evolution of the $\eltwo(\Omega \times \I^1)$ error measure for velocity and pressure as a function of the number of modes, $m$, for both the a priori and the a posteriori PGD approaches. 
The results in figure~\ref{fig:RadiusU} show that the a posteriori approach is able to provide highly accurate results, with an error below $10^{-3}$, with only five modes and using 11 snapshots. Increasing the number of snapshots to 21, the accuracy of the approximation is not improved and the error with five modes stagnates at approximately the same level achieved by the PGD with 11 snapshots. The additional ten snapshots introduced for $\nsnap {=} 21$ are responsible for slightly perturbing the PGD solution computed using 11 snapshots, without providing valuable information on the parametric solution. Of course, this result is strongly influenced by the choice of the sampling points and advanced sampling techniques are expected to improve the a posteriori approximation in this case. Nonetheless, since the a priori PGD algorithm constructs an approximation without any prior information, no specific sampling has been considered for the a posteriori PGD in order to present a fair comparison of the cost of the two solvers under similar working conditions. If a higher accuracy is required for the approximation in figure~\ref{fig:RadiusU}, the number of snapshots needs to be increased to 41. In this case, with seven computed modes, the a posteriori PGD is able to provide an error below $10^{-4}$.

For the a priori approach, the results show that with only one iteration in the AD scheme, the accuracy of the computed modes is limited and the error stagnates at a level almost two orders of magnitude higher than the corresponding results obtained performing two iterations. In addition, this example also shows that two iterations is the optimal value as higher values, for instance three or five iterations, provide results with almost the same accuracy but they require additional solutions of the spatial problem. With two iterations in the AD scheme, the a priori PGD approach requires up to 20 modes to reach an accuracy that is comparable to the accuracy obtained by the a posteriori approach with five modes and 41 snapshots.

%The results in figure~\ref{fig:RadiusUP} also show that the qualitative behaviour of the error for the velocity and pressure fields is almost identical despite, as previously observed when performing full order HDG computations~\cite{}, the error in the pressure field is higher than the error in the velocity field. 

\begin{figure}[!tb]
	\centering
	\subfigure[$\bL$] {\includegraphics[width=0.49\textwidth]{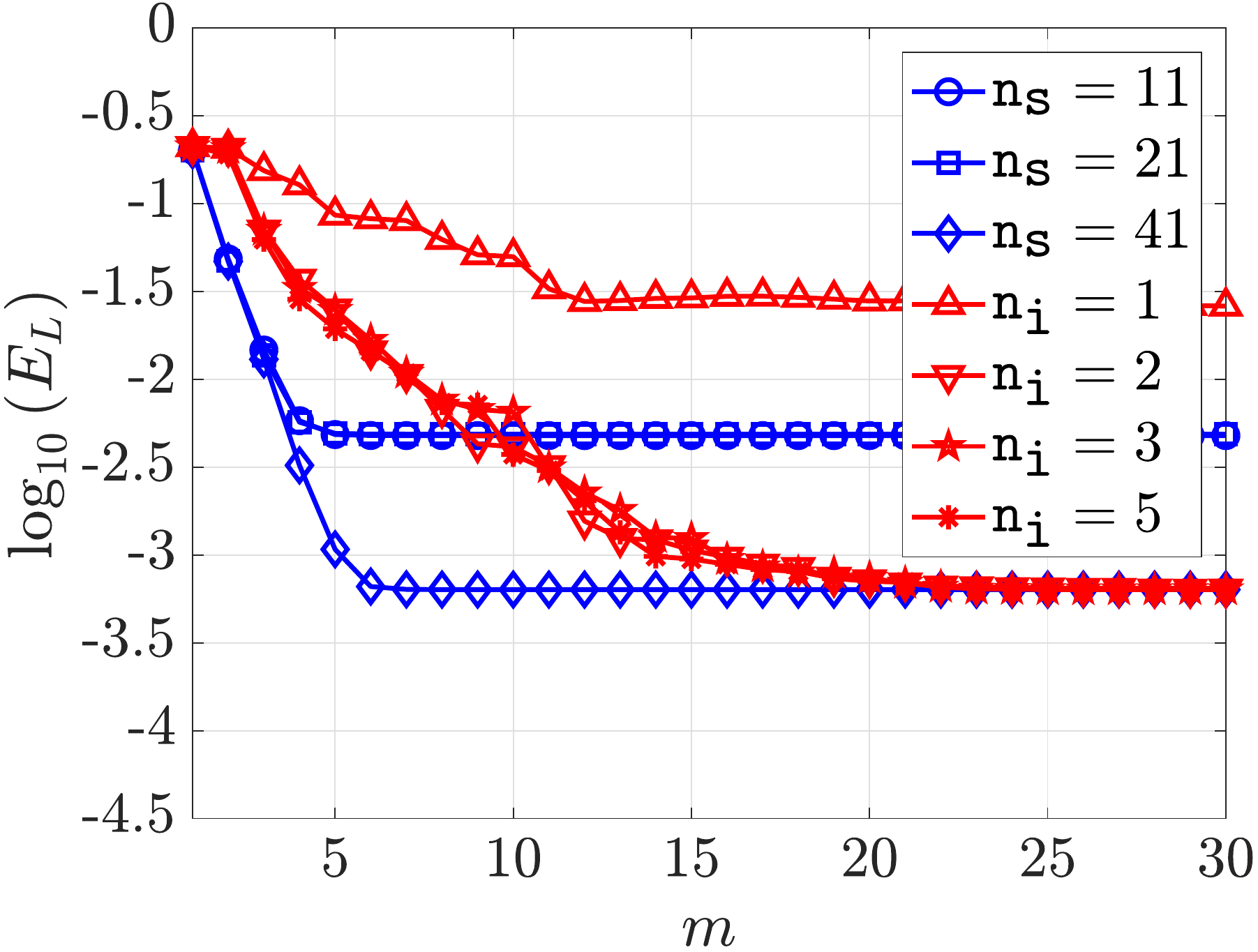}}
	\subfigure[$\FD$] {\includegraphics[width=0.49\textwidth]{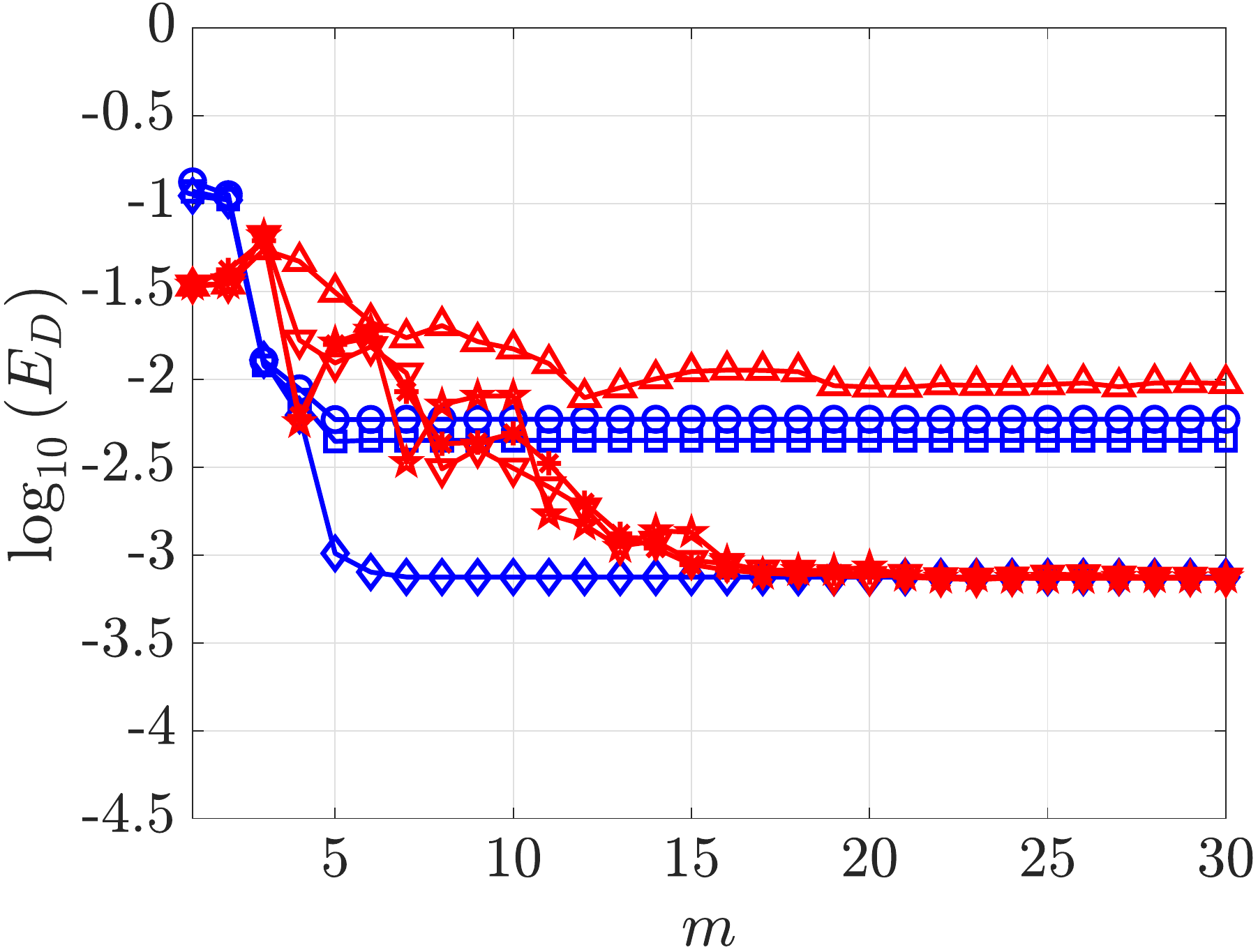}\label{fig:L2errDragRadius}}
	\caption{Evolution of (a) the $\eltwo(\Omega \times \I^1)$ error for the gradient of velocity and (b) the $\eltwo(\I^1)$ error for the drag force as a function of the number of PGD modes for the problem with one geometric parameter controlling the radius of the spherical bladders. The legend details the number $\nsnap$ of snapshots used by the a posteriori PGD approach (blue) and the number $\niter$ of nonlinear AD iterations used by the a priori PGD approach (red).}
	\label{fig:RadiusLFd}
\end{figure}

To further analyse the accuracy of the two PGD approaches, the evolution of the $\eltwo(\Omega \times \I^1)$ error for the gradient of velocity and the $\eltwo(\I^1)$ error for the drag force on the two spherical bladders is computed as a function of the number of modes (Fig.~\ref{fig:RadiusLFd}). 
It can be observed that the results for the gradient of velocity are very similar, qualitatively and quantitatively, to the ones presented in figure~\ref{fig:RadiusP} for the pressure field. This is due to the extra accuracy provided by the HDG formulation in the gradient of velocity, compared to other approaches based on primal formulations. Furthermore, this example also confirms that the accuracy that is obtained in the drag force is similar to the accuracy obtained in the pressure and in the gradient of velocity, from which it is computed. In all cases, the a posteriori PGD approach requires five modes and 41 snapshots to construct a solution with an error in the drag force below $10^{-3}$, whereas the a priori approach achieves a similar level of accuracy using two iterations and 15 modes.
Hence, the results show that the two PGD approaches require a similar computational cost to reach an error in the drag force below $10^{-3}$. The a priori approach requires the solution of 45 spatial problems (i.e. 15 modes, each computed with two iterations of the AD scheme plus the initial solve to perform the prediction of the mode, see algorithm~\ref{alg:PGDpriori}), whereas the a posteriori approach utilises 41 snapshots to reach the same level of accuracy.

\begin{figure}[!tb]
	\centering
	\subfigure[A priori PGD]{\includegraphics[width=0.48\textwidth]{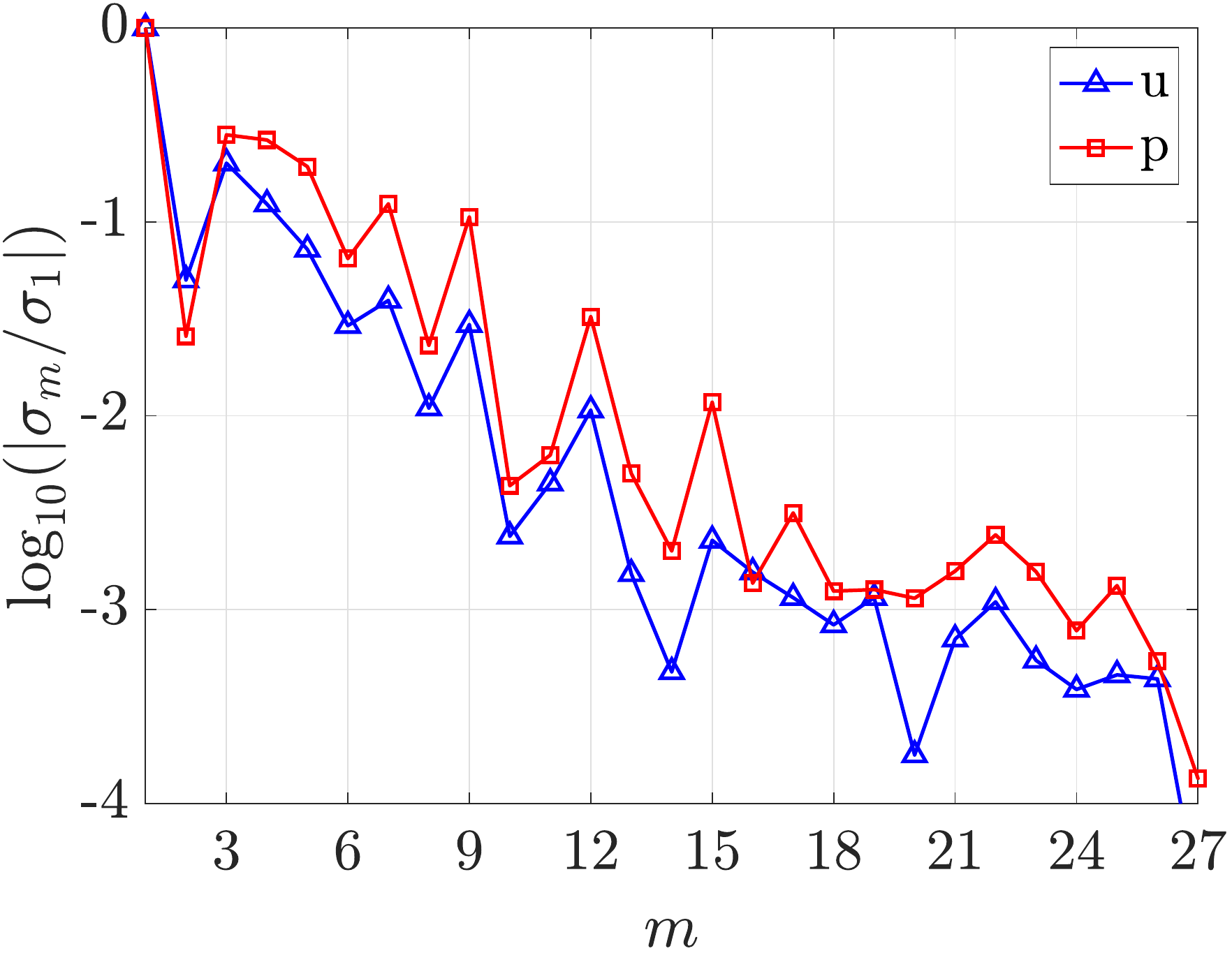}}
	\subfigure[A posteriori PGD]{\includegraphics[width=0.48\textwidth]{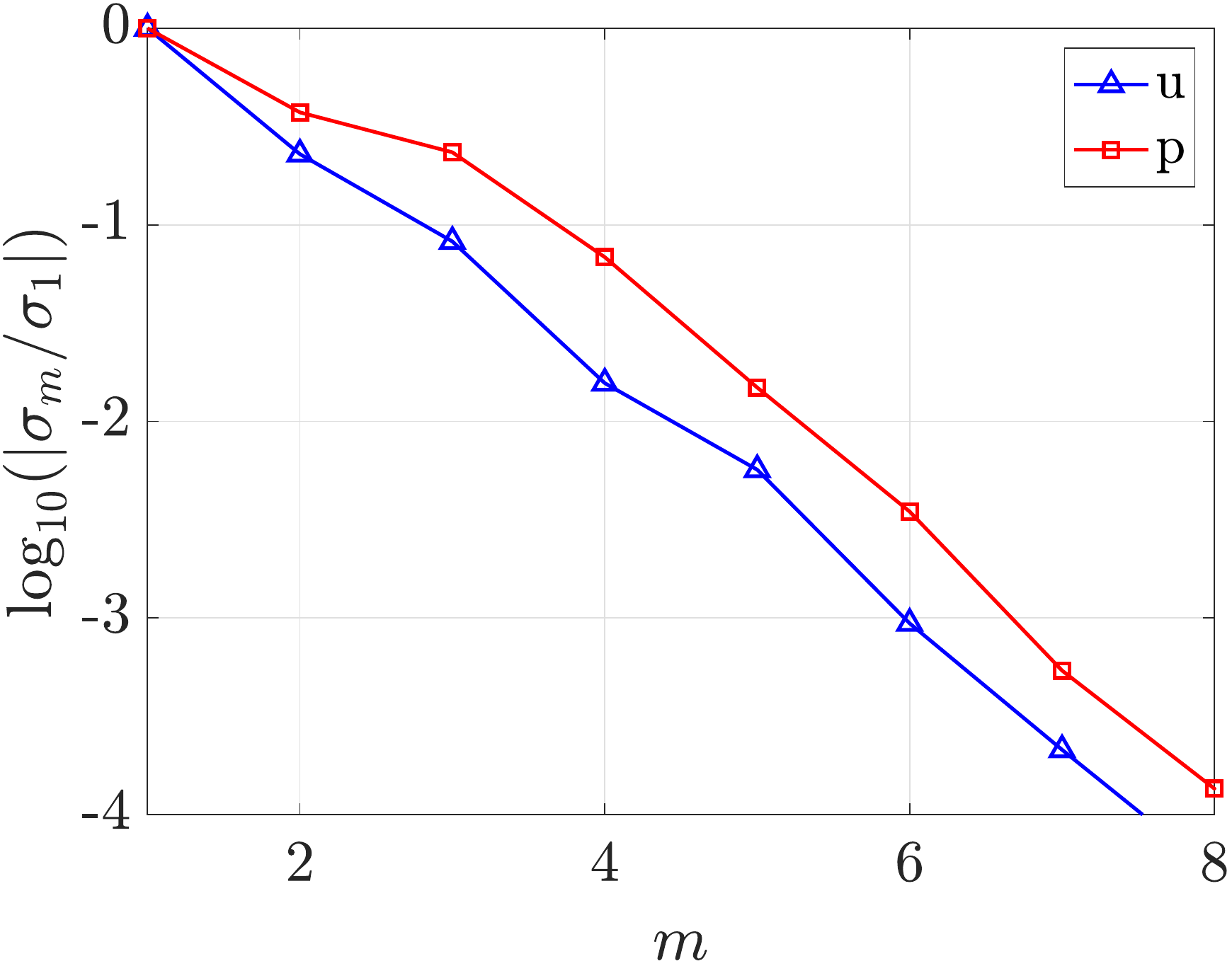}}

	\caption{Convergence of the mode amplitude computed using the (a) a priori and (b) a posteriori PGD algorithms for the parametric study of the radius.}
	\label{fig:RadiusAmp}
\end{figure}

Figure~\ref{fig:RadiusAmp} reports the evolution of the relative amplitude of the velocity and pressure modes computed using the a priori and a posteriori PGD algorithms. For the a priori PGD algorithm with $\niter {=} 2$, 24 modes are required to lower the relative amplitude of both the velocity and the pressure modes below $10^{-3}$, whereas the a posteriori algorithm with $\nsnap {=} 41$ only computes six modes. It is very important to emphasise that no compression of the modes obtained in the a priori approach, see~\cite{DM-MZH:15}, has been performed to enable the reader to clearly see the number of calls to the spatial solver required. However, the number of modes computed by the two approaches is expected to be the same when the PGD compression is performed. It is also worth noticing that both the a priori and the a posteriori algorithms stagnate at the same level of error as this is the error induced by the assumption of separability of the exact solution of the problem. 

\begin{figure}[!tb]
	\centering
	\subfigure[$m = 1$]{\includegraphics[width=0.48\textwidth]{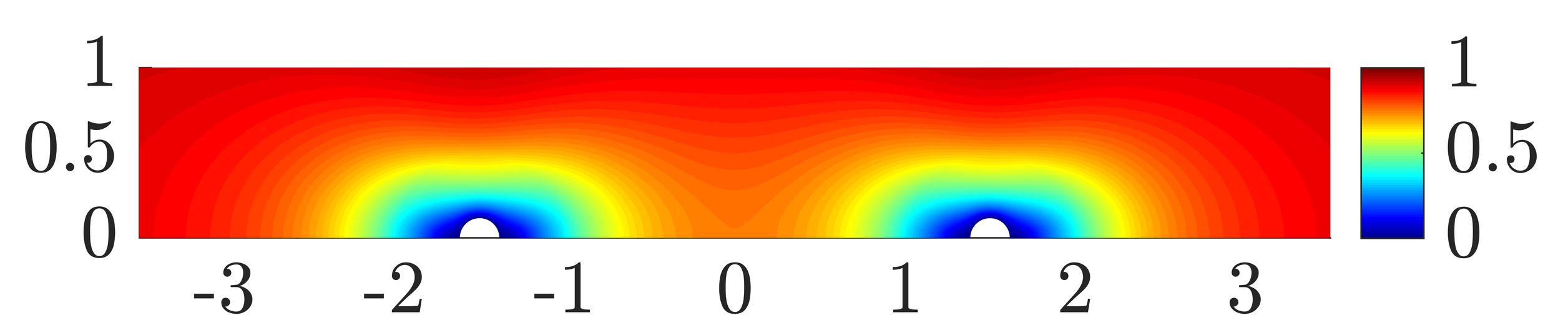}}
	\subfigure[$m = 2$]{\includegraphics[width=0.48\textwidth]{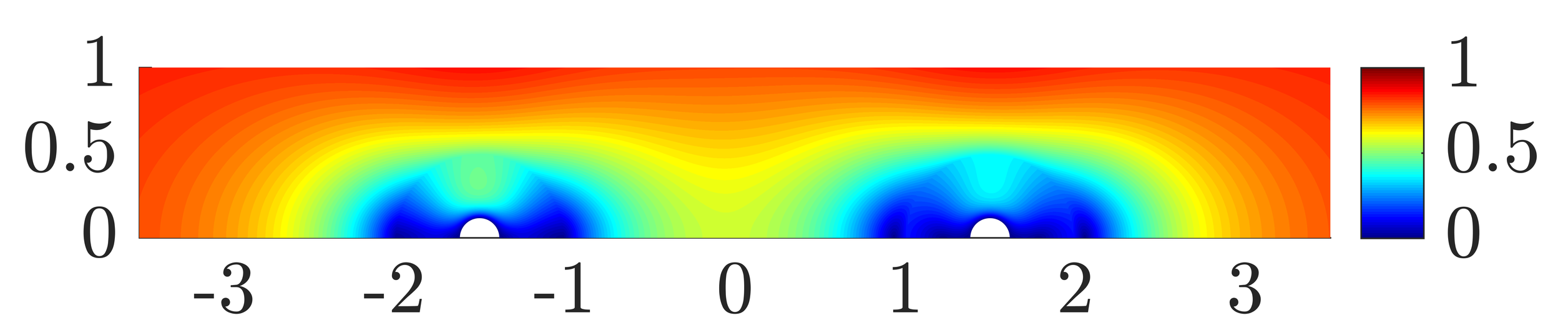}}

	\subfigure[$m = 3$] {\includegraphics[width=0.48\textwidth]{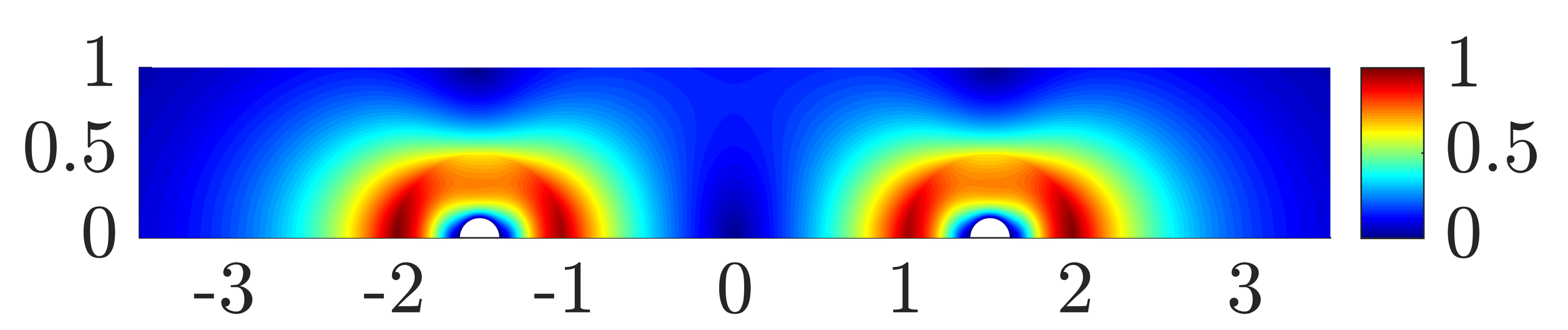}}
	\subfigure[$m = 4$] {\includegraphics[width=0.48\textwidth]{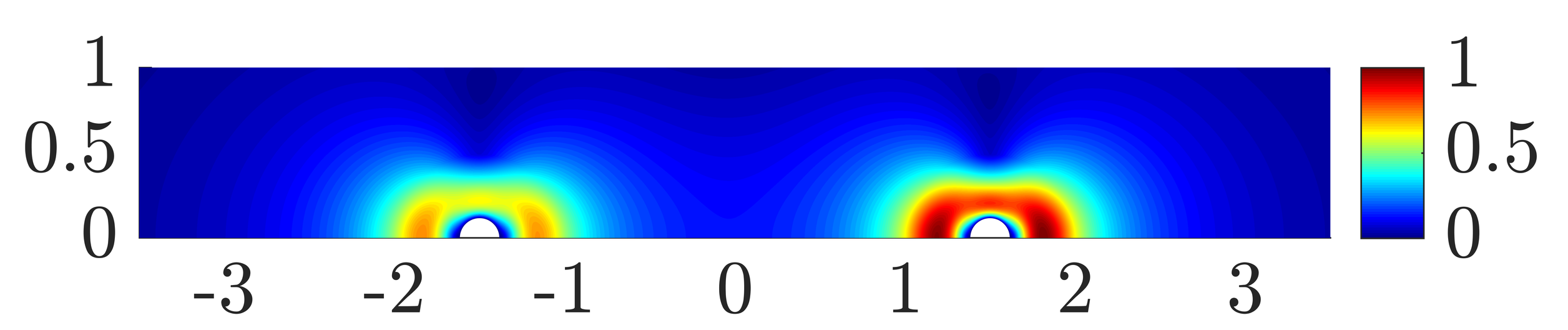}}
		
	\caption{First four normalised spatial modes of the module of the velocity computed using the a priori PGD algorithm for the parametric study of the radius.}
	\label{fig:RadiusModesUpriori}
\end{figure}

\begin{figure}[!tb]
	\centering
	\subfigure[$m = 1$]{\includegraphics[width=0.48\textwidth]{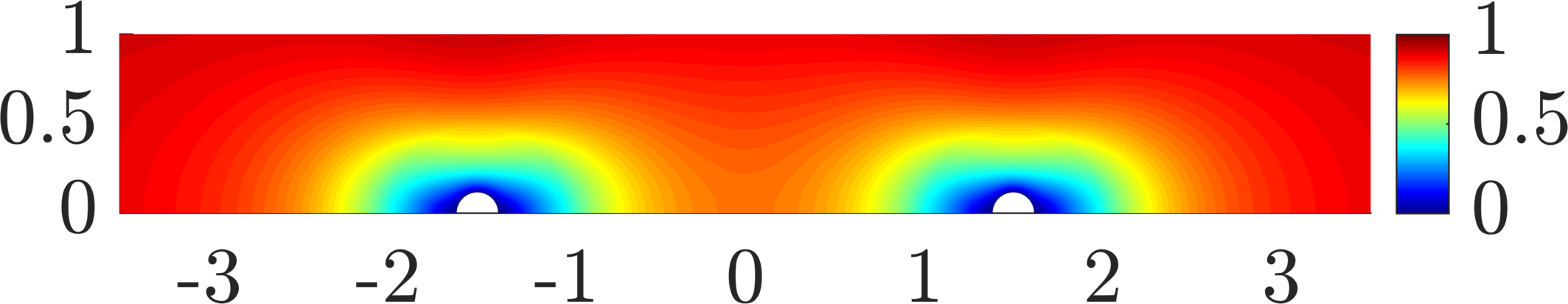}}
	\subfigure[$m = 2$]{\includegraphics[width=0.48\textwidth]{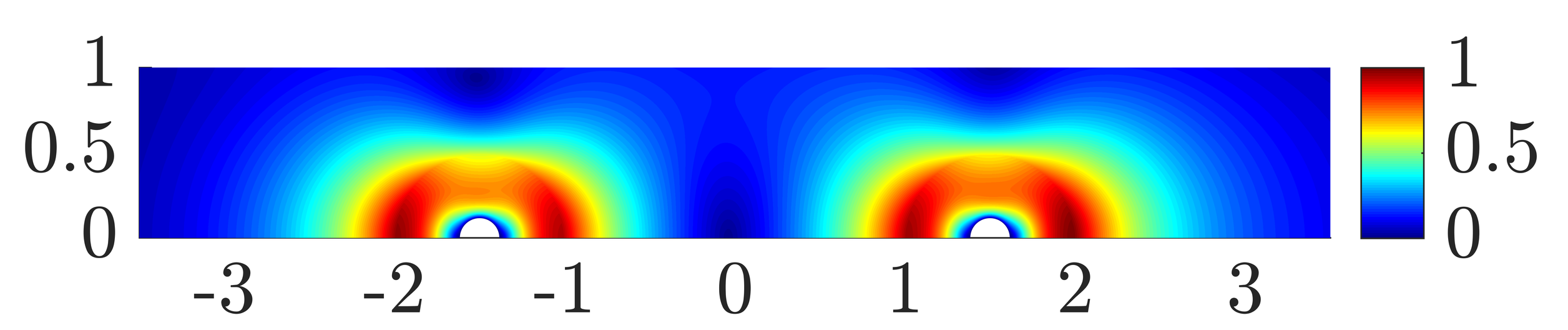}}

	\subfigure[$m = 3$] {\includegraphics[width=0.48\textwidth]{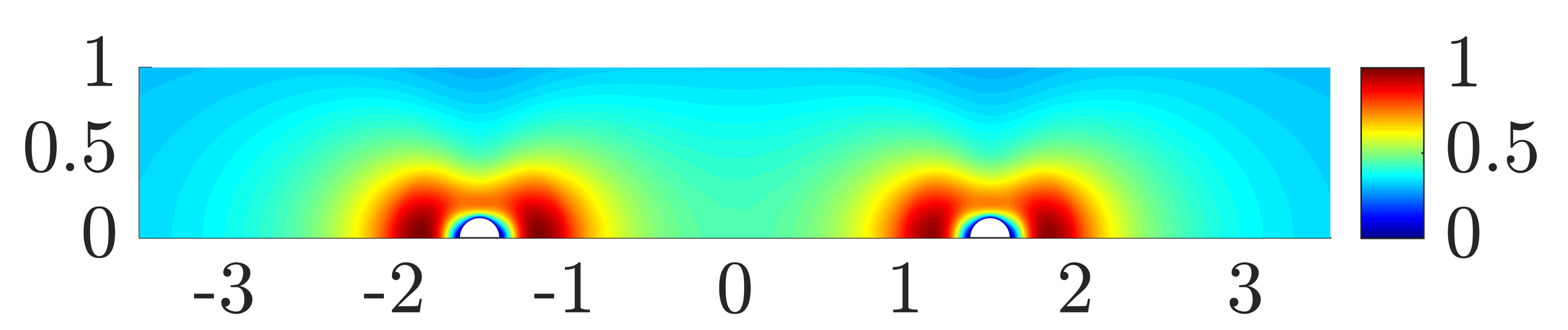}}
	\subfigure[$m = 4$] {\includegraphics[width=0.48\textwidth]{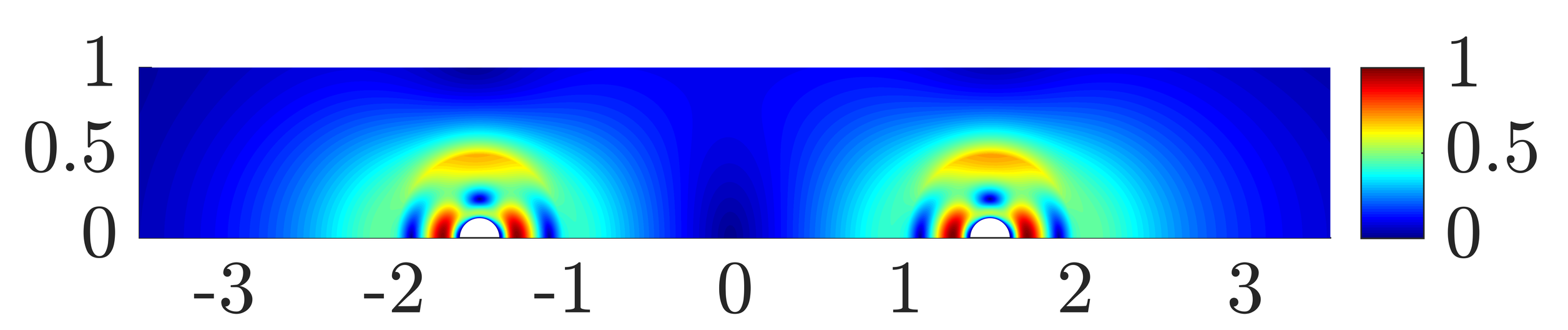}}
		
	\caption{First four normalised spatial modes of the module of the velocity computed using the a posteriori PGD algorithm for the parametric study of the radius.}
	\label{fig:RadiusModesUposteriori}
\end{figure}

\begin{figure}[!tb]
	\centering
	\subfigure[$m = 1$]{\includegraphics[width=0.48\textwidth]{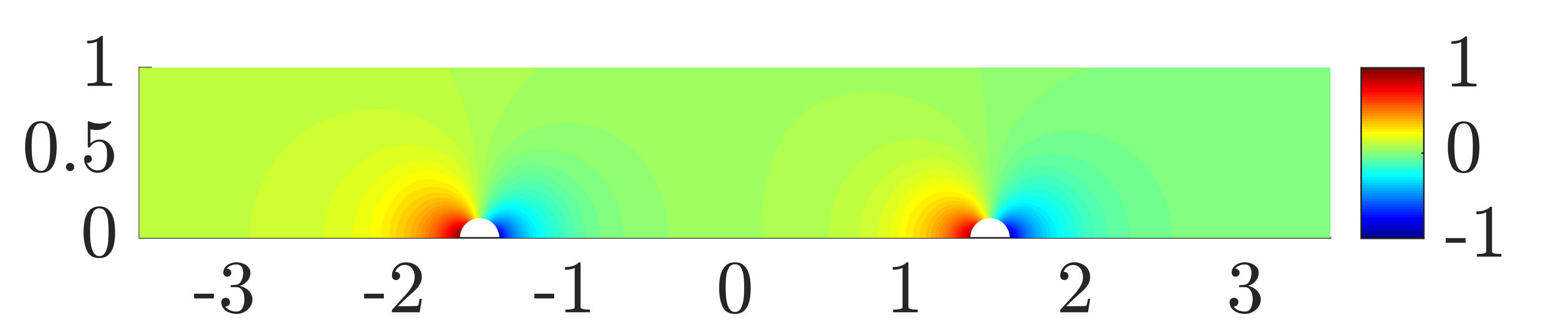}}
	\subfigure[$m = 2$]{\includegraphics[width=0.48\textwidth]{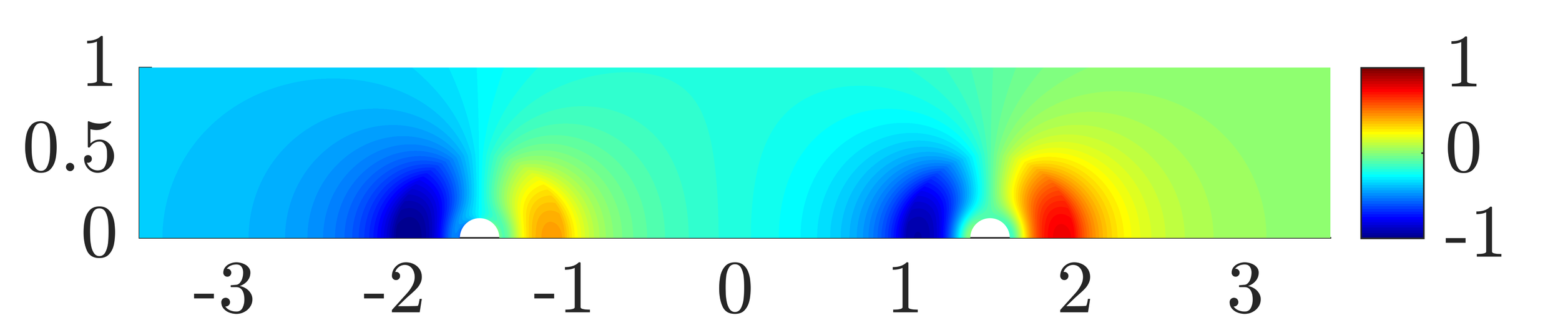}}

	\subfigure[$m = 3$] {\includegraphics[width=0.48\textwidth]{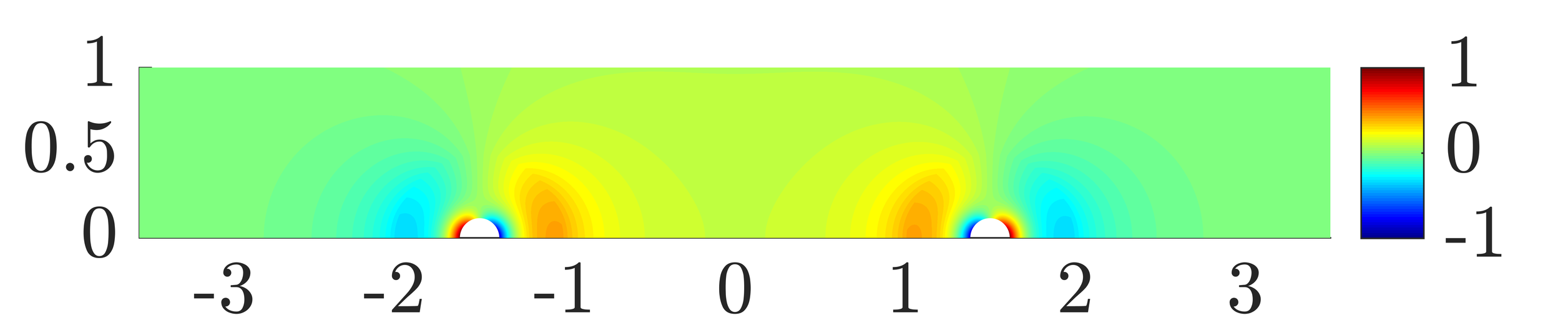}}
	\subfigure[$m = 4$] {\includegraphics[width=0.48\textwidth]{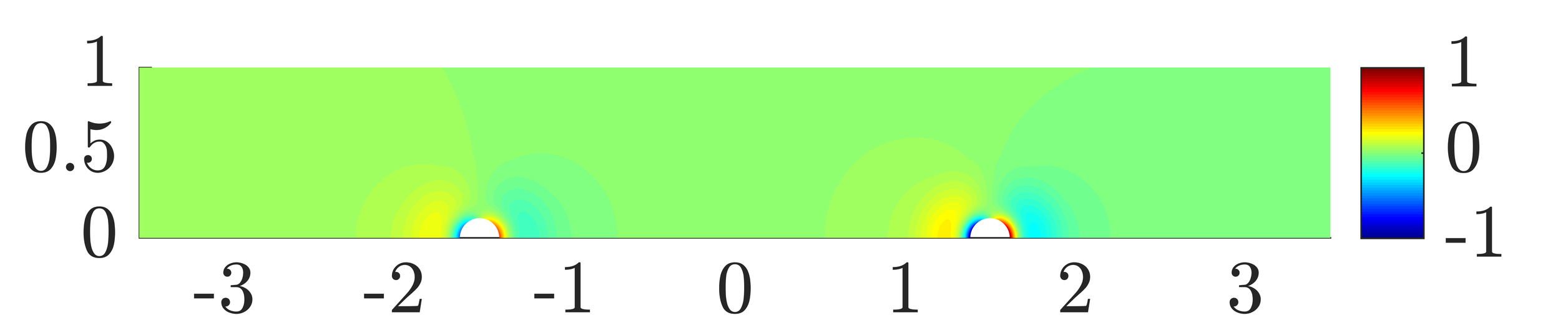}}
		
	\caption{First four normalised spatial modes of the pressure computed using the a priori PGD algorithm for the parametric study of the radius.}
	\label{fig:RadiusModesPpriori}
\end{figure}

\begin{figure}[!tb]
	\centering
	\subfigure[$m = 1$]{\includegraphics[width=0.48\textwidth]{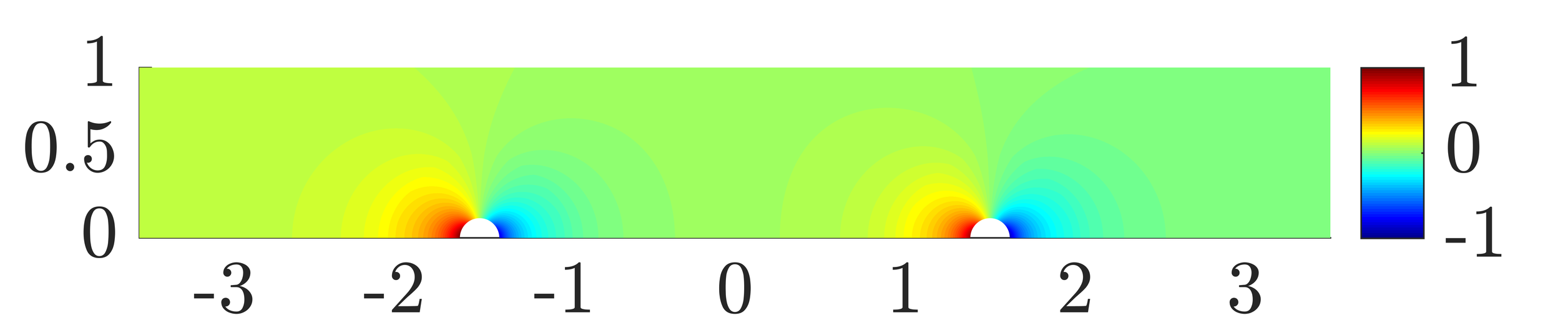}}
	\subfigure[$m = 2$]{\includegraphics[width=0.48\textwidth]{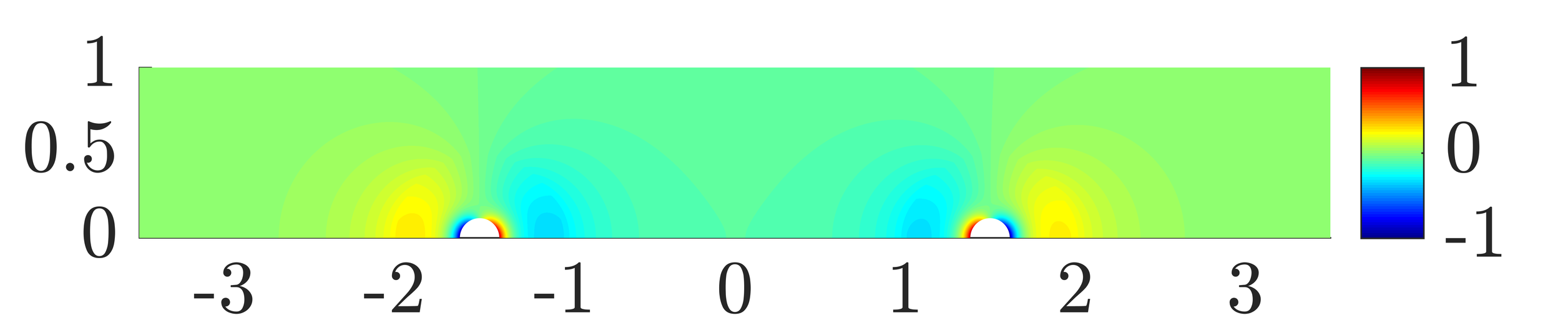}}

	\subfigure[$m = 3$] {\includegraphics[width=0.48\textwidth]{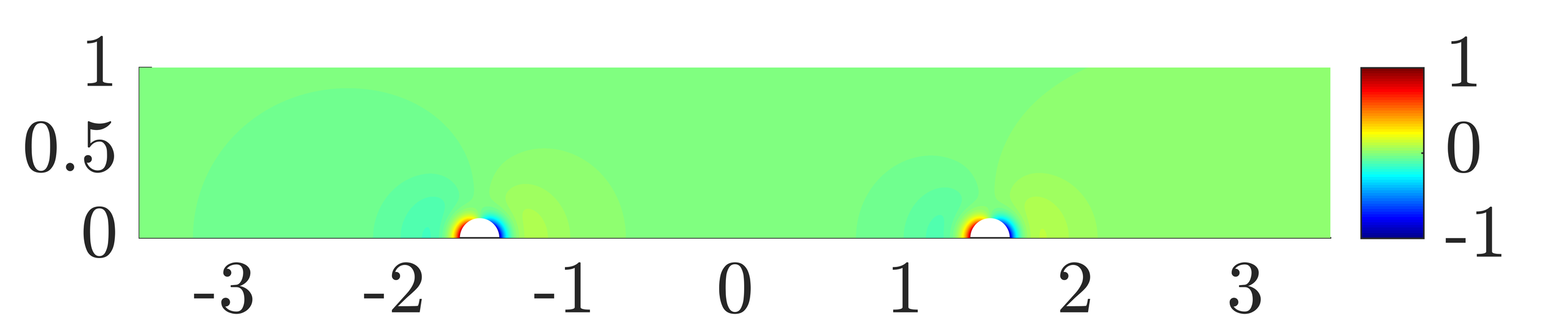}}
	\subfigure[$m = 4$] {\includegraphics[width=0.48\textwidth]{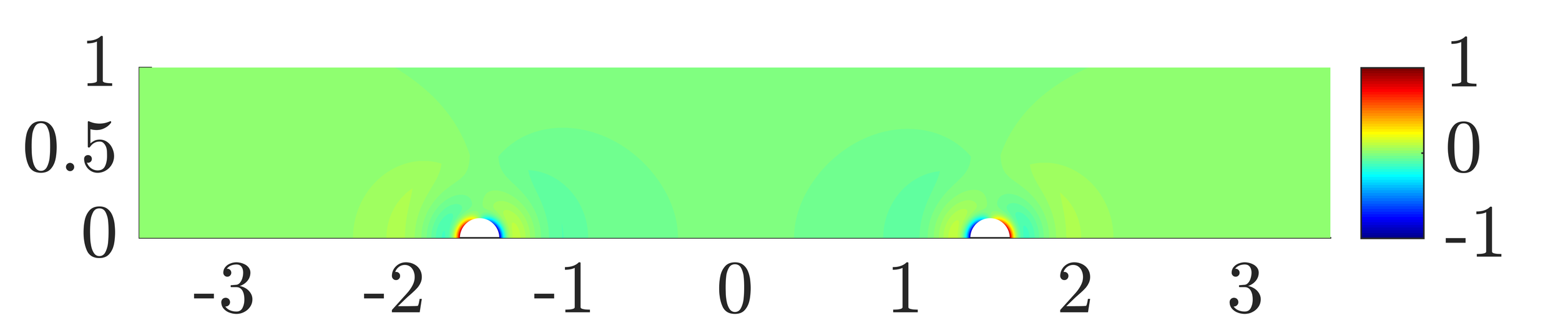}}
		
	\caption{First four normalised spatial modes of the pressure computed using the a posteriori PGD algorithm for the parametric study of the radius.}
	\label{fig:RadiusModesPposteriori}
\end{figure}

\begin{figure}[!tb]
	\centering
	\subfigure[A priori PGD]{\includegraphics[width=0.48\textwidth]{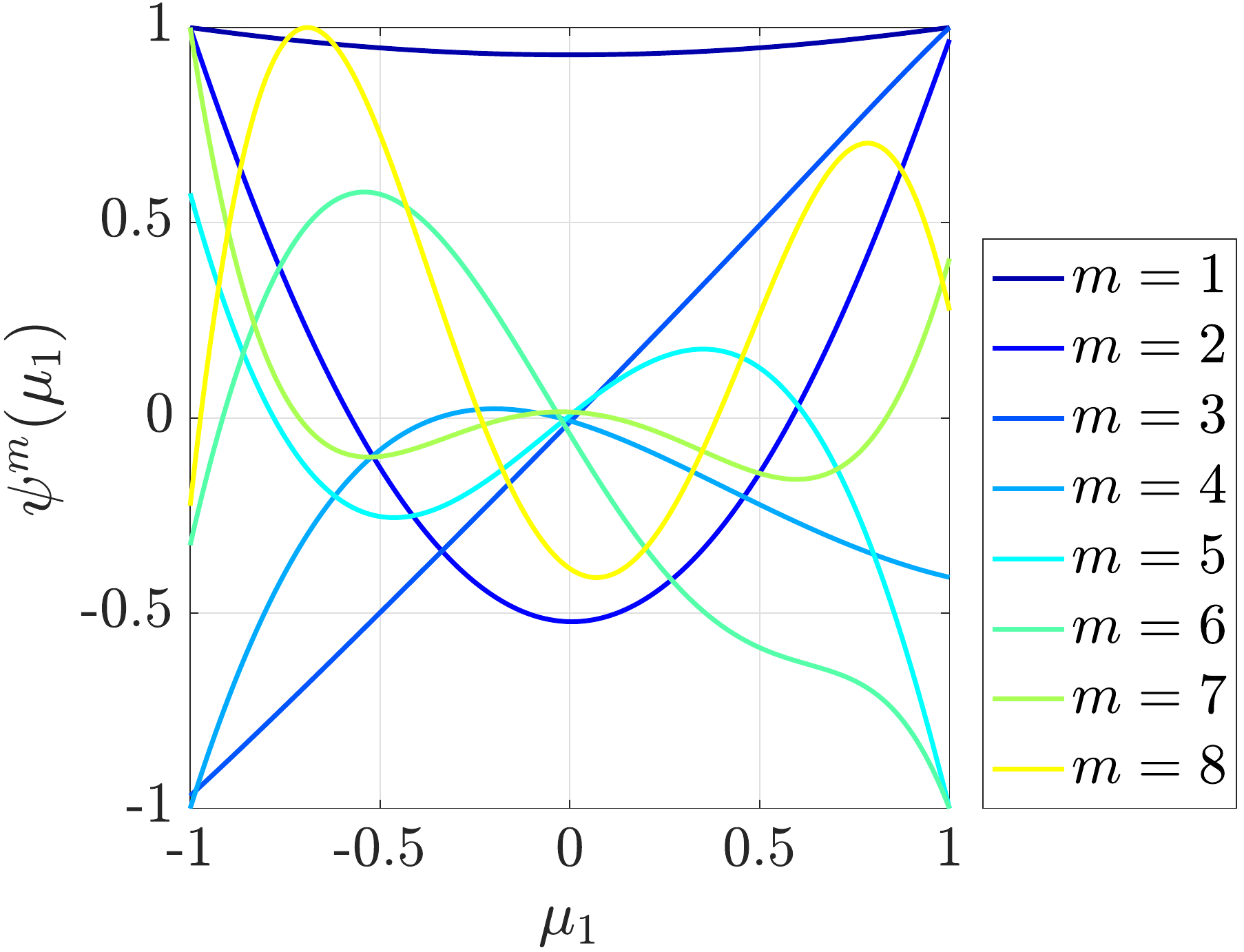}}
	\subfigure[A posteriori PGD]{\includegraphics[width=0.48\textwidth]{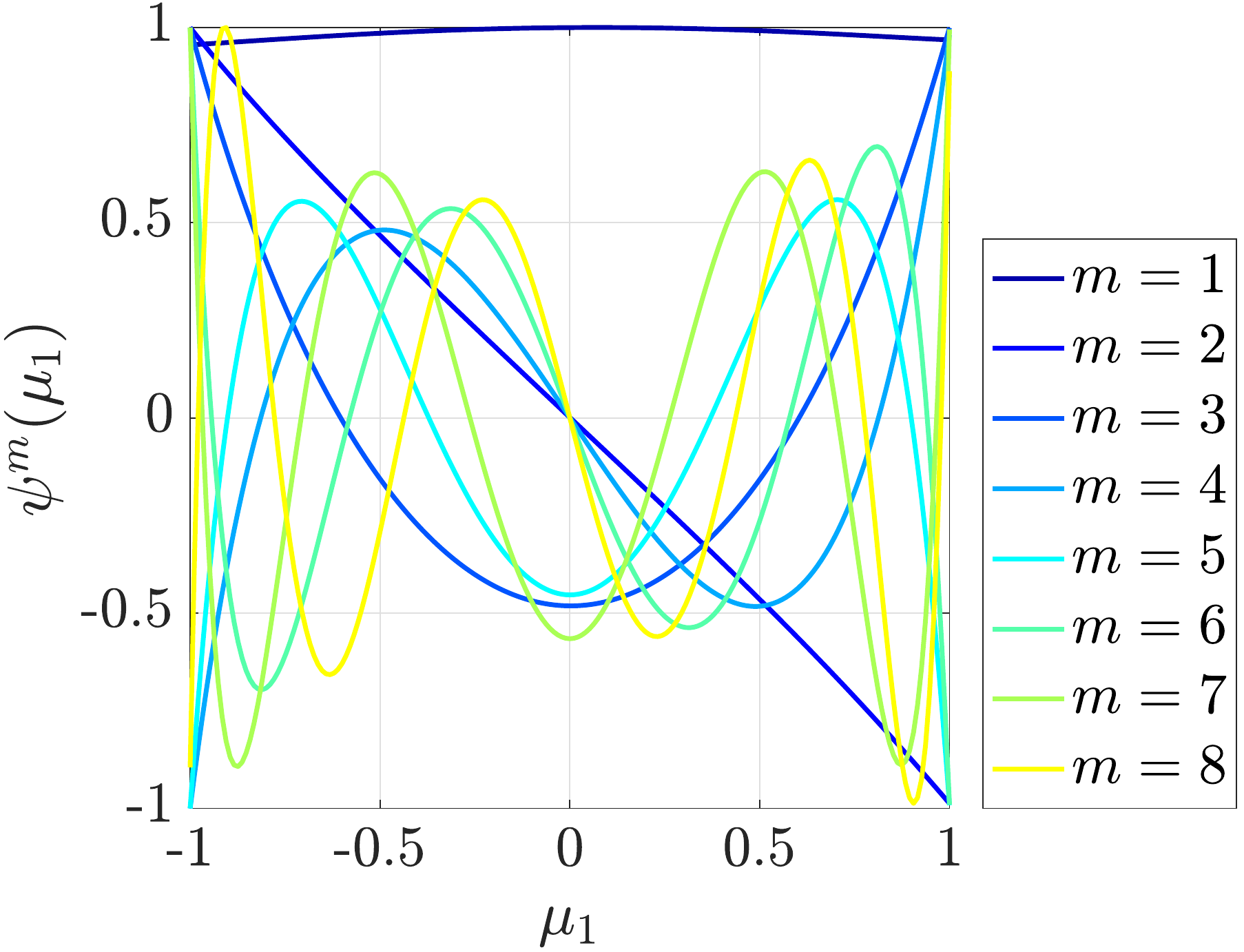}}

	\caption{First eight normalised parametric modes computed using the (a) a priori and (b) a posteriori PGD algorithms for the parametric study of the radius.}
	\label{fig:RadiusModesParam}
\end{figure}

As previously observed in figure~\ref{fig:RadiusUP}, the a posteriori PGD requires seven modes to achieve an error below $10^{-4}$, whereas 20 modes are computed by the a priori PGD algorithm. Hence, the first modes computed using the a posteriori approach capture a larger variability of the solution than the corresponding modes obtained by the a priori PGD. This result is also confirmed by figure~\ref{fig:RadiusModesUpriori} and~\ref{fig:RadiusModesUposteriori} where the first four modes of the module of the velocity computed using the a priori  and the a posteriori PGD , respectively, are reported. It is worth noticing that this behaviour stems from the orthogonality of the modes computed by the a posteriori PGD approach. On the contrary, the modes computed using the a priori PGD feature repeated information which can be eliminated through the above mentioned PGD compression strategy~\cite{DM-MZH:15}.
The corresponding normalised modes for pressure are displayed in figure~\ref{fig:RadiusModesPpriori} and~\ref{fig:RadiusModesPposteriori}. Finally, figure~\ref{fig:RadiusModesParam} shows the first eight parametric modes computed using the two algorithms. It is worth observing that the first two modes feature a comparable global behaviour, whereas, starting from the third one, the parametric functions computed by the a priori and a posteriori PGD approaches noticeably differ. On the one hand, the parametric modes provided by the a posteriori PGD present a regular structure in the interval $\I^1$, reminding the well-known hierarchical basis functions in 1D. As mentioned above, this follows from the orthogonal construction of the modes performed via the high-order PGD projection~\cite{DM-MZH:15}. On the other hand, the modes obtained using the a priori PGD approach display a less regular structure in the parametric domain. Nonetheless, it is expected that such a structure may be retrieved after eliminating redundant information in the PGD approximation by means of the PGD compression~\cite{DM-MZH:15}.

\begin{figure}[!tb]
	\centering
	\subfigure[Module of velocity, $\mu_1=-1$]{\includegraphics[width=0.48\textwidth]{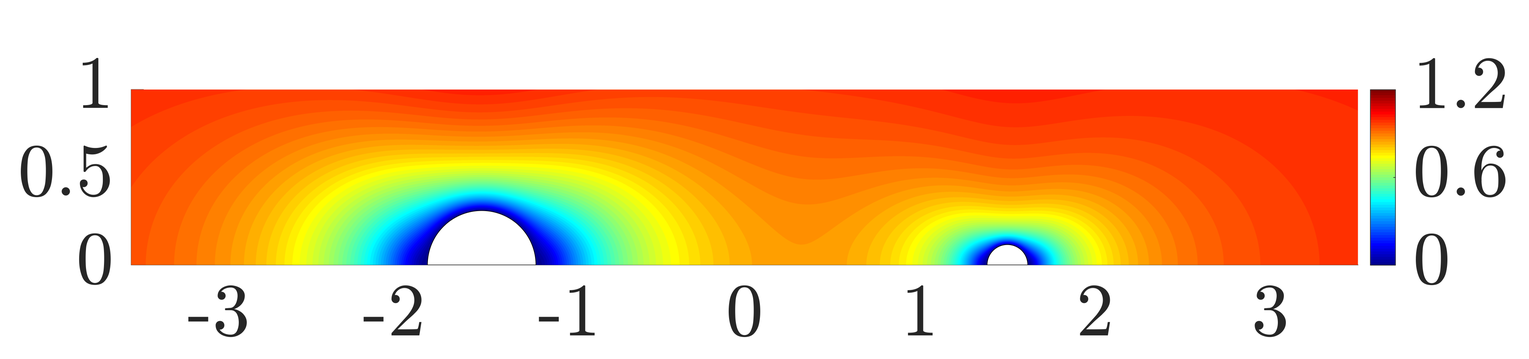}}
	\subfigure[Pressure, $\mu_1=-1$]{\includegraphics[width=0.48\textwidth]{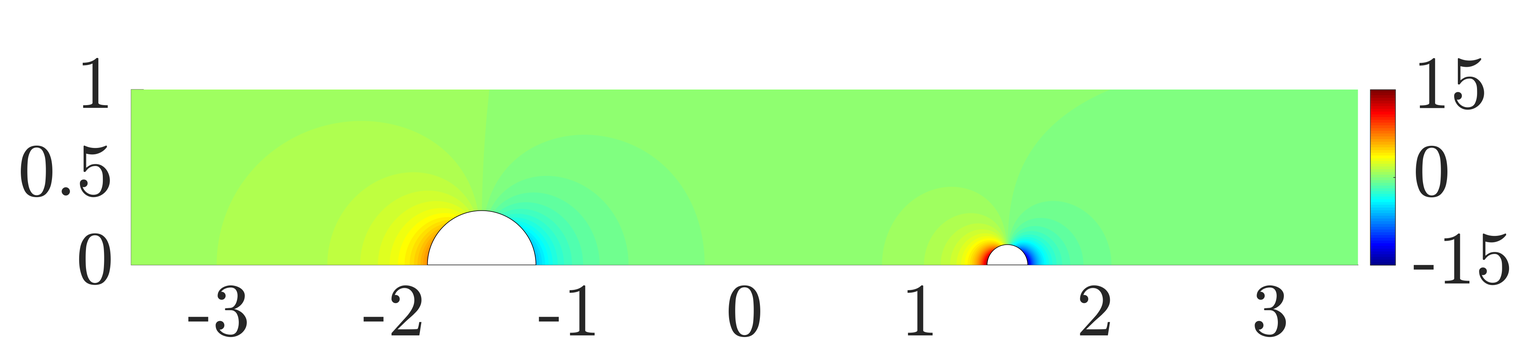}}

	\subfigure[Module of velocity, $\mu_1=0$] {\includegraphics[width=0.48\textwidth]{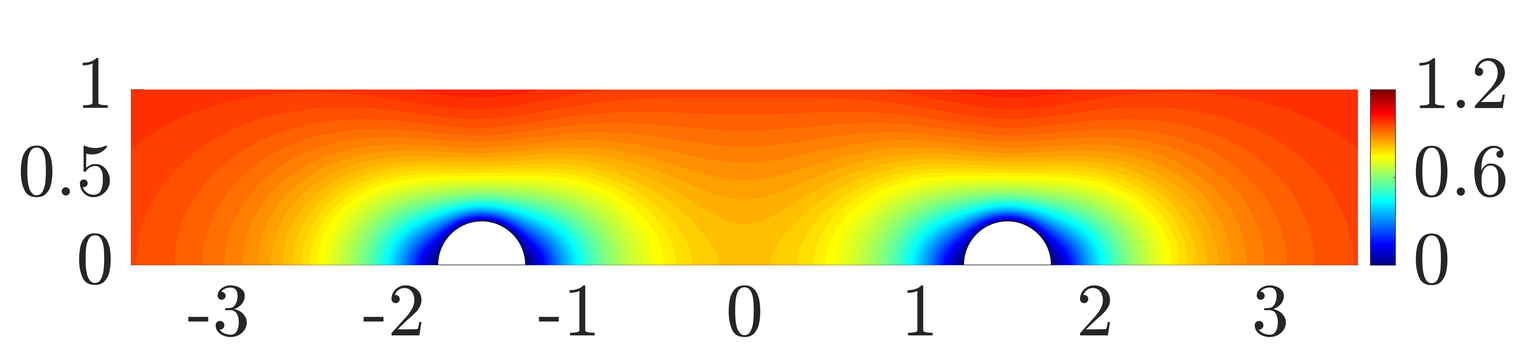}}
	\subfigure[Pressure, $\mu_1=0$] {\includegraphics[width=0.48\textwidth]{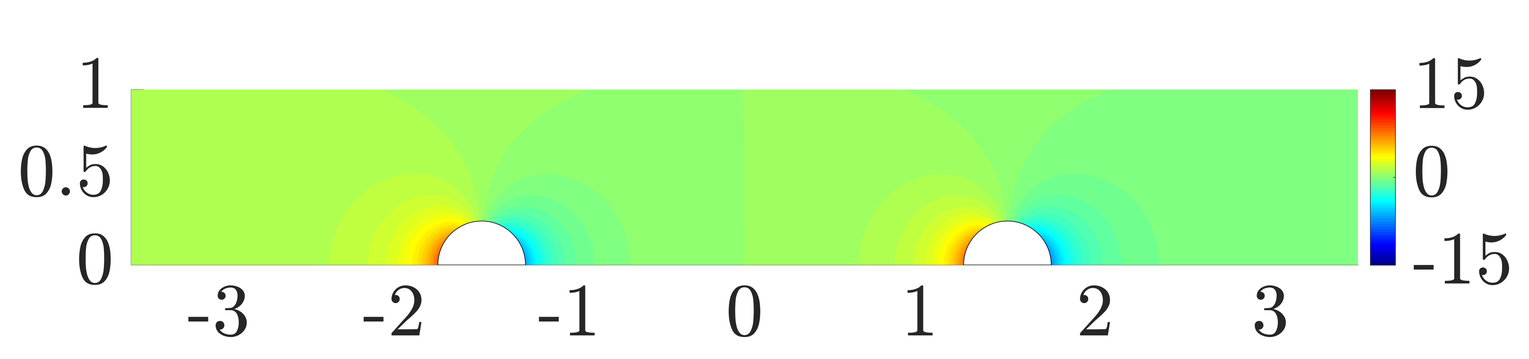}}
	
	\subfigure[Module of velocity, $\mu_1=1$] {\includegraphics[width=0.48\textwidth]{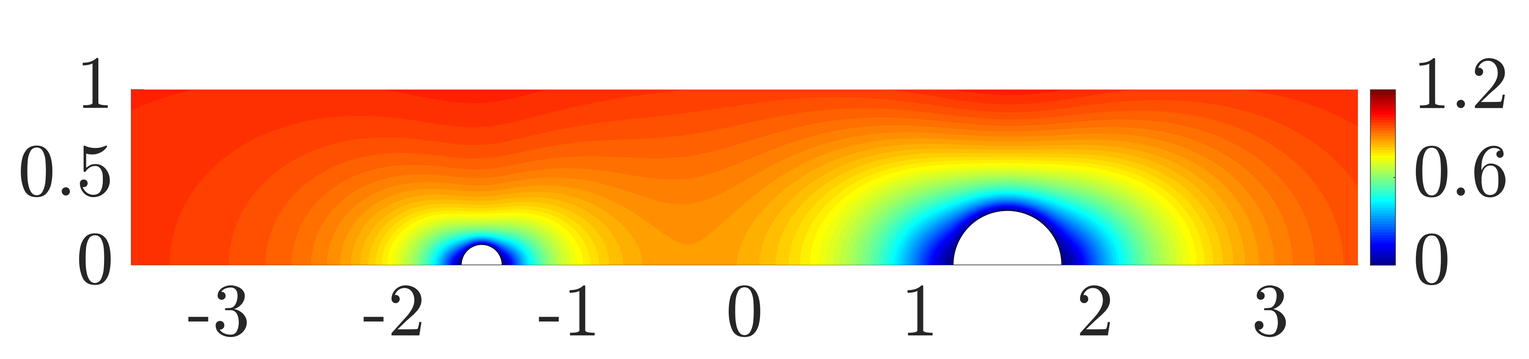}}
	\subfigure[Pressure, $\mu_1=1$] {\includegraphics[width=0.48\textwidth]{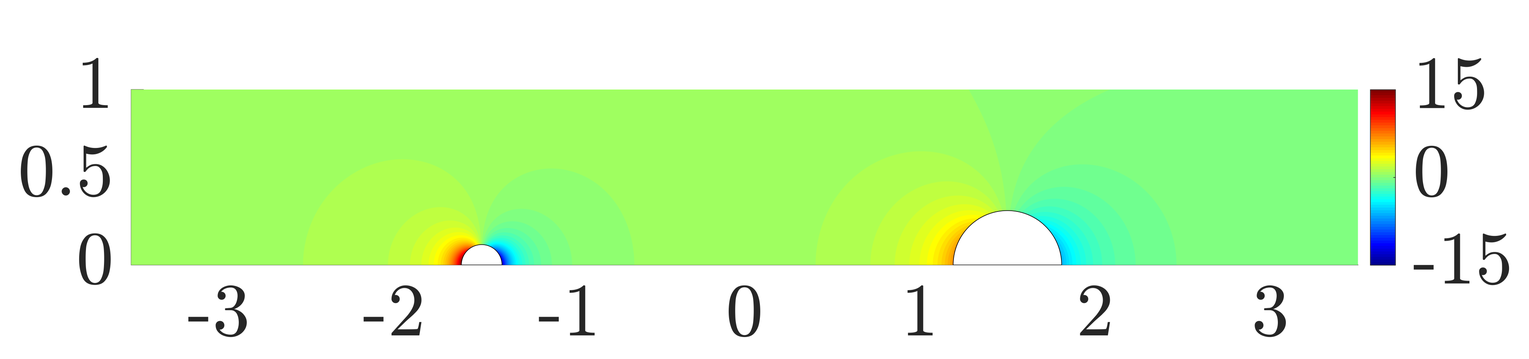}}
	
	\caption{Module of velocity and pressure field for three values of the parameter $\mu_1$ corresponding to the radius of the first sphere being maximum (top), equal to the radius of the second sphere (middle) and minimum (bottom).}
	\label{fig:RadiusOnline}
\end{figure}

Once the modes are computed, it is possible to perform queries in real-time by particularising the generalised velocity and pressure fields for a value of the parameter of interest. As an example, figure~\ref{fig:RadiusOnline} displays the velocity and pressure fields corresponding to three different values of the parameter $\mu_1$. 
Only the solution obtained using the a priori PGD is reported, since the velocity and pressure approximations computed by the two approaches are almost identical as shown in figure~\ref{fig:RadiusUP}. 
Of course, with the a priori algorithm, it is of interest to perform a compression of the PGD solution, see~\cite{DM-MZH:15,PD-DZGH-20}, before the online evaluation, in order to eliminate redundant information and reduce the global number of modes. 

%As an example, figure~\ref{fig:RadiusUOnline} displays the velocity field corresponding to three different values of the parameter $\mu_1$. 
%
%\begin{figure}[!tb]
%	\centering
%	\subfigure[$\mu_1=-1$]{\includegraphics[width=0.8\textwidth]{Param1_Velocity_Mu_M1}}
%	\subfigure[$\mu_1=0$] {\includegraphics[width=0.8\textwidth]{Param1_Velocity_Mu_0}}
%	\subfigure[$\mu_1=1$] {\includegraphics[width=0.8\textwidth]{Param1_Velocity_Mu_1}}
%	\caption{Velocity field for three values of the parameter $\mu_1$.}
%	\label{fig:RadiusUOnline}
%\end{figure}
%%
%Only the solution obtained by the a priori approach is reported since the velocity field obtained by the two approaches is almost identical, as shown in figure~\ref{fig:RadiusU}. 
%
%Similarly, figure~\ref{fig:RadiusPOnline} shows the pressure field corresponding to three different values of the parameter $\mu_1$ introduced above. 
%%
%\begin{figure}[!tb]
%	\centering
%	\subfigure[$\mu_1=-1$]{\includegraphics[width=0.8\textwidth]{Param1_Pressure_Mu_M1}}
%	\subfigure[$\mu_1=0$] {\includegraphics[width=0.8\textwidth]{Param1_Pressure_Mu_0}}
%	\subfigure[$\mu_1=1$] {\includegraphics[width=0.8\textwidth]{Param1_Pressure_Mu_1}}
%	\caption{Pressure field for three values of the parameter $\mu_1$.}
%	\label{fig:RadiusPOnline}
%\end{figure}

\begin{figure}[!tb]
	\centering
	\subfigure[Drag force on each sphere]{\includegraphics[width=0.48\textwidth]{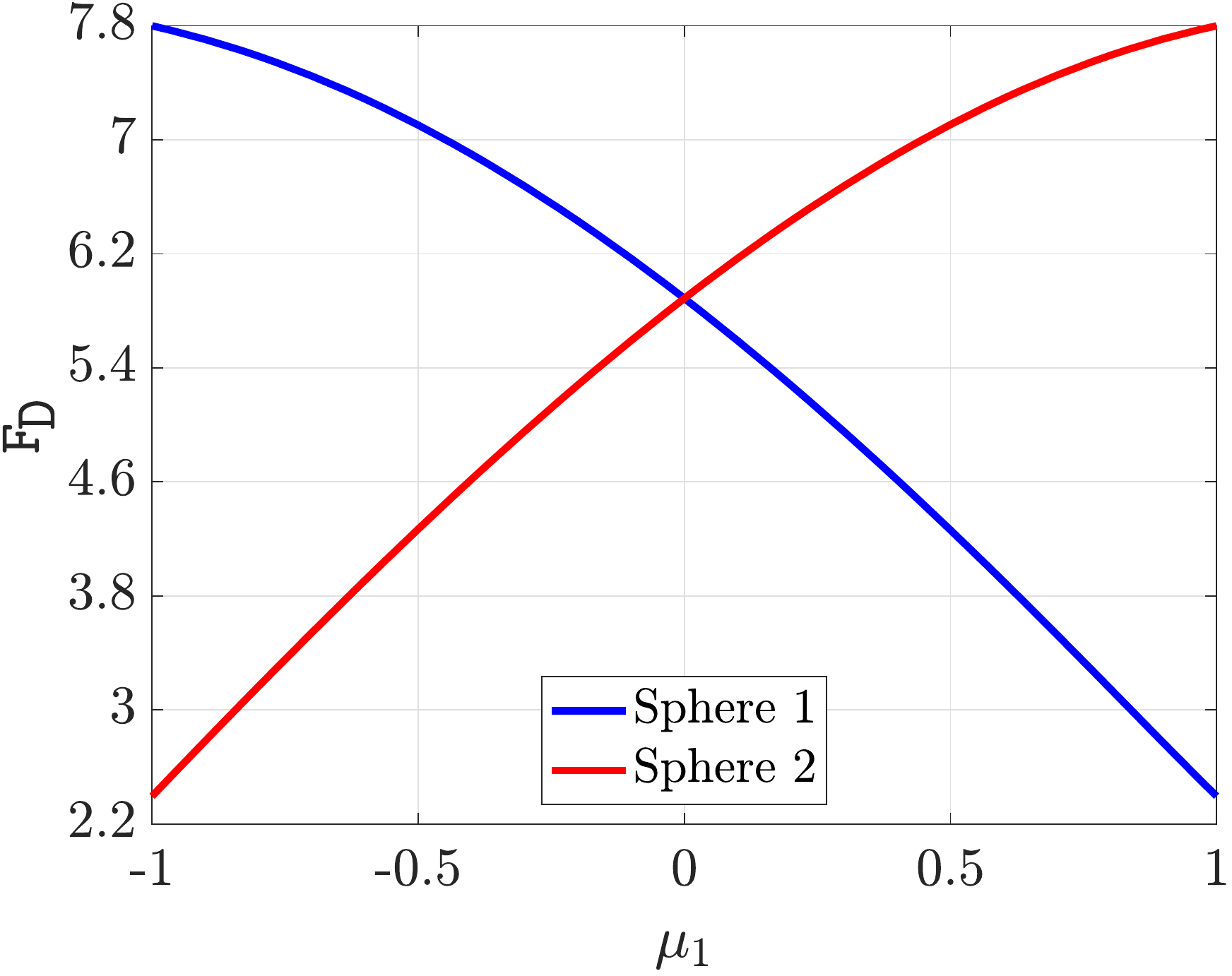}\label{fig:DragRadiusSpheres}}
	\subfigure[Total drag force]{\includegraphics[width=0.48\textwidth]{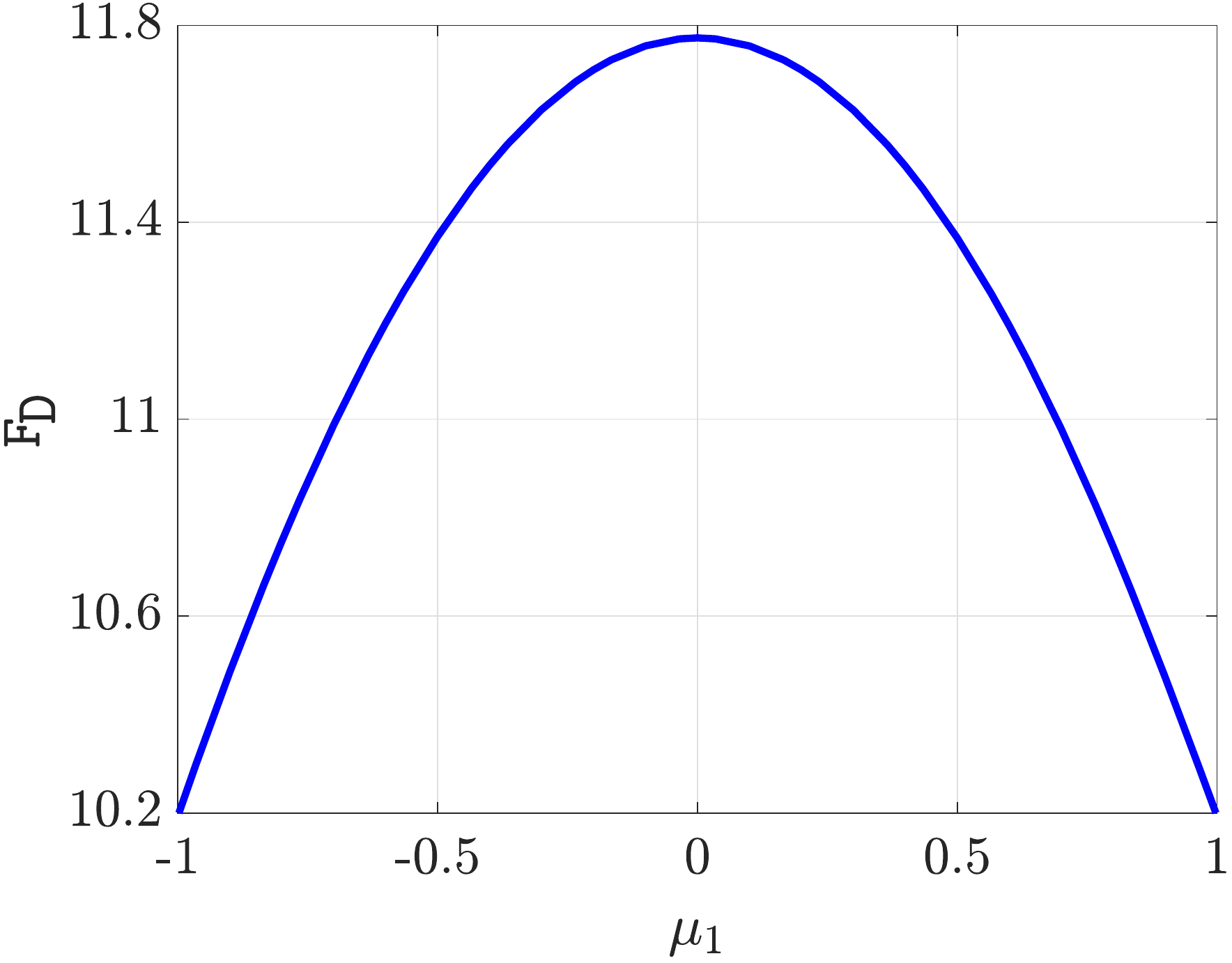}\label{fig:DragRadiusTot}}
	
	\caption{Response surfaces of the drag force as a function of the radius $\mu_1$ of the first sphere.}
	\label{fig:DragRadius}
\end{figure}

In a similar fashion, separated response surfaces for quantities of interest can be devised as explicit functions of the parameter. Figure~\ref{fig:DragRadius} reports the response surface of the drag force as a function of the radius of the two spherical bladders, computed using the a priori PGD. As expected, the drag is maximum on the first sphere when its radius is maximum (i.e. $\mu_1 {=} {-} 1$) and it decreases monotonically until reaching the configuration of minimum radius for $\mu_1 {=} 1$. An analogous behaviour is observed for the second sphere, with the drag force spanning from its minimum value at $\mu_1 {=} {-} 1$ to its maximum at $\mu_1 {=} 1$. Moreover, the forces on the two spheres are equal for the geometric configuration of $\mu_1 {=} 0$, that is, when the two bladders have the same volume (Fig.~\ref{fig:DragRadiusSpheres}). For the sake of completeness, figure~\ref{fig:DragRadiusTot} displays the total drag on the two spheres as a function of the geometric parameter $\mu_1$.
The results obtained with the a posteriori PGD are qualitatively and quantitatively similar, whence they are omitted for the sake of brevity. A detailed comparison of the accuracy of the two approaches is presented in section~\ref{sc:comparisonAccuracy}.

%-------------------------------------------------------
\subsubsection{Varying the distance between the spherical bladders} 
\label{sc:exMu2}
%-------------------------------------------------------

The second example considers a geometrically parametrised problem, where the parameter controls the distance between two equal spherical bladders with radius $0.25$. It is worth recalling that the reference geometry in figure~\ref{fig:swimmerGeometry} is characterised by two equal spheres of radius $\Rref {=} 0.116$. Hence, for the cases studied in this section, the geometric mapping accounts for both a parameter-dependent variation of the distance between the two bladders and an expansion of the spheres, independent of the parameter.

Two intervals $\I^2$ are considered to analyse the sensitivity of the PGD solutions to the range of variations of the parameter. 

\begin{figure}[!tb]
	\centering
	\subfigure[$\bu$] {\includegraphics[width=0.49\textwidth]{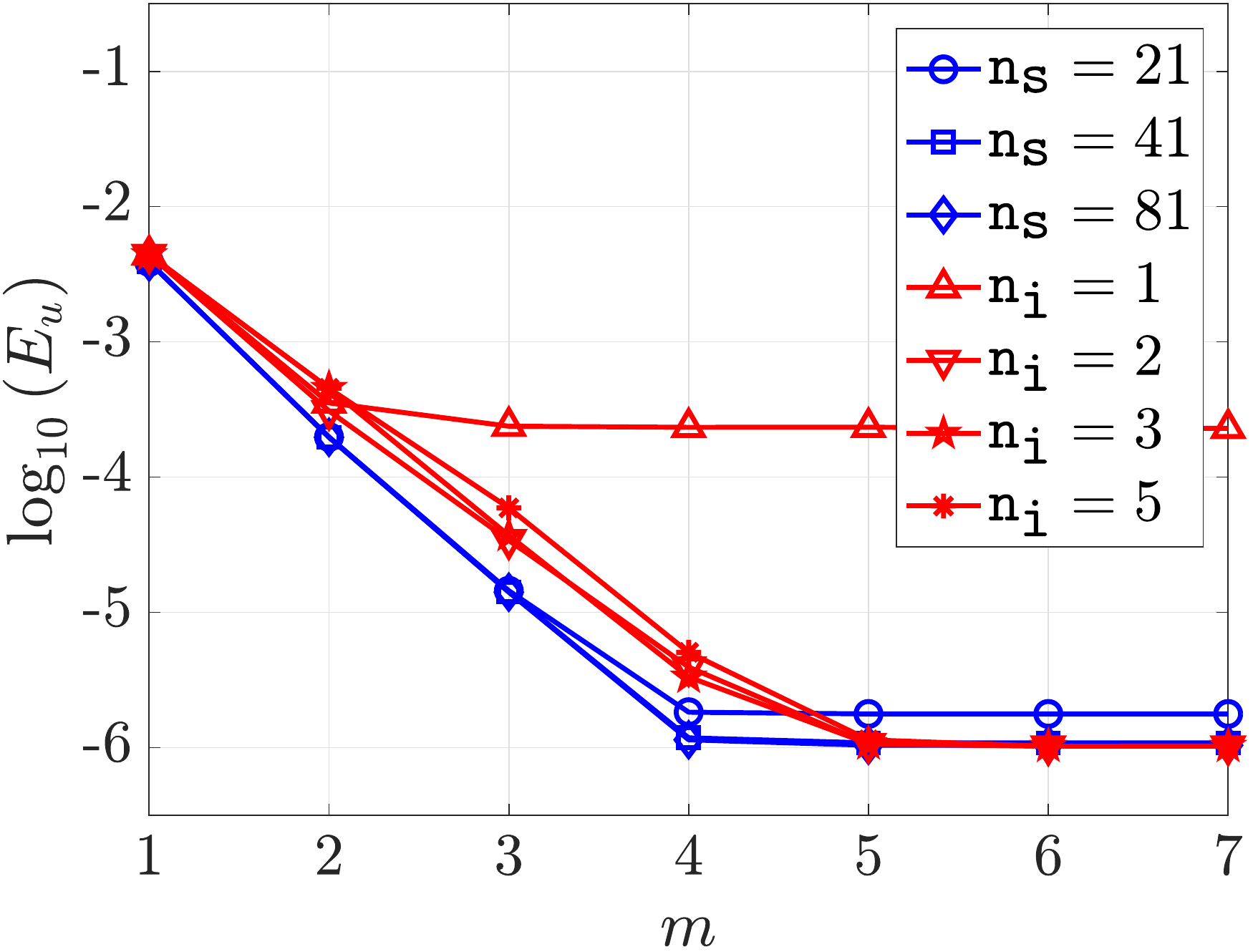}}
	\subfigure[$p$]   {\includegraphics[width=0.49\textwidth]{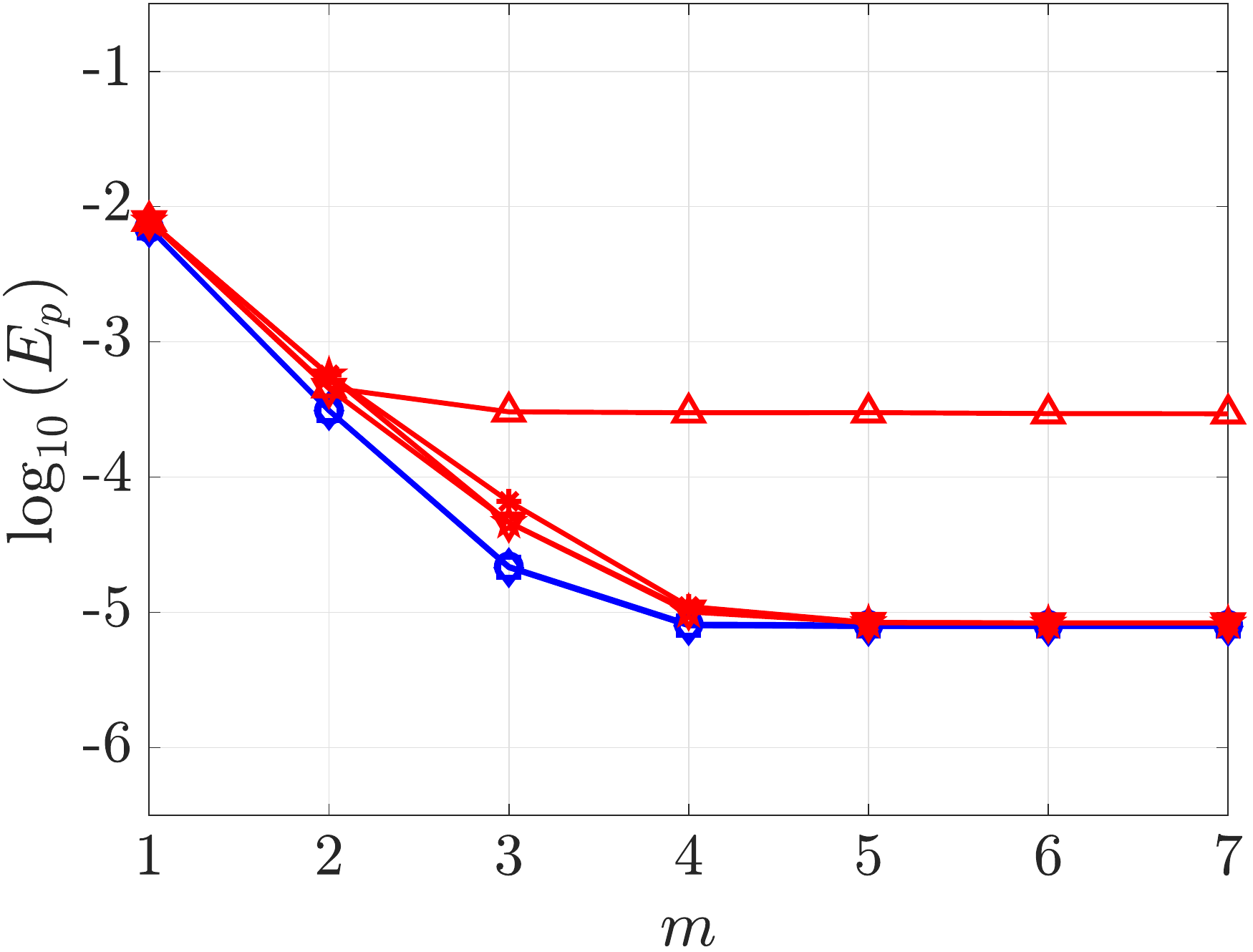}}
	\subfigure[$\bL$] {\includegraphics[width=0.49\textwidth]{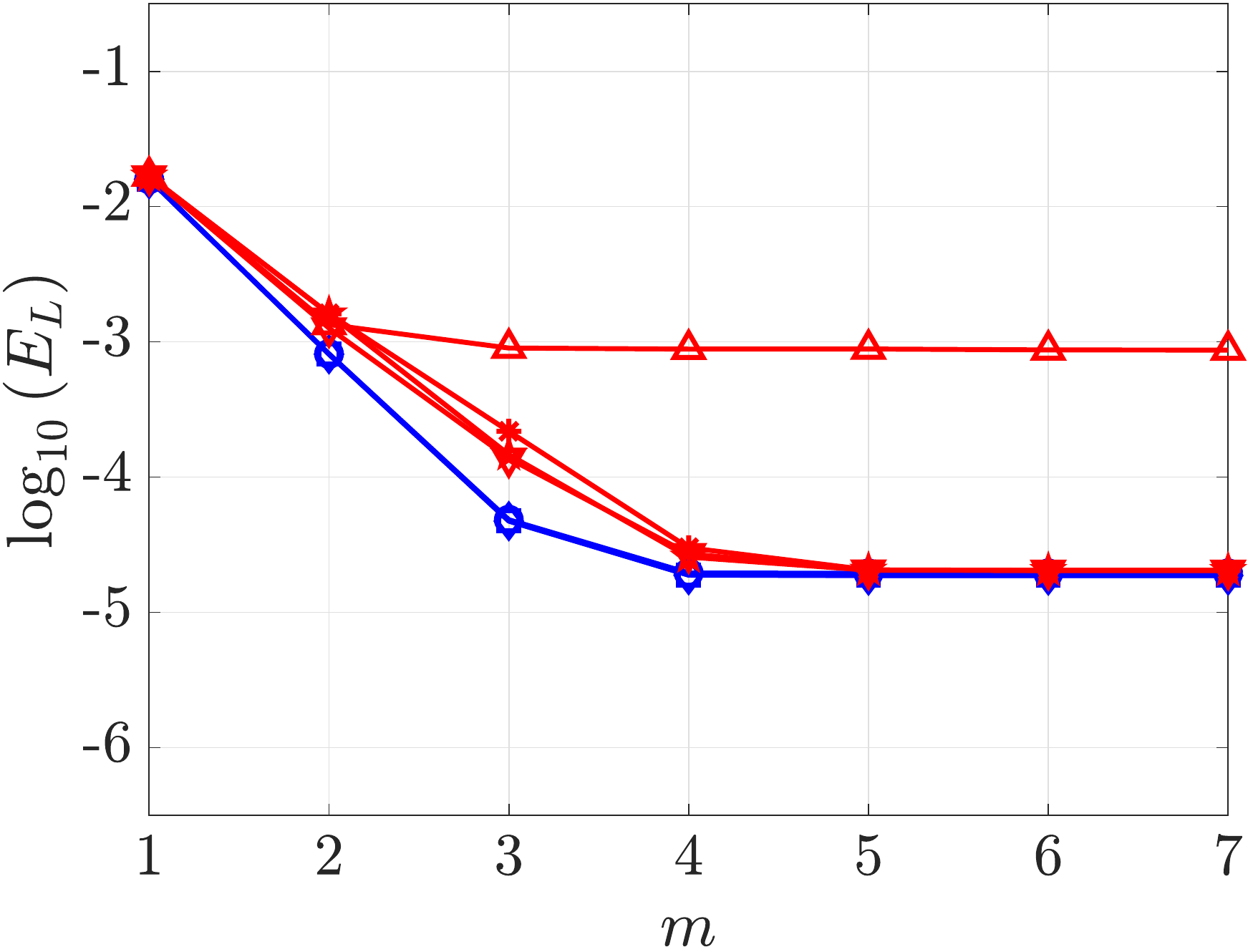}}
	\subfigure[$\FD$]   {\includegraphics[width=0.49\textwidth]{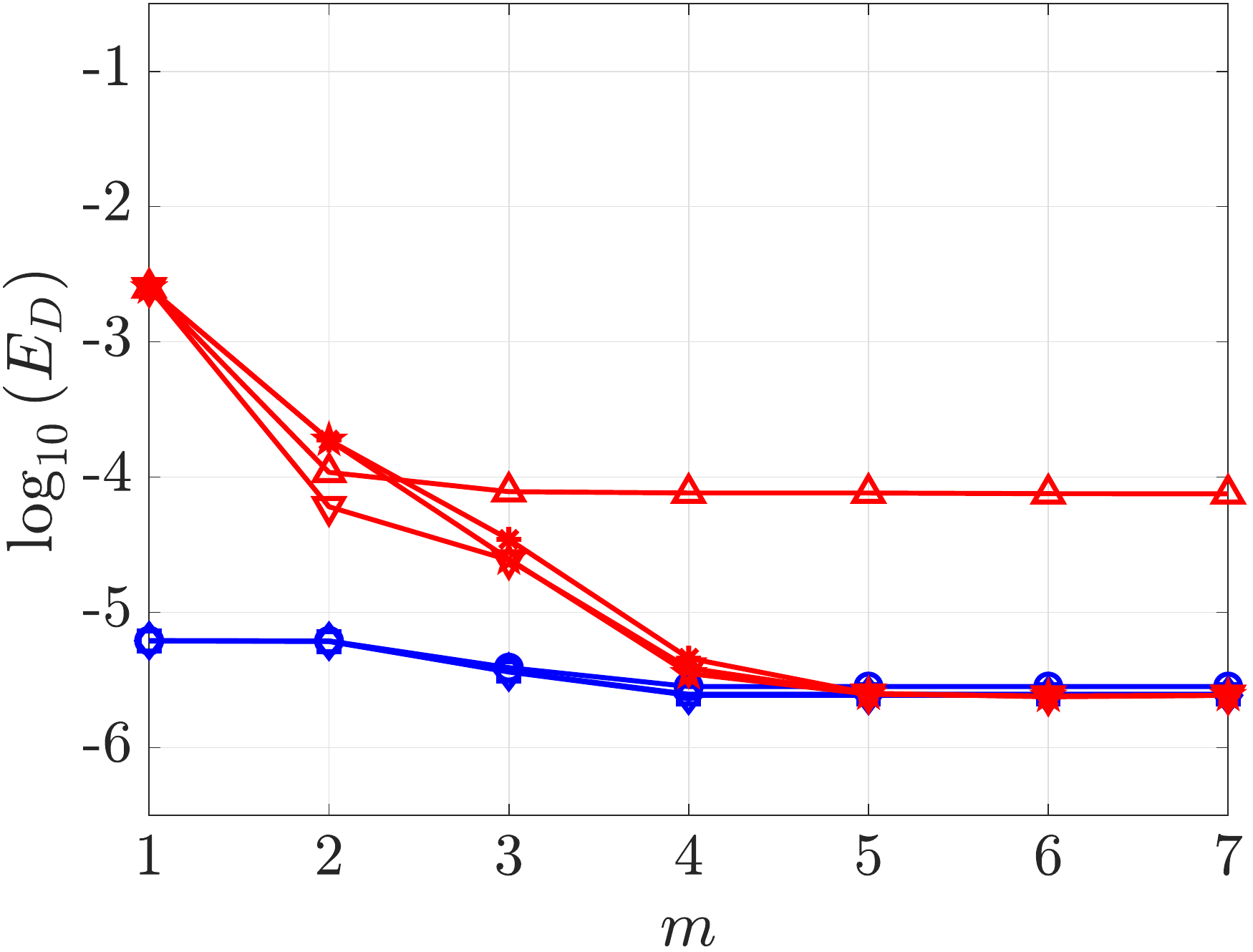}}
	\caption{Evolution of (a-c) the $\eltwo(\Omega \times \I^2)$ error for velocity, pressure and gradient of velocity and (d) the $\eltwo(\I^2)$ error for the drag force as a function of the number of PGD modes for the problem with one geometric parameter controlling the distance between the spherical bladders and $\I^2 {=} [-2,-1]$. The legend details the number $\nsnap$ of snapshots used by the a posteriori PGD approach (blue) and the number $\niter$ of nonlinear AD iterations used by the a priori PGD approach (red).}
	\label{fig:Distance1}
\end{figure}

The first interval is taken as $\I^2 {=} [-2,-1]$ and figure~\ref{fig:Distance1} reports the evolution of the $\eltwo(\Omega \times \I^2)$ error for velocity, pressure and gradient of velocity and the $\eltwo(\I^2)$ error for the drag force as a function of the number $m$ of modes. 
The results show that with only four modes, the a posteriori PGD approach is able to produce the most accurate results for all the variables, including the drag force computed from the pressure and the gradient of velocity. It is worth noticing that in this example the accuracy of the a posteriori approach in the drag force does not improve when increasing the number of snapshots and 21 snapshots are sufficient to provide a drag force with an error below $10^{-5}$. For the a priori approach, five modes computed with two AD iterations are required to obtain the maximum accuracy in all the variables. With one iteration, the error in the drag force is more than one order of magnitude higher than the one obtained with two iterations. Furthermore, a higher number of iterations does not produce any gain in accuracy despite the increased computational cost.  
In this case, the two PGD approaches show similar performance as the a priori algorithm provides an error in the drag force below $10^{-5}$ with 12 solutions of the spatial problem (i.e. four modes, each computed with two iterations of the AD scheme plus the initial solve to perform the prediction of the mode, see algorithm~\ref{alg:PGDpriori}), whereas the a posteriori approach requires 21 snapshots for a similar level of accuracy.

\begin{figure}[!tb]

	\subfigure[Mesh, $\mu_2=-3$]{\includegraphics[width=0.94\textwidth]{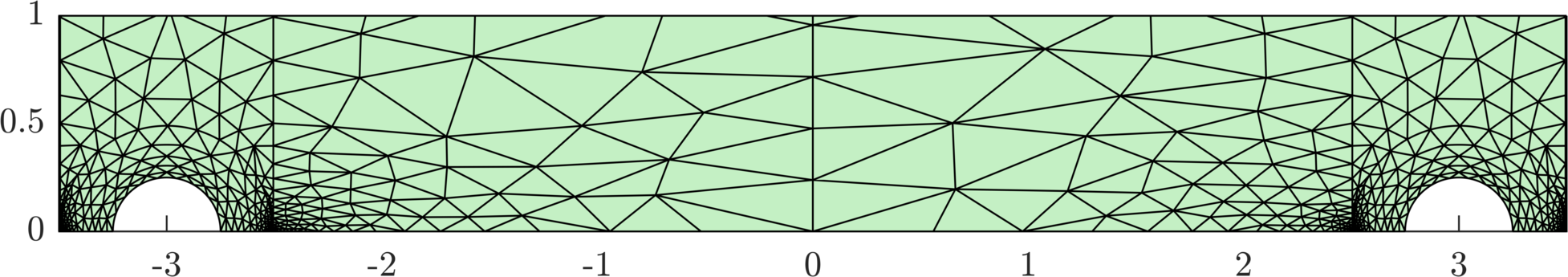}}
	
	{\centering		
	\subfigure[Quality, $\mu_2=-3$]{\includegraphics[width=\textwidth]{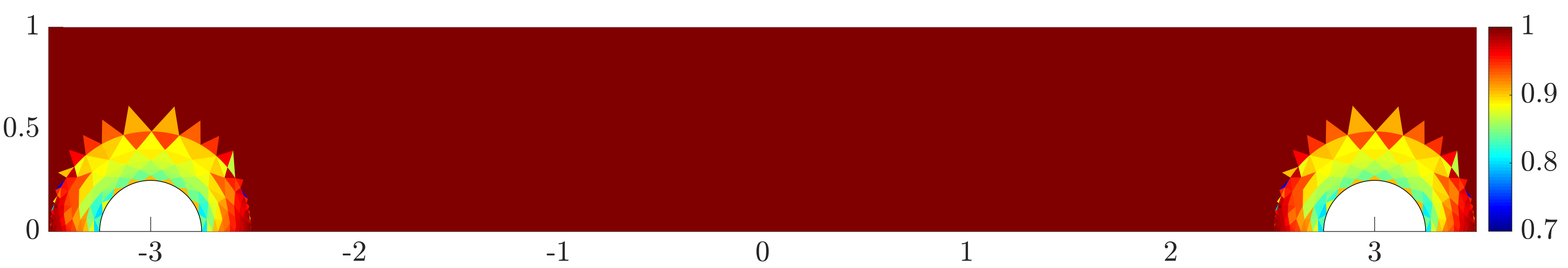}}
	}
	
	\subfigure[Mesh, $\mu_2= 2$]{\includegraphics[width=0.94\textwidth]{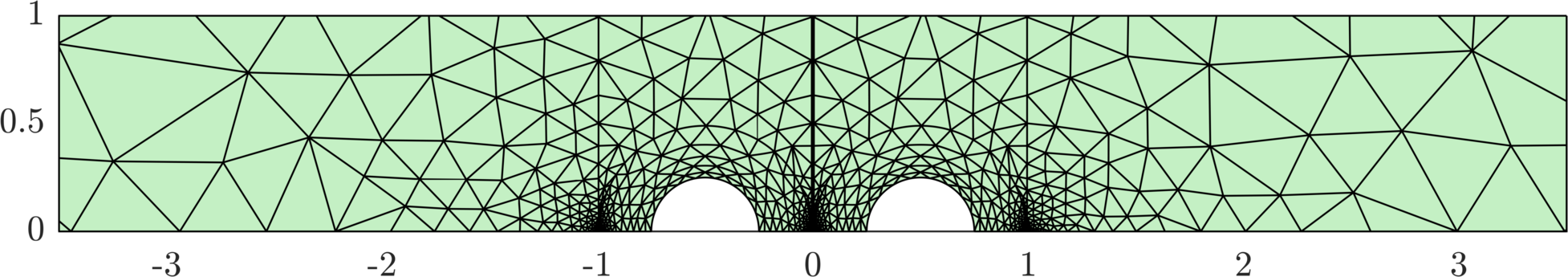}}
	
	{\centering	
	\subfigure[Quality, $\mu_2= 2$]{\includegraphics[width=\textwidth]{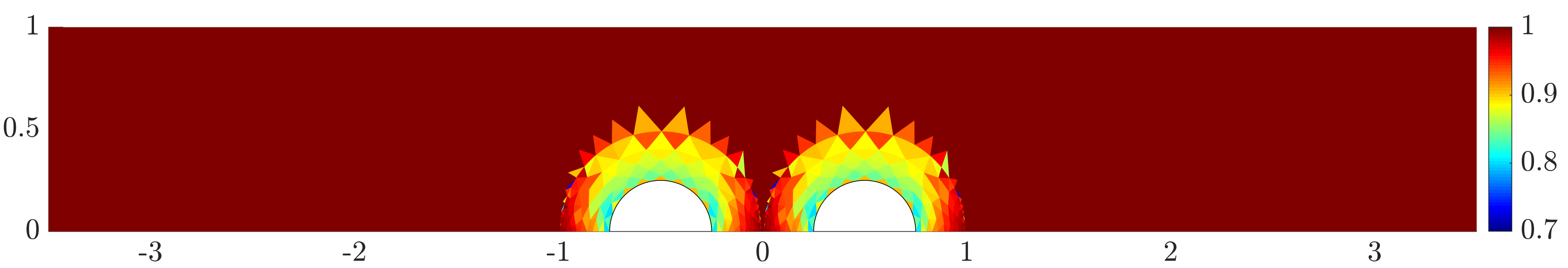}}
	}
	\caption{(a-c) Mesh configurations and (b-d) mesh quality of two deformed swimmers for the mapping with $\mu_2$ as a geometric parameter.}
	\label{fig:swimmerMeshMappingMu2}
\end{figure}

Second, the parametric interval is extended to $\I^2 {=} [-3,2]$. Figure~\ref{fig:swimmerMeshMappingMu2} presents the meshes of two microswimmer configurations obtained from the extreme values of the parameter $\mu_2$ describing the distance between the bladders. The figure also displays the mesh quality, measured as the scaled Jacobian of the isoparametric mapping, of the two deformed configurations.
The results report that the mesh quality is not affected by the mapping considered as the change in distance is piecewise linear and the lower mesh quality only concentrates in the vicinity of the spheres. This is due to the deformation required to transform the reference mesh with radius 0.116 into the geometric configuration under analysis, associated with the bladders of equal volume, in which the radius achieves the value 0.25. As previously observed for the case of the parametrised radius, only few elements present a mesh quality of 0.7 whereas most elements feature a scaled Jacobian of 0.9 or higher. Hence, the mesh in figure~\ref{fig:swimmerMesh} also provides a good approximation for the parametric study of the distance between the bladders.

\begin{figure}[!tb]
	\centering
	\subfigure[$\bu$] {\includegraphics[width=0.49\textwidth]{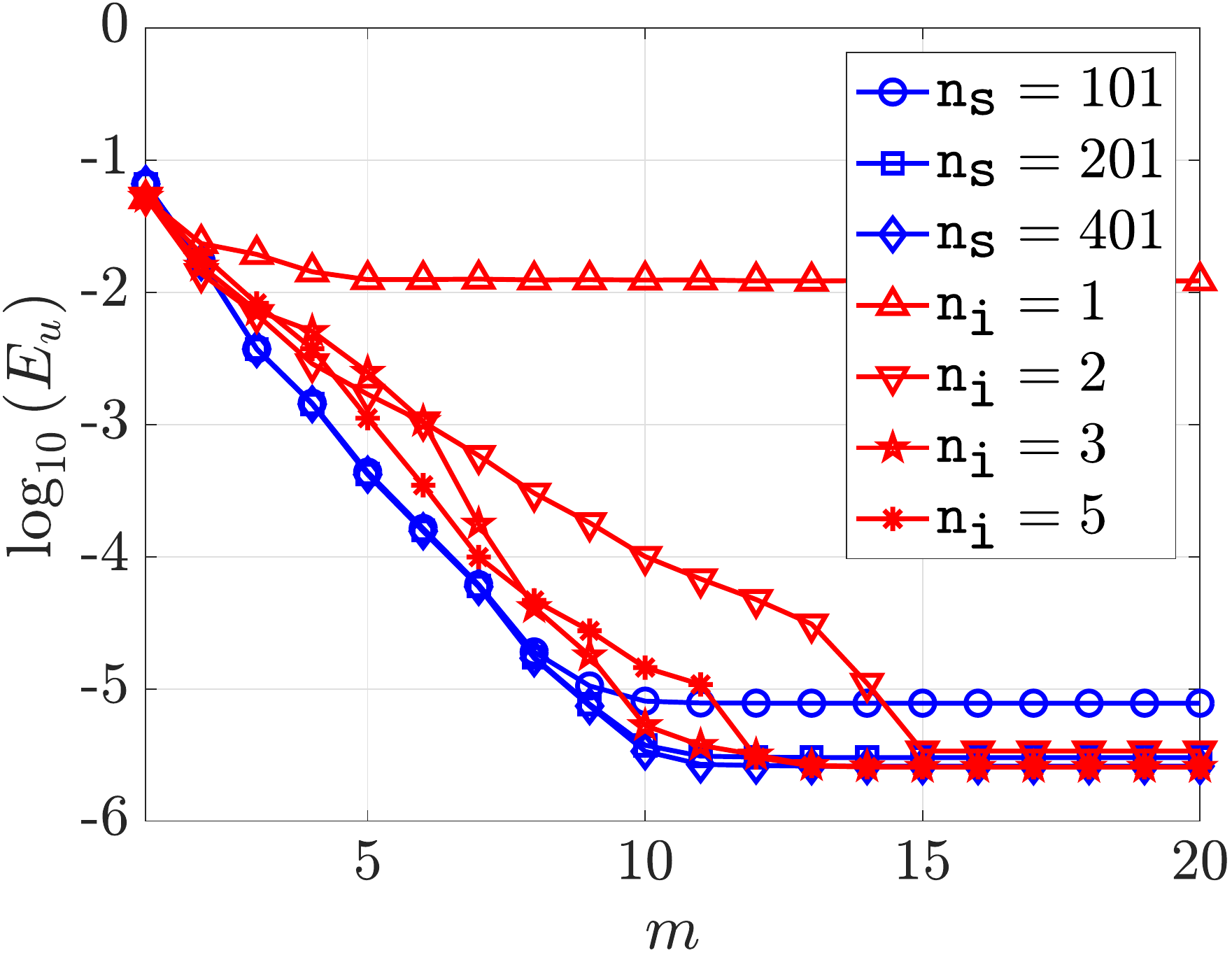}}
	\subfigure[$p$]   {\includegraphics[width=0.49\textwidth]{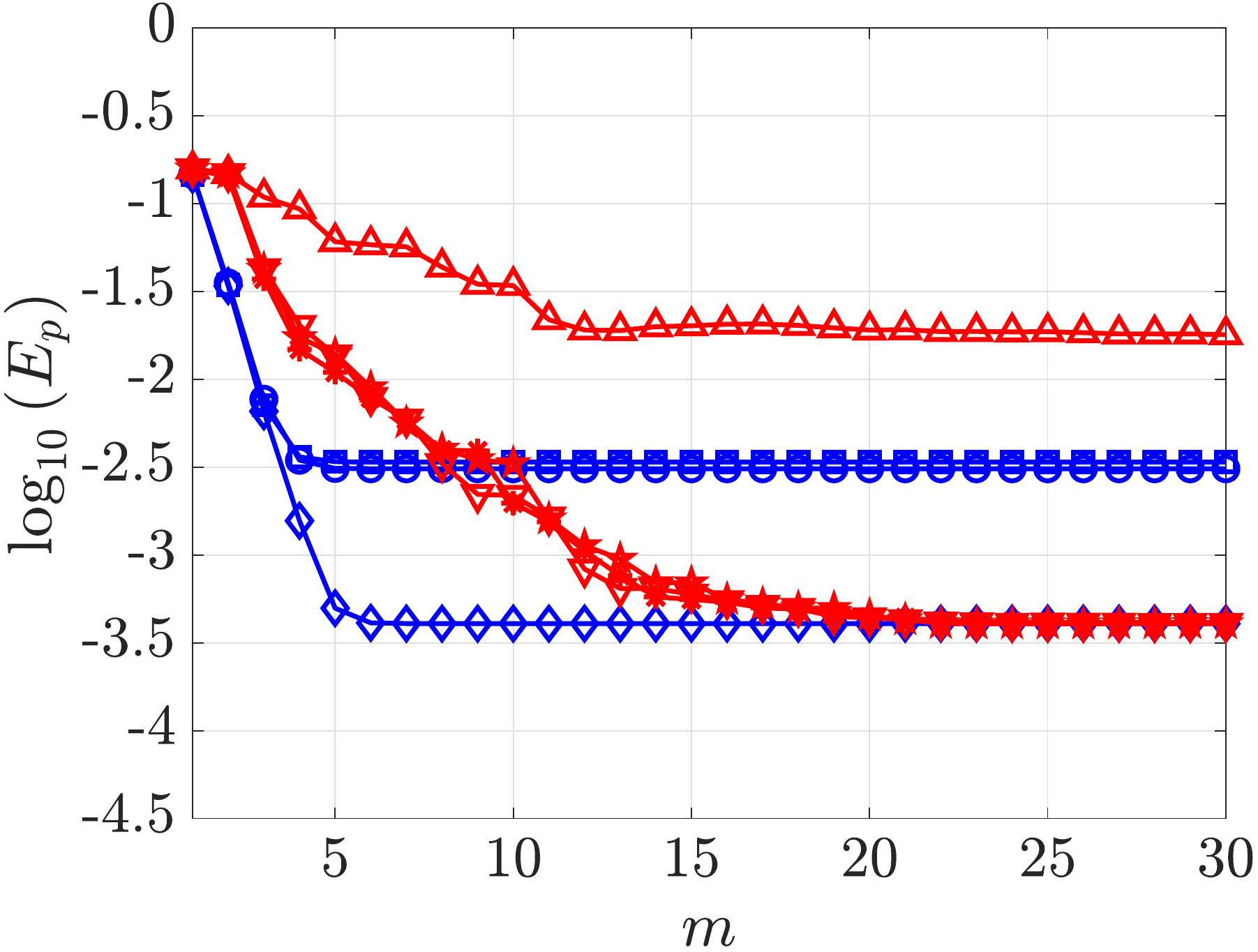}}
	\subfigure[$\bL$] {\includegraphics[width=0.49\textwidth]{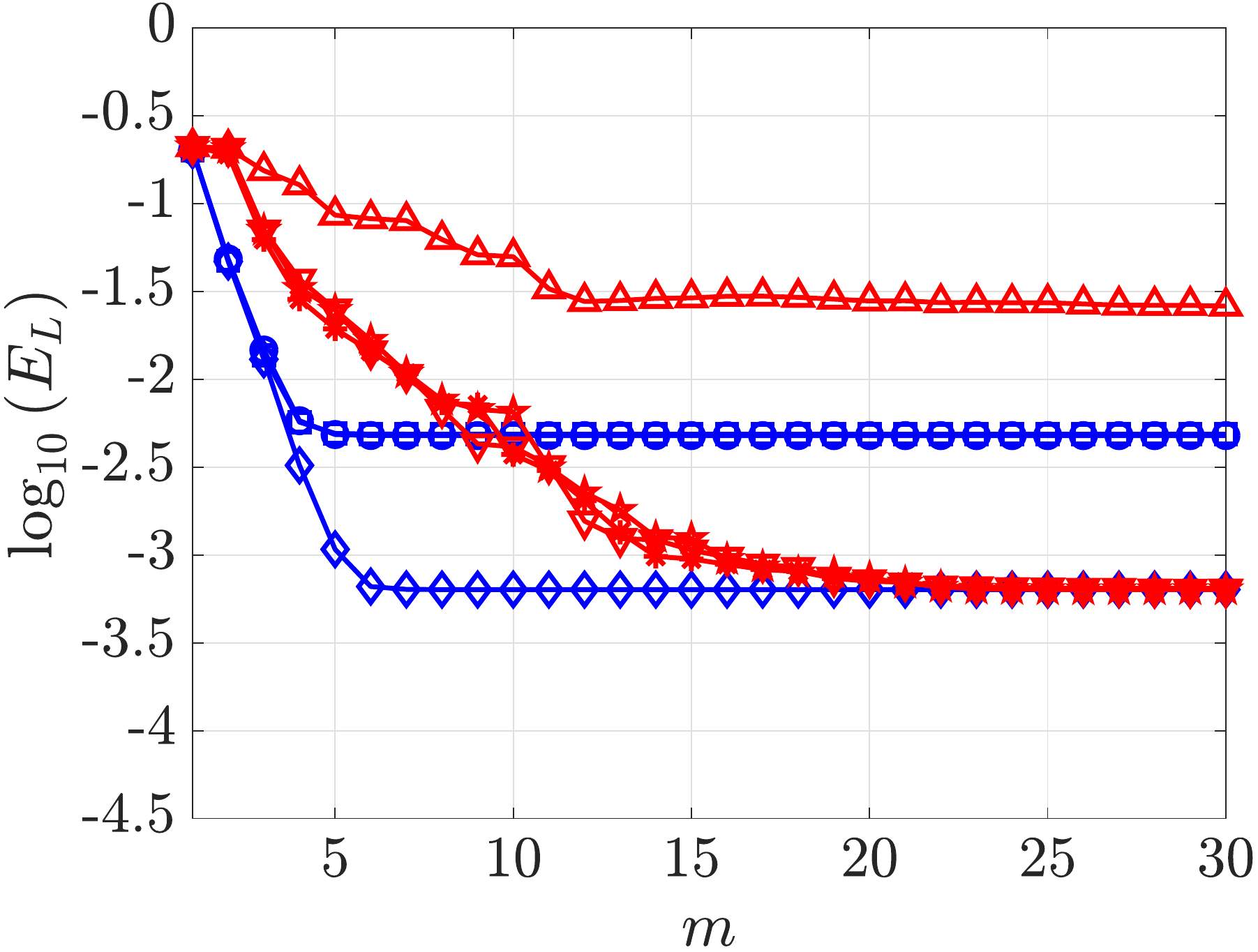}}
	\subfigure[$\FD$]   {\includegraphics[width=0.49\textwidth]{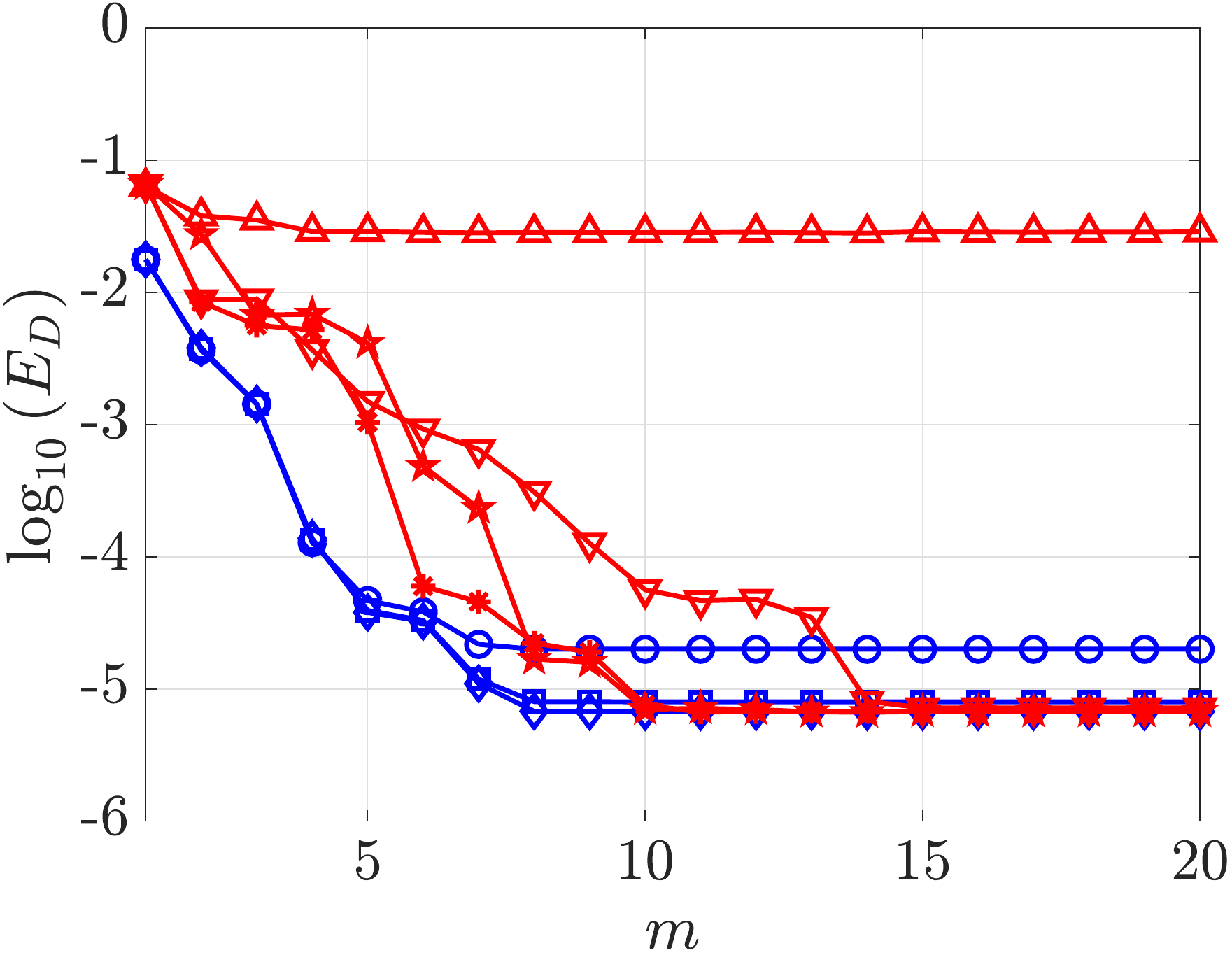}}
	\caption{Evolution of (a-c) the $\eltwo(\Omega \times \I^2)$ error for velocity, pressure and gradient of velocity and (d) the $\eltwo(\I^2)$ error for the drag force as a function of the number of PGD modes for the problem with one geometric parameter controlling the distance between the spherical bladders and $\I^2 {=} [-3,2]$. The legend details the number $\nsnap$ of snapshots used by the a posteriori PGD approach (blue) and the number $\niter$ of nonlinear AD iterations used by the a priori PGD approach (red).}
	\label{fig:Distance2}
\end{figure}

Figure~\ref{fig:Distance2} shows the evolution of the $\eltwo(\Omega \times \I^2)$ error for velocity, pressure and gradient of velocity and the $\eltwo(\I^2)$ error for the drag force as a function of the number $m$ of modes. 
It is worth empasising that a simple visual comparison of the results in figures~\ref{fig:Distance1} and~\ref{fig:Distance2} clearly illustrates the challenge that a larger interval of variation of the geometric parameter induces for both PGD approaches.

The results show that the a posteriori approach requires 10 modes in order to reach the maximum accuracy for velocity, pressure and gradient of velocity. In addition, it can be observed that the a posteriori algorithm requires eight modes and 201 snapshots to provide an error in the drag force below $10^{-5}$. A higher number of snapshots does not lead to a further reduction in the error, whereas a lower number of snapshots, 101, is responsible for a slight increase in the error. 
Concerning the a priori PGD, the AD scheme with one iteration leads to a stagnated error that is several orders of magnitude higher than the one obtained with two or more iterations. For two iterations, the number of modes required to reach the maximum accuracy is 15 and for three or five iterations the number of modes required varies between 10 and 12. When the drag force is considered, the a priori approach shows that an accuracy below $10^{-5}$ can be obtained with two iterations and 14 modes, three iterations and 10 modes or five iterations and 10 modes. The most efficient alternative thus consists of computing 10 modes with three AD iterations for a total of 40 spatial solves, requiring a marginally lower cost than the computation of 14 modes with two AD iterations, that is, 42 calls to the HDG spatial solver.

Comparing the performance of the a priori and a posteriori approaches for this more challenging problem, it is clear that the a priori approach is capable of achieving the same accuracy as the a posteriori approach with a significant lower computational cost. For instance, to reach an accuracy in the drag force below $10^{-5}$, the a priori approach requires 40 solutions of the spatial problem whereas the same level of accuracy cannot be reached by the a posteriori approach with 101 snapshots. In this case the a posteriori approach requires 201 snapshots, which is five times more than the a priori method.
Although this may seem a clear disadvantage of the a posteriori PGD, it is worth noticing that the 201 snapshots could be easily computed in parallel. On the contrary, the a priori PGD solves the spatial problems following a sequential approach. Hence, the resulting computational time of the a posteriori PGD may still be competitive despite the higher number of full-order HDG solves required. In addition, recall that the number of snapshots required by the a posteriori PGD algorithms strongly depends on the sampling strategy employed. Several advanced sampling methods proposed in the literature, see section~\ref{sc:intro}, can be utilised to reduce the number of snapshots and to improve the performance of the a posteriori PGD scheme. As previously mentioned, this is out of the scope of the current work: in order to perform an unbiased comparison of a priori and a posteriori PGD strategies, two versions of the algorithms exploiting neither prior information nor tailored improvements such as advanced sampling and error control techniques, are considered.

%As in the previous example, it is possible to construct separated response surfaces and to perform queries in real-time by particularising the solution for a value of the parameter of interest. Figures~\ref{fig:RadiusUOnline} and~\ref{fig:Distance2POnline} show the velocity and pressure fields, respectively, corresponding to two different values of the parameter $\mu_2$. 
%%
%\begin{figure}[!tb]
%	\centering
%	\subfigure[$\mu_2=-3$]{\includegraphics[width=0.8\textwidth]{Param2_Int2_Velocity_Mu_M3}}
%	\subfigure[$\mu_2=2$] {\includegraphics[width=0.8\textwidth]{Param2_Int2_Velocity_Mu_2}}
%	\caption{Velocity field for two values of the parameter $\mu_2$.}
%	\label{fig:Distance2UOnline}
%\end{figure}
%%
%\begin{figure}[!tb]
%	\centering
%	\subfigure[$\mu_2=-3$]{\includegraphics[width=0.8\textwidth]{Param2_Int2_Pressure_Mu_M3}}
%	\subfigure[$\mu_2=2$] {\includegraphics[width=0.8\textwidth]{Param2_Int2_Pressure_Mu_2}}
%	\caption{Pressure field for two values of the parameter $\mu_2$.}
%	\label{fig:Distance2POnline}
%\end{figure}
%%
%The results illustrate the increased difficulty of this problem as the flow pattern substantially changes as the distance between the two spherical bladders decrease.

\begin{figure}[!tb]
	\centering
	\subfigure[$m = 1$]{\includegraphics[width=0.48\textwidth]{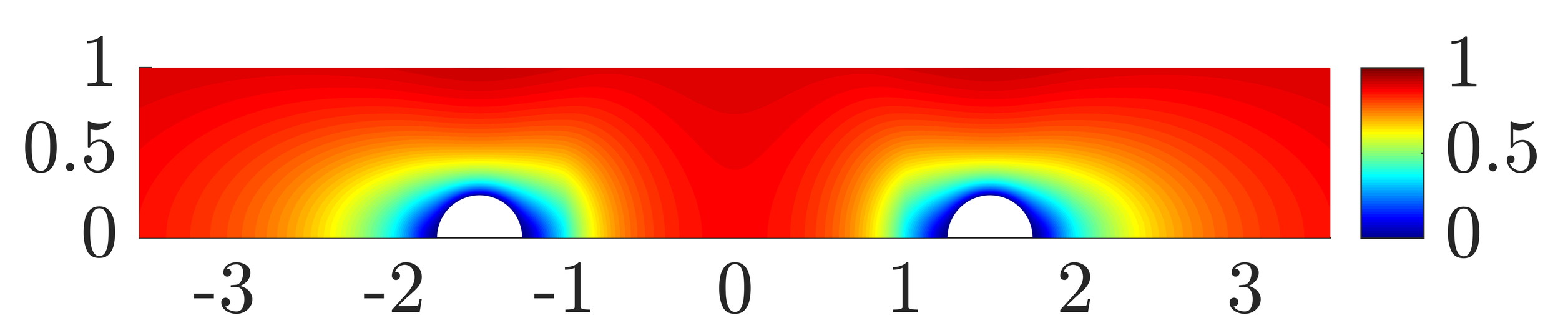}}
	\subfigure[$m = 2$]{\includegraphics[width=0.48\textwidth]{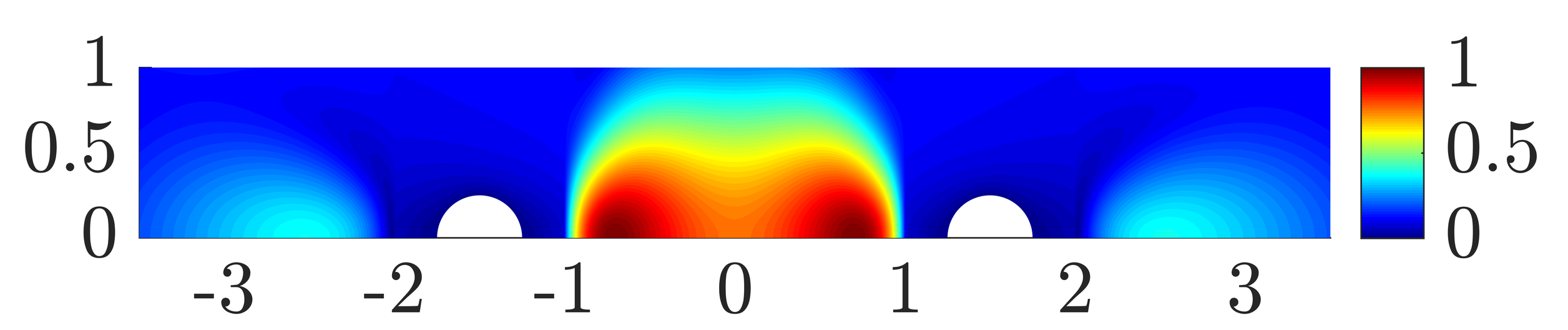}}

	\subfigure[$m = 3$] {\includegraphics[width=0.48\textwidth]{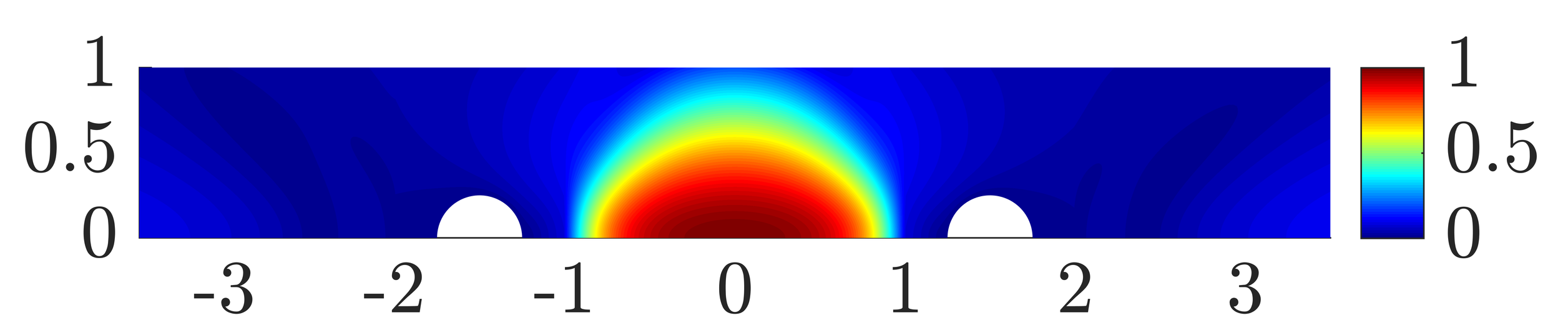}}
	\subfigure[$m = 4$] {\includegraphics[width=0.48\textwidth]{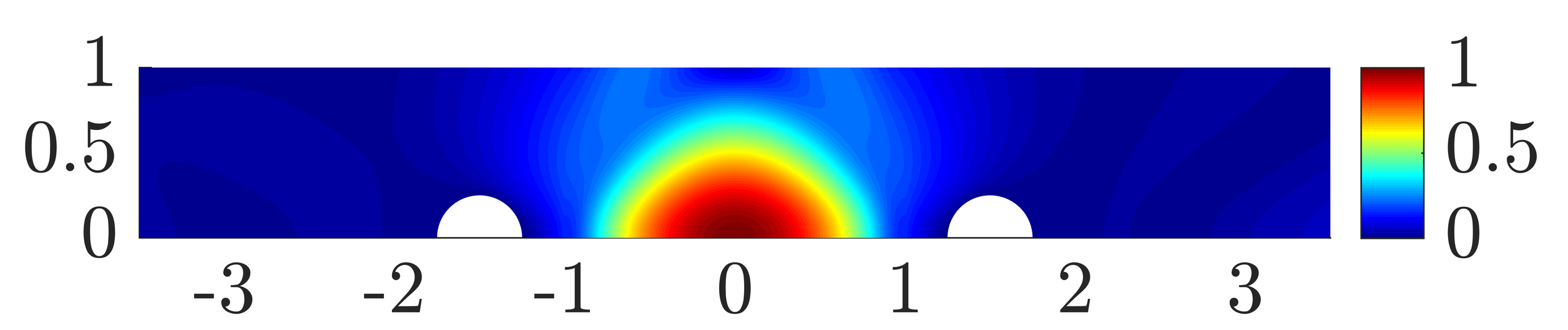}}
		
	\caption{First four normalised spatial modes of the module of the velocity computed using the a priori PGD algorithm for the interval $\I^2 = [-2, -1]$ of the parametric distance.}
	\label{fig:DistanceSmallModesU}
\end{figure}

\begin{figure}[!tb]
	\centering
	\subfigure[$m = 1$]{\includegraphics[width=0.48\textwidth]{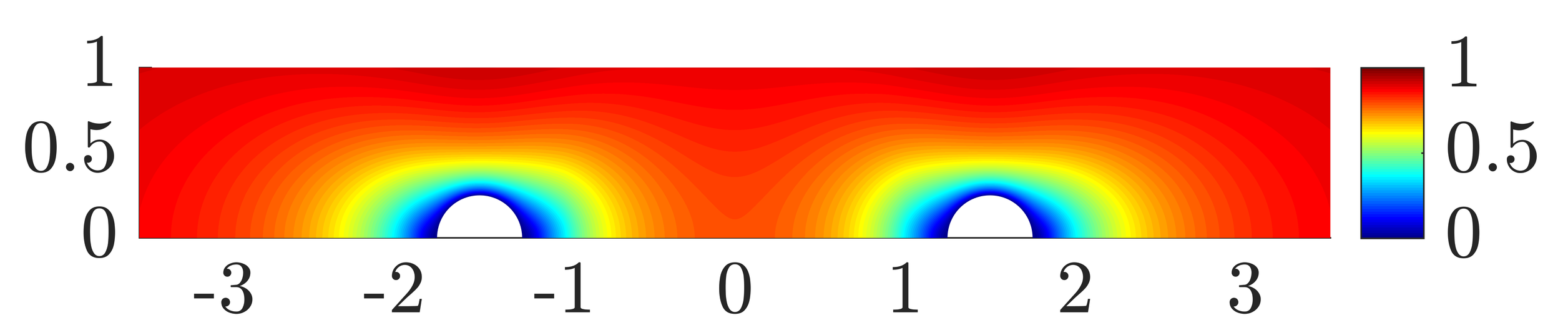}}
	\subfigure[$m = 2$]{\includegraphics[width=0.48\textwidth]{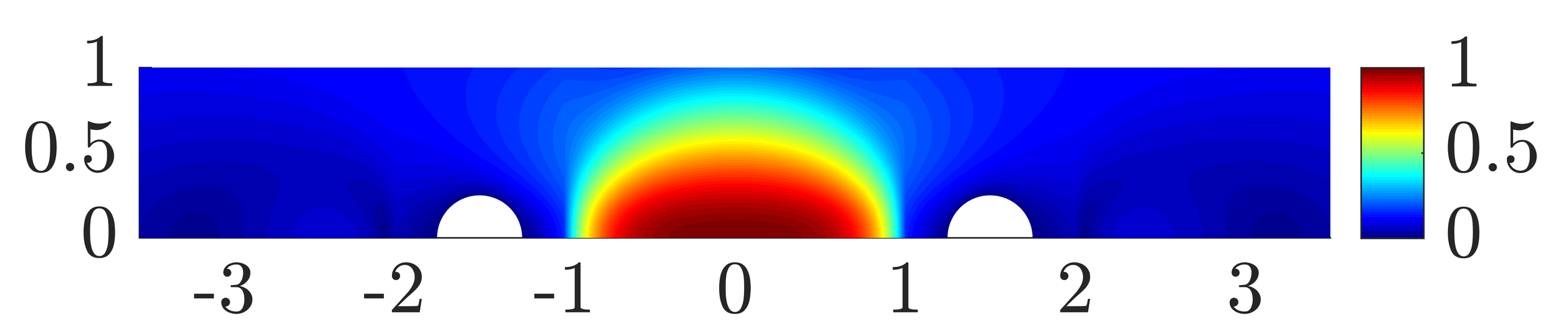}}
	
	\subfigure[$m = 3$] {\includegraphics[width=0.48\textwidth]{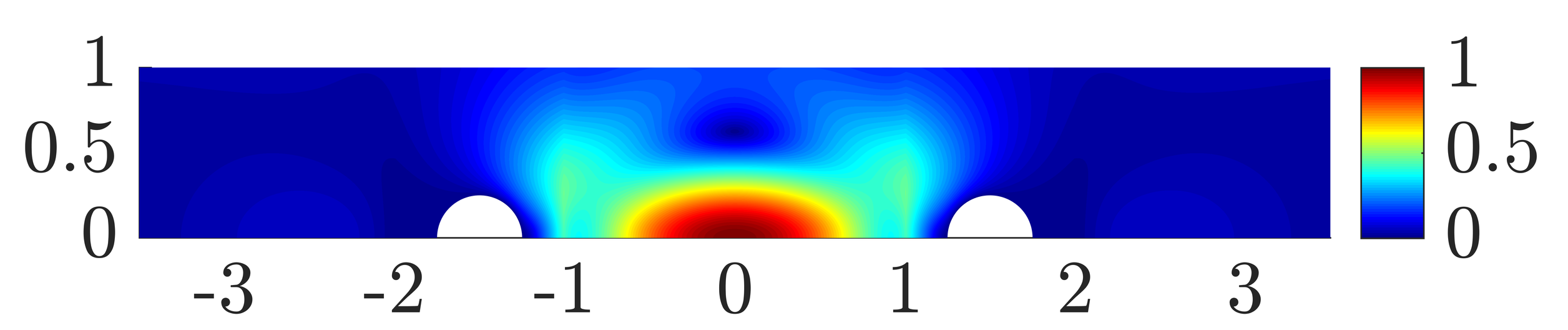}}
	\subfigure[$m = 4$] {\includegraphics[width=0.48\textwidth]{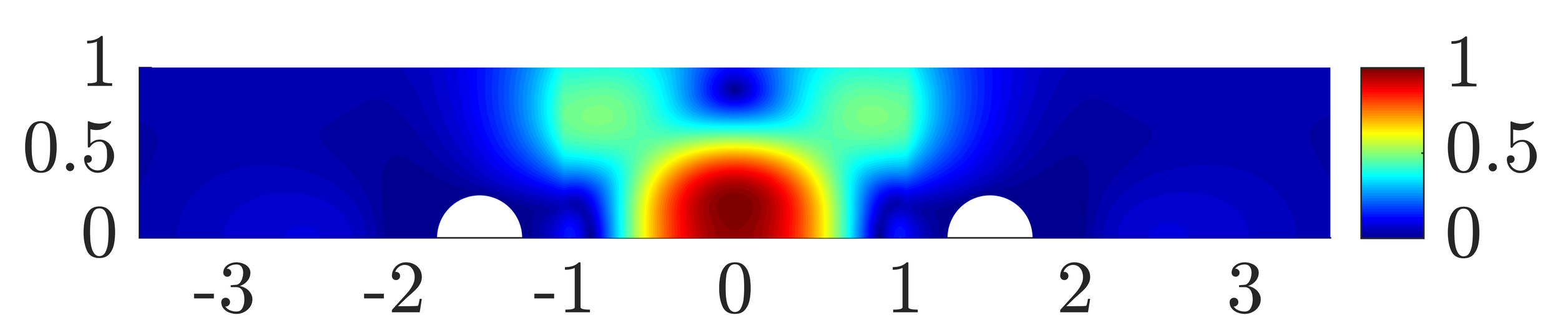}}
		
	\subfigure[$m = 5$]{\includegraphics[width=0.48\textwidth]{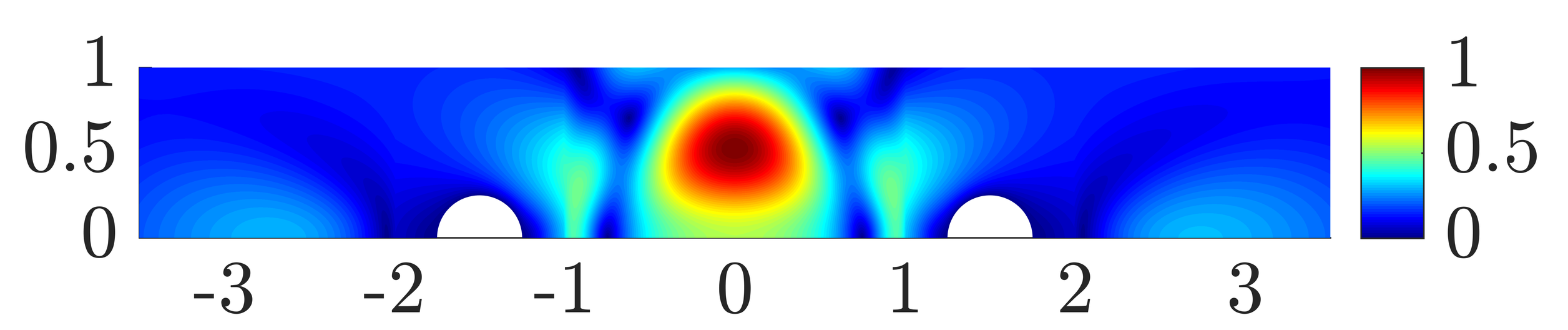}}
	\subfigure[$m = 6$]{\includegraphics[width=0.48\textwidth]{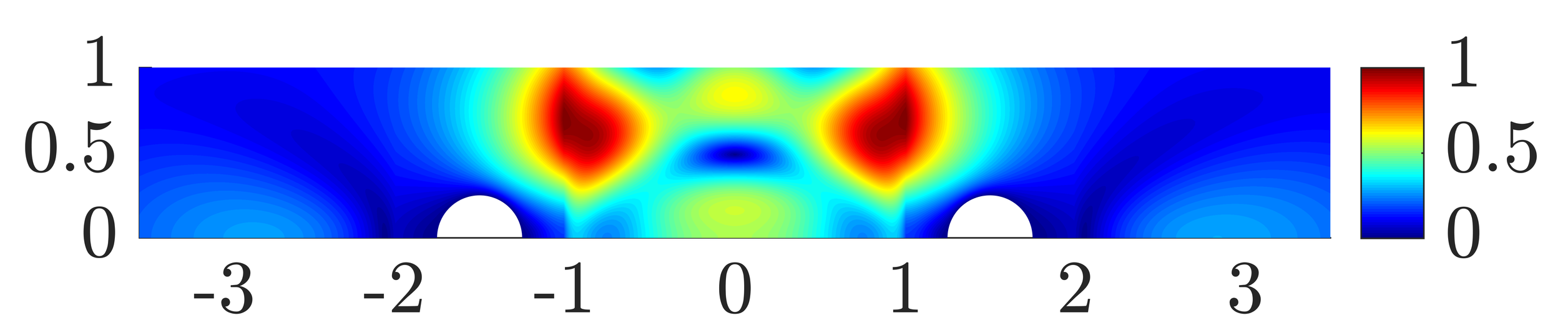}}
	
	\subfigure[$m = 7$] {\includegraphics[width=0.48\textwidth]{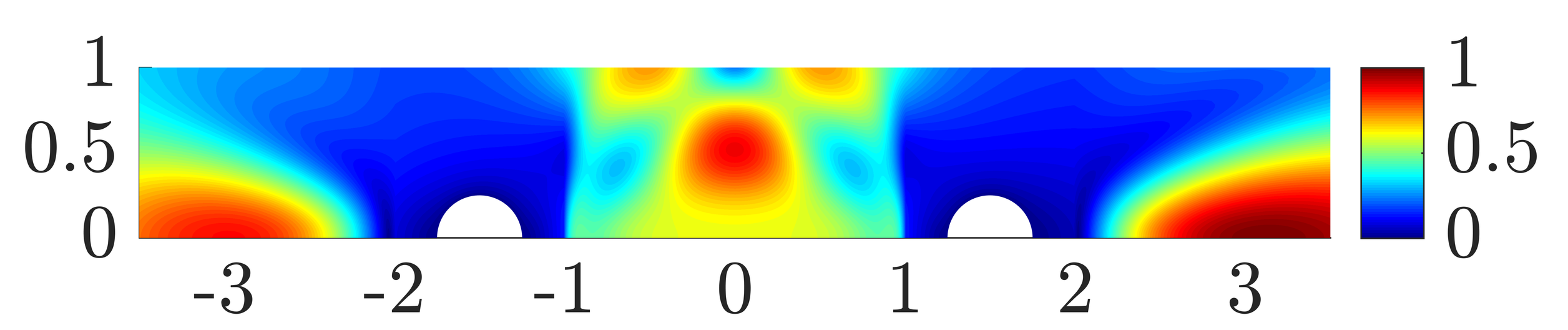}}
	\subfigure[$m = 8$] {\includegraphics[width=0.48\textwidth]{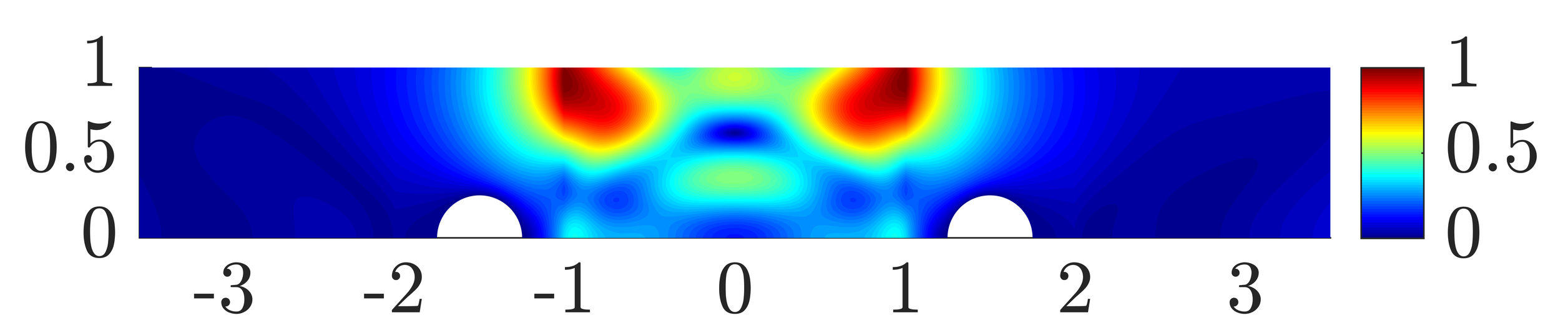}}
	
	\caption{First eight normalised spatial modes of the module of the velocity computed using the a priori PGD algorithm for the interval $\I^2 = [-3, 2]$ of the parametric distance.}
	\label{fig:DistanceLargeModesU}
\end{figure}

\begin{figure}[!tb]
	\centering
	\subfigure[{$\I^2 = [-2, -1]$}]{\includegraphics[width=0.48\textwidth]{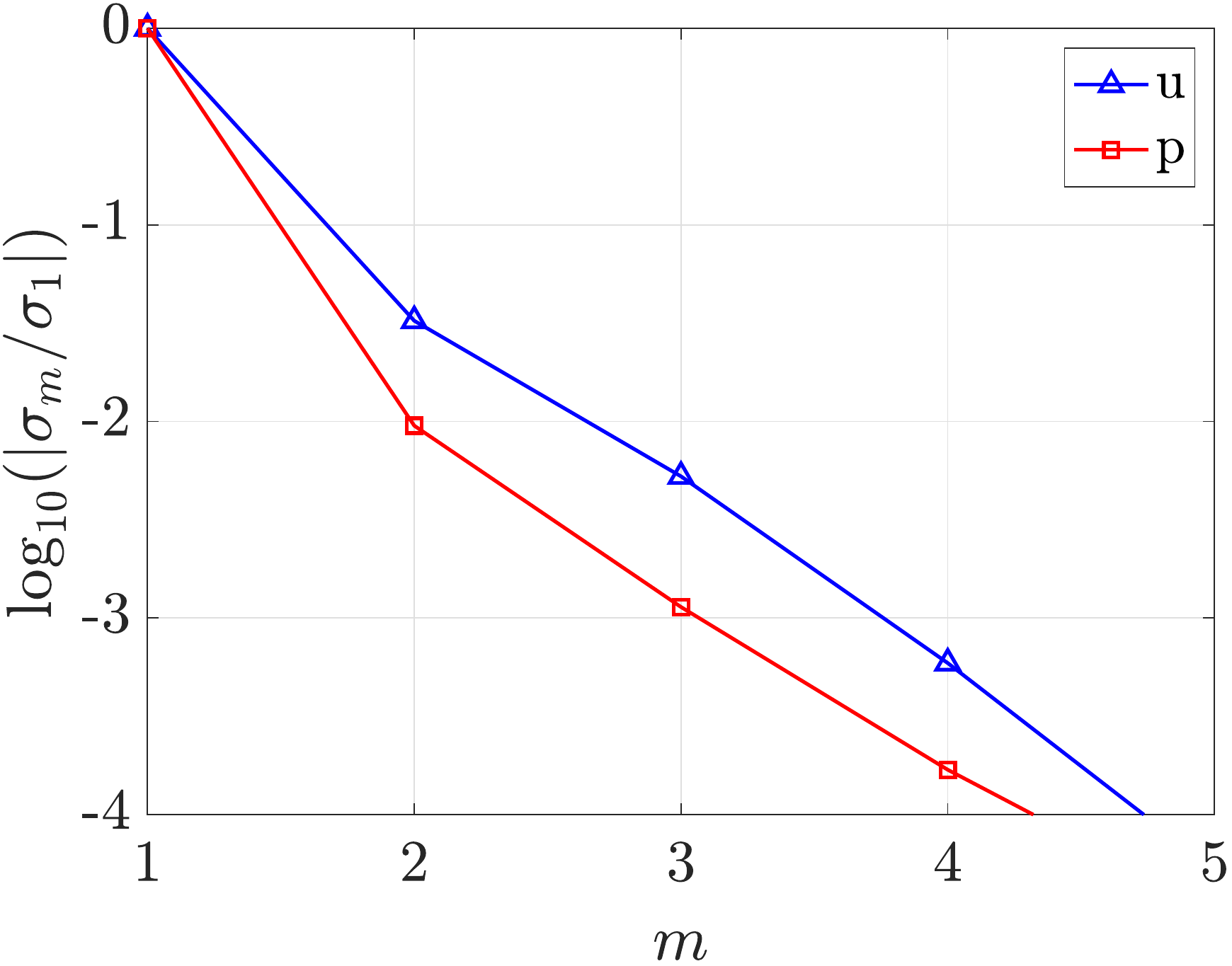}}
	\subfigure[{$\I^2 = [-3, 2]$}]{\includegraphics[width=0.48\textwidth]{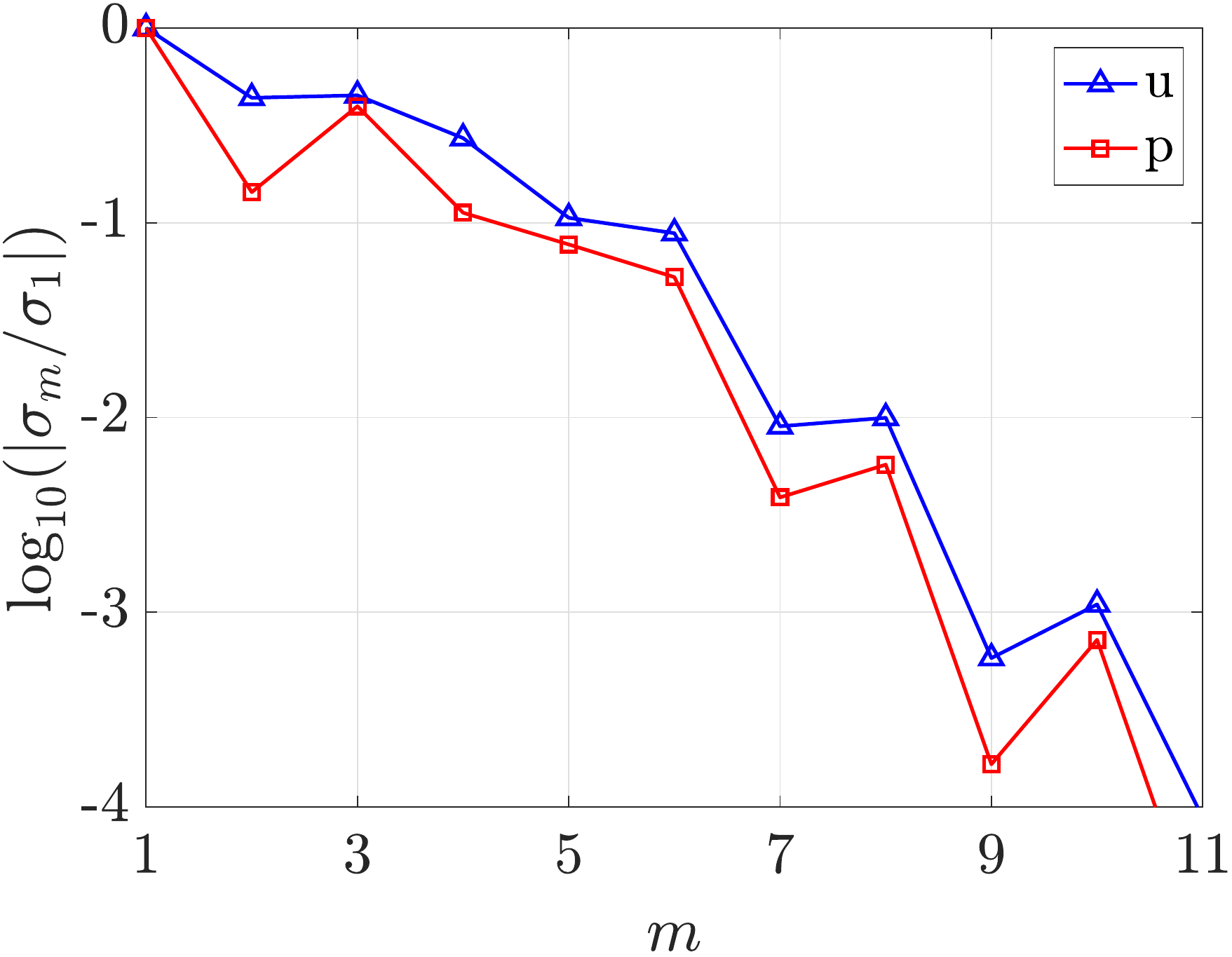}}

	\caption{Convergence of the mode amplitude computed using the a priori PGD algorithm for different intervals of the parametric distance.}
	\label{fig:DistanceAmp}
\end{figure}

\begin{figure}[!tb]
	\centering
	\subfigure[{$\I^2 = [-2, -1]$}]{\includegraphics[width=0.48\textwidth]{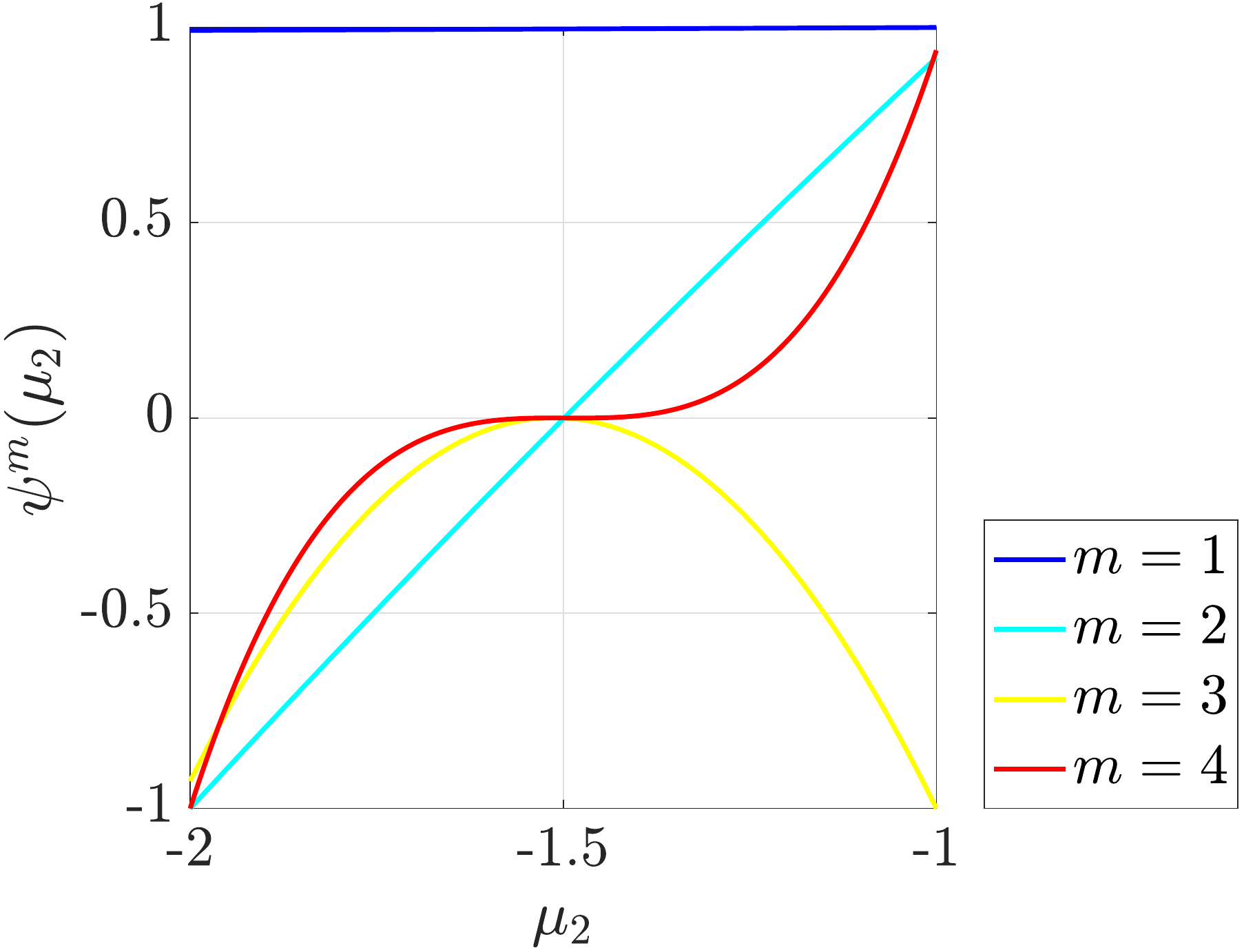}}
	\subfigure[{$\I^2 = [-3, 2]$}]{\includegraphics[width=0.48\textwidth]{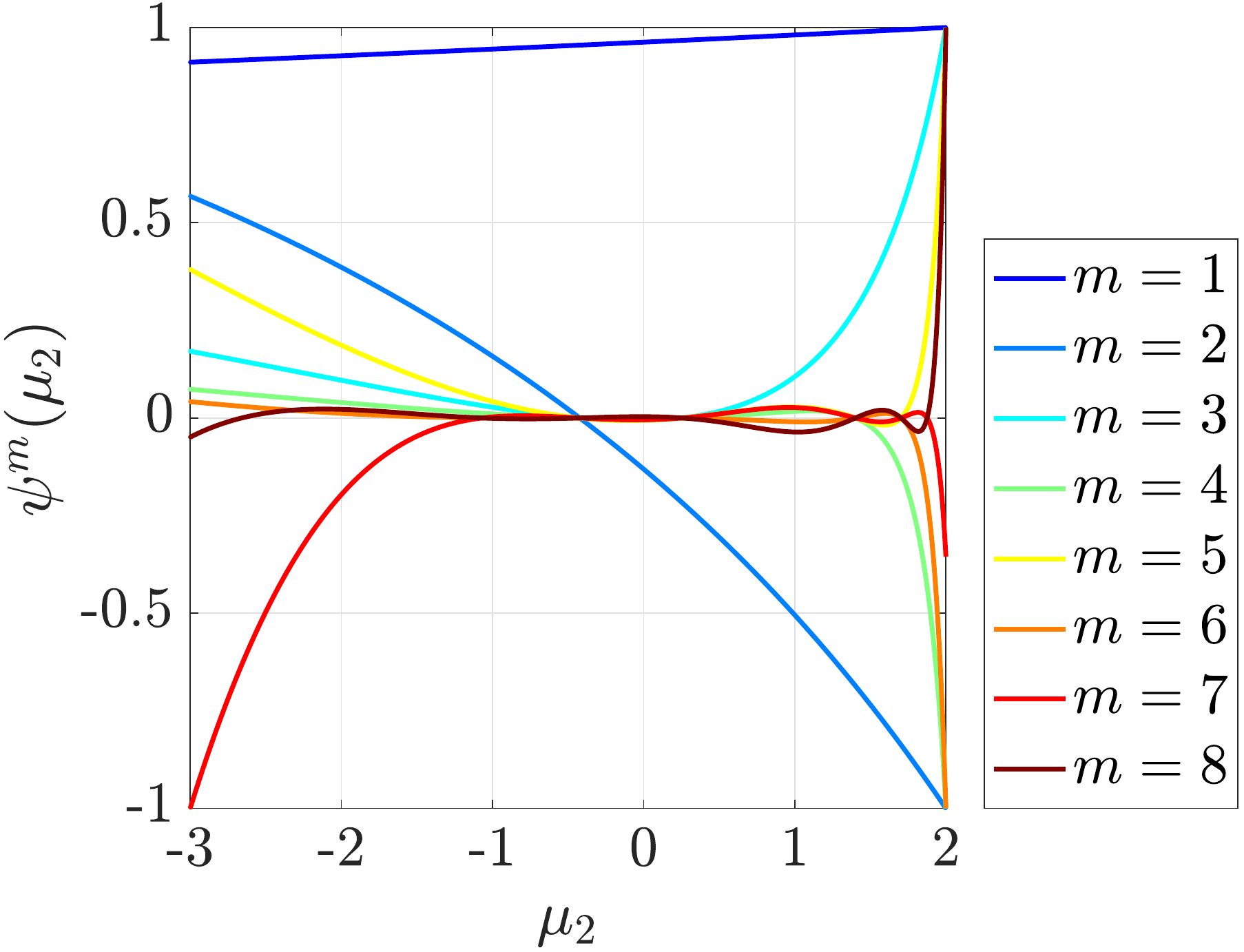}}

	\caption{First normalised parametric modes computed using the a priori PGD algorithm for different intervals of the parametric distance.}
	\label{fig:DistanceModesParam}
\end{figure}

To further highlight the additional difficulty of constructing a ROM for the extended range of values of the parameter $\mu_2$, the first normalised spatial modes of the module of the velocity are reported in figure~\ref{fig:DistanceSmallModesU} and~\ref{fig:DistanceLargeModesU} for the interval $\I^2 = [-2, -1]$ and $\I^2 = [-3, 2]$, respectively. The results, computed via the a priori PGD with $\niter {=} 3$, display an increased variability of the flow in the region between the two bladders, where modes accounting for localised spatial phenomena appear. A similar behaviour is experienced by the pressure modes, not reported here for brevity.
The importace of the presented modes in the final PGD approximation is analysed in figure~\ref{fig:DistanceAmp}. The evolution of the relative amplitude of the computed modes for the interval $\I^2 = [-3, 2]$ displays that the first four modes have a comparable importance in the PGD approximation. On the contrary, the fourth mode in the case of $\I^2 = [-2, -1]$ already has a relative amplitude below $10^{-3}$. To achieve this level of truncation when the extended interval $\I^2$ is considered, the a priori PGD needs to compute nine modes.
Finally, figure~\ref{fig:DistanceModesParam} shows the corresponding parametric modes computed using the a priori PGD strategy for the two parametric intervals studied above. The results display that the parametric functions for the small interval $\I^2 = [-2, -1]$ present a smooth transition between the extreme values of $\mu_2$ and no localised phenomena are identified. On the contrary, when the extended interval $\I^2 = [-3, 2]$ is considered, several parametric functions feature large variations in a localised region between $\mu_2 {=} 1$ and $\mu_2 {=} 2$. This region of the parametric domain is associated with the configurations of minimum distance between the bladders. In these scenarios, the influence of the two spheres on one another is maximum and small variations of the distance are expected to generate complex flow patterns.

\begin{figure}[!tb]
	\centering
	\subfigure[{$\I^2 = [-2,-1]$}]{\includegraphics[width=0.48\textwidth]{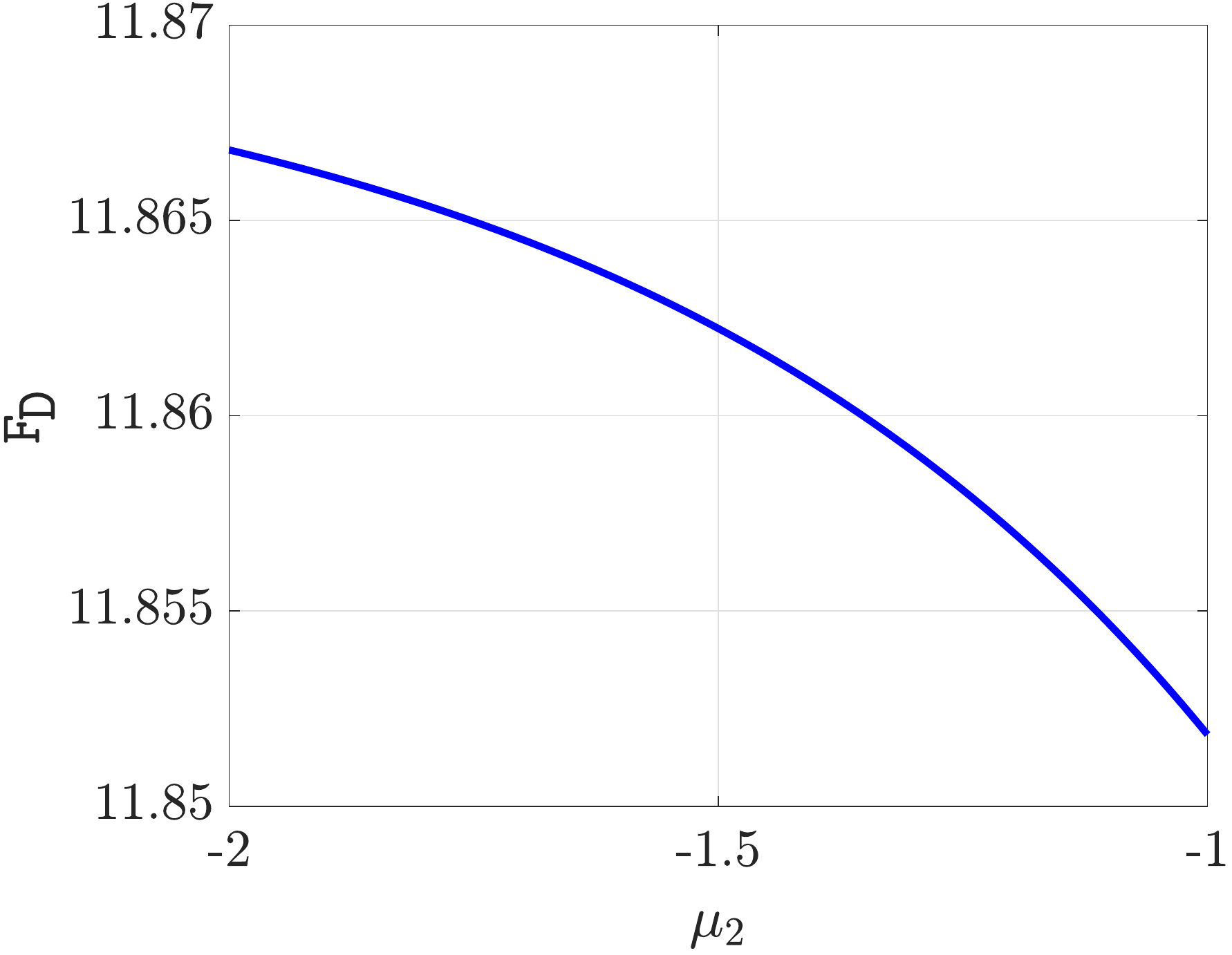}\label{fig:DragDistanceSmall}}
	\subfigure[{$\I^2 = [-3, 2]$}]{\includegraphics[width=0.48\textwidth]{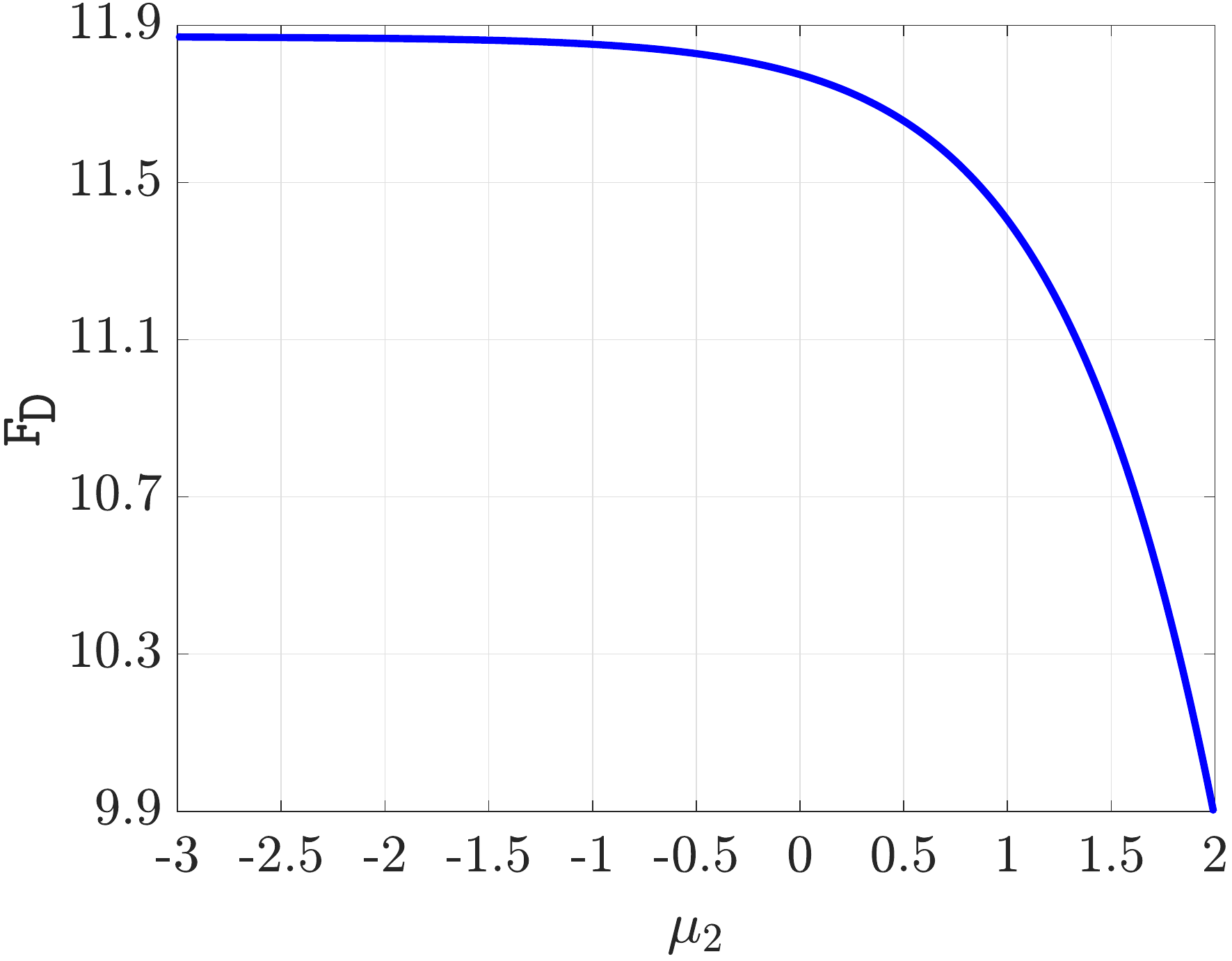}\label{fig:DragDistanceLarge}}
	
	\caption{Response surfaces of the total drag force as a function of the distance $\mu_2$ between the two spheres, for two different ranges of values of the parameter.}
	\label{fig:DragDistance}
\end{figure}

As in the previous example, the separated response surface of the total drag force on the two spheres is computed using the a priori PGD algorithm. It is worth noticing that the range of values of $\mu_2$ considered in figure~\ref{fig:DragDistanceSmall} is a subinterval of the one analysed in figure~\ref{fig:DragDistanceLarge}. The scales of the two figures confirm the higher variability of the flow quantities when larger parametric intervals are considered and the consequent additional difficulties faced by the PGD-ROM strategies to cope with the sensitivity to the range of values considered.

%-------------------------------------------------------
\subsubsection{Accuracy of a priori and a posteriori response surfaces} 
\label{sc:comparisonAccuracy}
%-------------------------------------------------------

\begin{figure}[!tb]
	\centering
	\subfigure[Error measure]{\includegraphics[width=0.49\textwidth]{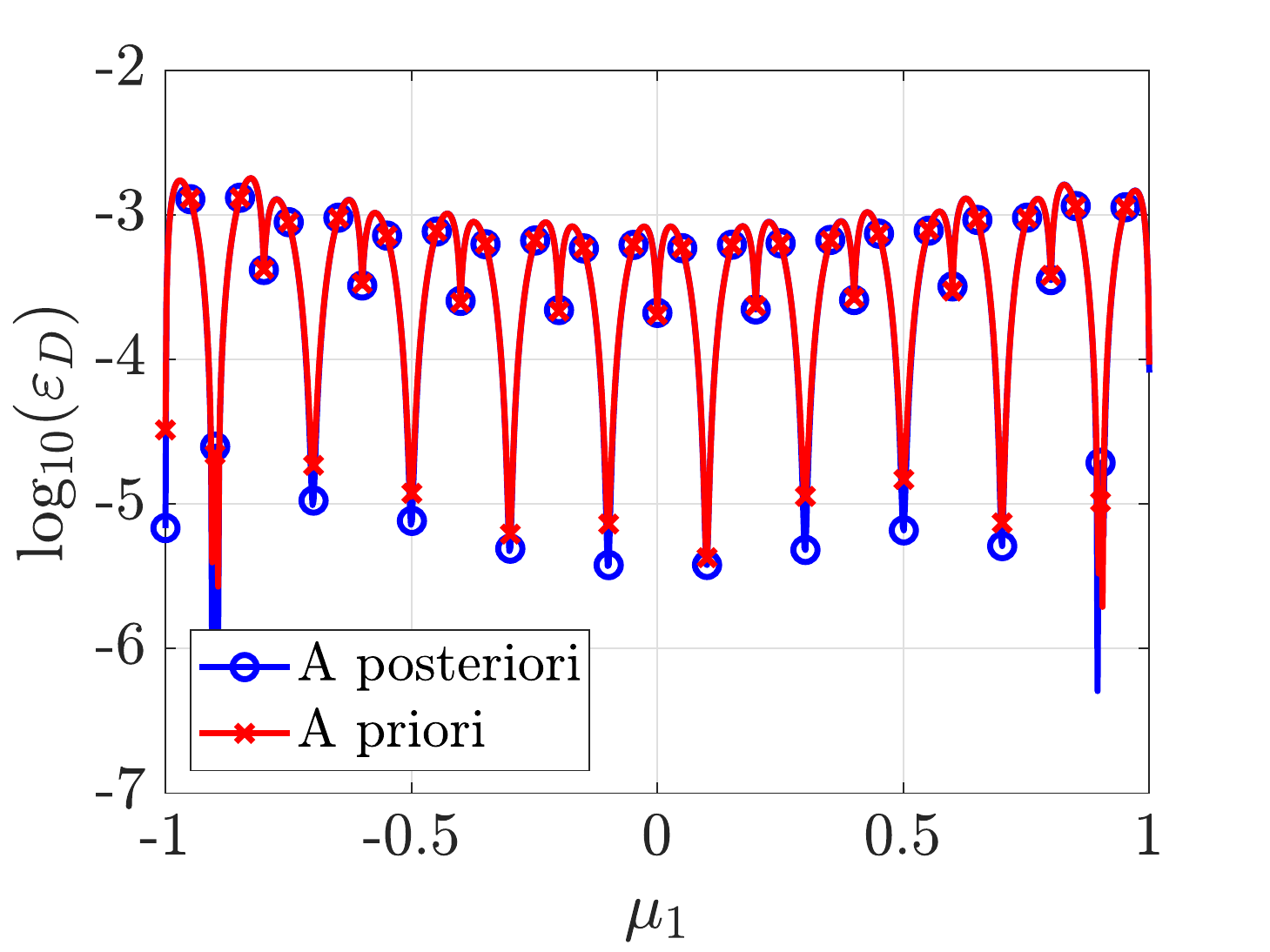}\label{fig:ErrD_param1Pointwise}}
	\subfigure[Smoothed error measure]{\includegraphics[width=0.49\textwidth]{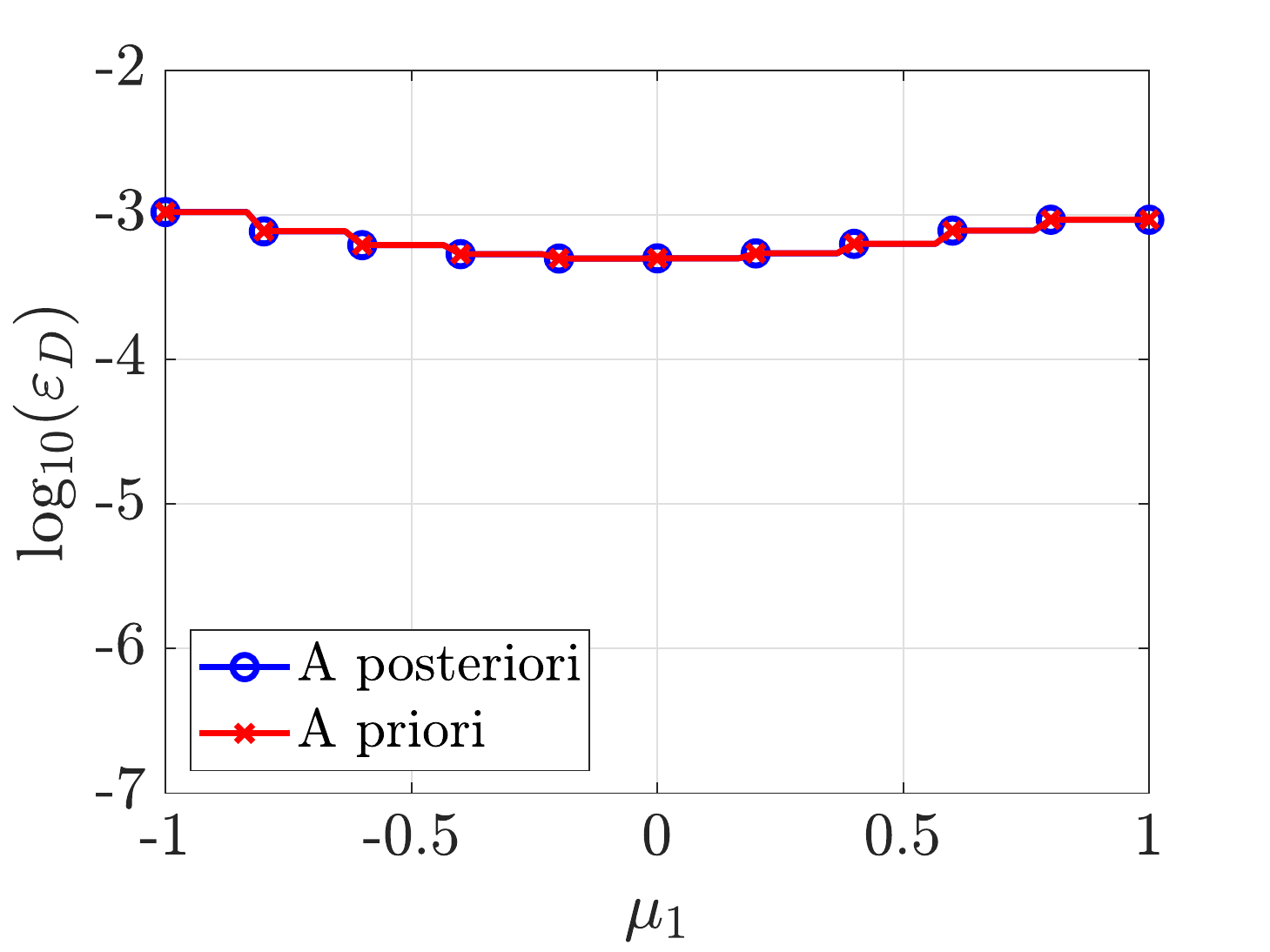}\label{fig:ErrD_param1Smooth}}
	\caption{Error in the drag, defined in equation~\eqref{eq:DragErrorPointwise}, as a function of the parameter $\mu_1$.}
	\label{fig:ErrD_param1}
\end{figure}

The previous examples with one geometric parameter have shown that, when the error in equation~\eqref{eq:DragErrorMulti} is considered, the computational cost of the a priori PGD is as competitive as the a posteriori one and outperforms it for a larger range of the parametric interval. However, this quantity measures the average accuracy over the whole parametric domain, without considering the worst case scenarios, that is, the cases where the maximum error is observed. 
To further compare the two approaches, figure~\ref{fig:ErrD_param1} displays the value of the error in equation~\eqref{eq:DragErrorPointwise} in the drag force, as a function of the parameter $\mu_1$ for the first example with one geometric parameter controlling the radius of the spherical bladders.
The minima observed for both the a priori and the a posteriori approaches in figure~\ref{fig:ErrD_param1Pointwise} coincide with the midpoints of the elements $\I_e^1, \ e {=} 1,\ldots,\numel^1$ as these locations correspond to both a high-order node and an integration point for the fourth-order polynomial approximation used in each element of the parametric space. More importantly, the results show that the accuracy of the a priori and the a posteriori approaches is almost identical, not only when measured in the $\eltwo(\I^1)$ norm (Fig.~\ref{fig:L2errDragRadius}), but also when the pointwise error in the drag force is displayed for every configuration in $\I^1$. To capture the qualitative behaviour of the error as a function of $\mu_1$, a smoothing is displayed in figure~\ref{fig:ErrD_param1Smooth}. The results clearly show that the error is slightly higher near the boundary of the parametric interval. The smoothing is performed by considering a single value for the error in each element, obtained as the average of the error at all integration points.

\begin{figure}[!tb]
	\centering
	\subfigure[{$\I^2 = [-2,-1]$}]{\includegraphics[width=0.49\textwidth]{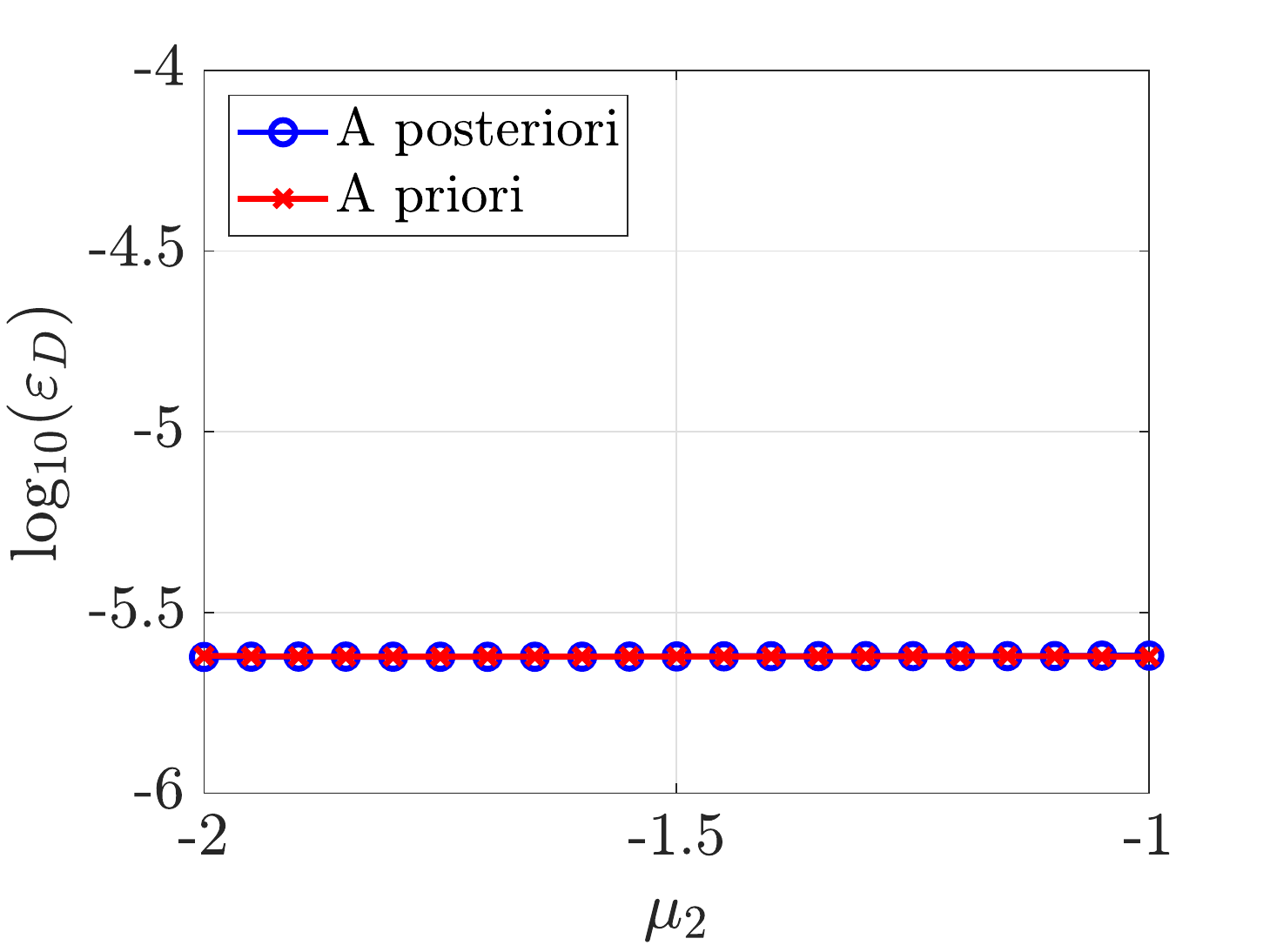}\label{fig:ErrD_param2A}}
	\subfigure[{$\I^2 = [-3, 2]$}]{\includegraphics[width=0.49\textwidth]{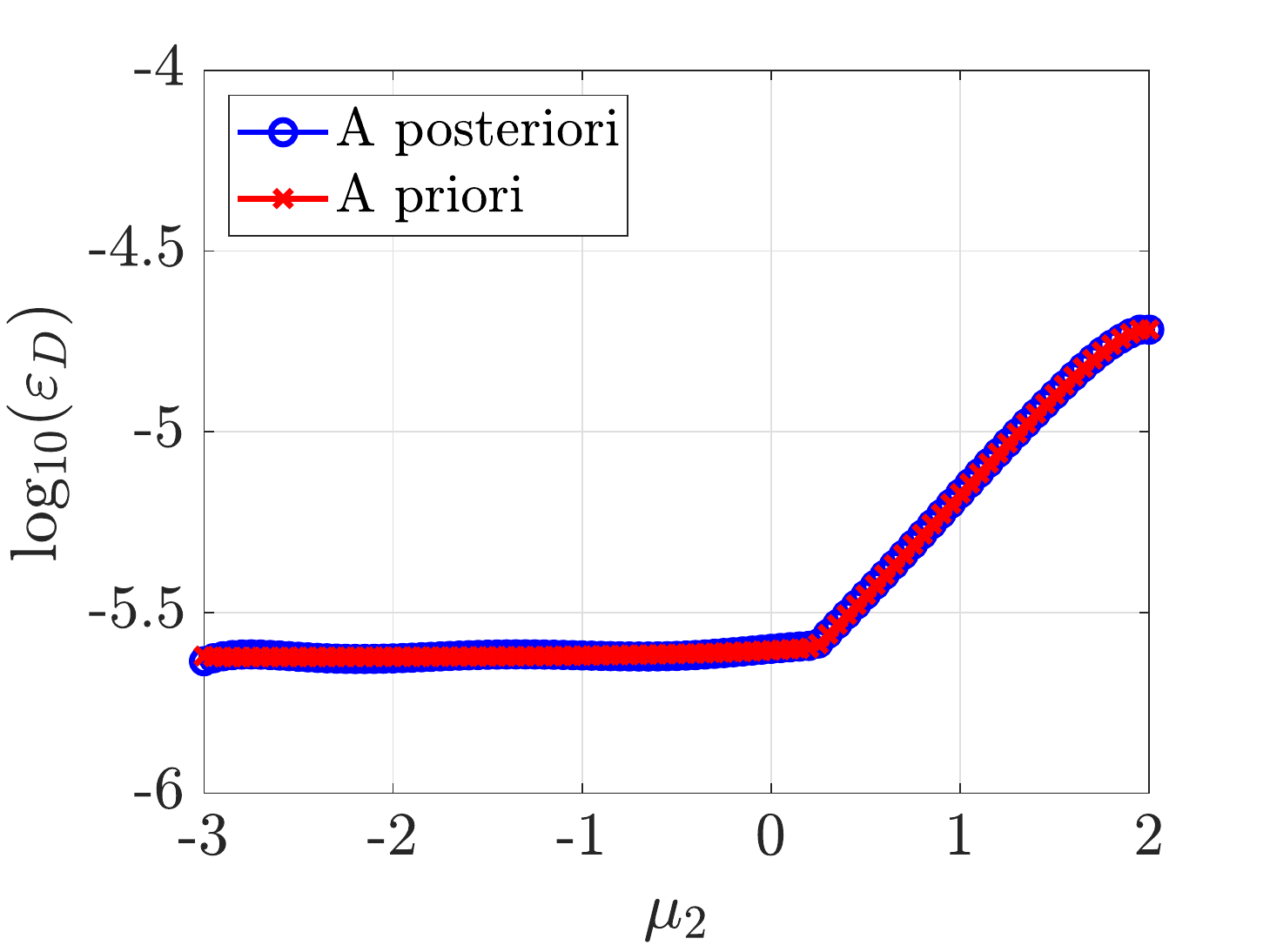}\label{fig:ErrD_param2B}}
	\caption{Smoothed error measure in the drag force, defined in equation~\eqref{eq:DragErrorPointwise}, as a function of the parameter $\mu_2$.}
	\label{fig:ErrD_param2}
\end{figure}

\begin{figure}[!tb]
	\centering
	\subfigure[Module of velocity, $\mu_2=-3$]{\includegraphics[width=0.48\textwidth]{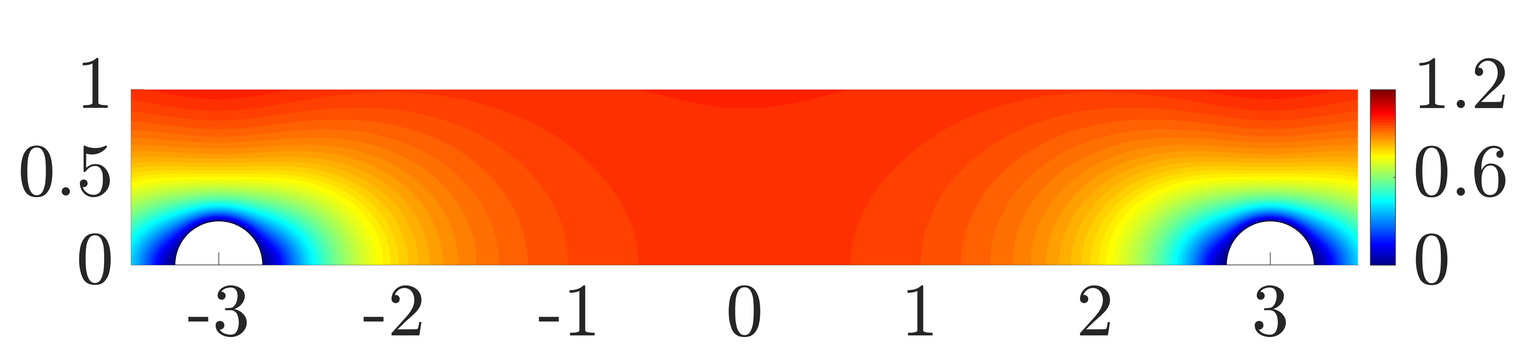}}
	\subfigure[Pressure, $\mu_2=-3$]{\includegraphics[width=0.48\textwidth]{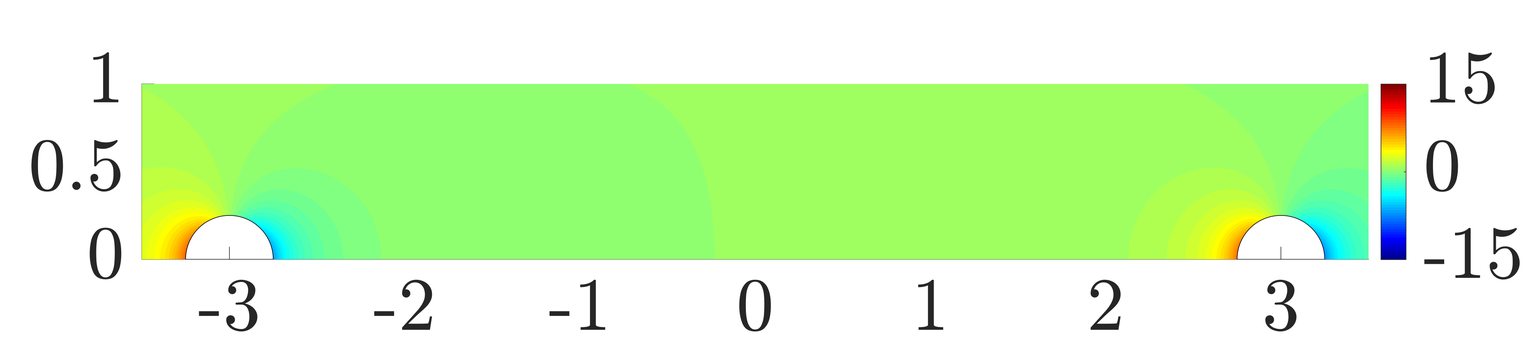}}
	
	\subfigure[Module of velocity, $\mu_2=2$] {\includegraphics[width=0.48\textwidth]{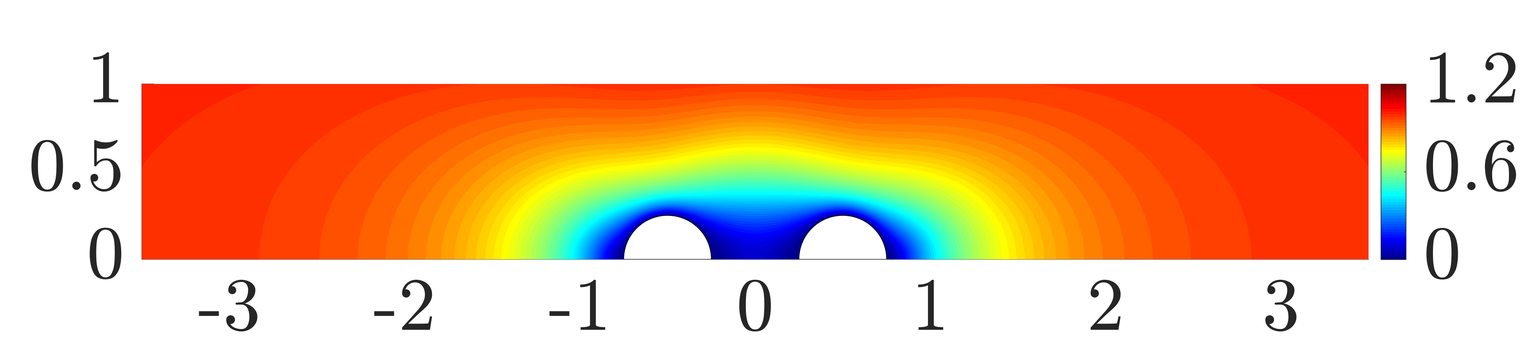}}	
	\subfigure[Pressure, $\mu_2=2$] {\includegraphics[width=0.48\textwidth]{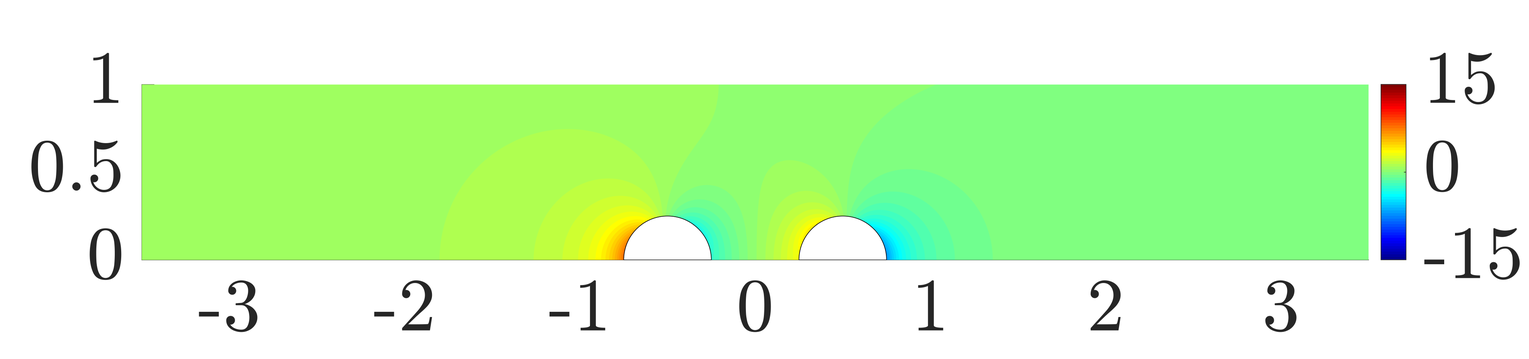}}

	\caption{Module of velocity and pressure field for two values of the parameter $\mu_2$ corresponding to maximum (top) and minimum (bottom) distance of the bladders.}
	\label{fig:Distance2Online}
\end{figure}

Similarly, figure~\ref{fig:ErrD_param2} compares the value of the smoothed error measure in the drag force as a function of the parameter $\mu_2$ for the second example, with one geometric parameter controlling the distance between the spherical bladders. 
The results of the two cases previously studied, with $\I^2 {=} [-2,-1]$ and $\I^2 {=} [-3,2]$, display the increased difficulty of computing an accurate response surface as the range of values in the parametric space increases. For $\I^2 {=} [-2,-1]$, the accuracy is almost independent of the value of the parameter, whereas for $\I^2 {=} [-3,2]$ a more significant dependence is observed, especially near $\mu_2 {=} 2$, that is, when the distance between the spherical bladders is minimum. It is clear that for large values of $\mu_2$, there is a strong influence in the flow impinging onto the second sphere caused by its proximity to the first sphere. Figure~\ref{fig:Distance2Online} reports the particularisation of the PGD solution computed for two values of the parameter of interest, highlighting the additional difficulties due to the small distance between the two bladders for $\mu_2 {=} 2$.

Finally, the comparison of the results in figures~\ref{fig:ErrD_param1} and~\ref{fig:ErrD_param2} clearly illustrates the rationale behind the choice of the resolution for the discretisation of the intervals $\I^1$ and $\I^2$. Given the limited variation of the solution in figure~\ref{fig:ErrD_param1Smooth}, only 10 elements were considered in the first parametric dimension, whereas the discretisation of the second parametric dimension contains 20 elements for the interval $\I^2 {=} [-2,-1]$ and 100 elements for the case of $\I^2 {=} [-3,2]$.

%-------------------------------------------------------
\subsection{Two geometric parameters} 
\label{sc:exTwo}
%-------------------------------------------------------

In this section, the two geometric parameters studied separately in the previous examples are considered in a single simulation. As observed in~\cite{RS-SBGH-20}, this problem is particularly challenging and an increased number of modes is required to capture the solution. To ease the visualisation, the figures in this section report the number of modes for the a priori approach after the PGD compression is performed, whereas the global number of computed modes is commented in the text.

\begin{figure}[!tb]
	\centering
	\subfigure[$\bu$]{\includegraphics[width=0.49\textwidth]{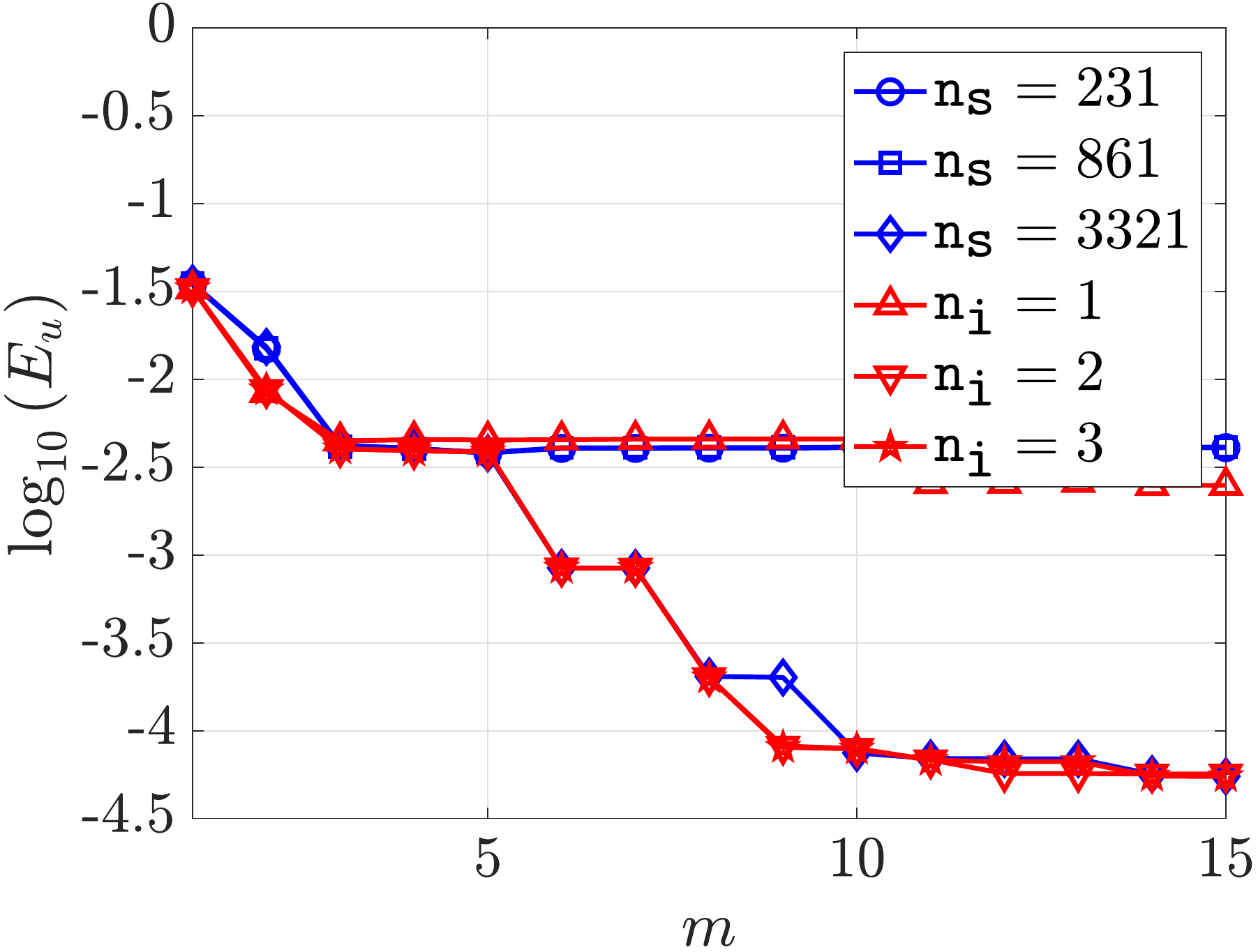}}
	\subfigure[$p$]  {\includegraphics[width=0.49\textwidth]{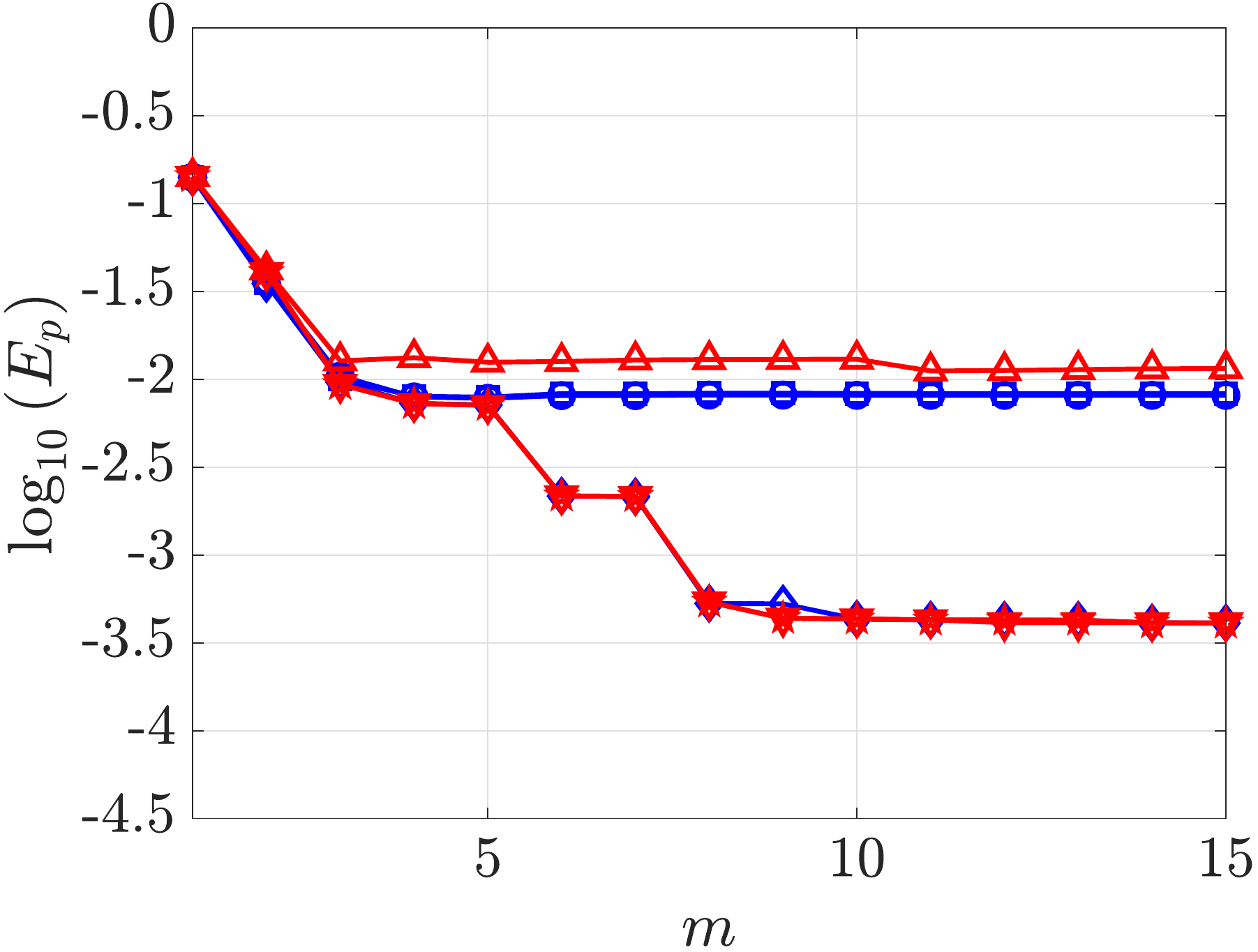}}
	\subfigure[$\bL$]{\includegraphics[width=0.49\textwidth]{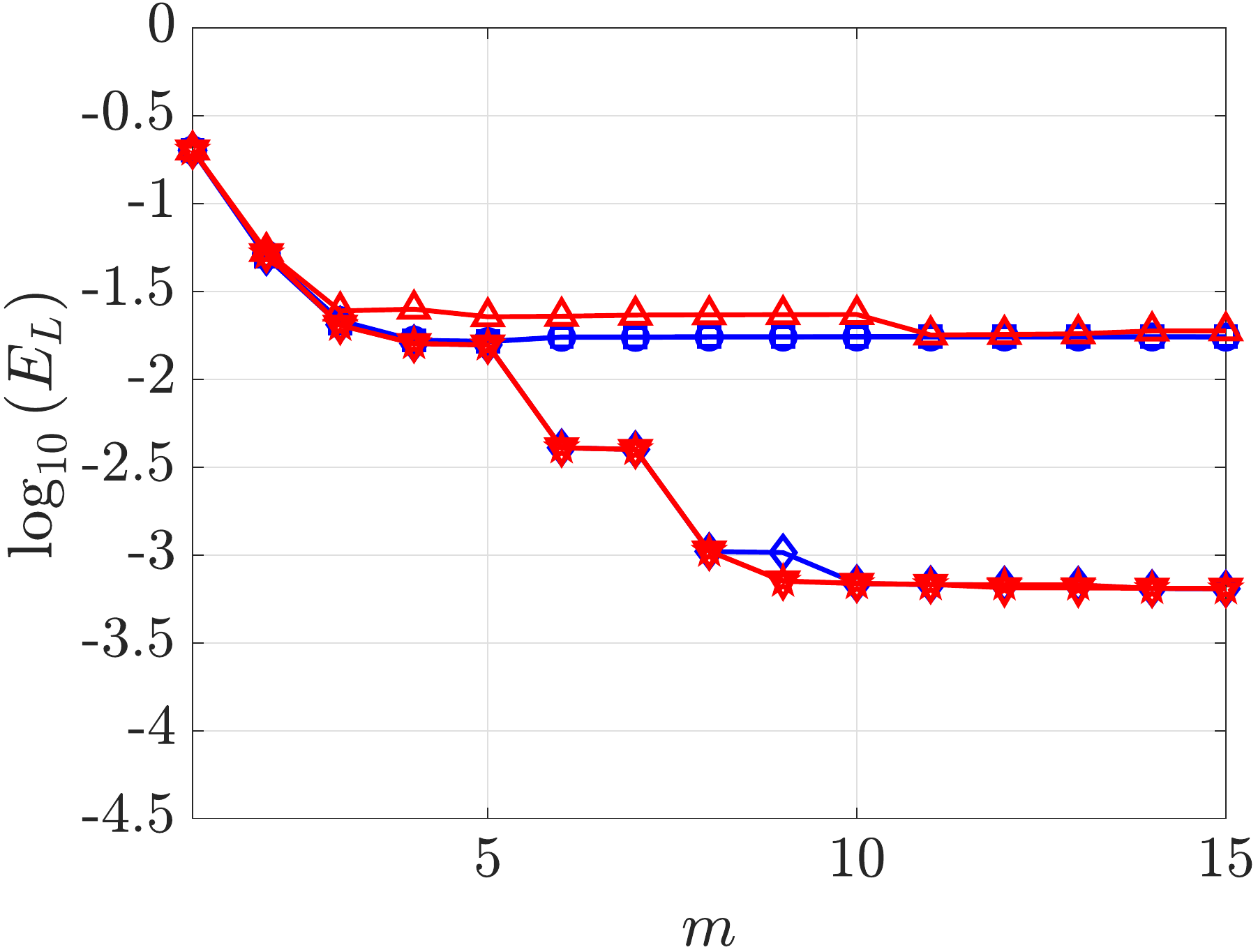}}
	\subfigure[$\FD$]{\includegraphics[width=0.49\textwidth]{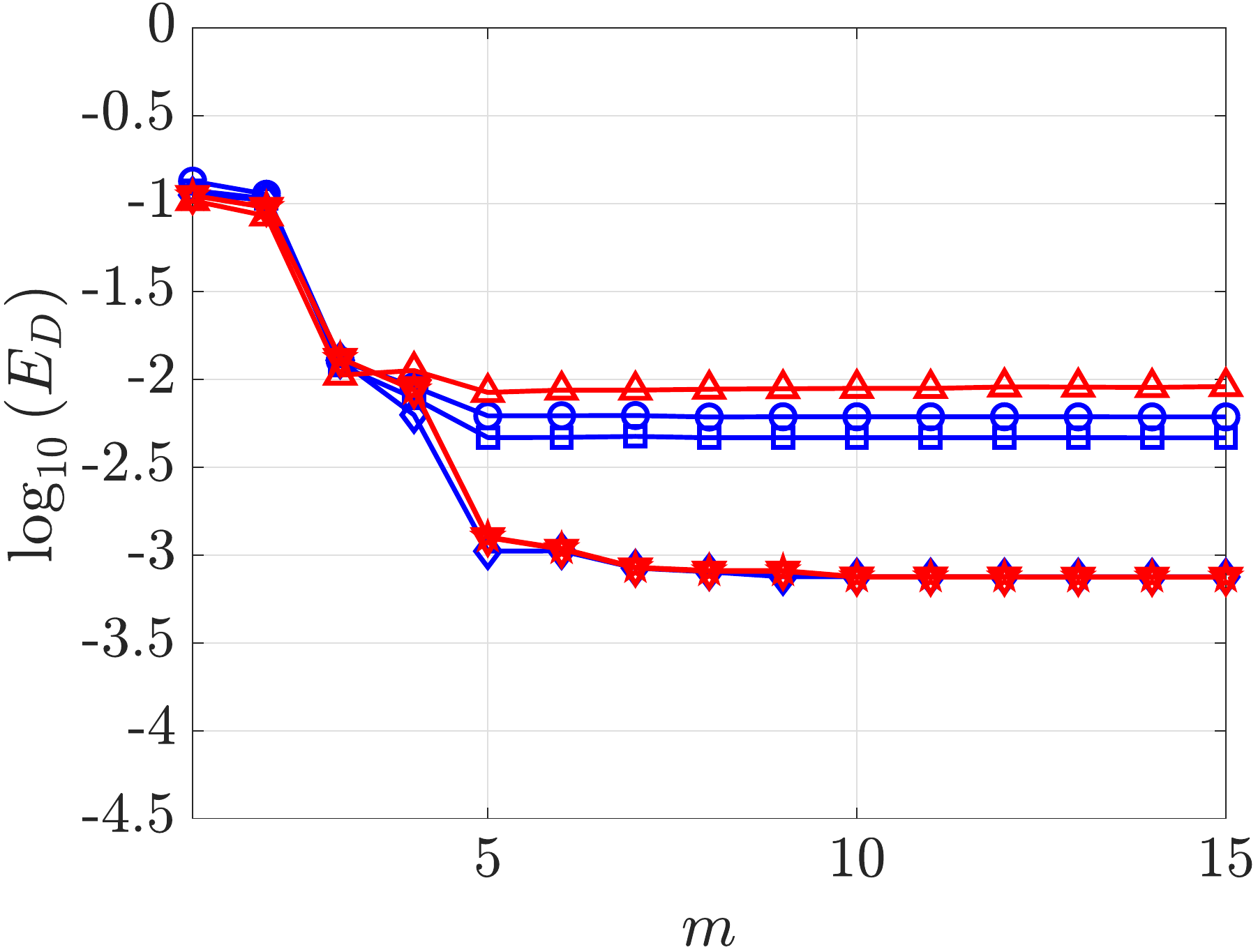}}
	\caption{Evolution of (a-c) the $\eltwo(\Omega \times \bI)$ error for velocity, pressure and gradient of velocity and (d) the $\eltwo(\bI)$ error for the drag force as a function of the number of PGD modes for the problem with two geometric parameters and $\I^2 {=} [-2,-1]$. The legend details the number $\nsnap$ of snapshots used by the a posteriori PGD approach (blue) and the number $\niter$ of nonlinear iterations used by the a priori PGD approach (red).}
	\label{fig:RadisDistance1}
\end{figure}

First, the interval for the parameter that controls the distance is set to $\I^2 {=} [-2,-1]$. Figure~\ref{fig:RadisDistance1} shows the evolution of the $\eltwo(\Omega \times \bI)$ error for velocity, pressure and gradient of velocity and the $\eltwo(\bI)$ error for the drag force as a function of the number $m$ of modes. 
The results reveal that the a posteriori approach provides almost identical accuracy using 231 and 861 snapshots. The errors in velocity and pressure are below $10^{-2}$ and the error in the gradient of  velocity is almost $10^{-2}$. In this case, four modes are sufficient to obtain the maximum accuracy in velocity, pressure and gradient of velocity, whereas an additional mode is required to achieve the most accurate results in the drag force. Using 3,321 snapshots, the a posteriori PGD computes 10 modes and provides much more accurate results, with an error one order of magnitude lower, compared to the computation with 861 snapshots. To obtain an error in the drag force below $10^{-2}$, the a posteriori approach requires four modes and 231 snapshots, whereas seven modes and 3,321 snapshots are required to achieve an error below $10^{-3}$.
When the a priori PGD algorithm is employed, one nonlinear iteration in the AD scheme is sufficient to obtain an accuracy almost identical to the one provided by the a posteriori approach with 231 and 861 snapshots. In addition, by considering only two nonlinear iterations, the a priori approach is capable of producing the same accuracy as the a posteriori PGD with 3,321 snapshots. In both cases, the number of modes required to obtain the maximum accuracy is the same for the a posteriori approach and the a priori one after PGD compression.

For this example, a fixed number of modes is computed for the a priori approach. This number is prescribed larger than what is really needed to achieve convergence, as can be seen by the saturation of the curves in figure~\ref{fig:RadisDistance1}. In this manner, the choice of the convergence criterion in the a priori PGD algorithm does not affect the number of computed modes, providing a fair comparison with the a posteriori PGD without tailored sampling strategies. More precisely, 200 modes are computed for the problem under analysis and this information is then compressed in the 15 modes reported in figure~\ref{fig:RadisDistance1}. The performance of the a priori PGD is therefore extremely competitive as, for an error in the drag force below $10^{-3}$, it requires the solution of 600 spatial problems (i.e. 200 modes, each computed with two iterations of the AD scheme plus the initial solve to perform the prediction of the mode, see algorithm~\ref{alg:PGDpriori}), whereas 3,321 snapshots are needed by the a posteriori approach. Hence, in this example, the a priori PGD method requires 18\% of the number of calls to the HDG solver performed by the a posteriori PGD algorithm. As mentioned in section~\ref{sc:comparison}, the a posteriori approach benefits from the possibility to compute the snapshots in parallel, but this example shows that the number of calls to the spatial solver required is significantly larger than the ones performed using the a priori algorithm.

\begin{figure}[!tb]
	\centering
	\subfigure[$\bu$]{\includegraphics[width=0.49\textwidth]{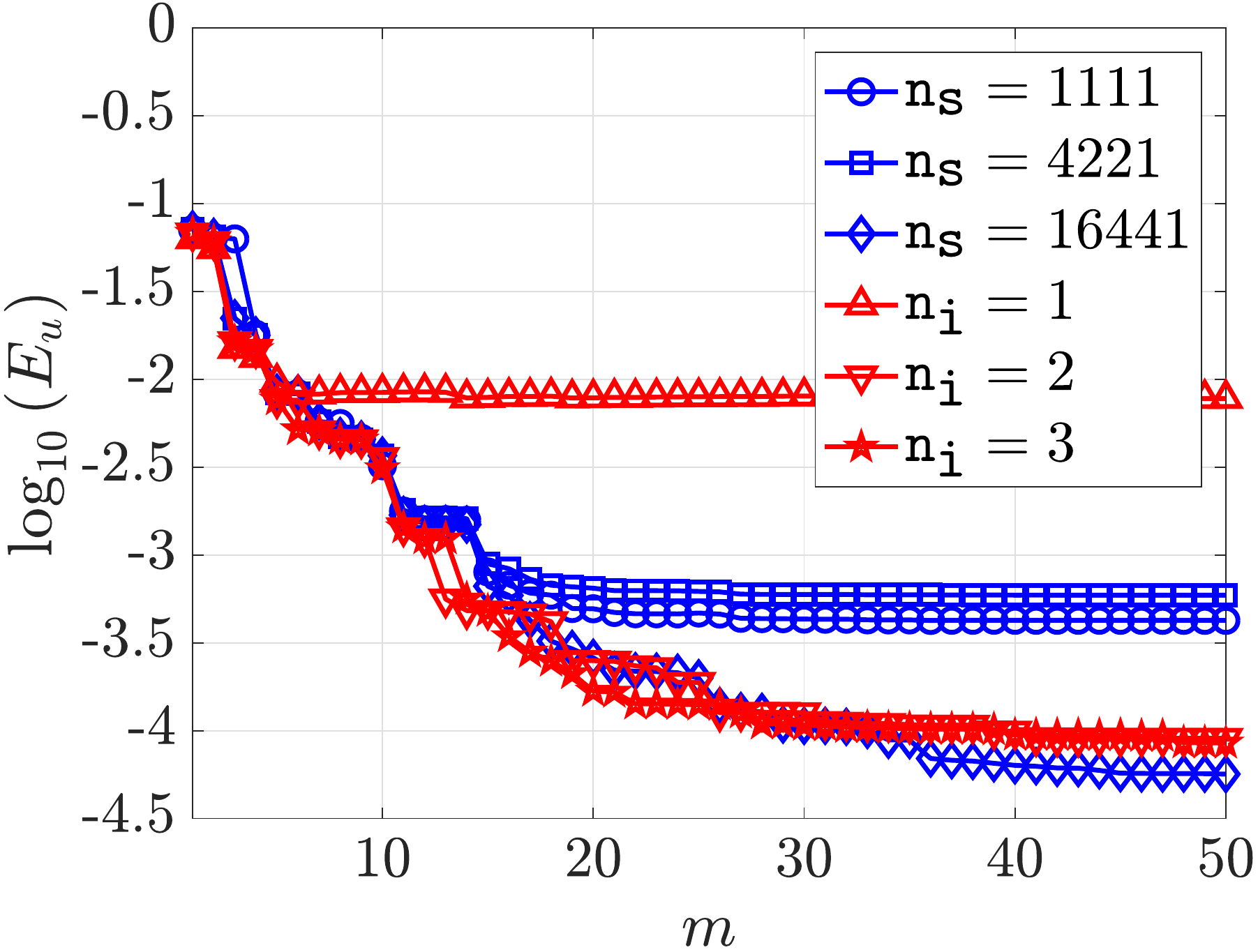}\label{fig:RadisDistance2U}}
	\subfigure[$p$]  {\includegraphics[width=0.49\textwidth]{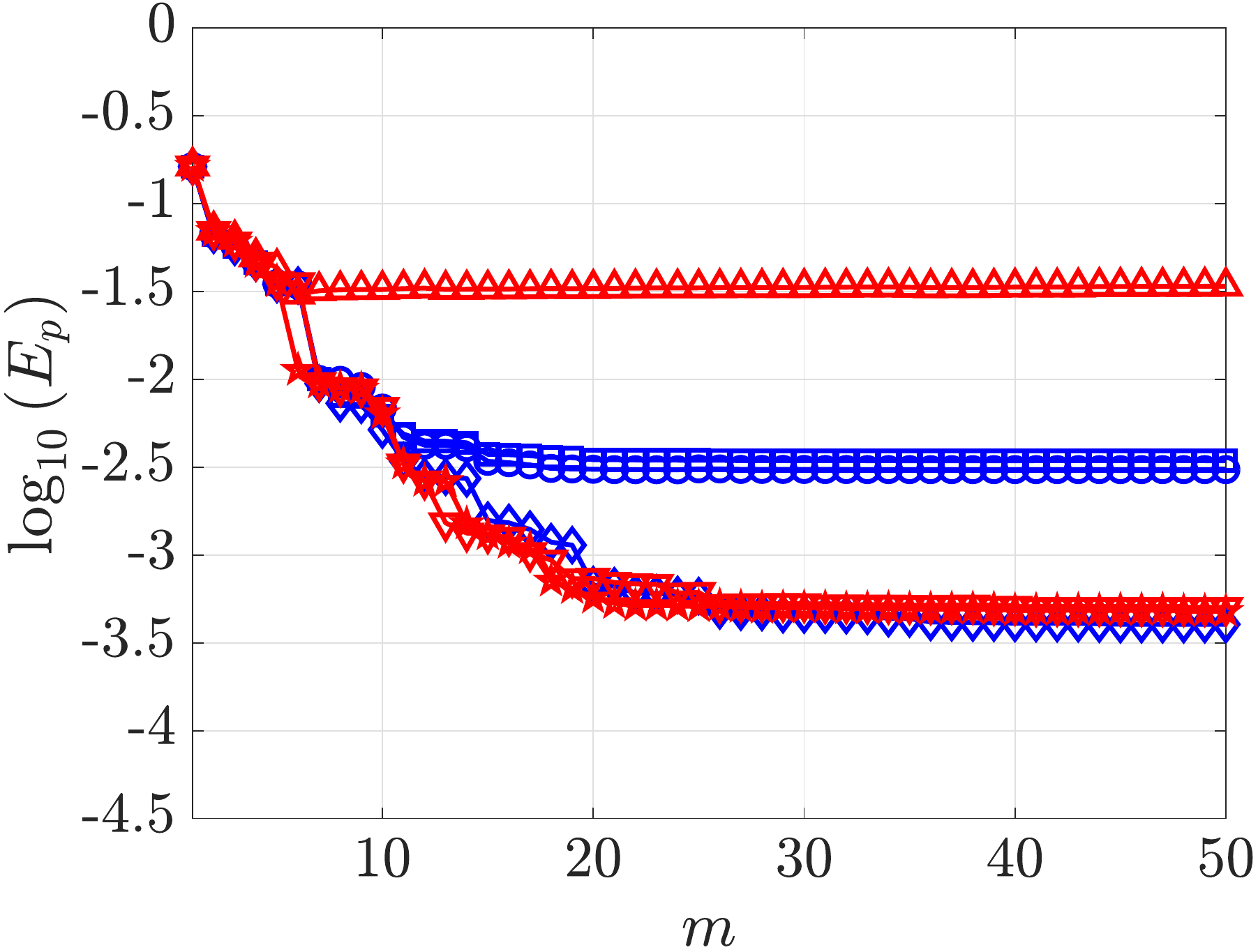}\label{fig:RadisDistance2P}}
	\subfigure[$\bL$]{\includegraphics[width=0.49\textwidth]{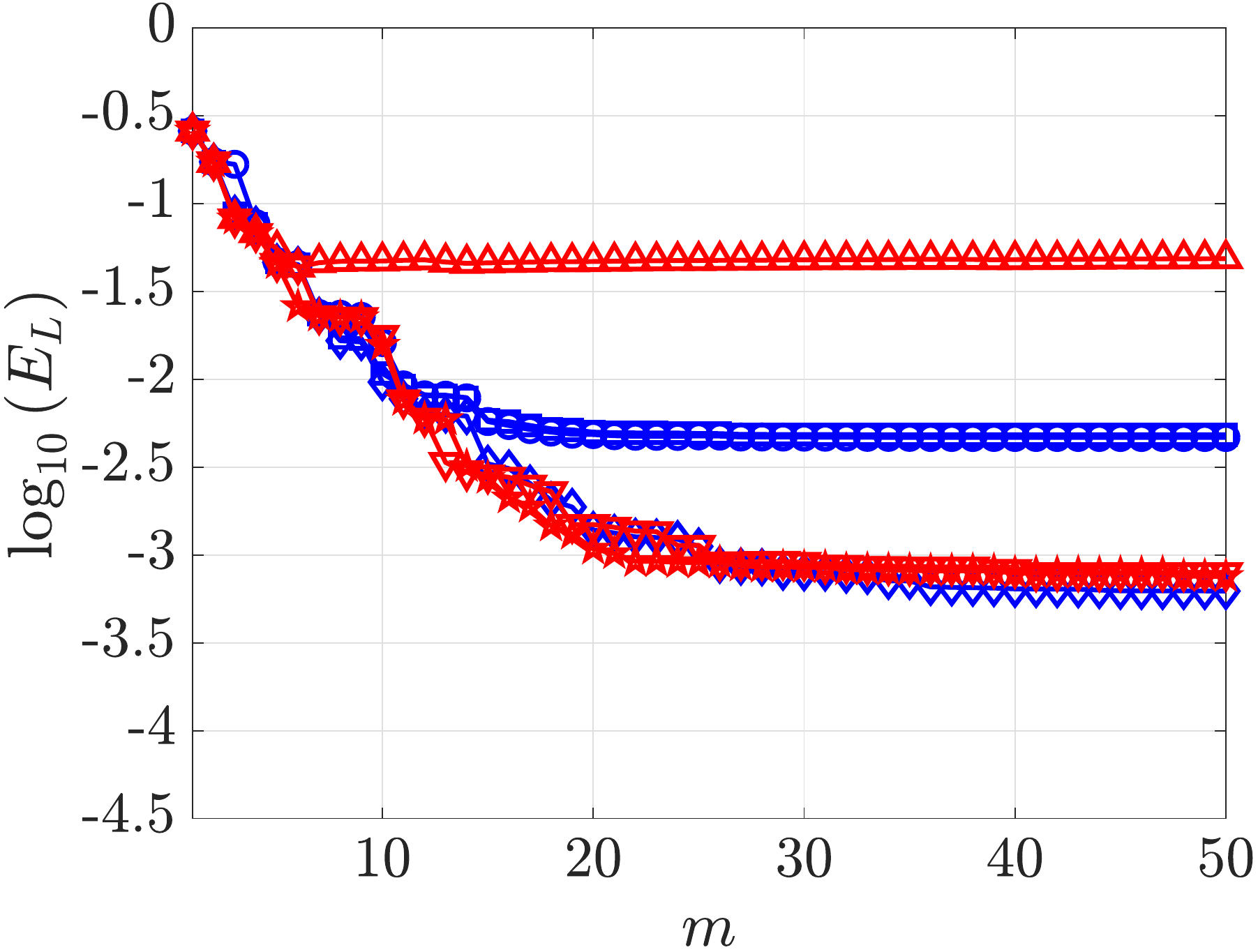}\label{fig:RadisDistance2L}}
	\subfigure[$\FD$]{\includegraphics[width=0.49\textwidth]{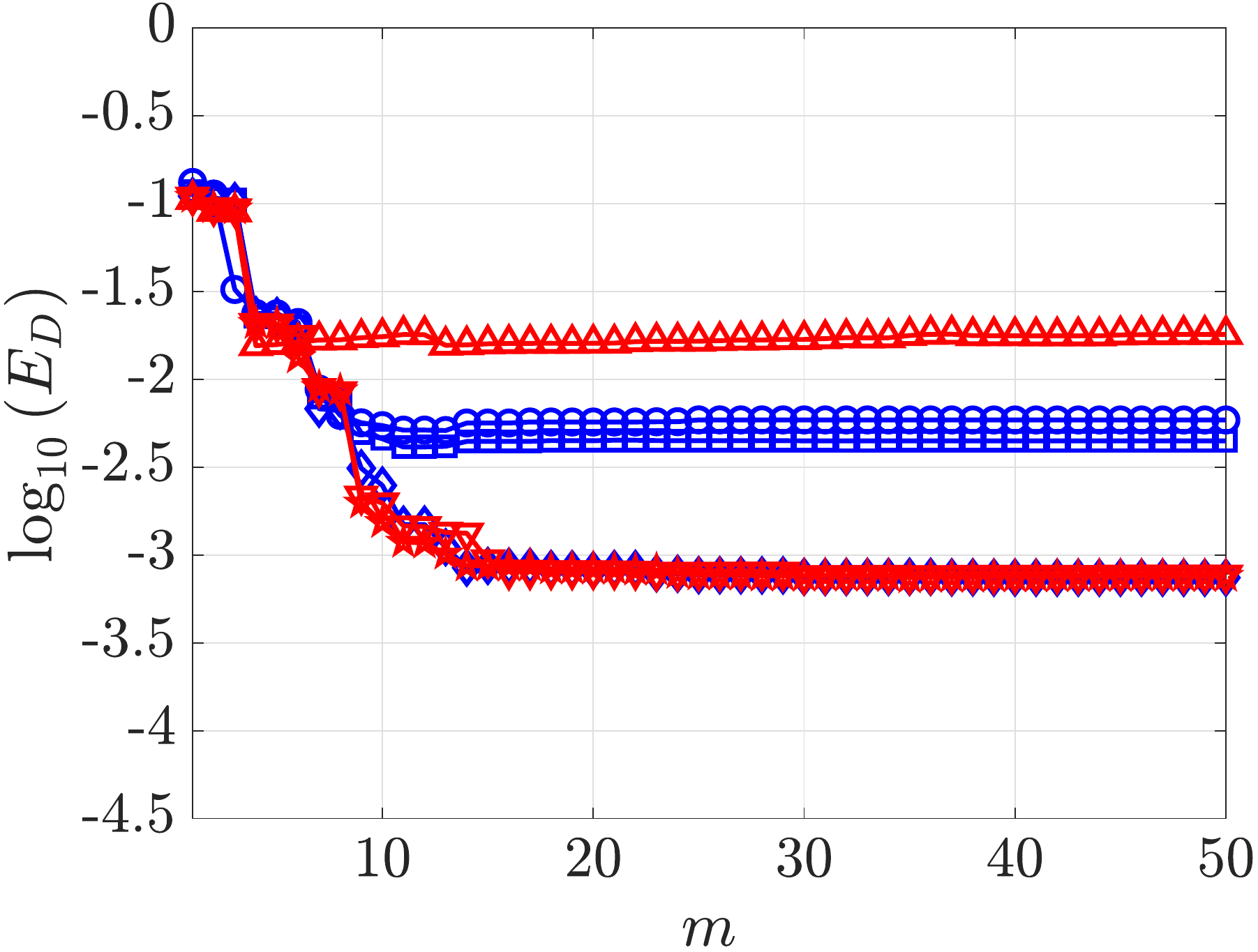}\label{fig:RadisDistance2D}}
	\caption{Evolution of (a-c) the $\eltwo(\Omega \times \bI)$ error for velocity, pressure and gradient of velocity and (d) the $\eltwo(\bI)$ error for the drag force as a function of the number of PGD modes for the problem with two geometric parameters and $\I^2 {=} [-3,2]$. The legend details the number $\nsnap$ of snapshots used by the a posteriori PGD approach (blue) and the number $\niter$ of nonlinear iterations used by the a priori PGD approach (red).}
	\label{fig:RadisDistance2}
\end{figure}

The last example considers the more challenging scenario with two geometric parameters and with the interval for the distance between the bladders equal to $\I^2 {=} [-3,2]$. Figure~\ref{fig:RadisDistance2} reports the evolution of the $\eltwo(\Omega \times \bI)$ error for velocity, pressure and gradient of velocity and the $\eltwo(\bI)$ error for the drag force as a function of the number $m$ of modes. 
The results are qualitatively similar to the previous example but the number of snapshots and modes required by the a posteriori and a priori PGD approaches changes significantly. The a posteriori approach with 1,111 and 4,221 snapshots provide almost identical accuracy in all the variables. In this case, 15 modes are sufficient to provide the maximum accuracy in velocity, pressure and gradient of velocity. One order of magnitude more accurate results are obtained if the number of snapshots is increased to 16,441. In terms of the drag force, with 1,111 snapshots and 10 modes the a posteriori approach is able to provide an accuray below $10^{-2}$. To obtain an accuracy below $10^{-3}$, the a posteriori approach requires 16,441 snapshots and 13 modes.
In this example, the a priori approach with only one nonlinear AD iteration is not able to produce results with an error in the drag force below $10^{-2}$. It is worth noticing that, despite an accurate velocity field is obtained, the error in both pressure and gradient of velocity is higher than $10^{-2}$. However, by performing only two nonlinear iterations in the AD scheme and computing enough modes, the error in the velocity field drops of two orders of magnitude and accurate results are obtained for both pressure and gradient of velocity, with an error below $10^{-3}$.

Similarly to the procedure described for the previous example, the a priori approach is here set to compute a fixed number of modes, namely 500 modes. The stagnation of the curves in figure~\ref{fig:RadisDistance2} confirms that this number is larger than what is required for convergence. The PGD approximation is then compressed in 50 modes, with the most relevant information concentrated in less than 30. To achieve an error in the drag force below $10^{-3}$, two nonlinear iterations are employed. To obtain the same accuracy, the a posteriori approach requires 14 modes but the number of snapshots needed for this challenging problem is 16,441. This means that the a posteriori approach requires 11 times extra spatial solutions to provide the same error as the a priori approach. The results illustrate again that the higher the accuracy required and the higher the variability in the solution introduced by the geometric parameters, the more beneficial is the use of the a priori approach. It is worth recalling that the presented results compare two basic versions of the a priori and a posteriori PGD algorithms which could be improved by introducing techniques to handle the space of parameters and error control strategies. While these techniques represent frontier research in the context of a priori ROMs, they are well established for a posteriori strategies. In particular, the employment of sampling methods is expected to reduce the number of snapshots required by the a posteriori PGD, whereas error control will provide more accurate information on the capability of the reduced basis to represent the multidimensional solution. These studies, which are out of the scope of the present work, represent promising lines of investigation to better understand advantages and disadvantages of different ROMs strategies.

\begin{figure}[!tb]
	\centering
	\subfigure[Module of velocity, $(\mu_1=-1,\mu_2=-3)$]{\includegraphics[width=0.49\textwidth]{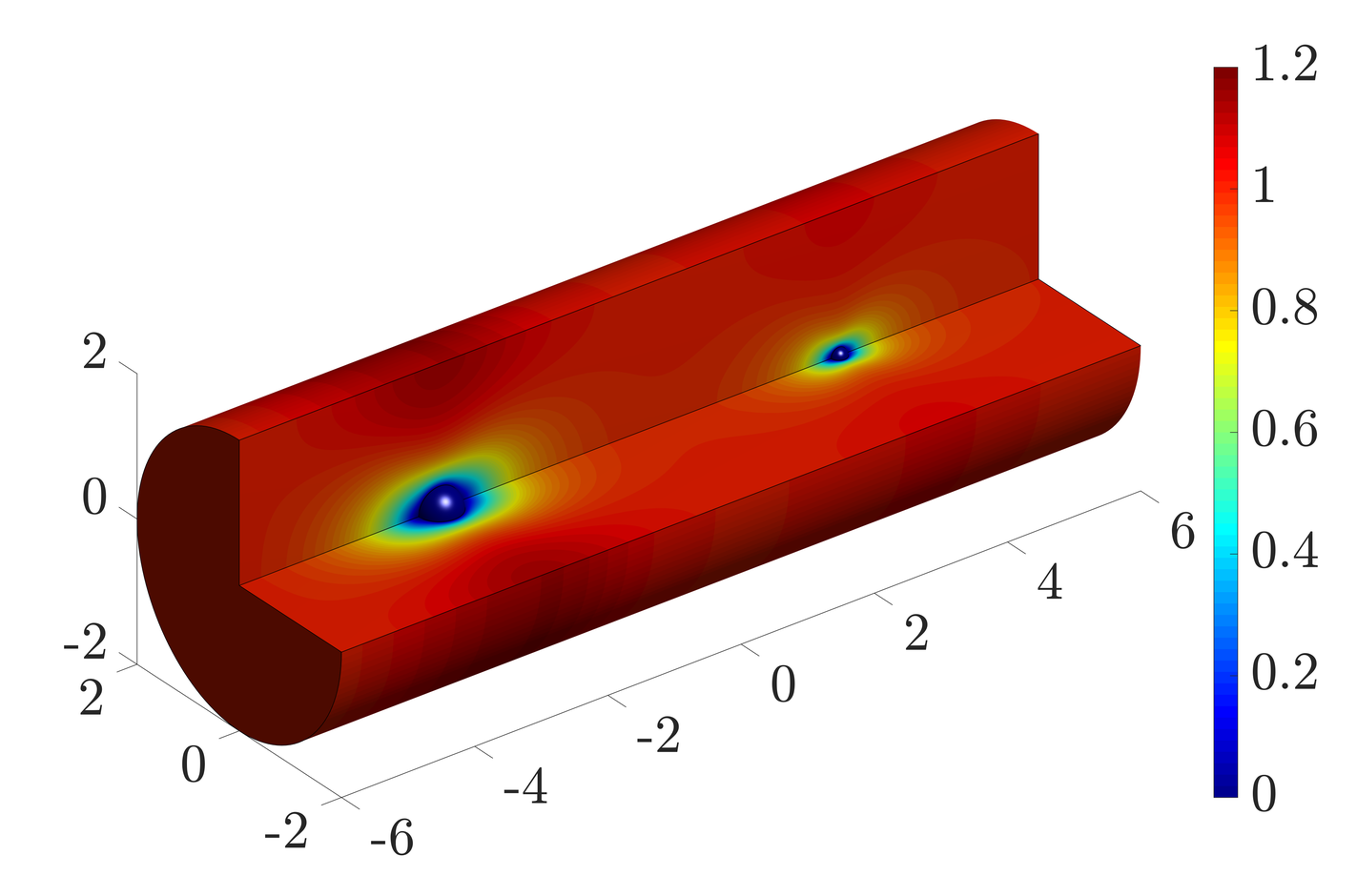}}
	\subfigure[Pressure, $(\mu_1=-1,\mu_2=-3)$]{\includegraphics[width=0.49\textwidth]{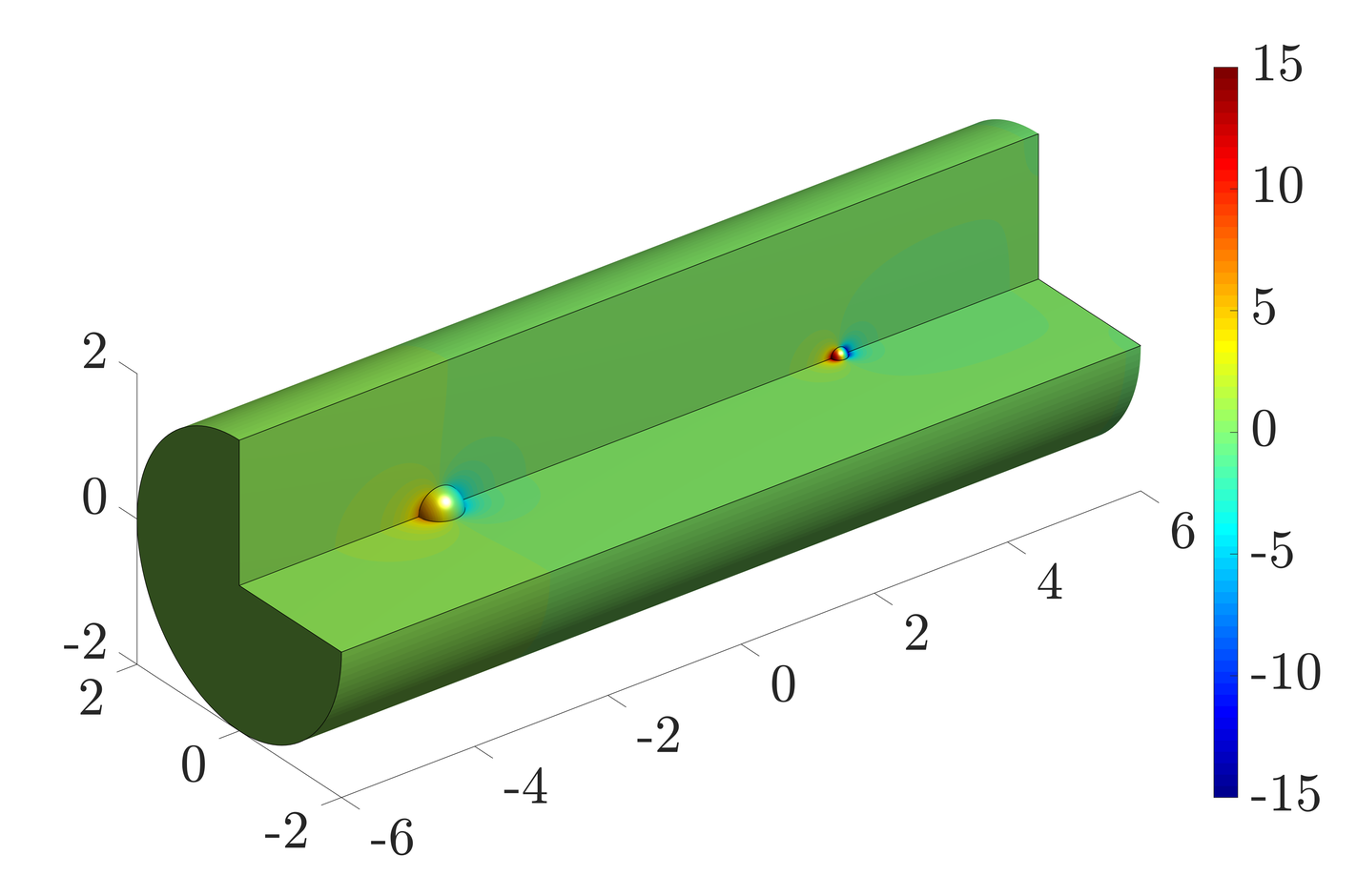}}
	
	\subfigure[Module of velocity, $(\mu_1=1,\mu_2=2)$]  {\includegraphics[width=0.49\textwidth]{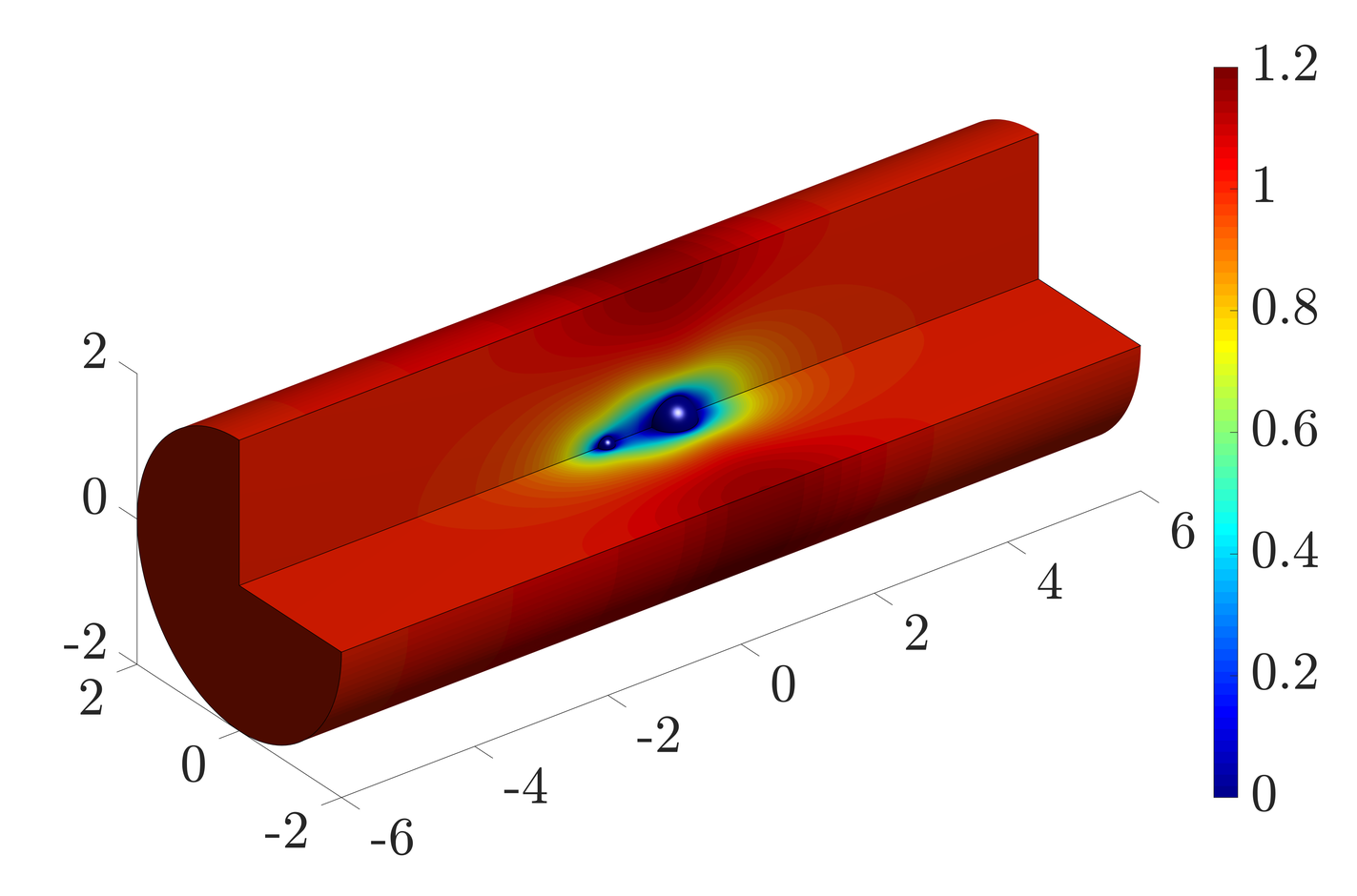}}
	\subfigure[Pressure, $(\mu_1=1,\mu_2=2)$]  {\includegraphics[width=0.49\textwidth]{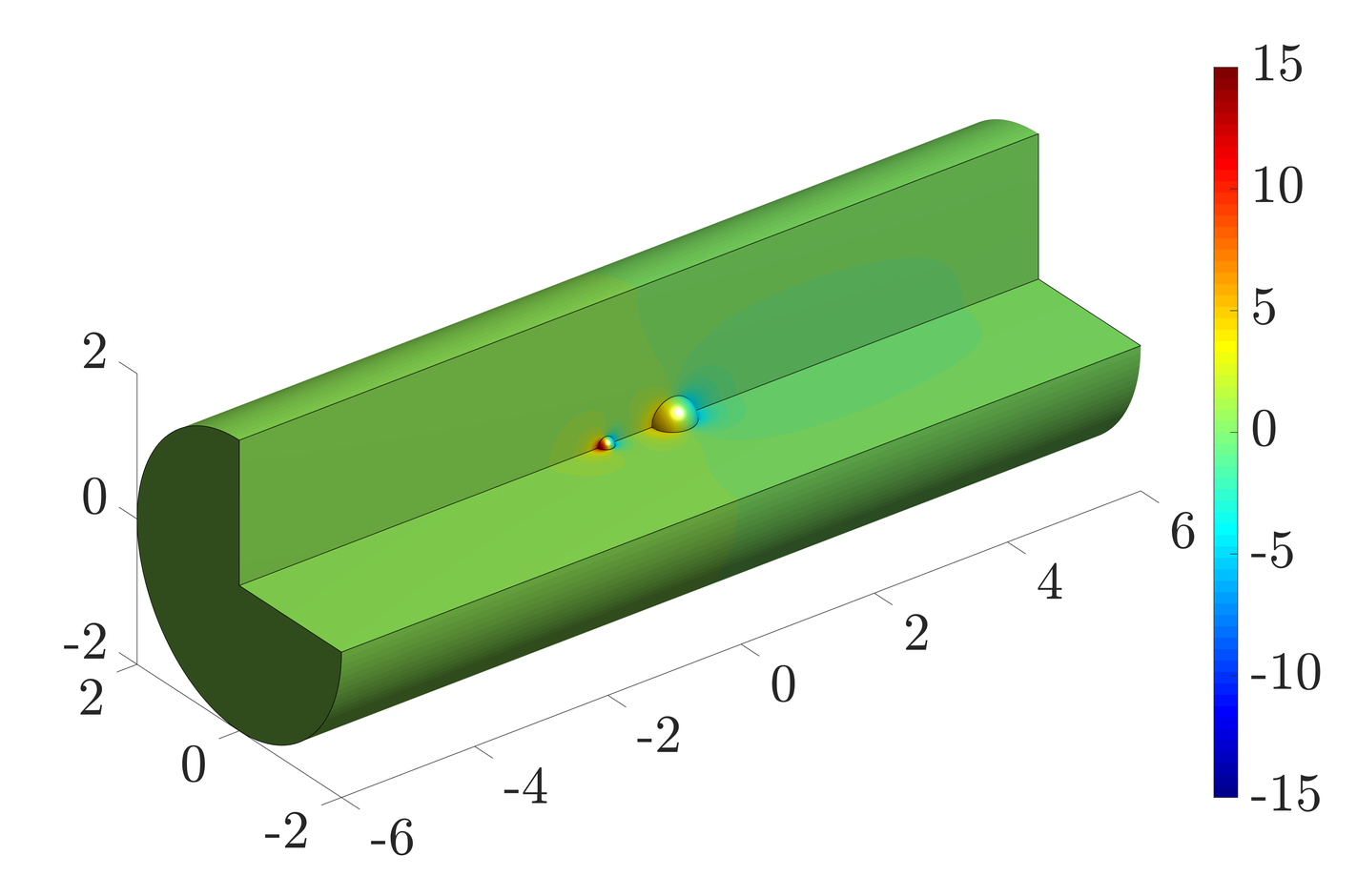}}
	
	\caption{Module of velocity and pressure field for two sets of parameters $\mu_1$ and $\mu_2$ describing the extremal configurations of the microswimmer.}
	\label{fig:RadisDistanceOnline}
\end{figure}

For the sake of brevity, only real-time evaluations of the velocity and pressure fields obtained from the a priori PGD computation are reported hereafter. Interested readers are referred to~\cite{RS-SBGH-20} for more in-depth presentation of the spatial and parametric modes for the problem under analysis.
To illustrate the online stage, figure~\ref{fig:RadisDistanceOnline} reports the velocity and pressure fields corresponding to the two extremal configurations of the \emph{push-me-pull-you} microswimmers described by the parameters $\mu_1$ and $\mu_2$. 

%-------------------------------------------------------
\subsubsection{Accuracy of a priori and a posteriori response surfaces} 
%-------------------------------------------------------

%To illustrate the online stage, figures~\ref{fig:RadisDistance2UOnline} and~\ref{fig:RadisDistance2POnline} show the velocity and pressure fields, respectively, corresponding to two different values of the parameters $\mu_1$ and $\mu_2$. 
%%
%\begin{figure}[!tb]
%	\centering
%	\subfigure[$\mu_1=-1$, $\mu_2=-3$]{\includegraphics[width=0.49\textwidth]{Param12_Int2_Velocity_Mu_M1M3}}
%	\subfigure[$\mu_1=1$, $\mu_2=2$]  {\includegraphics[width=0.49\textwidth]{Param12_Int2_Velocity_Mu_12}}
%	\caption{Velocity field for two values of the parameters $\mu_1$ and $\mu_2$.}
%	\label{fig:RadisDistance2UOnline}
%\end{figure}
%%
%\begin{figure}[!tb]
%	\centering
%	\subfigure[$\mu_1=-1$, $\mu_2=-3$]{\includegraphics[width=0.49\textwidth]{Param12_Int2_Pressure_Mu_M1M3}}
%	\subfigure[$\mu_1=1$, $\mu_2=2$]  {\includegraphics[width=0.49\textwidth]{Param12_Int2_Pressure_Mu_12}}
%	\caption{Pressure field for two values of the parameters $\mu_1$ and $\mu_2$.}
%	\label{fig:RadisDistance2POnline}
%\end{figure}

\begin{figure}[!tb]
	\centering
	\subfigure[{$\I^2 = [-2,-1]$}]{\includegraphics[width=0.48\textwidth]{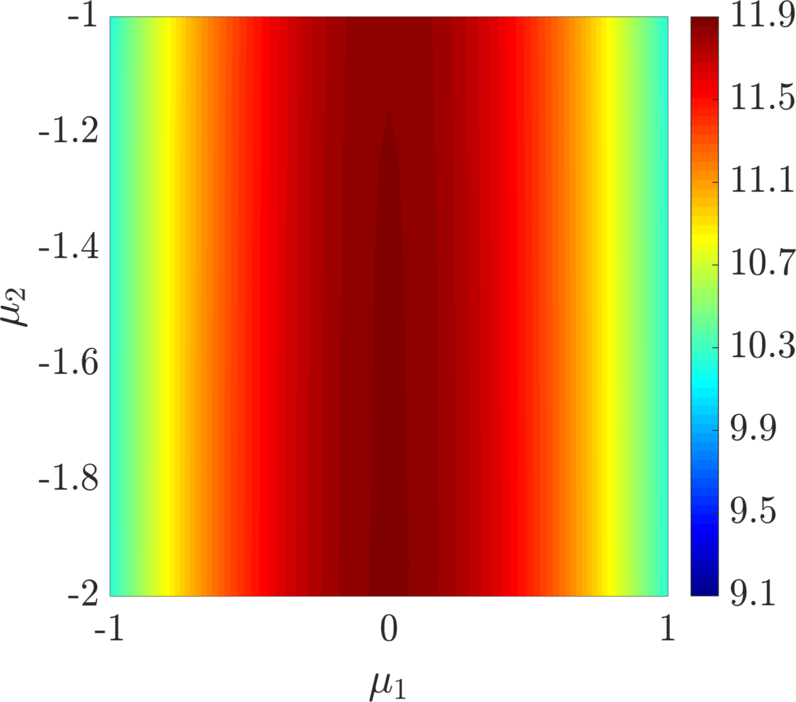}}
	\subfigure[{$\I^2 = [-3, 2]$}]{\includegraphics[width=0.48\textwidth]{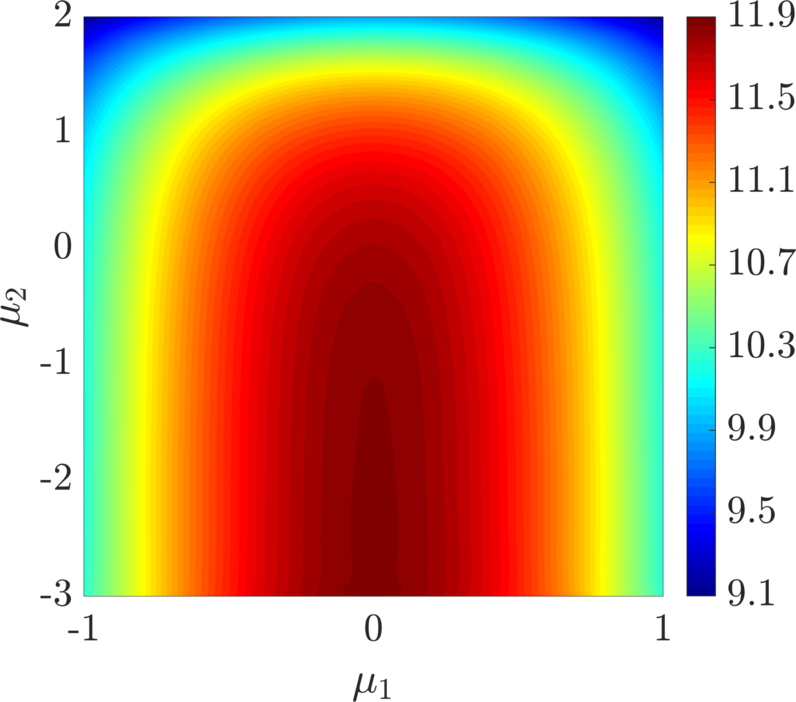}\label{fig:Drag2Dlarge}}
	
	\caption{Response surfaces of the total drag force as a function of the radius $\mu_1$ of the first sphere and the distance $\mu_2$ between the two bladders, for two different ranges of values of the parameter $\mu_2$.}
	\label{fig:Drag2D}
\end{figure}

The separated response surfaces for the total drag force on the spheres computed using the a priori PGD are presented in figure~\ref{fig:Drag2D}, as a function of the parameters $\mu_1$ and $\mu_2$.
The results confirm the increased sensitivity of the quantity of interest to the extended range of the parameter $\mu_2$, as already observed in figure~\ref{fig:DragDistance}, with the appearence of localised variations of the drag force in the vicinity of the value $\mu_2 {=} 2$ (Fig.~\ref{fig:Drag2Dlarge}). 

\begin{figure}[!tb]
	\centering
%	\subfigure[$\varepsilon_D$, a priori PGD]		{\includegraphics[width=0.49\textwidth]{Param12_Int1_ErrD_Prio}}
		\subfigure[$\varepsilon_D$, a priori PGD]	{\includegraphics[width=0.49\textwidth]{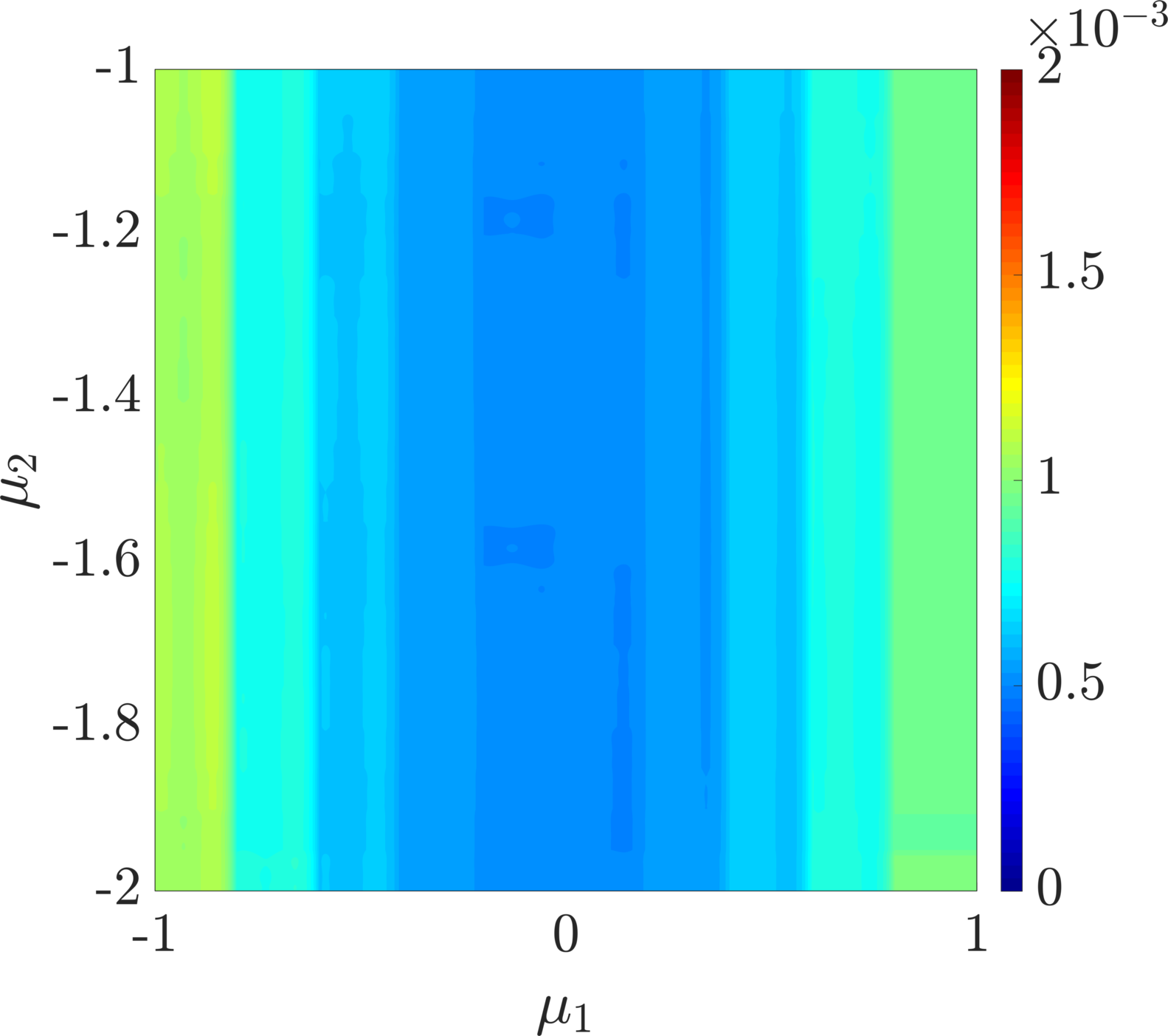}}
%	\subfigure[$\varepsilon_D$, a posteriori PGD]	{\includegraphics[width=0.49\textwidth]{Param12_Int1_ErrD_Post}}
	\subfigure[$\varepsilon_D$, a posteriori PGD]	{\includegraphics[width=0.49\textwidth]{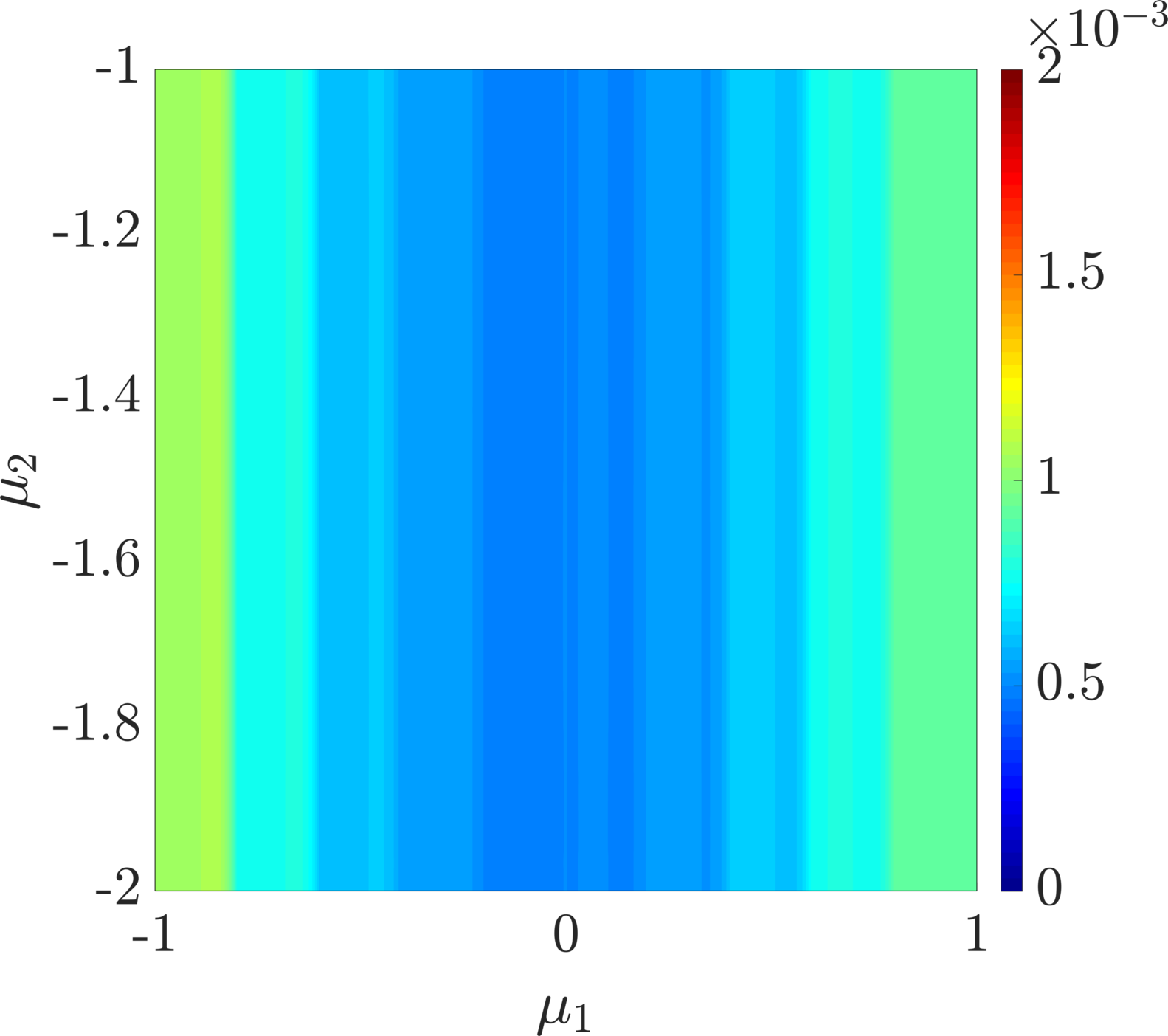}}
	\caption{Error map for the drag force as a function of the two parameters $\mu_1 \in [-1,1]$ and $\mu_2 \in [-2,-1]$.}
	\label{fig:ErrorRadisDistance1Drag}
\end{figure}

The previous two examples with two geometric parameters have shown that the a priori PGD approach is competitive when the multidimensional error measure in equation~\eqref{eq:DragErrorMulti} is considered. To further analyse the performance of both PGD approaches, figure~\ref{fig:ErrorRadisDistance1Drag} reports the smoothed pointwise error of the drag force as a function of the two parameters $\mu_1$ and $\mu_2$ for the first example in this section, when the second parameter belongs to the interval $\I^2 {=} [-2,-1]$.
The results show that both the a priori and the a posteriori approaches produce almost identical results for each value of the two geometric parameters. The behaviour is very similar to the one observed for the solution with only one parameter, as reported in figures~\ref{fig:ErrD_param1Smooth} and \ref{fig:ErrD_param2A}. A slightly higher error is observed for the a priori PGD near the left and right boundaries of the parametric domain, corresponding to the maximum and minimum radius of the first sphere, respectively. In addition, the accuracy obtained is almost independent of the value of the second parameter. This is attributed to the fact that, with the interval $\I^2 {=} [-2,-1]$ considered here, the minimum distance between the spheres does not induce a significant variation of the flow impinging onto the second sphere. 

\begin{figure}[!tb]
	\centering
%	\subfigure[$\varepsilon_D$, a priori PGD]		{\includegraphics[width=0.49\textwidth]{Param12_Int2_ErrD_Prio}}
	\subfigure[$\varepsilon_D$, a priori PGD]	{\includegraphics[width=0.49\textwidth]{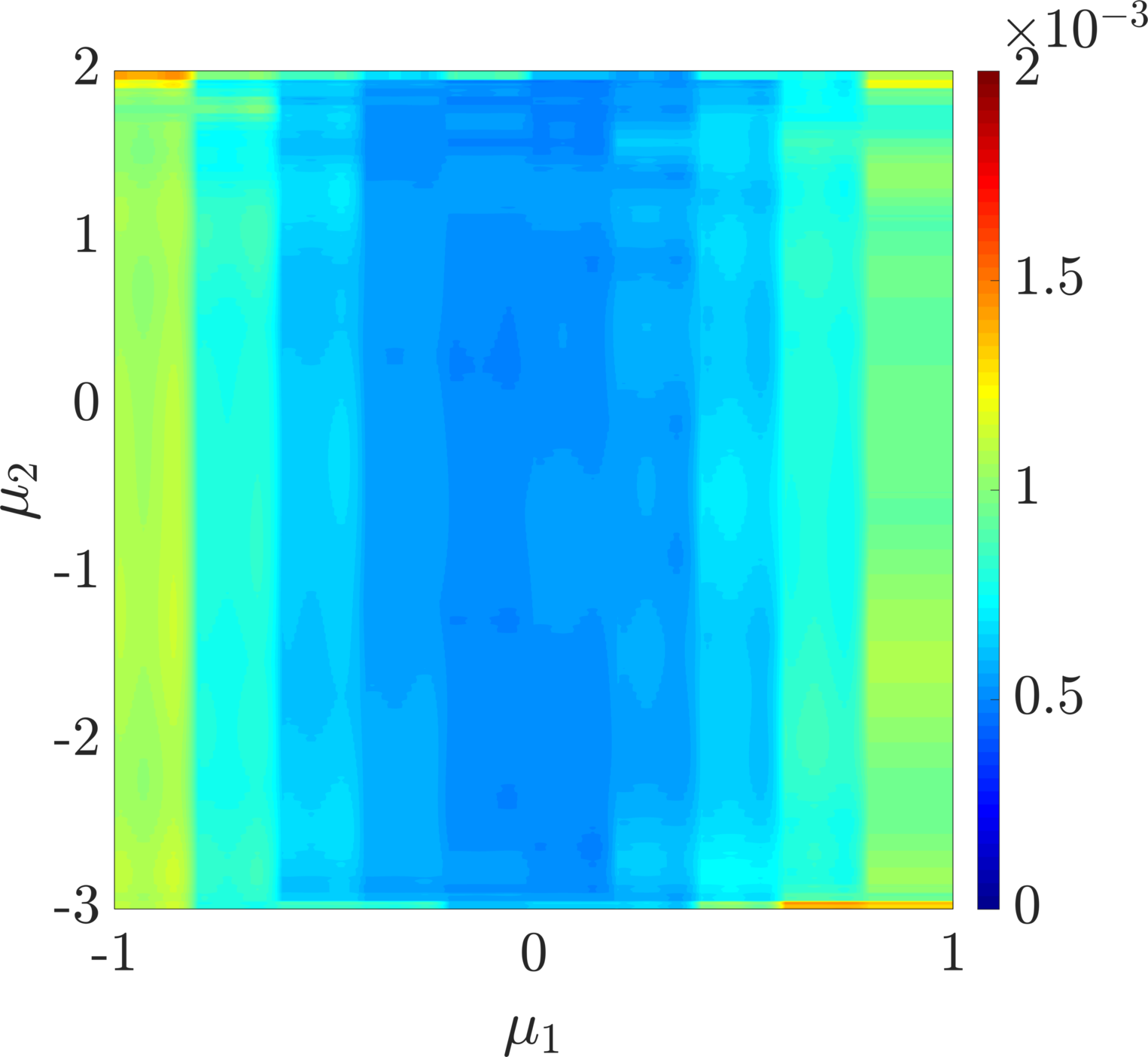}}
%	\subfigure[$\varepsilon_D$, a posteriori PGD]	{\includegraphics[width=0.49\textwidth]{Param12_Int2_ErrD_Post}}
	\subfigure[$\varepsilon_D$, a posteriori PGD]	{\includegraphics[width=0.49\textwidth]{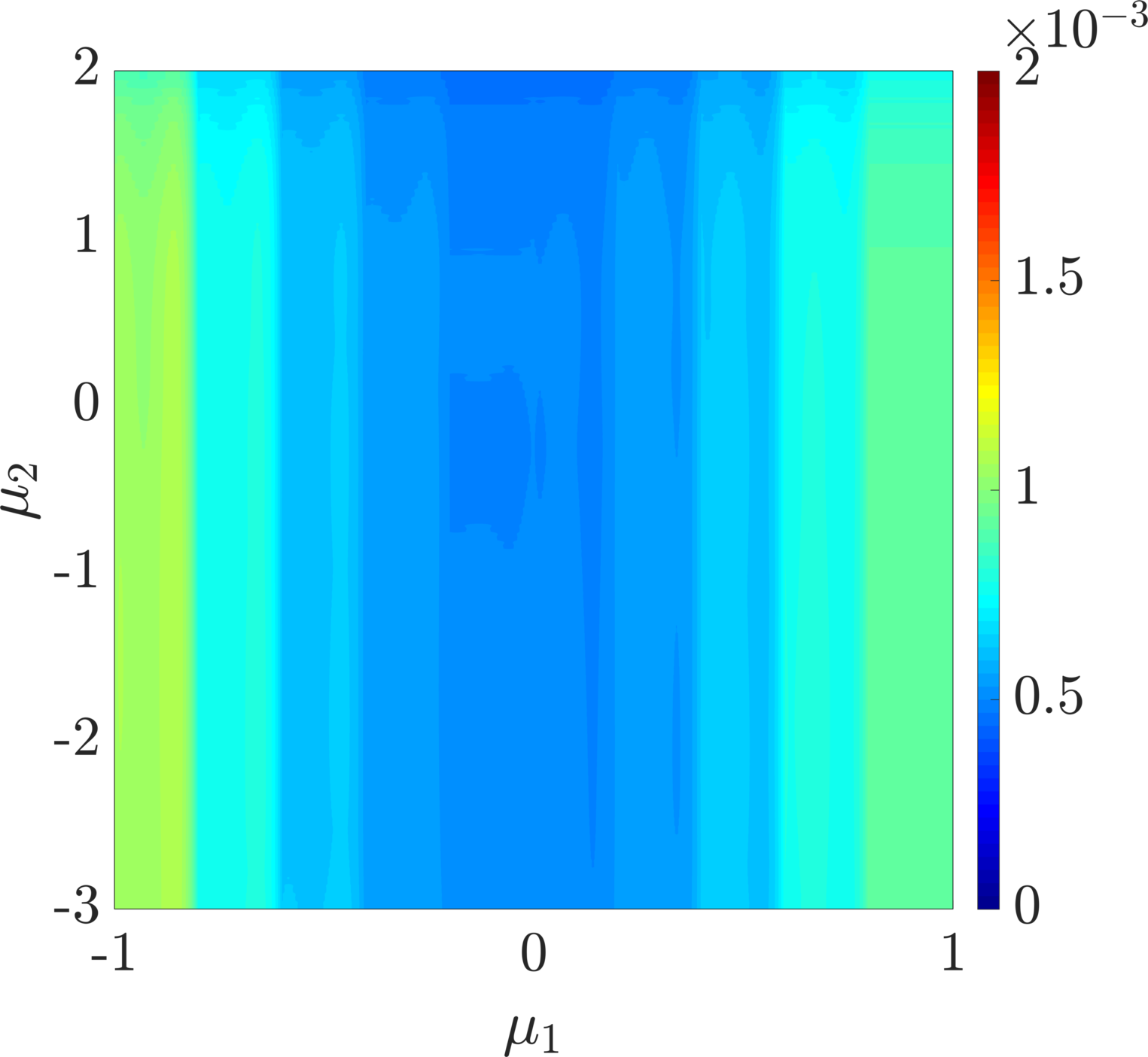}}
	\caption{Error map for the drag force as a function of the two parameters $\mu_1 \in [-1,1]$ and $\mu_2 \in [-3,2]$.}
	\label{fig:ErrorRadisDistance2Drag}
\end{figure}

The same study is repated for the case of $\I^2 {=} [-3,2]$. Figure~\ref{fig:ErrorRadisDistance2Drag} shows the smoothed error of the drag force as a function of the two parameters $\mu_1$ and $\mu_2$ for the second example, with $\mu_2 \in [-3,2]$.
Despite the $\eltwo(\bI)$ error measure is almost identical for the a priori and the a posteriori PGD approaches (Fig.~\ref{fig:RadisDistance2D}), the error map displays important differences between the two methods. More precisely, the error map of the a priori approach reveals higher error in the vicinity of the boundary of $\bI$, whereas the error of the a posteriori PGD does not show such increase near the boundary. It is worth noticing that the higher errors observed in the a priori approach are very localised and therefore they are not observed when the $\eltwo(\bI)$ error measure is computed. In addition, the higher errors are not only observed for the maximum value of the parameter $\mu_2$ but also for lower values of $\mu_2$. 

This result reveals the increased difficulty in addressing problems with more than one geometric parameter with the a priori approach. Furthermore, the study shows that the conclusions of independent studies with only one geometric parameter do not extend to problems with the same parameters considered in a single simulation.

%%%%%%%%%%%%%%%%%%%%%%%%%%%%%%%%%%%%%%%%%%%%%%
\section{Concluding remarks}
\label{sc:conclusion}
%%%%%%%%%%%%%%%%%%%%%%%%%%%%%%%%%%%%%%%%%%%%%%

A comparison of a priori and a posteriori PGD algorithms was presented for the challenging problem of an incompressible Stokes flow in geometrically parametrised domains. The full-order solver is based on a multidimensional HDG method which allows the use of equal order polynomial approximations for all the variables leading to an LBB-compliant discretisation with high-order isoparametric formulations. In addition, the HDG-PGD framework provides an exact separation of the terms appearing in the formulation of the geometrically parametrised PDE on a reference domain.

The a priori PGD algorithm, see~\cite{RS-SBGH-20}, is obtained devising a separated formulation of the multidimensional HDG solver and does not require any prior snapshot computation. The a posteriori PGD, also known as least-squares PGD~\cite{PD-DZGH-20}, constructs a separated approximation starting from a series of snapshots obtained as full-order solutions of the spatial problem.

The challenging problem of the flow around a geometrically parametrised \emph{push-me-pull-you} microswimmer is employed to test the performance of the two PGD approaches. More precisely, extensive numerical experiments are performed to test the sensitivity of the ROMs to the range of variation of the parameters and to the number of parameters considered. For problems with a unique geometric parameter inducing limited variations in the flow, accuracy and performance of the two approaches are comparable. When the range of variation of the parameter is extended and extreme geometric transformations are considered, the a priori approach requires a significantly lower computational cost, measured in terms of number of full-order HDG solves, with respect to the a posteriori PGD. Nonetheless, it is worth recalling that snapshots in the a posteriori PGD can be easily computed in parallel, whereas the computation of the modes in the a priori approach is sequential. The presented numerical results also highlight the additional difficulty introduced by the presence of multiple geometric parameters. In this case, the a priori PGD requires up to 11 times less calls to the full-order HDG solver than the a posteriori approach to achieve an accuracy in the drag force below $10^{-3}$. Hence, the numerical results display the superior performance of the a priori PGD when either the parametric solution features increased variability, due to the range of values of the geometric parameters or to their number, or higher accuracy is required by the user. Nonetheless, the employment of advanced sampling techniques in the a posteriori PGD is expected to reduce the number of required snapshots and competitive computing times may be achieved via their parallel computation.

It is worth emphasising that the conclusions of parametric studies considering only one geometric parameter at a time do not extend to parametric PDEs in which multiple parameters are concurrently considered in a unique problem. Hence, prompted by these results, further comparisons involving problems with more than two geometric parameters and a posteriori algorithms based on tailored sampling techniques are expected to provide additional insights on the applicability and limitations of PGD-based strategies in the context of industrial parametric studies.

%______________________________________________________________________
\section*{Acknowledgements}

This work was partially supported by the European Union's Horizon 2020 research and innovation programme under the Marie Sk\l odowska-Curie Actions (Grant number: 675919) that financed the Ph.D. fellowship of L.B. and by the Spanish Ministry of Economy and Competitiveness (Grant number: DPI2017-85139-C2-2-R). M.G. and A.H. are also grateful for the support provided by the Spanish Ministry of Economy and Competitiveness through the Severo Ochoa programme for centres of excellence in RTD (Grant number: CEX2018-000797-S) and the Generalitat de Catalunya (Grant number: 2017-SGR-1278). R.S. also acknowledges the support of the Engineering and Physical Sciences Research Council (Grant number: EP/P033997/1).

%________________________________________________________________________
%________________________________________________________________________
%\bibliographystyle{amsplain}
%\bibliographystyle{abbrv}
\bibliographystyle{elsarticle-num}
%\bibliography{Ref-PGD-Geo}
\bibliography{Ref-PGD-Comparison}

\begin{thebibliography}{10}
\expandafter\ifx\csname url\endcsname\relax
  \def\url#1{\texttt{#1}}\fi
\expandafter\ifx\csname urlprefix\endcsname\relax\def\urlprefix{URL }\fi
\expandafter\ifx\csname href\endcsname\relax
  \def\href#1#2{#2} \def\path#1{#1}\fi

\bibitem{Box-BW-51}
G.~Box, K.~Wilson, On the experimental attainment of optimum conditions,
  Journal of the royal statistical society: Series b (Methodological) 13~(1)
  (1951) 1--38.

\bibitem{Breitkopf-BNRV-05}
P.~Breitkopf, H.~Naceur, A.~Rassineux, P.~Villon, Moving least squares response
  surface approximation: formulation and metal forming applications, Computers
  \& Structures 83~(17-18) (2005) 1411--1428.

\bibitem{Breitkopf-ZBKZ-11}
P.~Zhang, P.~Breitkopf, C.~Knopf-Lenoir, W.~Zhang, Diffuse response surface
  model based on moving {L}atin hypercube patterns for reliability-based design
  optimization of ultrahigh strength steel {NC} milling parameters, Structural
  and Multidisciplinary Optimization 44~(5) (2011) 613--628.

\bibitem{AH-CHRW:17}
F.~Chinesta, A.~Huerta, G.~Rozza, K.~Willcox, Model {R}eduction {M}ethods, in:
  E.~Stein, R.~{de Borst}, T.~Hughes (Eds.), Encyclopedia of Computational
  Mechanics Second Edition, Vol. Part 1 Solids and Structures, John Wiley \&
  Sons, Ltd., Chichester, 2017, Ch.~3, pp. 1--36.

\bibitem{Gunzburger-PWG-18}
B.~Peherstorfer, K.~Willcox, M.~Gunzburger, Survey of multifidelity methods in
  uncertainty propagation, inference, and optimization, SIAM Review 60~(3)
  (2018) 550--591.

\bibitem{Alonso-LA-00}
P.~LeGresley, J.~Alonso, Airfoil design optimization using reduced order models
  based on proper orthogonal decomposition, in: AIAA Fluids 2000 conference and
  exhibit, 2000, p. 2545.

\bibitem{Maute-WEM-09}
G.~Weickum, M.~Eldred, K.~Maute, A multi-point reduced-order modeling approach
  of transient structural dynamics with application to robust design
  optimization, Structural and Multidisciplinary Optimization 38~(6) (2009)
  599.

\bibitem{Breitkopf-XBCKSV-10}
M.~Xiao, P.~Breitkopf, R.~Coelho, C.~Knopf-Lenoir, M.~Sidorkiewicz, P.~Villon,
  Model reduction by {CPOD} and {K}riging, Structural and Multidisciplinary
  Optimization 41~(4) (2010) 555--574.

\bibitem{Carlberg-CF-11}
K.~Carlberg, C.~Farhat, A low-cost, goal-oriented ‘compact proper orthogonal
  decomposition’ basis for model reduction of static systems, International
  Journal for Numerical Methods in Engineering 86~(3) (2011) 381--402.

\bibitem{Rozza-MQR-12}
A.~Manzoni, A.~Quarteroni, G.~Rozza, Shape optimization for viscous flows by
  reduced basis methods and free-form deformation, International Journal for
  Numerical Methods in Fluids 70~(5) (2012) 646--670.

\bibitem{Passieux-GP-13}
C.~Gogu, J.-C. Passieux, Efficient surrogate construction by combining response
  surface methodology and reduced order modeling, Structural and
  Multidisciplinary Optimization 47~(6) (2013) 821--837.

\bibitem{Zahr-ZF-15}
M.~J. Zahr, C.~Farhat, Progressive construction of a parametric reduced-order
  model for {PDE}-constrained optimization, International Journal for Numerical
  Methods in Engineering 102~(5) (2015) 1111--1135.

\bibitem{Breitkopf-Coelho}
P.~Breitkopf, R.~Coelho, Multidisciplinary design optimization in computational
  mechanics, John Wiley \& Sons, 2013.

\bibitem{Conover-79-Latin}
M.~McKay, R.~Beckman, W.~Conover, A comparison of three methods for selecting
  values of input variables in the analysis of output from a computer code,
  Technometrics 21~(1) (1979) 239--245.

\bibitem{Gunzburger-DFG-99}
Q.~Du, V.~Faber, M.~Gunzburger, Centroidal {V}oronoi tessellations:
  {A}pplications and algorithms, SIAM Review 41~(4) (1999) 637--676.

\bibitem{Patera-GP-05}
M.~Grepl, A.~Patera, A posteriori error bounds for reduced-basis approximations
  of parametrized parabolic partial differential equations, ESAIM: Mathematical
  Modelling and Numerical Analysis 39~(1) (2005) 157--181.

\bibitem{Patera-VP-05}
K.~Veroy, A.~Patera, Certified real-time solution of the parametrized steady
  incompressible {N}avier--{S}tokes equations: rigorous reduced-basis a
  posteriori error bounds, International Journal for Numerical Methods in
  Fluids 47~(8-9) (2005) 773--788.

\bibitem{Willcox-BWG-08}
T.~Bui-Thanh, K.~Willcox, O.~Ghattas, Model reduction for large-scale systems
  with high-dimensional parametric input space, SIAM Journal on Scientific
  Computing 30~(6) (2008) 3270--3288.

\bibitem{Ryckelynck-05}
D.~Ryckelynck, A priori hyperreduction method: an adaptive approach, Journal of
  Computational Physics 202~(1) (2005) 346--366.

\bibitem{Sorensen-CS-10}
S.~Chaturantabut, D.~Sorensen, Nonlinear model reduction via discrete empirical
  interpolation, SIAM Journal on Scientific Computing 32~(5) (2010) 2737--2764.

\bibitem{Carlberg-CBF-11}
K.~Carlberg, C.~Bou-Mosleh, C.~Farhat, Efficient non-linear model reduction via
  a least-squares {P}etrov--{G}alerkin projection and compressive tensor
  approximations, International Journal for Numerical Methods in Engineering
  86~(2) (2011) 155--181.

\bibitem{Amsallem-AZF-12}
D.~Amsallem, M.~Zahr, C.~Farhat, Nonlinear model order reduction based on local
  reduced-order bases, International Journal for Numerical Methods in
  Engineering 92~(10) (2012) 891--916.

\bibitem{Carlberg-CFCA-13}
K.~Carlberg, C.~Farhat, J.~Cortial, D.~Amsallem, The {GNAT} method for
  nonlinear model reduction: effective implementation and application to
  computational fluid dynamics and turbulent flows, Journal of Computational
  Physics 242 (2013) 623--647.

\bibitem{Breitkopf-PBBVZ-20}
P.~Phalippou, S.~Bouabdallah, P.~Breitkopf, P.~Villon, M.~Zarroug,
  `{O}n-the-fly' snapshots selection for {P}roper {O}rthogonal {D}ecomposition
  with application to nonlinear dynamics, Computer Methods in Applied Mechanics
  and Engineering 367 (2020) 113120.

\bibitem{Chinesta-Keunings-Leygue}
F.~Chinesta, R.~Keunings, A.~Leygue, The proper generalized decomposition for
  advanced numerical simulations. {A} primer, Springer Briefs in Applied
  Sciences and Technology, Springer, Cham, 2014.

\bibitem{AH-CLBACGAAH-13}
F.~Chinesta, A.~Leygue, F.~Bordeu, J.~Aguado, E.~Cueto, D.~Gonz{\'a}lez,
  I.~Alfaro, A.~Ammar, A.~Huerta, {PGD}-based computational vademecum for
  efficient design, optimization and control, Archives of Computational Methods
  in Engineering 20~(1) (2013) 31--59.

\bibitem{Maday-BMNP:04}
M.~Barrault, Y.~Maday, N.~Nguyen, A.~Patera, An `empirical interpolation'
  method: application to efficient reduced-basis discretization of partial
  differential equations, Comptes Rendus Mathematique 339~(9) (2004) 667 --
  672.

\bibitem{Iollo-ILD:00}
A.~Iollo, S.~Lanteri, J.-A. D{\'e}sid{\'e}ri, {S}tability {P}roperties of
  {POD}-{G}alerkin {A}pproximations for the {C}ompressible {N}avier-{S}tokes
  {E}quations, Theoretical and Computational Fluid Dynamics 13~(6) (2000)
  377--396.

\bibitem{Volkwein-KV:02}
K.~Kunisch, S.~Volkwein, Galerkin proper orthogonal decomposition methods for a
  general equation in fluid dynamics, SIAM Journal on Numerical Analysis 40~(2)
  (2002) 492--515.

\bibitem{Amsallem-AF-08}
D.~Amsallem, C.~Farhat, Interpolation method for adapting reduced-order models
  and application to aeroelasticity, AIAA Journal 46~(7) (2008) 1803--1813.

\bibitem{Chinesta-GACAC-18}
D.~Gonz{\'a}lez, J.~Aguado, E.~Cueto, E.~Abisset-Chavanne, F.~Chinesta,
  k{PCA}-based parametric solutions within the {PGD} framework, Archives of
  Computational Methods in Engineering 25~(1) (2018) 69--86.

\bibitem{Arroyo-MA-13}
D.~Mill{\'a}n, M.~Arroyo, Nonlinear manifold learning for model reduction in
  finite elastodynamics, Computer Methods in Applied Mechanics and Engineering
  261 (2013) 118--131.

\bibitem{Breitkopf-LRB-15}
G.~Le~Quilliec, B.~Raghavan, P.~Breitkopf, A manifold learning-based reduced
  order model for springback shape characterization and optimization in sheet
  metal forming, Computer Methods in Applied Mechanics and Engineering 285
  (2015) 621--638.

\bibitem{Breitkopf-MBLRV-18}
L.~Meng, P.~Breitkopf, G.~Le~Quilliec, B.~Raghavan, P.~Villon, Nonlinear
  shape-manifold learning approach: concepts, tools and applications, Archives
  of Computational Methods in Engineering 25~(1) (2018) 1--21.

\bibitem{Breitkopf-RXBRV-13}
B.~Raghavan, L.~Xia, P.~Breitkopf, A.~Rassineux, P.~Villon, Towards
  simultaneous reduction of both input and output spaces for interactive
  simulation-based structural design, Computer Methods in Applied Mechanics and
  Engineering 265 (2013) 174--185.

\bibitem{Breitkopf-XZBVZ-18}
M.~Xiao, G.~Zhang, P.~Breitkopf, P.~Villon, W.~Zhang, Extended {C}o-{K}riging
  interpolation method based on multi-fidelity data, Applied Mathematics and
  Computation 323 (2018) 120--131.

\bibitem{Bungartz-Griebel}
H.-J. Bungartz, M.~Griebel, Sparse grids, Acta Numerica 13~(1) (2004) 147--269.

\bibitem{Ibanez-IBAACLC-17}
R.~{Ib\'{a}\~nez}, D.~Borzacchiello, J.~Aguado, E.~Abisset-Chavanne, E.~Cueto,
  P.~Ladev{\`e}ze, F.~Chinesta, Data-driven non-linear elasticity: constitutive
  manifold construction and problem discretization, Computational Mechanics
  60~(5) (2017) 813--826.

\bibitem{Ibanez-IAAGCC-18}
R.~{Ib\'{a}\~nez}, E.~Abisset-Chavanne, J.~Aguado, D.~Gonzalez, E.~Cueto,
  F.~Chinesta, A manifold learning approach to data-driven computational
  elasticity and inelasticity, Archives of Computational Methods in Engineering
  25~(1) (2018) 47--57.

\bibitem{Borzacchiello-BAC-19}
D.~Borzacchiello, J.~Aguado, F.~Chinesta, Non-intrusive sparse subspace
  learning for parametrized problems, Archives of Computational Methods in
  Engineering 26~(2) (2019) 303--326.

\bibitem{Ibanez-IAAGCHDC-18}
R.~{Ib\'{a}\~nez}, E.~Abisset-Chavanne, A.~Ammar, D.~Gonz\'alez, E.~Cueto,
  A.~Huerta, J.-L. Duval, F.~Chinesta, A multi-dimensional data-driven sparse
  identification technique: the sparse {P}roper {G}eneralized {D}ecomposition,
  Complexity (5608286) (2018) 1--11.

\bibitem{Chinesta-CCADE-20}
F.~Chinesta, E.~Cueto, E.~Abisset-Chavanne, J.-L. Duval, F.~El~Khaldi, Virtual,
  digital and hybrid twins: a new paradigm in data-based engineering and
  engineered data, Archives of Computational Methods in Engineering 27~(1)
  (2020) 105--134.

\bibitem{Modesto-MYZH-20}
D.~Modesto, B.~Ye, S.~Zlotnik, A.~Huerta, Fast solution of elliptic harbor
  agitation problems under frequency-direction input spectra by model order
  reduction and {NURBS}-enhanced {FEM}, Coastal Engineering 156 (2020) 103618.

\bibitem{Chamoin-KCLP-19}
K.~Kergrene, L.~Chamoin, M.~Laforest, S.~Prudhomme, On a goal-oriented version
  of the proper generalized decomposition method, Journal of Scientific
  Computing 81~(1) (2019) 92--111.

\bibitem{Moitinho-RMDZ-20}
J.~Reis, J.~{Moitinho de Almeida}, P.~D{\'\i}ez, S.~Zlotnik, Error estimation
  for proper generalized decomposition solutions: {D}ual analysis and
  adaptivity for quantities of interest, International Journal for Numerical
  Methods in Engineering 121~(23) (2020) 5275--5294.

\bibitem{Smetana-SZ-20}
K.~Smetana, O.~Zahm, Randomized residual-based error estimators for the proper
  generalized decomposition approximation of parametrized problems,
  International Journal for Numerical Methods in Engineering 121~(23) (2020)
  5153--5177.

\bibitem{Aguado-AHCC-15}
J.~Aguado, A.~Huerta, F.~Chinesta, E.~Cueto, Real-time monitoring of thermal
  processes by reduced-order modeling, International Journal for Numerical
  Methods in Engineering 102~(5) (2015) 991--1017.

\bibitem{DM-MZH:15}
D.~Modesto, S.~Zlotnik, A.~Huerta, {P}roper {G}eneralized {D}ecomposition for
  parameterized {H}elmholtz problems in heterogeneous and unbounded domains:
  application to harbor agitation, Computer Methods in Applied Mechanics and
  Engineering 295 (2015) 127--149.

\bibitem{diez2017generalized}
P.~D{\'\i}ez, S.~Zlotnik, A.~Huerta, Generalized parametric solutions in
  {S}tokes flow, Computer Methods in Applied Mechanics and Engineering 326
  (2017) 223--240.

\bibitem{Sibileau-SGAMD-18}
A.~Sibileau, A.~Garc{\'\i}a-Gonz{\'a}lez, F.~Auricchio, S.~Morganti,
  P.~D{\'\i}ez, Explicit parametric solutions of lattice structures with proper
  generalized decomposition ({PGD}). {A}pplications to the design of
  3{D}-printed architectured materials, Computational Mechanics 62~(4) (2018)
  871--891.

\bibitem{Barroso-BGLMH-20}
G.~Barroso, A.~Gil, P.~Ledger, M.~Mallett, A.~Huerta, A regularised-adaptive
  {P}roper {G}eneralised {D}ecomposition implementation for coupled
  magneto-mechanical problems with application to {MRI} scanners, Computer
  Methods in Applied Mechanics and Engineering 358 (2020) 112640.

\bibitem{AH-AHCCL:14}
A.~Ammar, A.~Huerta, F.~Chinesta, E.~Cueto, A.~Leygue, Parametric solutions
  involving geometry: a step towards efficient shape optimization, Computer
  Methods in Applied Mechanics and Engineering 268 (2014) 178--193.

\bibitem{SZ-ZDMH:15}
S.~Zlotnik, P.~D\'{\i}ez, D.~Modesto, A.~Huerta, Proper generalized
  decomposition of a geometrically parametrized heat problem with geophysical
  applications, International Journal for Numerical Methods in Engineering
  103~(10) (2015) 737--758.

\bibitem{sevilla2020solution}
R.~Sevilla, S.~Zlotnik, A.~Huerta, Solution of geometrically parametrised
  problems within a {CAD} environment via model order reduction, Computer
  Methods in Applied Mechanics and Engineering 358 (2020) 112631.

\bibitem{RS-SBGH-20}
R.~Sevilla, L.~Borchini, M.~Giacomini, A.~Huerta, Hybridisable discontinuous
  {G}alerkin solution of geometrically parametrised {S}tokes flows, Computer
  Methods in Applied Mechanics and Engineering 372 (2020) 113397.

\bibitem{MG-GSH-20}
M.~Giacomini, R.~Sevilla, A.~Huerta, Tutorial on {H}ybridizable {D}iscontinuous
  {G}alerkin ({HDG}) formulation for incompressible flow problems, in: L.~{De
  Lorenzis}, A.~D\"{u}ster (Eds.), Modeling in Engineering Using Innovative
  Numerical Methods for Solids and Fluids, Vol. 599 of CISM International
  Centre for Mechanical Sciences, Springer International Publishing, 2020, pp.
  163--201.

\bibitem{Jay-CG:09}
B.~Cockburn, J.~Gopalakrishnan, The derivation of hybridizable discontinuous
  {G}alerkin methods for {S}tokes flow, SIAM Journal on Numerical Analysis
  47~(2) (2009) 1092--1125.

\bibitem{Cockburn-CS:14}
B.~Cockburn, K.~Shi, Devising {HDG} methods for {S}tokes flow: an overview,
  Computers \& Fluids 98 (2014) 221--229.

\bibitem{Nguyen-NPC:10}
N.~Nguyen, J.~Peraire, B.~Cockburn, A hybridizable discontinuous {G}alerkin
  method for {S}tokes flow, Computer Methods in Applied Mechanics and
  Engineering 199~(9-12) (2010) 582--597.

\bibitem{Jay-CGL:09}
B.~Cockburn, J.~Gopalakrishnan, R.~Lazarov, Unified hybridization of
  discontinuous {G}alerkin, mixed, and continuous {G}alerkin methods for second
  order elliptic problems, SIAM Journal on Numerical Analysis 47~(2) (2009)
  1319--1365.

\bibitem{Jay-CGNPS:11}
B.~Cockburn, J.~Gopalakrishnan, N.~Nguyen, J.~Peraire, F.-J. Sayas, Analysis of
  {HDG} methods for {S}tokes flow, Mathematics of Computation 80~(274) (2011)
  723--760.

\bibitem{HDG-NEFEM}
R.~Sevilla, A.~Huerta, {HDG-NEFEM} with degree adaptivity for {S}tokes flows,
  Journal of Scientific Computing 77~(3) (2018) 1953--1980.

\bibitem{giacomini2018superconvergent}
M.~Giacomini, A.~Karkoulias, R.~Sevilla, A.~Huerta, A superconvergent {HDG}
  method for {S}tokes flow with strongly enforced symmetry of the stress
  tensor, Journal of Scientific Computing 77~(3) (2018) 1679--1702.

\bibitem{Nguyen-CNP:10}
B.~Cockburn, N.~Nguyen, J.~Peraire, A comparison of {HDG} methods for {S}tokes
  flow, Journal of Scientific Computing 45~(1-3) (2010) 215--237.

\bibitem{Nguyen-NPC:09}
N.~Nguyen, J.~Peraire, B.~Cockburn, An implicit high-order hybridizable
  discontinuous {G}alerkin method for linear convection-diffusion equations,
  Journal of Computational Physics 228~(9) (2009) 3232--3254.

\bibitem{Nguyen-NPC:09b}
N.~Nguyen, J.~Peraire, B.~Cockburn, An implicit high-order hybridizable
  discontinuous {G}alerkin method for nonlinear convection-diffusion equations,
  Journal of Computational Physics 228~(23) (2009) 8841--8855.

\bibitem{Patera-Rozza:07}
A.~Patera, G.~Rozza, Reduced {B}asis {A}pproximation and {A-P}osteriori {E}rror
  {E}stimation for {P}arametrized {P}artial {D}ifferential {E}quations, {MIT
  Pappalardo Graduate Monographs in Mechanical Engineering}, Massachusetts
  Institute of Technology, Cambridge, MA, USA (2007).

\bibitem{Rozza:14}
G.~Rozza, Fundamentals of reduced basis method for problems governed by
  parametrized {PDE}s and applications, in: Separated representations and
  {PGD}-based model reduction, Vol. 554 of CISM Courses and Lectures, Springer,
  Vienna, 2014, pp. 153--227.

\bibitem{tsiolakis2020nonintrusive}
V.~Tsiolakis, M.~Giacomini, R.~Sevilla, C.~Othmer, A.~Huerta, Nonintrusive
  proper generalised decomposition for parametrised incompressible flow
  problems in {OpenFOAM}, Computer Physics Communications 249 (2020) 107013.

\bibitem{poya2016unified}
R.~Poya, R.~Sevilla, A.~Gil, A unified approach for a posteriori high-order
  curved mesh generation using solid mechanics, Computational Mechanics 58~(3)
  (2016) 457--490.

\bibitem{HO-Meshing}
Z.~Xie, R.~Sevilla, O.~Hassan, K.~Morgan, The generation of arbitrary order
  curved meshes for {3D} finite element analysis, Computational Mechanics
  51~(3) (2013) 361--374.

\bibitem{PD-GBCD-17}
R.~Garc{\'\i}a-Blanco, D.~Borzacchiello, F.~Chinesta, P.~D\'{i}ez, Monitoring a
  {PGD} solver for parametric power flow problems with goal-oriented error
  assessment, International Journal for Numerical Methods in Engineering
  111~(6) (2017) 529--552.

\bibitem{PD-GDBC-18}
R.~Garc{\'\i}a-Blanco, P.~D{\'\i}ez, D.~Borzacchiello, F.~Chinesta, Algebraic
  and parametric solvers for the power flow problem: towards real-time and
  accuracy-guaranteed simulation of electric systems, Archives of Computational
  Methods in Engineering 25~(4) (2018) 1003--1026.

\bibitem{PD-DZGH-18}
P.~D\'{i}ez, S.~Zlotnik, A.~Garc\'{i}a-Gonz\'{a}lez, A.~Huerta, Algebraic {PGD}
  for tensor separation and compression: an algorithmic approach, Comptes
  Rendus M\'{e}canique 346~(7) (2018) 501--514.

\bibitem{PD-DZGH-20}
P.~D\'{i}ez, S.~Zlotnik, A.~Garc\'ia-Gonz\'alez, A.~Huerta, Encapsulated {PGD}
  {A}lgebraic {T}oolbox {O}perating with {H}igh-{D}imensional {D}ata, Archives
  of Computational Methods in Engineering 27 (2020) 1321--1336.

\bibitem{Ladeveze-CNLB-16}
A.~Courard, D.~N{\'e}ron, P.~Ladev{\`e}ze, L.~Ballere, Integration of
  {PGD}-virtual charts into an engineering design process, Computational
  Mechanics 57~(4) (2016) 637--651.

\bibitem{zou2018nonintrusive}
X.~Zou, M.~Conti, P.~D{\'\i}ez, F.~Auricchio, A nonintrusive proper generalized
  decomposition scheme with application in biomechanics, International Journal
  for Numerical Methods in Engineering 113~(2) (2018) 230--251.

\bibitem{avron2005pushmepullyou}
J.~Avron, O.~Kenneth, D.~Oaknin, Pushmepullyou: an efficient micro-swimmer, New
  Journal of Physics 7~(1) (2005) 234.

\bibitem{Alouges-ADL-09}
F.~Alouges, A.~DeSimone, A.~Lefebvre, Optimal strokes for axisymmetric
  microswimmers, The European Physical Journal E 28~(3) (2009) 279--284.

\bibitem{lovgren2009global}
A.~L{\o}vgren, Y.~Maday, E.~R{\o}nquist, Global $\mathcal{C}^1$ maps on general
  domains, Mathematical Models and Methods in Applied Sciences 19~(05) (2009)
  803--832.

\end{thebibliography}

%==========================================================================
\appendix

%==========================================================================
\section{Separated expressions of the bilinear and linear forms}  
\label{app:separatedPGD}
%==========================================================================

In this appendix, the separated form of the PGD approximation of the HDG local~\eqref{eq:localWeak} and global~\eqref{eq:globalWeak} problems is briefly reported. For its detailed derivation, interested readers are referred to~\cite{RS-SBGH-20}.

From the separated form~\eqref{eq:PGDapprox}, the computation of the $m$-th mode is performed in two steps, corresponding to the HDG local and global problems. The PGD spatial equation arising from the local problem~\eqref{eq:localWeak} is: find $(\ampL^m \De \fL, \ampU^m \De \fU, \ampP^m \De \fP) \in \Wh \times \Vh \times \sVh$ such that
\begin{equation}\label{eq:localWeakPGDSpatial}
	\begin{aligned}
	\sum_{k=1}^{\ndet} \beta_{\theta}^k \mathcal{A}^k_{LL}(\de \fL, \ampL^m \De \fL) & + \sum_{k=1}^{\nadj} \beta^k_{\vartheta} \mathcal{A}^k_{Lu}(\de \fL, \ampU^m \De \fU)	\\
	= & \mathcal{R}_L^m (\de \fL \psi^m) 
	+ \sum_{k=1}^{\nadj} \beta^k_{\vartheta} \mathcal{A}^k_{L\hu}(\de \fL, \ampHU^m \De \fHU),  
	\\	
	\sum_{k=1}^{\nadj} \beta^k_{\vartheta} \mathcal{A}^k_{uL}(\de \fU, \ampL^m \De \fL)  + \beta & \mathcal{A}_{uu} (\de \fU, \ampU^m \De \fU) \\
	+ \sum_{k=1}^{\nadj} \beta^k_{\vartheta} \mathcal{A}^k_{up}(\de \fU, \ampP^m \De \fP)  = & \mathcal{R}_u^m (\de \fU \psi^m) + \beta \mathcal{A}_{u\hu}(\de \fU, \ampHU^m \De \fHU),
	\\
	\sum_{k=1}^{\nadj} \beta^k_{\vartheta} \mathcal{A}^k_{pu}(\de \fP, \ampU^m \De \fU)  = & \mathcal{R}_p^m (\de \fP \psi^m) + \sum_{k=1}^{\nadj} \beta^k_{\vartheta} \mathcal{A}^k_{p\hu}(\de \fP, \ampHU^m \De \fHU) 
	\\
	\beta \mathcal{A}_{\rho p}(1, \ampP^m \De \fP)  = & \mathcal{R}_{\overline{p}}^m (\psi^m) + \beta \mathcal{A}_{\rho \rho}(1, \ampR^m \De \fR),
	\end{aligned}
\end{equation}
for all $(\de \fL,\de \fU,\de \fP) \in \Wh \times \Vh \times \sVh$.  \\
The separated bilinear forms in equation~\eqref{eq:localWeakPGDSpatial} are given by
\begin{equation} \label{eq:bilinearPGDSpatialLocal}
	\begin{aligned}	
	\mathcal{A}^k_{LL}(\de \fL, \fL)	 & := - \intE{\de \fL}{\nu^{-1} D^k \fL }, 	&	
	\mathcal{A}^k_{Lu}(\de \fL, \fU)	 & :=  \intE{\bA^k \DivX \de \fL}{ \fU },  		\\
	\mathcal{A}^k_{L\hu}(\de \fL, \fHU)& := \intBNoD{\bA^k \bn \cdot \de \fL}{ \fHU }, &
	\mathcal{A}^k_{uL}(\de \fU, \fL)	 & :=  \intE{\de \fU}{\bA^k \DivX \fL },  		\\
	\mathcal{A}_{uu}(\de \fU, \fU)	 & :=  \intBE{\de \fU}{\btau \fU },				&
	\mathcal{A}^k_{up}(\de \fU, \fP)	 & :=  \intE{\de \fU}{\bA^k \GradX \fP },  		\\
	\mathcal{A}_{u\hu}(\de \fU, \fHU)	 & := \intBNoD{\de \fU}{\btau \fHU }			&
	\mathcal{A}^k_{pu}(\de \fP, \fU)	 & :=   \intE{\bA^k \GradX \de \fP}{\fU},  		\\
	\mathcal{A}^k_{p\hu}(\de \fP, \fHU)& :=  \intBNoD{\de \fP}{\fHU \cdot \bA^k \bn},  &
	\mathcal{A}_{\rho p}(\de\fR, \fP)  & :=  \intE{\de\fR}{|\Omega_e|^{-1} \fP },						\\
	\mathcal{A}_{\rho \rho}(\de\fR, \fR) & :=  \de\fR \,\fR, 
\end{aligned}	
\end{equation}
whereas the corresponding linear forms are
\begin{equation} \label{eq:linearLocalPGDSpatial}
\begin{aligned}	
\mathcal{R}_L^m (\de \fL \psi)	:= & \sum_{k=1}^{\nadj} \sum_{l=1}^{\nDir} \intBD{\bA^k \bn \cdot \de \fL}{\bgD^l} \mathcal{A}^k_{\vartheta}(\psi, \lambda_D^l) 	\\
& {-} \sum_{i=1}^{m} \sum_{k=1}^{\ndet} \mathcal{A}^k_{LL}(\de \fL, \ampL^i \fL^i) \mathcal{A}^k_{\theta}(\psi, \psi^i)   \\
& {-} \sum_{i=1}^{m} \sum_{k=1}^{\nadj} \left\{ \mathcal{A}^k_{Lu}(\de \fL, \ampU^i \fU^i) - \mathcal{A}^k_{L\hu}(\de \fL, \ampHU^i \fHU^i)\right\} \mathcal{A}^k_{\vartheta}(\psi, \psi^i) 
\\
\mathcal{R}_u^m (\de \fU \psi)	 := & \sum_{k=1}^{\ndet} \sum_{l=1}^{\nSou} \intE{\de \fU}{D^k \bgS^l} \mathcal{A}^k_{\theta}(\psi, \lambda_S^l) 	\\
& {+}  \sum_{l=1}^{\nDir} \intBD{\de \fU}{\btau \bgD^l} \mathcal{A}(\psi, \lambda_D^l)  \\
& {-}  \sum_{i=1}^{m} \sum_{k=1}^{\nadj} \left\{ \mathcal{A}^k_{uL}(\de \fU, \ampL^i \fL^i) + \mathcal{A}^k_{up}(\de \fU, \ampP^i \fP^i) \right\} \mathcal{A}^k_{\vartheta}(\psi, \psi^i) \\ 
& {-}  \sum_{i=1}^{m} \left\{ \mathcal{A}_{uu}(\de \fU, \ampU^i \fU^i) - \mathcal{A}_{u\hu}(\de \fU, \ampHU^i \fHU^i) \right\} \mathcal{A}(\psi, \psi^i)
\\
\mathcal{R}_p^m (\de \fP \psi)	 := & \sum_{k=1}^{\nadj} \sum_{l=1}^{\nDir} \intBD{\de \fP}{\bgD^l \cdot \bA^k \bn} \mathcal{A}^k_{\vartheta}(\psi, \lambda_D^l)	\\
& {-}  \sum_{i=1}^{m} \sum_{k=1}^{\nadj} \left\{ \mathcal{A}^k_{pu}(\de \fP, \ampU^i \fU^i) - \mathcal{A}^k_{p\hu}(\de \fP, \ampHU^i \fHU^i) \right\} \mathcal{A}^k_{\vartheta}(\psi, \psi^i)
\\
\mathcal{R}_{\overline{p}}^m (\de\fR\psi) := & -\sum_{i=1}^{m} \left\{ \mathcal{A}_{\rho p}(\de\fR, \ampP^i \fP^i) - \mathcal{A}_{\rho \rho}(\de\fR, \ampR^i \fR^i) \right\} \mathcal{A}(\psi, \psi^i) .
\end{aligned}	
\end{equation}
Finally, the parametric constants appearing in equation~\eqref{eq:localWeakPGDSpatial} are defined as
\begin{equation} \label{eq:cttLocalPGDSpatial}
\beta^k_{\theta}    :=  \mathcal{A}^k_{\theta}(\psi^m,\psi^m) 		\qquad
\beta^k_{\vartheta} :=  \mathcal{A}^k_{\vartheta}(\psi^m,\psi^m),  	\qquad
\beta    			:=  \mathcal{A}(\psi^m,\psi^m) ,
\end{equation}
where the bilinear forms in the parametric space are given by
\begin{equation} \label{eq:bilinearPGDParam}
\begin{aligned}	
\mathcal{A}^k_{\theta}(\de \psi, \psi)	& :=  \intI{\de \psi}{\theta^k \psi},  		\\
\mathcal{A}^k_{\vartheta}(\de \psi, \psi)	& :=  \intI{\de \psi}{\vartheta^k \psi}, 		\\
\mathcal{A}(\de \psi, \psi)				& :=  \intI{\de \psi}{\psi}.	
\end{aligned}	
\end{equation}

The separated approximation of the trial and test functions, see equations~\eqref{eq:PGDapprox} and~\eqref{eq:tangentU}, is also exploited to construct the separated form of the HDG global problem~\eqref{eq:globalWeak}: find $\ampHU^m \De \fHU \in \hVh$ and $\ampR^m \De \fR \in \mathbb{R}^{\numel}$ such that, for all $\de \fHU \in \hVh$, it holds
\begin{equation}\label{eq:globalWeakPGDSpatialTransmission}
	\begin{aligned}	
		\sum_{e=1}^{\numel} \Biggl\{ \sum_{k=1}^{\nadj}\beta^k_{\vartheta} \mathcal{A}^k_{\hu L}  (\de \fHU, \ampL^m \De \fL) 	
	 + \beta \mathcal{A}_{\hu u}(\de \fHU, \ampU^m & \De \fU) \\
	 + \sum_{k=1}^{\nadj} \beta^k_{\vartheta} \mathcal{A}^k_{\hu p}(\de \fHU, \ampP^m \De \fP)  
     + \beta \mathcal{A}_{\hu \hu}(& \de \fHU, \ampHU^m \De \fHU) \\
	  + \sum_{k=1}^{\nadj} \beta^k_{\vartheta}\mathcal{A}^k_{\hu \hu} (\de \fHU, \ampHU^m \De \fHU) \Biggr\}
	& = \sum_{e=1}^{\numel} \mathcal{R}_{\hu}^m (\de \fHU \psi^m),  \\[1ex]
	\sum_{k=1}^{\nadj} \beta^k_{\vartheta} \mathcal{A}^k_{p \hu} (1,\ampHU^m \De \fHU)  &= \mathcal{R}_{\rho}^m (\psi^m), \\[-2ex]
	& \quad \ e=1,\dots,\numel .
	\end{aligned}
\end{equation}
The bilinear forms in equation~\eqref{eq:globalWeakPGDSpatialTransmission} are defined as
\begin{equation} \label{eq:bilinearPGDSpatialGlobal}
\begin{aligned}	
	\mathcal{A}^k_{\hu L}(\de \fHU, \fL)	 	& := \intBNoDS{\de \fHU}{\bA^k \bn \cdot \fL } - \intBS{\de \fHU}{\bA^k \bn \cdot \fL \bE },  		\\
	\mathcal{A}_{\hu u}(\de \fHU, \fU)     	& :=  \intBNoDS{\de \fHU}{\btau \fU} - \intBS{\de \fHU}{(\btau \fU) \!\cdot\! \bE} , 	\\
	\mathcal{A}^k_{\hu p}(\de \fHU, \fP)   	& :=  \intBNoDS{\de \fHU}{\fP \bA^k \bn }, \\
	\mathcal{A}_{\hu \hu}(\de \fHU, \fHU)		& := -\intBNoDS{\de \fHU}{\btau \fHU}  + \intBS{\de \fHU}{(\btau \fHU) \!\cdot\! \bE} ,  \\
	\mathcal{A}^k_{\hu \hu}(\de \fHU, \fHU)	& := \intBS{\de \fHU}{\fHU \cdot \bA^k \bD},
\end{aligned}	
\end{equation}
whereas the corresponding linear form are given by
\begin{equation} \label{eq:linearGlobalPGDSpatial}
\begin{aligned}	
\mathcal{R}_{\hu}^m (\de \fHU \psi)	:=  & -  \sum_{l=1}^{\nNeu} \intBN{\de \fHU}{\bgN^l} \mathcal{A}(\psi, \lambda_N^l) \\
& {-}  \sum_{i=1}^{m} \left\{ \mathcal{A}_{\hu u}(\de \fHU, \ampU^i \fU^i)  +  \mathcal{A}_{\hu \hu}(\de \fHU, \ampHU^i \fHU^i) \right\} \mathcal{A}(\psi, \psi^i) \\
& {-}  \sum_{i=1}^{m} \sum_{k=1}^{\nadj}  \Bigl\{ \mathcal{A}^k_{\hu L}(\de \fHU, \ampL^i \fL^i) \mathcal{A}^k_{\vartheta}(\psi, \psi^i) \\
& {+} \left[ \mathcal{A}^k_{\hu p}(\de \fHU, \ampP^i \fP^i) + \mathcal{A}^k_{\hu \hu}(\de \fHU, \ampHU^i \fHU^i) \right] \mathcal{A}^k_{\vartheta}(\psi, \psi^i) \Bigr\} ,
\\
\mathcal{R}_{\rho}^m (\de \fR \psi)	 := & -\sum_{k=1}^{\nadj} \sum_{l=1}^{\nDir} \intBD{\de \fR}{\bgD^l \cdot \bA^k \bn} \mathcal{A}^k_{\vartheta}(\psi, \lambda_D^l)	\\
& {-}  \sum_{i=1}^{m} \sum_{k=1}^{\nadj} \mathcal{A}^k_{p\hu}(\de \fR, \ampHU^i \fHU^i) \mathcal{A}^k_{\vartheta}(\psi, \psi^i) .
\end{aligned}	
\end{equation}

Following from remark~\ref{rmrk:paramFunc}, a unique parametric function is considered for all the variables in the PGD approximation~\eqref{eq:PGDapprox}. Hence, the PGD parametric problem is: find $\De\psi \in \Lh(\bI)$ such that 
\begin{equation}\label{eq:WeakPGDParam}
\sum_{k=1}^{\ndet} \gamma_{LL}^k \mathcal{A}^k_{\theta}(\de\psi,\De\psi) +  \sum_{k=1}^{\nadj} \gamma_{\vartheta}^k \mathcal{A}_{\vartheta}^k(\de\psi,\De\psi) + \gamma \mathcal{A}(\de \psi,\De\psi) =
\mathcal{R}^m(\de\psi),
\end{equation}
for all $\de \psi \in \Lh(\bI)$, where
\begin{equation}\label{eq:WeakPGDParamCtt}
\begin{aligned}
\gamma_{\vartheta}^k 	 :=& \gamma_{Lu}^k - \gamma_{L\hu}^k + \gamma_{uL}^k + \gamma_{up}^k + \gamma_{pu}^k - \gamma_{p \hu}^k + \gamma^k_{\hu L} + \gamma^k_{\hu p} + \gamma^k_{\hu \hu} + \gamma_{\rho \hu}^k , \\
\gamma   				 :=& \gamma_{uu} - \gamma_{u\hu} + \gamma_{\rho p} - \gamma_{\rho \rho} + \gamma_{\hu u} + \gamma_{\hu \hu}, \\
\mathcal{R}^m(\de\psi) 			  :=& \mathcal{R}_L^m (\ampL^m \fL^m \de \psi) + \mathcal{R}_u^m (\ampU^m \fU^m \de \psi) + \mathcal{R}_p^m (\ampP^m \fP^m \de \psi)\\
&+ \mathcal{R}_{\overline{p}}^m (\de \psi) + \mathcal{R}_{\hu}^m (\ampHU^m \fHU^m \de \psi) + \mathcal{R}_{\rho}^m (\de \psi).
\end{aligned}
\end{equation}

The spatial constants appearing in equation~\eqref{eq:WeakPGDParam} are defined as
\begin{equation} \label{eq:cttPGDParam}
\begin{aligned}	
\gamma_{LL}^k 	& := \mathcal{A}^k_{LL}		(\ampL^m \fL^m, \ampL^m \fL^m),		&
\gamma_{Lu}^k 	& := \mathcal{A}^k_{Lu}		(\ampL^m \fL^m, \ampU^m \fU^m), 	\\
\gamma_{L\hu}^k & := \mathcal{A}^k_{L\hu}	(\ampL^m \fL^m, \ampHU^m \fHU^m), 	&
\gamma_{uL}^k 	& := \mathcal{A}^k_{uL}		(\ampU^m\fU^m, \ampL^m \fL^m),		\\
\gamma_{uu} 	& := \mathcal{A}_{uu}		(\ampU^m\fU^m, \ampU^m\fU^m),		&
\gamma_{up}^k	& := \mathcal{A}^k_{up}		(\ampU^m\fU^m, \ampP^m\fP^m),		\\
\gamma_{u\hu}	& := \mathcal{A}_{u\hu}		(\ampU^m\fU^m, \ampHU^m\fHU^m),	&
\gamma_{pu}^k	& := \mathcal{A}^k_{pu}		(\ampP^m \fP^m, \ampU^m\fU^m),		\\
\gamma_{p\hu}^k	& := \mathcal{A}^k_{p\hu}	(\ampP^m\fP^m, \ampHU^m\fHU^m),  	&
\gamma_{\rho p}	& := \mathcal{A}_{\rho p}	(1, \ampP^m\fP^m), 	\\
\gamma_{\rho \rho}& := \mathcal{A}_{\rho \rho}(1, \ampR^m\fR^m), & & \\
\gamma_{\hu L}^k 	& := \mathcal{A}^k_{\hu L}		(\ampHU^m\fHU^m, \ampL^m \fL^m),		&
\gamma_{\hu u} 	& := \mathcal{A}_{\hu u}		(\ampHU^m\fHU^m, \ampU^m\fU^m),		\\
\gamma_{\hu p}^k	& := \mathcal{A}^k_{\hu p}		(\ampHU^m\fHU^m, \ampP^m\fP^m),		&
\gamma_{\hu \hu}	& := \mathcal{A}_{\hu \hu}		(\ampHU^m\fHU^m, \ampHU^m\fHU^m),	\\
\gamma_{\hu \hu}^k	& := \mathcal{A}_{\hu\hu}^k		(\ampHU^m\fHU^m, \ampHU^m\fHU^m),	&
\gamma_{\rho\hu}^k	& := \mathcal{A}_{p\hu}^k		(1, \ampHU^m\fHU^m).
\end{aligned}	
\end{equation}

\end{document}